\documentclass[a4paper,11pt,french]{article} 
        \usepackage[french,english]{babel}
\usepackage[T1]{fontenc}
\usepackage[latin1]{inputenc}
\usepackage[margin=1in]{geometry}

        \usepackage{amsfonts}
        \usepackage{amstext}
        \usepackage{mathtools}
        \usepackage{appendix}
        \usepackage{enumerate}
        \usepackage{exscale}
        \usepackage{relsize}
        \usepackage{url}
        \usepackage[all]{xy}
        \usepackage{paralist}
        \usepackage{dsfont}
\usepackage{latexsym}
\usepackage{upgreek}
\usepackage{mathabx}
\usepackage{amssymb,amsmath,amsthm,amscd}
\usepackage{mathrsfs}
\usepackage{hyperref}
\usepackage[nottoc]{tocbibind}

\numberwithin{equation}{section}
\usepackage{makeidx}
\makeindex

\newtheorem{theo}{Th\'eor\`eme}[section]
\newtheorem{lem}[theo]{Lemme}
\newtheorem{prop}[theo]{Proposition}
\newtheorem{cor}[theo]{Corollaire}
\newtheorem{rem}[theo]{Remarque}

\newtheorem{defi}[theo]{D\'efinition}

\newcommand{\beqyn}{\begin{eqnarray}}
\newcommand{\enqyn}{\end{eqnarray}}
\newcommand{\blem}{\begin{lem}}
\newcommand{\elem}{\end{lem}}
\newcommand{\brop}{\begin{prop}}
\newcommand{\erop}{\end{prop}}
\newcommand{\bcor}{\begin{cor}}
\newcommand{\ecor}{\end{cor}}
\newcommand{\brem}{\begin{rem}}
\newcommand{\erem}{\end{rem}}
\newcommand{\benum}{\begin{enumerate}}
\newcommand{\enum}{\end{enumerate}}

\newcommand{\btheo}{\begin{theo}}
\newcommand{\etheo}{\end{theo}}
\newcommand{\bdem}{\begin{proof}}
\newcommand{\edem}{\end{proof}}
\newcommand{\bdefi}{\begin{defi}}
\newcommand{\edefi}{\end{defi}}

\newcommand{\rmE}{\mathrm{E}}
\newcommand{\rmF}{\mathrm{F}}

\newcommand{\rmK}{\mathrm{K}}

\newcommand{\rmN}{\mathrm{N}}

\newcommand{\Hom}{\mathrm{Hom}}
\newcommand{\Gal}{\mathrm{Gal}}
\newcommand{\Lie}{\mathrm{Lie}}
\newcommand{\Id}{\mathrm{Id}}

\newcommand{\Ad}{\mathrm{Ad}}

\newcommand{\Tr}{\mathrm{Tr}}
\newcommand{\End}{\mathrm{End}}

\newcommand{\Gl}{\mathrm{GL}}

\newcommand{\Res}{\mathrm{Res}}

\newcommand{\ind}{\mathds{1}}

\newcommand{\upla}{\underline{\rho}}
\newcommand{\uphi}{\upphi}
\newcommand{\al}{\alpha}

\newcommand{\la}{\lambda}
\newcommand{\La}{\Lambda}

\newcommand{\bilif}{\langle \cdot,\cdot \rangle}

\newcommand{\dsl}{\displaystyle \left(}
\newcommand{\rb}{\right)}

\newcommand{\bsl}{\backslash}

\newcommand{\rar}{\rightarrow}
\newcommand{\irar}{\xrightarrow{\sim}}
\newcommand{\hrar}{\hookrightarrow}
\newcommand{\smin}{\smallsetminus}
\newcommand{\sbs}{\subseteq}
\newcommand{\sbn}{\subsetneq}

\newcommand{\sps}{\supseteq}

\newcommand{\ps}[1]{\prescript{#1}{}}

\newcommand{\Gf}{G(\rmF)}

\newcommand{\Pf}{P(\rmF)}
\newcommand{\Qf}{Q(\rmF)}

\newcommand{\Uf}{U(\rmF)}

\newcommand{\PGf}{\Pf \bsl \Gf}

\newcommand{\PQf}{\Pf \bsl \Qf}

\newcommand{\PUf}{\Pf \bsl \Uf}

\newcommand{\gll}{\mathfrak{gl}}

\newcommand{\all}{\mathfrak{a}}

\newcommand{\gl}{\mathfrak{g}}

\newcommand{\ml}{\mathfrak{m}}
\newcommand{\nl}{\mathfrak{n}}

\newcommand{\sn}{\mathfrak{s}}

\newcommand{\ul}{\mathfrak{u}}
\newcommand{\zl}{\mathfrak{z}}

\newcommand{\ol}{\mathfrak{o}}
\newcommand{\Sgl}{\mathfrak{S}}

\newcommand{\N}{\mathbb{N}}
\newcommand{\Z}{\mathbb{Z}}
\newcommand{\A}{\mathbb{A}}
\newcommand{\C}{\mathbb{C}}

\newcommand{\Q}{\mathbb{Q}}
\newcommand{\R}{\mathbb{R}}

\newcommand{\Gm}{\mathbb{G}_m}

\newcommand{\calB}{\mathcal{B}}
\newcommand{\calC}{\mathcal{C}}

\newcommand{\calF}{\mathcal{F}}

\newcommand{\calH}{\mathcal{H}}

\newcommand{\calM}{\mathcal{M}}

\newcommand{\calO}{\mathcal{O}}
\newcommand{\calP}{\mathcal{P}}

\newcommand{\calS}{\mathcal{S}}

\newcommand{\calU}{\mathcal{U}}
\newcommand{\calV}{\mathcal{V}}

\newcommand{\calX}{\mathcal{X}}

\newcommand{\cad}{c'est-\`a-dire }

\newcommand{\nin}{\notin}
\newcommand{\resp}[1]{(resp. #1)}

\newcommand{\tlal}{\widetilde{\alpha}}
\newcommand{\tlgam}{\tilde \gamma}
\newcommand{\tldlt}{\tilde \delta}

\newcommand{\tlxi}{\tilde \xi}

\newcommand{\tlvpi}{\widetilde{\varpi}}

\newcommand{\tlPhi}{\tilde \Phi}

\newcommand{\tlB}{\widetilde{B}}

\newcommand{\tlE}{\widetilde{E}}

\newcommand{\tlG}{\widetilde{G}}
\newcommand{\tlH}{\widetilde{H}}

\newcommand{\tlK}{\widetilde{K}}

\newcommand{\tlM}{\widetilde{M}}

\newcommand{\tlP}{\widetilde{P}}
\newcommand{\tlR}{\widetilde{R}}
\newcommand{\tlQ}{\widetilde{Q}}
\newcommand{\tlS}{\widetilde{S}}
\newcommand{\tlT}{\widetilde{T}}
\newcommand{\tlU}{\widetilde{U}}

\newcommand{\tlX}{\widetilde{X}}

\newcommand{\tlb}{\widetilde{b}}

\newcommand{\tlg}{\widetilde{g}}
\newcommand{\tlh}{\widetilde{h}}

\newcommand{\tlk}{\widetilde{k}}

\newcommand{\tlm}{\widetilde{m}}
\newcommand{\tln}{\widetilde{n}}

\newcommand{\tls}{\widetilde{s}}

\newcommand{\tlx}{\widetilde{x}}
\newcommand{\tly}{\widetilde{y}}

\newcommand{\tlul}{\widetilde{\ul}}

\newcommand{\tlgl}{\widetilde{\gl}}

\newcommand{\tlzero}{\widetilde{0}}
\newcommand{\tlone}{\widetilde{1}}
\newcommand{\tltwo}{\widetilde{2}}
\newcommand{\tltr}{\widetilde{3}}
\newcommand{\tlqt}{\widetilde{4}}
\newcommand{\tlcq}{\widetilde{5}}
\newcommand{\tlsix}{\widetilde{6}}

\newcommand{\htau}{\hat \tau}

\newcommand{\hDelta}{\widehat \Delta}

\newcommand{\hiota}{\hat \iota}

\newcommand{\brg}{\bar g}

\newcommand{\brn}{\bar n}

\newcommand{\brvpi}{\overline{\varpi}}

\newcommand{\brI}{\overline{I}}

\newcommand{\brX}{\overline{X}}

\newcommand{\brdel}{\overline{\delta}}

\newcommand{\alltz}{\all_{\tlzero}}

\newcommand{\tlBmin}{\tlB_{0}}
\newcommand{\relPz}{\calF(M_{\tlzero}, P_{0})}
\newcommand{\relPh}{\calF(M_{\tlH}, P_{H})}
\newcommand{\relPb}{\calF(M_{\tlzero}, B)}
\newcommand{\relPbe}{\calF(M_{\tlzero, \rmE}, B_{\rmE})}

\newcommand{\Laz}{\Lambda_{m}}
\newcommand{\Lazt}{\Lambda_{m,3}}

\newcommand{\Laiyd}{\La_{IY,2}}

\newcommand{\Lajlr}{\La_{JLR}}
\newcommand{\Lajlrq}{\La_{JLR,4}}
\newcommand{\Lad}{\La_{d}}
\newcommand{\Ladud}{\La_{d,12}}
\newcommand{\Latqd}{\La_{d,34}}

\newcommand{\relPr}{\calF(M_{\tlzero}, B')}
\newcommand{\hGp}{(H_{\rmE} \times G'_{\rmE}) \times (H_{\rmE} \times \tlG'_{\rmE})}
\newcommand{\hGpa}{ (H_{\rmE}(\A)^{1} \times G'_{\rmE}(\A)) \times (H_{\rmE}(\A)^{1} \times \tlG'_{\rmE}(\A)) }

\newcommand{\GetlG}{G_{\rmE} \times \tlG_{\rmE}}
\newcommand{\GetlGa}{(G_{\rmE} \times \tlG_{\rmE})(\A)}

\newcommand{\relP}{\relPb}

\newcommand{\bX}{\mathbb{X}}

\numberwithin{equation}{section}
\begin{document}
\selectlanguage{french}
\title{Les formules des traces relatives de {J}acquet-{R}allis 
    grossières}
\date{}
\author{Micha\l \ Zydor
\footnote{Faculty of Mathematics and Computer Science, 
Weizmann Institute of Science
P.O. Box 26, 
Rehovot 76100,
Israel.
mail: \texttt{michalz@weizmann.ac.il}
}}
\maketitle 
\selectlanguage{english}
\begin{abstract} 
We establish the coarse relative trace formulae of Jacquet-Rallis for linear and unitary groups. 
Both formulae are of the form: a sum of spectral distributions equals a sum of geometric distributions. 
In order to obtain the spectral decompositions we introduce new truncation operators and we investigate their properties. 
On the geometric side, by means of the Cayley transform, 
the decompositions are derived from a procedure of descent to the tangent spaces for which the formulae 
are known thanks to our previous work. 
\end{abstract}
\selectlanguage{french}
\begin{abstract}
On établit les formules des traces relatives de Jacques-Rallis grossières pour les groupes 
linéaires et unitaires. 
Les deux formules sont sous la forme suivante: une somme des distributions spectrales est égale à une somme des distributions géométriques. 
Pour établir les développements spectraux on introduit de nouveaux opérateurs de troncature et on étudie 
leur propriétés. 
Du côté géométrique, en utilisant les applications de Cayley, 
les développements s'obtiennent par un argument de descente vers les espaces 
tangents pour lesquels les formules sont connues grâce à nos travaux précédents. 
\end{abstract}

\section*{Introduction}

\subsection{Le contexte}

La motivation de ce travail est l'article \cite{jacqrall} de Jacquet et Rallis
qui propose une approche 
à la conjecture globale de Gan-Gross-Prasad (GGP) pour le produit de groupes unitaires $U(n) \times U(n+1)$ via une formule des traces relative .
Dans cet article on présente les versions grossières de ces formules, \cad une identité entre une somme des distributions spectrales 
indexée par les données cuspidales et une somme des distributions géométriques suivants les éléments d'un quotient catégorique.

Décrivons d'abord brièvement la conjecture GGP. 
 Soit $\rmE/\rmF$
une extension quadratique de corps de nombres et soit $\sigma$ 
le générateur du groupe de Galois de cette extension. 
Soit $W$ un $\rmF$-espace vectoriel de dimension finie $n + 1$ 
muni d'une décomposition en somme directe $W = V \oplus D$ 
avec $V$ de dimension $n$ et $D$ de dimension $1$, où  $n \in \N$. 
On note $G = \Gl(V)$ que l'on voit comme un sous-groupe 
de $\tlG := \Gl(W)$ qui agit trivialement sur $D$. 
On note $G_{\rmE} = \Res_{\rmE/\rmF}\Gl(V_{\rmE})$ 
et $\tlG_{\rmE} := \Res_{\rmE/\rmF}\Gl(W_{\rmE})$ 
où $V_{\rmE} = V \otimes_{\rmF} \rmE$ etc.
On a donc des inclusions naturelles $G \hrar G_{\rmE} \hrar \tlG_{\rmE} \hookleftarrow \tlG$.  
On suppose que $W_{\rmE}$ est muni d'une forme $\sigma$-hermitienne 
non-dégénérée $\tlPhi$ de façon que $V_{\rmE}$ soit orthogonal à $D_{\rmE}$. 
On note $U = U(V_{\rmE}, \tlPhi|_{V_{\rmE}})$ le groupe unitaire associé qu'on voit comme un sous-groupe de 
$\tlU := U(W_{\rmE}, \tlPhi)$ qui agit trivialement sur $D_{\rmE}$. 
Le groupe $\tlU$ est aussi naturellement 
un sous-groupe de $\tlG_{\rmE}$.

Soit $\A$ l'anneau des adèles de $\rmF$. 
Pour un $\rmF$-groupe algébrique $H$ quelconque on note $[H] = H(\rmF) \bsl H(\A)$. 
Pour une représentation automorphe cuspidale $\pi$ de $U \times \tlU$ définissons
$\calP_{U,\pi} : \pi \rar \C$ par
\[
\pi \ni \upphi \mapsto \int_{[U]}\upphi(x,x)dx.
\]
On appelle $\calP_{U,\pi}$ une période. 
La conjecture GGP prédit alors, en gros, que  
la valeur centrale $L(\Pi, 1/2)$ de la fonction $L$ du changement de base $\Pi$ de $\pi$
à $G_{\rmE} \times \tlG_{\rmE}$ est non-nulle si est seulement si
la période $\calP_{U,\pi}$, 
quitte à changer $(U \times \tlU, \pi)$ 
par un autre membre du $L$-paquet de Vogan,
est non-triviale. On renvoie à \cite{ggp} pour les 
énoncés précis. De plus, Ichino et Ikeda \cite{ichinoIkeda} ont formulé un raffinement 
de la conjecture GGP dans le cas des groupes orthogonaux et le cas des groupes unitaires discuté 
ici a été traité par N. Harris \cite{harris}. Ce raffinement 
exprime la valeur $L(\Pi, 1/2)$ en fonction de la période $\calP_{U,\pi}$.


On décrit maintenant  l'approche de Jacquet et Rallis \cite{jacqrall} à la conjecture GGP pour les groupes unitaires. 
Soient $F_{U} \in C_{c}^{\infty}((U \times \tlU)(\A))$ et $F_{G} \in C_{c}^{\infty}((G_{\rmE} \times \tlG_{\rmE})(\A))$ des fonctions 
lisses à support compact. Soient $k_{F_{U}}$ et $k_{F_{G}}$ leur noyaux automorphes.
Jacquet et Rallis proposent d'étudier les intégrales suivantes:
\begin{equation*}
J(F_{U}) := \int_{[U]}\int_{[U]}k_{F_{U}}(x,x,y,y)dxdy, \quad
I(F_{G}) := \int_{[G_{\rmE}]}\int_{[G]}\int_{[\tlG]}
k_{F_{G}}(g,g,h,\tlh)\eta(h,\tlh)d\tlh dh dg,
\end{equation*}
où $\eta(h,\tlh) = \eta_{\rmE/\rmF}(\det h)^{n+1} \eta_{\rmE/\rmF}(\det \tlh)^{n}$ et 
$\eta_{\rmE/\rmF}$ c'est le caractère quadratique sur le groupe des idèles de $\A$ associé à 
l'extension $\rmE/\rmF$ par la théorie du corps de classes. 
 On note aussitôt que les intégrales ci-dessus ne sont pas en général convergentes. 
Pourtant, si $F_{U}$ et $F_{G}$ vérifient certaines conditions locales supplémentaires les intégrales
convergent et 
elles ont des développements spectraux, avec des contributions des représentations cuspidales seulement, 
du type $J(F_{U}) = \sum_{\pi}J_{\pi}(F_{U})$ et $I(F_{G}) = \sum_{\Pi} I_{\Pi}(F_{G})$ 
où les sommes portent sur les représentations automorphes cuspidales de $U \times \tlU$ 
et $G_{\rmE} \times \tlG_{\rmE}$ respectivement. De plus, la distribution
$J_{\pi}$ est non-nulle si et seulement si $\calP_{U,\pi} \neq 0$. De l'autre côté  
la distribution $I_{\Pi}$ est liée à la période de Jacquet, Piatetski-Schapiro et Shalika \cite{jpps} 
et si elle est non-nulle on a $L(\Pi, 1/2) \neq 0$. De surcroît, $I_{\Pi}$ fait intervenir les périodes de Flicker-Rallis \cite{flick}, 
qui permettent (au moins conjecturalement) de sélectionner seulement les représentations $\Pi$ qui 
proviennent par changement de base d'un produit des groupes unitaires $U' \times \tlU'$ comme ci-dessus.

Pour décrire les décompositions géométriques de 
$I(F_{G})$ et $J(F_{U})$ notons qu'on a des isomorphismes de quotients suivants: 
\[
\Delta U \bsl U \times \tlU / \Delta U \cong \tlU // U, \quad  (x, \tlx) \mapsto x^{-1}\tlx, \quad 
\Delta G_{\rmE} \bsl G_{\rmE} \times \tlG_{\rmE} / G \times \tlG \cong S_{W}// G, \quad 
(g, \tlg) \mapsto \brg^{-1}  \overline{\tlg} \tlg^{-1} g
\]
où $\Delta U$ (resp. $\Delta G_{\rmE}$) 
est l'image de $U$ (resp. de $G_{\rmE} $) par l'inclusion diagonale $U \hrar U \times \tlU$ 
(resp. $G_{\rmE} \hrar G_{\rmE} \times \tlG_{\rmE}$),
$S_{W} \sbs \tlG_{\rmE}$ est la $\rmF$-variété des $y \in \tlG_{\rmE}$ tels que $y \sigma(y) = 1$, 
 $U$ agit par conjugaison sur $\tlU$ et 
$G$ agit sur $S_{W}$ par restriction de son action adjointe à $\tlG_{\rmE}$.

En utilisant les isomorphismes des quotients ci-dessus 
on associe à $F_{U}$ une fonction $\Phi_{U} \in C_{c}^{\infty}(\tlU(\A))$ et 
à $F_{G}$ une fonction $\Phi_{G} \in C_{c}^{\infty}(S_{W}(\A))$ de sorte qu'on a formellement:
\begin{equation}\label{eq:geomIntU}
J(F_{U}) = \int_{[U]}k_{\Phi_{U}}(x)dx, \quad I(F_{G}) = \int_{[G]}k_{\Phi_{G}}(h)\eta_{\rmE/\rmF}(\det h)dh
\end{equation}
où $k_{\Phi_{U}}(x) = \sum_{\tlgam\in \tlU(\rmF)} \Phi(x^{-1}\tlgam x)$ 
pour $x \in U(\A)$ et $k_{\Phi_{G}}(h) = \sum_{\gamma \in S_{W}(\rmF)}\Phi_{G}(h \gamma h)$ pour $h \in G(\A)$.

De nouveau, les intégrales ci-dessus ne convergent pas en général, mais sous certains conditions locales elles admettent 
des décompositions en somme d'intégrales orbitales relatives indexée par les orbites dites \textit{semi-simples régulières} 
(qui sont fermées de centralisateurs triviaux). 
On peut comparer ces intégrales orbitales par un transfert des fonctions dual de transfert de classes de conjugaison 
semi-simples régulières ce qui mène 
aux identités entre les distributions spectrales sur les deux groupes et permet de déduire la conjecture GGP.
Cette stratégie a été utilisée par Wei Zhang \cite{zhang2}. Zhang démontre une partie substantielle 
de la conjecture globale de Gan-Gross-Prasad pour les groupes unitaires. 
 L'un des ingrédients importants dans sa preuve est le lemme fondamental 
 de Jacquet-Rallis démontré par Yun \cite{yun}. Notons aussi que Zhang \cite{zhang1}, 
 en utilisant toujours les formules des traces décrites ci-dessus, démontre certains cas du raffinement de la 
 conjecture GGP qui précise la valeur 
 centrale $L(\Pi, 1/2)$ du à Ichino-Ikeda \cite{ichinoIkeda} et N. Harris \cite{harris}.
 
Afin d'étendre les résultats de Zhang, il faut des formules des traces valables 
 pour toutes les fonctions lisses à support compact. 
En général, du côté spectral il y a d'autres contributions que celles provenant des représentations 
cuspidales et du côté géométrique il y a des contributions des orbites qui ne sont pas semi-simples régulières. 
Ce sont ces contributions qui rendent les intégrales définissant $J(F_{U})$, $I(F_{G})$, $J(\Phi_{U})$ et $I(\Phi_{G})$
divergentes en général. 
 Dans cet article on 
 donne des formules des traces 
 pour les groupes unitaires et linéaires valables pour toutes les fonctions lisses à support compact 
 qui prennent en compte toutes les contributions des côtés spectraux et géométriques.  
 
  \subsection{Nos résultats}
  
  Notre approche de fait est par un processus de troncature à la Arthur. 
  Soient $\chi'$ une donnée cuspidale de   $\GetlG$ et $\chi$ une donnée cuspidale de $U \times \tlU$. 
  On a alors les $\chi'$ et $\chi$-parties des noyaux $k_{F_{G}}$ et $k_{F_{U}}$ notées 
  respectivement $k_{F_{G}, \chi'}$ et $k_{F_{U}, \chi}$ (voir les paragraphes \ref{par:decompSpectr} et 
  \ref{par:repRegul}). Au début de la section \ref{sec:SpecG} on définit le noyau modifié 
  $k_{F_{G}, \chi'}^{T'}$ où $T'$ est un paramètre de troncature qui appartient à un cône aigu, noté $\all_{\tlzero}^{+}$ 
  (voir le paragraphe \ref{par:prelimstraceSp}),
engendré par les copoids associés à un sous-groupe parabolique minimal de $\tlG_{\rmE}$. 
De même, au début de la section \ref{sec:SpecU} on définit le noyau modifié $k_{F_{U}, \chi}^{T}$ où $T \in \all_{0}^{+}$. 
Nos résultats du côté spectral sont alors:
  
  \btheo 
  \begin{enumerate}
  \item[(cf. théorèmes \ref{thm:noyauCVGG} et \ref{thm:noyauCVGU}).]
   Pour tous $T' \in \all_{\tlzero}^{+}$ et $T \in \all_{0}^{+}$ assez réguliers ainsi que pour tout $\sigma \in \R$ on a:
  \[
\mkern-18mu \mkern-18mu 
\sum_{\chi'} \int\limits_{[G_{\rmE}]}\int\limits_{[G]}\int\limits_{[\tlG]}|k_{F_{G}, \chi'}^{T'}(g,g,h,\tlh)||\rmN_{\rmE/\rmF}\det g|_{\A}^{\sigma} d\tlh dh dg < \infty, 
\quad
\sum_{\chi} \int\limits_{[U]}\int\limits_{[U]}|k_{F_{U},\chi}^{T}(x,x,y,y)|dxdy < \infty,
  \]
  où $\rmN_{\rmE/\rmF} :(\A \otimes_{\rmF} \rmE)^{*} \rar \A^{*}$ c'est la norme, $\det$ c'est le determinant 
  et $|\cdot|_{\A}$ c'est la valeur absolue standard sur les idèles de $\A$. 
  \item[(cf. théorèmes \ref{thm:mainQualitThmG} et \ref{thm:mainQualitThmU}).] Soit $s \in \C \smin\{ -1,1\}$. 
  Les intégrales 
  \[
    \int_{[G_{\rmE}]}\int_{[G]}\int_{[\tlG]}k_{F_{G}, \chi'}^{T'}(g,g,h,\tlh)|\rmN_{\rmE/\rmF}\det g|_{\A}^{s} \eta(h,\tlh) d\tlh dh dg, \quad
  \int_{[U]}\int_{[U]}k_{F_{U},\chi}^{T}(x,x,y,y)dxdy
  \]
  sont des polynômes-exponentielles en leur paramètre de troncature dont les termes purement 
  polynomiaux sont constants, notés ensuite $I_{\chi'}(s, F_{G})$ et $J_{\chi}(F_{U})$ respectivement. 
  \item[(cf. théorèmes \ref{thm:invarianceTheoG} et \ref{thm:invarianceTheoU}).] 
  La distribution $I_{\chi'}(s, \cdot)$ 
  est $|\rmN_{\rmE/\rmF}\det \cdot |_{\A}^{s}$-équivariante pour l'action de $G_{\rmE}(\A)$ à gauche
  et $\eta(\cdot, \cdot)$-équivariante pour l'action de $G(\A) \times \tlG(\A)$ à droite
  et $J_{\chi}(\cdot)$ est invariante pour les actions de $U(\A)$ à gauche et à droite.
  \end{enumerate}
  \etheo

Disons quelques mots sur la preuve de la convergence. 
Dans la section \ref{sec:opsTronc} on introduit
certains opérateurs de troncature et on étudie leur propriétés. Ils jouent le rôle crucial dans les preuves des théorèmes 
\ref{thm:noyauCVGG} et \ref{thm:noyauCVGU}. 
Plus précisément, l'un des opérateurs qu'on utilise c'est l'opérateur de Jacquet-Lapid-Rogawski \cite{jlr} 
qui tronque les fonctions définies sur $G_{\rmE}(\A)$.
Dans le paragraphe \ref{par:opTronqm} on introduit une modification de leur opérateur qui prend en compte 
le centre de $G_{\rmE}$. Finalement, dans le paragraphe \ref{par:opTronqD} on introduit des opérateurs  
qui tronquent les fonctions sur $(G_{\rmE} \times \tlG_{\rmE})(\A)$ et $(U \times \tlU)(\A)$. Leur construction 
est similaire à celle des opérateurs considérés par Ichino et Yamana \cite{ichYamGl, ichYamUn}. 


Décrivons maintenant les côtés géométriques. 
On définit $\calO$ l'ensemble de classes d'équivalence sur $\tlG_{\rmE}(\rmF)$ pour la relation d'équivalence 
définie à partir des invariants géométriques pour l'action par conjugaison de $G_{\rmE}(\rmF)$ sur $\tlG_{\rmE}(\rmF)$ 
(voir le paragraphe \ref{par:invpar}). En particulier, toute orbite semi-simple régulière définit une unique classe dans $\calO$ 
et on obtient des décompositions des ensembles $\tlU(\rmF)$ et $S_{W}(\rmF)$.
Soient $\Phi_{G} \in C_{c}^{\infty}(S_{W}(\A))$ et $\Phi_{U} \in C_{c}^{\infty}(\tlU(\A))$. Pour 
$\ol \in \calO$ on pose $k_{\Phi_{G}, \ol}(x) = \sum_{\gamma \in S_{W}(\rmF) \cap \ol}\Phi_{G}(x^{-1}\gamma x)$ et 
$k_{\Phi_{U}, \ol}(x) = \sum_{\tlgam\in \tlU(\rmF) \cap \ol}\Phi_{U}(x^{-1}\tlgam x)$. 
Les intégrales \eqref{eq:geomIntU} de ces fonctions ne convergent pas en général. Comme dans le cas spectral, 
on définit alors les noyaux modifiés $k_{\Phi_{G}, \ol}^{T'}$ (paragraphe \ref{par:geomGroupesG}) 
et $k_{\Phi_{U}, \ol}^{T}$ (paragraphe \ref{par:geomGroupesU})
où $T' \in \all_{\tlzero}^{+}$ et $T \in \all_{0}^{+}$. On obtient alors:

\btheo[cf. théorèmes \ref{thm:noyauGeomCVGG} et \ref{thm:noyauGeomCVGU}] 
\begin{enumerate}
\item  Pour tous $T' \in \all_{\tlzero}^{+}$ et $T \in \all_{0}^{+}$ assez réguliers ainsi que pour tout $\sigma \in \R$ on a: 
\[
\sum_{\ol \in \calO} \int_{[G]}|k_{\Phi_{G}, \ol}^{T'}(x)||\det x|_{\A}^{\sigma}dx < \infty, \quad 
\sum_{\ol \in \calO} \int_{[U]}|k_{\Phi_{U}, \ol}^{T}(x)|dx < \infty.
\]
\item Soit $s \in \C \smin \{-1,1\}$ et $\ol \in \calO$.
  Les intégrales 
  $\int_{[G]}k_{\Phi_{G}, \ol}^{T'}(x)|\det x|_{\A}^{s} \eta_{\rmE/\rmF}(\det x)dx$ et \\ $\int_{[U]}k_{\Phi_{U}, \ol}^{T}(x)dx$
  sont des polynômes-exponentielles en leur paramètre de troncature dont les termes purement 
  polynomiaux sont constants, notés ensuite $I_{\ol}(s, \Phi_{G})$ et $J_{\ol}(\Phi_{U})$ respectivement. 
  \item La distribution $I_{\ol}(s, \cdot)$ est $\eta$-équivariante pour l'action de $G(\A)$ sur $S_{W}(\A)$ 
  et la distribution $J_{\ol}(\cdot)$ est invariante pour l'action de $U(\A)$ sur $\tlU(\A)$.
\end{enumerate}
\etheo

En fait, dans \cite{leMoi, leMoi2} on démontre ces résultats pour les versions infinitésimales des distributions 
$I_{\ol}(s, \cdot)$ et $J_{\ol}(\cdot)$. Ce qu'on démontre alors dans les sections \ref{sec:RTFG} et \ref{sec:RTFU} 
c'est qu'on peut se ramener aux cas infinitésimaux. L'outil principal est ici l'application de Cayley (voir le paragraphe \ref{par:invpar}) 
qui a été 
introduite dans le contexte de la formule des traces relative de Jacquet-Rallis 
dans \cite{zhang2}, section 3.
Elle permet
de passer des algèbres de Lie vers les groupes et vice-versa 
et elle envoie les classes $\ol \in \calO$ en classes dans $\calO$.
En faisant alors un argument
de descente, détaillé dans le paragraphe \ref{par:Liedescente}, on réduit la preuve du théorème ci-dessus 
aux résultats démontrés dans \cite{leMoi, leMoi2}.

Voici le résultat principal de cet article qui à ce stade est une conséquence directe des théorèmes ci-dessus.

\btheo[Formules des traces relatives de Jacquet-Rallis] 
\begin{enumerate}
\item[]
\item[(cf. théorème \ref{thm:RTFG}).]
Soit $F_{G} \in C_{c}^{\infty}(\GetlGa)$ et soit
$\Phi_{G} \in C_{c}^{\infty}(S_{W}(\A))$ définie par
\begin{equation*}
\Phi_{G}(y) := \int_{G_{\rmE}(\A)} \int_{\tlG(\A)} F_{G}(x, x \tlg \tlh) d\tlh dx \quad 
y = \overline{\tlg} \tlg^{-1} \in S_{W}(\A).
\end{equation*}
On a alors pour tout $s \in \C \smin \{-1,1\}$:
\[
\sum_{\chi'}
I_{\chi'}(s, F_{G}) = \sum_{\ol \in \calO} I_{\ol}(s, \Phi_{G}).
\]
\item[(cf. théorème \ref{thm:RTFU}).]
Soit $F_{U} \in C_{c}^{\infty}((U \times \tlU)(\A))$ et soit
$\Phi_{U} \in C_{c}^{\infty}(\tlU(\A))$ définie par
\begin{equation*}
\Phi_{U}(x) := \int_{U(\A)}F_{U}(y, x y) dy., \quad x \in \tlU(\A).
\end{equation*}
On a alors:
\[
\sum_{\chi}
J_{\chi}(F_{U}) = \sum_{\ol \in \calO} J_{\ol}(\Phi_{U}).
\]
\end{enumerate} 
\etheo

Commentons finalement l'apparition du déterminant à la puissance complexe dans les définitions des distributions 
$I_{\chi}(s, \cdot)$ et $I_{\ol}(s, \cdot)$. On ajoute ce terme en vu de possibles applications 
aux dérivées des intégrales orbitales relatives de Jacquet-Rallis. 
Pour les classes semi-simples régulières, les dérivées des analogues locaux des distributions 
$I_{\ol}(s, \cdot)$ définies sur la variété $S_{W}$ ont d'abord été introduites dans 
 \cite{zhang3} et ensuite étudiées dans \cite{zhangRapTers, zhangRapSmith} 
 pour des fonctions test particulières. 
Ces travaux lient ces dérivées à des nombres d'intersection sur certains espaces de Rapoport-Zink 
et s'inscrivent dans le cadre de la variante arithmétique de la conjecture de Gan-Gross-Prasad.
 Récemment, dans \cite{mihatsch}, 
 les dérivées des analogues locaux des distributions 
$I_{\ol}(s, \cdot)$ (leurs versions infinitésimales étudiées dans \cite{leMoi2}) ont aussi été étudiées. 
Peut-être, que notre formule pourrait avoir un intérêt dans ces questions.

\textbf{Remerciements}. 
Je remercie mon directeur de thèse, Pierre-Henri Chaudouard, pour m'avoir introduit à ce projet de recherche 
et au domaine de formes automorphes. Il m'a guidé depuis mon mémoire de Master 2 et j'ai beaucoup appris
pendant cette période. Je lui suis très reconnaissant pour ses conseils éclairés et la confiance qu'il m'a apporté. 
Ce travail a été partiellement soutenu par l'Institut Universitaire de France et son projet Ferplay ANR-13-BS01-0012.
\section{Prolégomènes}\label{sec:prolegoSpec}

\subsection{Préliminaires pour la formule des traces}\label{par:prelimstraceSp}

Soient $\rmF$ un corps de nombres et 
$G$ un $\rmF$-groupe algébrique 
réductif.
Pour tout $\rmF$-sous-groupe de Levi $M$ 
de $G$ (\cad un facteur de Levi d'un $\rmF$-sous-groupe parabolique 
de $G$) soit $\calF(M)$ l'ensemble de $\rmF$-sous-groupes 
paraboliques de $G$ contenant $M$ et $\calP(M)$ le sous-ensemble 
de $\calF(M)$ composé de sous-groupes paraboliques 
admettant $M$ comme facteur de Levi.
On fixe 
un sous-groupe de Levi minimal $M_{0}$ 
de $G$. 
On appelle les éléments de $\calF(M_{0})$ 
les sous-groupes paraboliques semi-standards de $G$
et les éléments de $\calP(M_{0})$ les sous-groupes paraboliques minimaux.
On utilisera toujours le symbole $P$, avec des indices éventuellement, 
pour noter un sous-groupe parabolique semi-standard.
Pour tout $P \in \calF(M_{0})$ soit $N_{P}$ le radical unipotent 
de $P$ et $M_{P}$ le facteur de Levi de $P$ contenant $M_{0}$. 
On a alors $P = M_{P}N_{P}$. On note $A_{P}$ 
le tore central de $M_{P}$ déployé sur $\rmF$ 
et maximal pour cette propriété.
Pour $P_{1} \in \calF(M_{0})$, quand il n'y aura pas d'ambiguïté, on écrit 
$N_{1}$ au lieu de $N_{P_{1}}$, $M_{1}$ au lieu de $M_{P_{1}}$ etc.

Soit  $P \in \calF(M_{0})$.
On définit le $\R$-espace vectoriel 
$\mathfrak{a}_{P} := \Hom_{\Z}(\Hom_{\rmF}(M_{P}, \Gm), \R)$, 
isomorphe à $\Hom_{\Z}(\Hom_{\rmF}(A_{P}, \Gm), \R)$ grâce à l'inclusion $A_{P} \hrar M_{P}$, 
ainsi que son espace dual $\all_{P}^{*} = \Hom_{\rmF}(M_{P}, \Gm) \otimes_{\Z} \R$ et on pose
 \begin{equation}\label{eq:dPDefSp}
 d_{P} = \dim_{\R} \all_{P}, \quad d_{Q}^{P} = d_{Q} - d_{P}, \ 
 Q \sbs P.
\end{equation}

Si $P_{1}\subseteq P_{2}$, on a un
homomorphisme injectif canonique
$\mathfrak{a}_2^{*} \hookrightarrow \mathfrak{a}_1^{*}$ 
qui donne la projection 
$\mathfrak{a}_1  \twoheadrightarrow \mathfrak{a}_2$, 
dont on note 
$\mathfrak{a}_{1}^{2} =\mathfrak{a}_{P_{2}}^{P_{2}}$ le noyau.
On a aussi l'inclusion 
$\mathfrak{a}_{{2}}\hookrightarrow \mathfrak{a}_{{1}}$, 
qui est une section de
$\mathfrak{a}_1 \twoheadrightarrow \mathfrak{a}_2$, 
grâce à la restriction des caractères de $A_{{1}}$ à
$A_{{2}}$. 
Il s'ensuit  que si 
$P_{1} \subseteq P_{2}$ alors
\begin{equation}\label{eq:decompSp}
\mathfrak{a}_1 = \mathfrak{a}_1^2\oplus \mathfrak{a}_2.
\end{equation}
Conformément à cette décomposition, on pose aussi 
$(\all_{1}^{2})^{*} = 
\{\la \in \all_{1}^{*}| \la(H) = 0 \ \forall H \in \all_{2}\}$. 
On aura besoin aussi de 
$(\all_{1,\C}^{2})^{*} := (\all_{1}^{2})^{*} \otimes_{\R} \C$ 
et de $\all_{1,\C}^{*} = \all_{1}^{*} \otimes_{\R} \C$.

Si $P \sbs Q$ sont des sous-groupes paraboliques semi-standards 
où $P$ est un sous-groupe parabolique minimal on note simplement 
$\all_{0} = \all_{P}$,
$\all_{0}^{Q} = \all_{P}^{Q}$, $\all_{0}^{*} = \all_{P}^{*}$ etc. Cela ne dépend pas du choix 
de $P$. 
En général donc, si $P_{1} \sbs P_{2}$, 
grâce à la décomposition (\ref{eq:decompSp}) ci-dessus, on considère les espaces
$\all_{1}$ et $\all_{1}^{2}$
(resp. $\all_{1}^{*}$ et $(\all_{1}^{2})^{*}$) comme des sous-espaces de $\all_{0}$ (resp. de $\all_{0}^{*}$). 
En particulier, quand on parle de la projection d'un élément de $\all_{0}$ à $\all_{1}^{2}$ (resp. de $\all_{0}^{*}$ à $(\all_{1}^{2})^{*}$)
ce sera toujours par rapport à la somme directe $\all_{0} = \all_{0}^{1} \oplus \all_{2} \oplus \all_{1}^{2}$ 
(resp. $(\all_{0})^{*} = (\all_{0}^{1})^{*} \oplus \all_{2}^{*} \oplus (\all_{1}^{2})^{*}$).

Notons $\Delta_{P}^{G} = \Delta_{P}$ 
l'ensemble de racines 
simples pour l'action de $A_{P}$ sur $N_{P}$. 
Il y a une correspondance bijective entre les sous-groupes paraboliques $P_{2}$
contenant $P_{1}$ et les sous-ensembles 
$\Delta_{1}^{2} = \Delta_{P_{1}}^{P_{2}}$ de 
$\Delta_1 = \Delta_{P_{1}}$. En fait, 
$\Delta_1^2$
est l'ensemble de racines simples 
pour l'action de $A_1$ sur $N_1 \cap M_2$ et 
l'on a
\begin{displaymath}
\mathfrak{a}_2 = \{H \in \mathfrak{a}_1| \al(H) = 0 
\ \forall \al \in \Delta_1^2\}. 
\end{displaymath}
En plus $\Delta_{1}^{2}$ (les restrictions de ses éléments 
à $\all_{1}^{2}$) est une base de $(\all_{1}^{2})^{*}$.

Fixons $P_{1} \subseteq P_{2}$ et soit $B \in \calP(M_{0})$ 
contenu dans $P_{1}$. On a alors l'ensemble 
$\Delta_{B}^{\vee} = \{\al^{\vee} \in \all_{0}| \al  \in \Delta_{B}\}$ 
de coracines simples associées aux racines simples $\Delta_{0}$. 
Soit $(\Delta_{{1}}^{{2}})^{\vee}$ l'ensemble de projections d'éléments de 
$\Delta_{B}^{\vee}$ à $\all_{1}^{2}$ privé de zéro. Cela ne dépend pas du choix de $B$. 
L'ensemble 
$\Delta_{{1}}^{{2}}$ est en bijection canonique avec $(\Delta_{{1}}^{{2}})^{\vee}$, 
la bijection étant: à $\al \in \Delta_1^2$ on associe 
l'unique $\al^{\vee} \in (\Delta_{{1}}^{{2}})^{\vee}$ tel 
que $\al(\al^{\vee}) >0$.  
Notons également $\hDelta_1^2$ et 
$(\hDelta_1^2)^{\vee}$ 
les bases de $(\mathfrak{a}_{{1}}^{{2}})^{*}$ et 
$\mathfrak{a}_{{1}}^{{2}}$ duales \`a 
$(\Delta_{{1}}^{{2}})^{\vee}$ et
$\Delta_{{1}}^{{2}}$ 
 respectivement. 
Si $P_{2} = G$ on note simplement $\Delta_{{1}}, 
\Delta_{{1}}^{\vee}$ etc. 
Dans ce cadre, on fixe aussi la bijection naturelle entre $\hDelta_{1}$ et $\Delta_{1}$ 
qui à un poids $\varpi \in \hDelta_{1}$ associe l'unique racine $\al \in \Delta_{1}$ 
tel que si l'on note $\al_{0} \in \Delta_{0} \smin \Delta_{0}^{1}$ la racine simple qui se projette sur $\al$, 
alors $\varpi(\al_{0}^{\vee}) = 1$.

Soient $P, P_{1}, P_{2} \in \calF(M_{0})$, on note
\begin{displaymath}
\mathfrak{a}_{P}^{+} = 
\{H \in \mathfrak{a}_{P}| \alpha(H) > 0 \
\forall \alpha \in \Delta_{P}\}
\end{displaymath}
et si $P_{1} \subseteq P_{2} $ 
notons $\tau_{1}^{2}$, $\htau_1^2$
les fonction caractéristiques de
\begin{equation*}
\{H \in \mathfrak{a}_0| \alpha(H) > 0 \
\forall \alpha \in \Delta_1^2\}, \quad  
\{H \in \mathfrak{a}_{0}| \varpi(H) > 0 \
\forall \varpi\ \in \hDelta_{{1}}^{{2}}\}
\end{equation*}
respectivement. 
On note $\tau_{P}$ pour $\tau_{P}^{G}$ et 
$\htau_{P}$ pour $\htau_{P}^{G}$. 

Soit $\A = \A_{\rmF}$ l'anneau des adèles de $\rmF$ 
et soit $|\cdot |_{\A}$ la valeur absolue standard sur le groupe 
des idèles $\A^{*}$.
 Pour tout $P \in \calF(M_{0})$, posons 
 $H_{P} : M_{P}(\A) \rightarrow \mathfrak{a}_{P}$
défini comme
\begin{equation*}
\langle H_{P}(m),\chi \rangle = \log (|\chi (m)|_{\A}), \quad \chi \in 
\Hom_{\rmF}(M_{P}, \Gm), \ m \in M_{P}(\A).
\end{equation*}
C'est un homomorphisme continu et surjectif, 
donc si l'on note $M_{P}(\A_{})^{1}$ 
son noyau, on obtient la suite exacte suivante
\begin{displaymath}
1 \rightarrow M_{P}(\A_{})^{1} \rightarrow M_{P}(\A_{}) 
\rightarrow \mathfrak{a}_{P} \rightarrow 0.
\end{displaymath}
Soit $A_{P}^{\infty}$ la composante neutre du groupe des $\R$-points du 
tore déployé et défini sur $\Q$ maximal pour cette propriété dans 
le $\Q$-tore 
$\Res_{\rmF/\Q}A_{P}$. Alors, comme $\rmF \otimes_{\Q} \R$ 
s'injecte dans $\A$, on a 
naturellement $A_{P}^{\infty}\hookrightarrow A_{P}(\A_{}) \hookrightarrow M_{P}(\A_{})$. En plus, 
la restriction de $H_{P}$ \`a $A_{P}^{\infty}$ est 
un isomorphisme donc 
$M_{P}(\A)$ est un produit direct de 
$M_{P}(\A)^{1}$ et $A_{P}^{\infty}$. 
Pour $Q \in \calF(M_{0})$ contenant $P$ on pose $A_{P}^{Q,\infty} = A_{P}^{\infty} \cap M_{Q}(\A)^{1}$. 
L'application $H_{P}$ induit alors un isomorphisme entre $A_{P}^{Q,\infty}$ et $\all_{P}^{Q}$.

Fixons $K$ un sous-groupe compact maximal admissible de 
$G(\A)$ par rapport à $M_{0}$ (voir le paragraphe 1 de \cite{arthur2} pour la définition).
On a donc, que pour tout sous-groupe parabolique semi-standard $P$, 
$K \cap M_{P}(\A)$ est admissible dans $M_{P}(\A)$ par rapport à $M_{0}$
et on obtient aussi la décomposition d'Iwasawa
$G(\A) = P(\A)K = N_{P}(\A)M_{P}(\A) K$ ce 
qui nous permet d'étendre $H_{P}$ à $G(\A)$ en posant 
$H_{P}(x) = H_{P}(m)$ où $x = nmk$ avec 
$m \in M_{P}(\A), n \in N_{P}(\A), k \in K$. Dans ce cas $H_{P}(x)$ 
ne dépend pas du choix de $m$.

On note $\Omega = \Omega^{G}$ 
le groupe de Weyl de $(G,M_{{0}})$. 
Pour tout $s \in \Omega$, on choisit un représentant $w_{s}$ de $s$ 
dans $G(\rmF)$. 
Pour un $\rmF$-sous-groupe $H$ de $G$ et $s \in \Omega$ 
on note $sH$ le $\rmF$-sous-groupe $w_{s}Hw_{s}^{-1}$. On n'utilisera cette notation que si cela ne dépend pas du choix de $w_{s}$.
Le groupe $\Omega$ agit donc ainsi sur $\calF(M_{0})$.
Pour tout $B \in \calP(M_{0})$ l'application 
$\Omega \ni s \mapsto sB \in \calP(M_{0})$ est une bijection 
et si $B' = sB$ on a $N_{B'} = sN_{B}$.
Pour tout $P \in \calF(M_{0})$ 
soit $\Omega^{P}$ le sous-groupe de 
$\Omega$ stabilisant $P$. On a donc 
$\Omega^{P} = \{s \in \Omega| w_{s} \in M_{P}(\rmF)\}$.

Soit $P \in \calF(M_{0})$. 
 Suivant le paragraphe 2 de  \cite{arthur2} on introduit 
la fonction $\Gamma_{P}'$, qu'on utilisera dans
 les sections \ref{sec:SpecG} et \ref{sec:SpecU}, comme étant
\begin{equation}\label{eq:GammaQDefSp}
\Gamma_{P}'(H,X) = 
\sum_{R \sps P}(-1)^{d_{R}^{G}}
\htau_{R}(H-X)\tau_{P}^{R}(H), \quad H,X \in \all_{P}.
\end{equation}
On a l'égalité suivante
\begin{equation}\label{eq:GammaRecurr}
\htau_{P}(H - X) = 
\sum_{R \supseteq P}
(-1)^{d_{R}^{G}}
\htau_{P}^{R}(H)
\Gamma_{R}'(H,X), \quad H,X \in \all_{P}.
\end{equation}

\blem[cf. \cite{arthur2}, lemme 2.1]\label{lem:GammaIsComp}
Pour un $X \in \all_{P}$ fixé, la fonction $\all_{P}^{G} \ni H \mapsto \Gamma_{P}'(H,X)$ 
est une fonction mesurable à support compact dans $\all_{P}^{G}$. 
\elem

Notons finalement, que parfois, pour économiser l'espace,
 on utilisera la notation $[H]$ et $[H]^{1}$ 
 pour noter $H(\rmF) \bsl H(\A)$ et $H(\rmF) \bsl H(\A)^{1}$ respectivement.

 \subsection{Le domaine de Siegel}\label{par:SiegelSp}

Soient $B, P \in \calF(M_{0})$ tels que $P \sps B$ et $B \in \calP(M_{0})$.
Pour un réel négatif $c$
posons
$A_{B}^{\infty}(P,c) = \{a \in A_{B}^{\infty}|
\al(H_{B}(a)) > c, \forall \al \in \Delta_{B}^{P}\}$
et pour un compact $\omega_{B} \sbs M_{0}(\A)^{1}N_{B}(\A)$ notons 
\begin{displaymath}
\mathfrak{S}_{B}^{P}(\omega,c) = 
\{mak \in G(\A)|  
m \in \omega, k \in K, a \in A_{B}^{\infty}(P,c)\}.
\end{displaymath}

Le résultat classique de la théorie de réduction, qu'on peut trouver, par exemple 
dans \cite{godement} est qu'il existe un $c_{0} <0$ et pour tout $B \in \calP(M_{0})$ 
un compact $\omega_{B} 
 \subseteq M_{0}(\A)^{1}N_{B}(\A)$ tels que pour 
tout sous-groupe parabolique semi-standard
$P$ contenant  $B$ l'on a:
\begin{equation}\label{eq:siegelReal}
G(\A) = P(\rmF) \mathfrak{S}_{B}^{P}(\omega_{B},c_{0}).
\end{equation}

Fixons la constante $c_{0}$ comme ci-dessus. 
Pour tout $B \in \calP(M_{0})$ 
on fixe aussi un $\omega_{B}$ comme dans l'équation
\eqref{eq:siegelReal} de façon que si 
$B' \in \calP(M_{0})$ est tel que $sB = B'$ 
on a $\omega_{B'} = w_{s}\omega_{B}w_{s}^{-1}$.
Les définitions de ce paragraphe sont valables en particulier 
pour les sous-groupes de Levi de $G$. On voit donc qu'on peut 
fixer un $\omega_{B}$ de façon que 
pour tout sous-groupe parabolique semi-standard $P$ et tout $B \in \calP(M_{0})$ 
le contenant
le compact $M_{P}(\A) \cap \omega_{B}$ 
ainsi que le sous-groupe $K_{P} := M_{P}(\A) \cap K$ 
jouissent des rôles de $\omega_{B}$ et $K$ ci-dessus 
par rapport au groupe réductif $M_{P}$ et son 
sous-groupe de Borel $B \cap M_{P}$, la constante $c_{0}$ 
restant la même.  

Soient $B \in \calP(M_{0})$, $P \sps B$ et $T \in \all_{0}$. On définit 
$F_{B}^{P}(x,T)$ comme la fonction 
caractéristique de l'ensemble:
\begin{displaymath}
\{x \in G(\A)| \ \exists \delta \in P(\rmF) \
\delta x \in \mathfrak{S}_{B}^{P}(\omega_{B},c_{0})
, \ \varpi(H_{B}(\delta  x)-T) < 0 \
\forall \varpi \in \hDelta_{B}^{P}\}.
\end{displaymath}
Visiblement, la fonction $G(\A) \ni x \mapsto F^{P}_{B}(x,T)$ est $P(\rmF)$-invariante.

Une fois les compacts $\omega_{B}$ et la constante $c_{0}$ choisis, 
on appelle les ensembles $\Sgl_{B}^{P}(\omega_{B}, c_{0})$ les domaines de Siegel et on les notera simplement par
$\Sgl_{B}^{P}$.

\subsection{Les mesures de Haar}\label{par:haarMes}

Soit $P$ un sous-groupe parabolique semi-standard de $G$.
On fixe $dx$ une mesure de Haar sur $G(\A)$, 
ainsi que 
pour tout sous-groupe connexe $V$ de $N_{P}$ 
(resp. toute sous-alg\`ebre $\mathfrak{h}$ de $\nl_{P}$)
l'unique mesure de Haar sur $V(\A)$ (resp. $\mathfrak{h}(\A)$)
pour laquelle le volume de $V(\rmF)\backslash V(\A)$ 
(resp. $\mathfrak{h}(\rmF)\backslash \mathfrak{h}(\A)$)
soit $1$. Choisissons les mesures de Haar sur $\rmK$ et tout ses sous-groupes fermés
normalisées de même façon. 

On fixe aussi une norme 
euclidienne $\|\cdot \|$ sur $\mathfrak{a}_{0}$ 
invariante par le groupe de Weyl
$\Omega$
et sur tout sous-espace de $\all_{0}$ la mesure de Haar compatible 
avec cette norme. 
Pour tout $Q \in \calF(M_{0})$ tel que $Q \sps P$, 
on en déduit
les mesures de Haar sur 
$A_{P}^{Q,\infty}$ et $A_{P}^{\infty}$
via l'isomorphisme 
$H_{P}$. 

Soit $dp$ la mesure de Haar sur $P(\A)$ invariante 
\`a gauche normalisé de façon que $dx = dpdk$ 
(grâce à la décomposition d'Iwasawa).
Notons $\rho_{P}^{G} = \rho_{P}$ 
l'élément de $(\all_{P}^{G})^{*}$ 
tel que $d(\Ad(m)n) = e^{2\rho_{P}(H_{P}(m))}dn$ 
pour $m \in M_{P}(\A)$ et $n \in N_{P}(\A)$.
Il s'ensuit 
qu'il existe une unique 
mesure de Haar $dm$ sur 
$M_{P}(\A)$ telle 
que si l'on écrit 
$p = nm$ 
où $p \in P(\A)$, $n \in N_{P}(A)$ et
$m \in M_{P}(\A)$ 
alors $dp = e^{-2\rho_{P}(H_{P}(m))}dndm$. Les mesures de Haar sur 
$M_{P}(\A)$ et $A_{P}^{\infty}$ induisent alors une unique mesure de Haar sur 
$M_{P}(\A)^{1}$, que l'on fixe, telle que la mesure de Haar sur $M_{P}(\A)$ 
soit le produit de mesures sur $A_{P}^{\infty}$ et sur $M_{P}(\A)^{1}$.

On introduit au passage certaines fonctions utiles dans les paragraphes \ref{par:fonsPolExpSp} 
et \ref{par:fonsPolExpU}.
Soit 
 $v_{P}$ 
 le volume dans $\all_{P}^{G}$ du 
 parallélotope engendré par $(\hDelta_{P})^{\vee}$.
 Suivant le paragraphe 2 de  \cite{arthur2}, posons
\begin{equation}\label{eq:thetaHatDefSp}
\hat \theta_{P}(\mu) = 
v_{P}^{-1}
\prod_{\tlvpi \in \hDelta_{P}^{\vee}} \mu(\tlvpi^{\vee}), 
\quad \mu \in \all_{P,\C}^{*}.
\end{equation}

\subsection{Fonctions lisses à support compact}

Soit $\A_{f}$ l'anneau des adèles finis et $\rmF_{\infty}$ le produit 
des toutes les complétions de $\rmF$ en places archimédiennes 
de sorte que $\A \cong \rmF_{\infty} \times \A_{f}$. 
On pose $C_{c}^{\infty}(G(\A)) = C_{c}^{\infty}(G(\rmF_{\infty})) \otimes_{\C} C_{c}^{\infty}(G(\A_{f}))$ 
et $C_{c}^{\infty}(G(\A)^{1}) = C_{c}^{\infty}(G(\rmF_{\infty}) \cap G(\A)^{1}) \otimes_{\C} C_{c}^{\infty}(G(\A_{f}))$, 
où $C_{c}^{\infty}(G(\rmF_{\infty}))$ et $C_{c}^{\infty}(G(\rmF_{\infty}) \cap G(\A)^{1})$ sont des espaces de fonctions 
lisses à support compact sur les groupes de Lie correspondants et 
$C_{c}^{\infty}(G(\A_{f}))$ c'est l'espace des fonctions localement constantes à support compact à valeurs complexes 
sur $G(\A_{f})$.

\subsection{Hauteurs}\label{par:hauteurs}
Soit $V$ un $\rmF$-espace vectoriel de dimension finie $n \in \N$. Fixons une $\rmF$-base 
$\xi_{1}, \ldots, \xi_{n}$ de $V$. 
Pour une place $v$ de $\rmF$ notons $\rmF_{v}$ le complété de $\rmF$ à $v$. 
Pour tout $x_{v} = \sum_{i=1}^{n}x_{v, i}\xi_{i} \in V \otimes_{\rmF} \rmF_{v}$ 
on définit sa norme par:
\[
|x_{v}|_{v} := \sup_{i}|x_{i,v}|_{v}
\]
où $|\cdot|_{v}$ c'est la valeur absolue standard sur $\rmF_{v}$. 
Pour tout $x = (x_{v})_{v} \in V \otimes_{\rmF} \A$ 
on définit la hauteur de $x$ par:
\[
\|x\| := \prod_{v}|x_{v}|_{v}.
\]

On fixe un plongement $\rho$ de $G$ dans $\Gl_{n}$.
Pour tout $g_{v} \in G(\rmF_{v})$ on définit la norme de $g_{v}$, notée $|g_{v}|_{v}$, 
comme la 
norme de  $(\rho(g_{v}), \ps{t}\rho(g_{v})^{-1})$ vu comme un élément 
de $\gll_{n}(\rmF_{v}) \oplus \gll_{n}(\rmF_{v})$ où $\gll_{n} = \Lie(\Gl_{n})$.
Pour tout $g = \prod_{v} g_{v} \in G(\A)$ on définit la hauteur de $g$, notée $\|g\|$, 
comme $\|g\| := \prod_{v}|g_{v}|_{v}$.

Soit $B \in \calF(M_{0})$, on a alors les propriétés suivantes:
\begin{gather}
\exists c>0, \ \forall g \in G(\A) \ \|g\| > c, \label{eq:haut1} \\
\exists c>0, \ \forall g_{1},g_{2} \in G(\A) \ \|g_{1}g_{2}\| <c \|g_{1}\|\|g_{2}\| , \label{eq:haut2} \\
\forall g \in G(\A) \ \|g^{-1}\| = \|g\|, \label{eq:haut3} \\
\exists c, c''>0, c' \in \R, \ \forall a \in A_{B}^{\infty} \ c\|H_{B}(a)\| \le \log \|a\| + c' \le c''\|H_{B}(a)\|, \label{eq:haut4} \\
\exists c>0,  \ \forall g \in \Sgl_{B}, \gamma \in G(\rmF) \ \|g\| \le c  \|\gamma g\|, \label{eq:haut5} \\
\exists c, t, t'>0, \ \forall a \in A_{G}^{\infty}, g \in G(\A)^{1} \ c\|a\|^{t'}\|g\|^{t} \le \|ag\|, \label{eq:haut6} \\
\exists c>0, \ \forall g \in G(\A) \ \|H_{B}(x)\| \le c(1+ \log \|x\|) \label{eq:haut7}. 
\end{gather}

On note une conséquence immédiate des propriétés (\ref{eq:haut1}) et (\ref{eq:haut6}):
\begin{equation}
\exists c, t >0, \ \forall a \in A_{G}^{\infty}, g \in G(\A)^{1} \ c\|g\|^{t} \le \|ag\|. \label{eq:haut8}
\end{equation}

\subsection{Décomposition spectrale}\label{par:decompSpectr}

Dans ce paragraphe on rappelle les résultats de la théorie spectrale de formes automorphes qu'on 
peut trouver dans \cite{moeWald2}.

Fixons un $P_{0} \in \calP(M_{0})$ et notons 
$\calF(P_{0})$ l'ensemble des sous-groupes paraboliques de $G$ contenant $P_{0}$. 
Pour $Q \in \calF(P_{0})$, on écrira 
$\bX_{Q} = \bX_{Q}^{G} = N_{Q}(\A)M_{Q}(\rmF) \bsl G(\A)$ et 
$\bX_{Q}^{1} = \bX_{Q}^{G,1} = N_{Q}(\A)M_{Q}(\rmF) \bsl G(\A)^{1}$.

Soit $P \in \calF(M_{0}, P_{0})$ et 
$\sigma$ une représentation automorphe cuspidale de $M_{P}(\A)^{1}$. 
On définit $\calH_{P,\sigma}^{G} = \calH_{P,\sigma}$ comme l'espace de fonctions 
$\phi : N_{P}(\A)M_{P}(\rmF)A_{P}^{\infty} \bsl G(\A) \rar \C$ qui vérifient
\[
\forall \ x \in G(\A), \ 
([M_{P}]^{1} \ni m \mapsto \phi(mx)) \in L^{2}_{\sigma}([M_{P}]^{1}), \quad \text{et} \quad
\int_{K}\int_{[M_{P}]^{1}}|\phi(mk)|^{2}dmdk < \infty
\]
où $L^{2}_{\sigma}([M_{P}]^{1}$ c'est la partie $\sigma$-isotypique de $L^{2}([M_{P}])^{1}$.
Notons $\calH_{P,\sigma}^{G, 0} = \calH_{P,\sigma}^{0} \sbs \calH_{P,\sigma}$ 
le sous-espace des fonctions 
lisses en places infinies, localement constantes en 
places finies, $K$ et $\zl$-finies, où $\zl$ c'est le centre de l'algèbre enveloppante 
de la complexification d'algèbre de Lie de $G(\A^{\infty})$.

Soient $V \sbs \calH_{P,\sigma}^{0}$ un sous-espace de dimension finie, 
 $\Psi^{G} : (\all_{P,\C}^{G})^{*} \rar V$ et 
  $\Psi: \all_{P,\C}^{*} \rar V$ 
des fonctions de Paley-Wiener. 
On pose
\begin{equation*}
\begin{split}
\psi^{G}(x) & = \int_{i(\all_{P}^{G})^{*}}
e^{(\la + \rho_{P})(H_{P}(x))}\Psi^{G}(\la,x)d\la, \quad x \in \bX_{P}^{1},\\
\psi(x) & = \int_{i\all_{P}^{*}}
e^{(\la + \rho_{P})(H_{P}(x))}\Psi(\la,x)d\la, \quad 
x \in \bX_{P},
\end{split}
\end{equation*}
où $i(\all_{P}^{G})^{*} \sbs (\all_{P,\C}^{G})^{*} $ (resp. $i\all_{P}^{*} \sbs \all_{P,\C}^{*}$) est 
muni de la mesure de Haar duale à celle sur $\all_{P}^{G}$ (resp. sur $\all_{P}$). 
Pour tout $Q \in \calF(P_{0})$ tel que $Q \sps P$ on note:
\begin{equation*}
\begin{split}
E_{P}^{Q}\psi^{G}(x) & = \sum_{\delta \in P(\rmF) \bsl Q(\rmF)}
\psi^{G}(\delta x), \quad x  \in \bX_{Q}^{1}, \\
E_{P}^{Q}\psi(x) & = \sum_{\delta \in \PQf}
\psi(\delta x), \quad x \in \bX_{Q}.
\end{split}
\end{equation*}
On appelle $E_{P}^{Q}\psi^{G}$ et $E_{P}^{Q}\psi$ des pseudo-séries d'Eisenstein.
Si $Q = G$ on note simplement $E_{P} = E_{P}^{G}$.
On a alors $E_{P}^{Q}\psi^{G} \in L^{2}(\bX_{Q}^{1})$ 
et $E_{P}^{Q}\psi \in L^{2}(\bX_{Q})$.

On appelle donnée cuspidale de $G$ un couple $(M, \sigma)$ 
où $M$ est un $\rmF$-sous-groupe de Levi de $G$ contenant $M_{0}$ 
et $\sigma$ est une représentation cuspidale de $M(\A)^{1}$.
On dit que deux données cuspidales $(M, \sigma)$ et $(M', \sigma')$ 
sont équivalentes s'il existe un $s \in \Omega$ tel que $sM = M'$ et $\sigma' \circ \Ad(w_{s})$ est équivalente 
à $\sigma$ en tant qu'une représentation de $M$.
On note $\calX^{G}$ l'ensemble de classes d'équivalences de données cuspidales de $G$.

Pour tout $\chi \in \calX^{G}$ et $Q \in \calF(P_{0})$
on note $L^{2}_{\chi}(\bX_{Q}^{1})$ (resp. $L^{2}_{\chi}(\bX_{Q})$), le complété $L^{2}$ de l'espace engendré par 
$E_{P}^{Q}\psi^{G}$ (resp. $E_{P}^{Q}\psi$) 
comme ci-dessus, où $P_{0} \sbs P \sbs Q$ 
et $(M_{P},\sigma) \in \chi$. La définition ne dépend pas du choix de $P_{0} \in \calP(M_{0})$ contenant $Q$.
 On a alors des décompositions orthogonales suivantes:
\begin{equation*}
L^{2}(\bX_{Q}^{1}) = \bigobot_{\chi \in \calX^{G}} L^{2}_{\chi}(\bX_{Q}^{1}), \quad
L^{2}(\bX_{Q}) = \bigobot_{\chi \in \calX^{G}} L^{2}_{\chi}(\bX_{Q}).
\end{equation*}

On a alors aussi pour tout $P \sbs Q$ les décompositions analogues 
des espaces $L^{2}(\bX_{M_{Q} \cap P}^{M_{Q},1})$ et $L^{2}(\bX_{M_{Q} \cap P}^{M_{Q}})$ 
par rapport à l'ensemble $\calX^{M_{Q}}$.
Pour tout $\chi \in \calX^{G}$ et tout $P \sbs Q$ on définit l'espace
$L_{\chi}^{2}(\bX_{M_{Q} \cap P}^{M_{Q},1})$ comme la somme directe des espaces
$L_{\chi'}^{2}(\bX_{M_{Q} \cap P}^{M_{Q},1})$ où $\chi'$ parcourt 
la pré-image de $\chi$ par l'application naturelle à fibre finie $\calX^{M_{Q}} \rar \calX^{G}$. 
On définit $L_{\chi}^{2}(\bX_{M_{Q} \cap P}^{M_{Q}})$ de même façon.
On a alors
\[
L^{2}(\bX_{M_{Q} \cap P}^{M_{Q},1}) = \bigobot_{\chi \in \calX^{G}} 
L_{\chi}^{2}(\bX_{M_{Q} \cap P}^{M_{Q},1}), \quad 
L^{2}(\bX_{M_{Q} \cap P}^{M_{Q}}) = \bigobot_{\chi \in \calX^{G}} 
L_{\chi}^{2}(\bX_{M_{Q} \cap P}^{M_{Q}}).
\]

\blem\label{lem:pseudoDescente}
 Soient $\chi \in \calX^{G}$, $P , Q \in \calF(M_{0})$ tels que 
 $P \sbs Q$ et $E_{P}^{Q}\psi \in L^{2}_{\chi}(\bX_{Q})$ une pseudo-série d'Eisenstein. 
 Alors:
\begin{enumerate}[1)]
\item Soient $\la \in \all_{G,\C}^{*}$
et $A^{\infty}$ un sous-groupe de $A_{G}^{\infty}$.
Alors, la fonction
\[
\bX_{Q}^{1} \ni x \mapsto \int_{A^{\infty}}e^{\la(H_{G}(a))}E_{P}^{Q}\psi(ax)da 
\]
est bien définie et appartient à $L^{2}_{\chi}(\bX_{Q}^{1})$.
\item 
Soient 
$K' \sbs K$ un sous-groupe compact 
et $\upphi : K' \rar \C$ une fonction continue.
La fonction
\begin{equation*}
\bX_{M_{Q} \cap P}^{M_{Q}} \ni x \mapsto \int_{K'}E_{P}^{Q}\psi(xk')e^{-\rho_{Q}(H_{Q}(x))}\upphi(k')dk'.
\end{equation*}
est bien définie et appartient
à $L^{2}_{\chi}(\bX_{M_{Q} \cap P}^{M_{Q}})$. 
\end{enumerate}
\elem

\subsection{Représentations régulières}\label{par:repRegul}

On note $R_{Q}^{G,1}$ (resp. $R_{Q}^{G}$) la représentation régulière 
à droite (par multiplication à droite) du
groupe $G(\A)^{1}$  (resp. $G(\A)$) sur 
$L^{2}(\bX_{Q}^{1})$ (resp. $L^{2}(\bX_{Q}^{1})$). 
Les espaces $L^{2}_{\chi}(\bX_{Q}^{1})$ et $L^{2}_{\chi}(\bX_{Q})$ 
introduits dans le paragraphe précédent sont alors invariants. 

La représentation $R_{Q}^{G,1}$ (resp. $R_{Q}^{G}$)
induit une représentation
de l'algèbre $C_{c}^{\infty}(G(\A)^{1})$ 
(resp. $C_{c}^{\infty}(G(\A))$), notée aussi $R_{Q}^{G,1}$ (resp. $R_{Q}^{G}$),  
sur l'espace $L^{2}(\bX_{Q}^{1})$ (resp. $L^{2}(\bX_{Q})$)
Pour tout $f \in C_{c}^{\infty}(G(\A)^{1})$ et $\Phi \in C_{c}^{\infty}(G(\A))$ les opérateurs 
$R_{Q}^{G,1}(f)$ et $R_{Q}^{G}(\Phi)$ sont intégrales, donnés par les noyaux
\[
k_{f,Q}(x,y) = \sum_{\gamma \in M_{Q}(\rmF)}\int_{N_{Q}(\A)}k(x^{-1}\gamma ny)dn, \quad 
k_{\Phi,Q}(x,y) = \sum_{\gamma \in M_{Q}(\rmF)}\int_{N_{Q}(\A)}k(x^{-1}\gamma ny)dn. 
\]
Pour tout $\chi \in \calX^{G}$, soit $\Pi_{Q, \chi}^{G,1}$ (resp. $\Pi_{Q, \chi}^{G}$) la projection orthogonale 
de $L^{2}(\bX_{Q}^{1})$ (resp. $L^{2}(\bX_{Q})$)
 à 
$L_{\chi}^{2}(\bX_{Q}^{1})$ (resp. $L_{\chi}^{2}(\bX_{Q})$). 
L'opérateur $R_{Q}^{G,1}(f)\Pi_{Q,\chi}^{G,1}$ 
(resp. $R_{Q}^{G}(\Phi)\Pi_{Q,\chi}^{G}$) est alors intégrale et 
l'on note $k_{f,Q, \chi} \in L^{2}(\bX_{Q}^{1} \times \bX_{Q}^{1})$ (resp. $k_{\Phi,Q, \chi} \in L^{2}(\bX_{Q} \times \bX_{Q})$) son noyau.
Si $Q = G$ on écrit $k_{f} = k_{f,G}$, $k_{f,\chi} = k_{f,G,\chi}$ etc.
On a donc:
\[
k_{f,Q} = \sum_{\chi \in \calX^{G}}k_{f,Q,\chi}, \quad 
k_{\Phi,Q} = \sum_{\chi \in \calX^{G}}k_{\Phi,Q,\chi}.
\]

\blem\label{lem:kChiDescente}
 Soient $\Phi \in C_{c}^{\infty}(G(\A))$, $P, Q \in \calF(M_{0})$ 
tels que $P \sbs Q$ et $\chi \in \calX^{G}$.
\begin{enumerate}[1)]
\item Soient $\la \in \all_{G,\C}^{*}$, $A^{\infty}$ un sous-groupe de $A_{G}^{\infty}$ et
\[
\bar \Phi(x) = \int_{A^{\infty}}e^{\la(H_{G}(a))}\Phi(ax)da, \quad x \in G(\A)^{1}.
\]
Alors $\bar \Phi \in C_{c}^{\infty}(G(\A)^{1})$ et 
pour tout $x,y \in \bX_{Q}^{1}$ on a
\[
k_{\bar \Phi, Q, \chi} (x,y)= \int_{A^{\infty}}e^{\la(H_{G}(a))}k_{\Phi, Q, \chi} (x,ay)da.
\]
\item Soient $K_{1}, K_{2}$ deux sous-groupes compacts de $K$ 
et $\Psi : K_{1} \times K_{2} \rar \C$ une fonction continue.
Posons 
\[
\Phi_{Q}(x) = 
\int_{K_{1}}\int_{K_{2}}\int_{N_{Q}(\A)}
e^{\rho_{Q}(H_{Q}(x))}
\Phi(k_{1}^{-1}xnk_{2})
\Psi_{1}(k_{1},k_{2})
dndk_{2}dk_{1}, \quad 
x \in M_{Q}(\A).
\]
Alors $\Phi_{Q} \in C_{c}^{\infty}(M_{Q}(\A))$. 
Soient $\{\chi_{Q}\} \in \calX^{M_{Q}}$ qui s'envoient sur $\chi$ par l'application naturelle 
$\calX^{M_{Q}} \rar \calX^{G}$. Alors, pour tout $x,y \in \bX_{P \cap M_{Q}}^{M_{Q}}$ on a
\[
\sum_{\chi_{Q}}k_{\Phi_{Q}, M_{Q} \cap P, \chi_{Q}}(x,y) = 
\int_{K_{1}}\int_{K_{2}}k_{\Phi, P, \chi}(x k_{1}, yk_{2})e^{-\rho_{Q}(xy)}
\Psi_{1}(k_{1},k_{2})
dk_{1}dk_{2}.
\]
\end{enumerate}
\bdem 
On expliquera le point 2), le point 1) étant analogue.
Le résultat est claire sans $\chi$. 
Pour démontrer le résultat avec $\chi$, il suffit de montrer que 
l'opérateur intégrale 
sur $L^{2}(\bX_{P \cap M_{Q}}^{M_{Q}})$ défini par le noyau
$(x,y) \mapsto \int_{K_{1}}\int_{K_{2}}k_{\Phi, P, \chi}(x k_{1}, yk_{2})e^{-\rho_{Q}(xy)}dk_{1}dk_{2}$
agit trivialement sur $L_{\chi'}^{2}(\bX_{P \cap M_{Q}}^{M_{Q}})$ pour tout $\chi' \neq \chi$. 
Pour cela, on décompose $k_{\Phi, P, \chi}(x,y)$ dans la base hilbertienne de 
$L_{\chi}^{2}(\bX_{P})$, que l'on peut supposer être composée de pseudo-séries d'Eisenstein, 
et l'on applique le lemme \ref{lem:pseudoDescente}. La seule chose 
qui est à vérifier c'est le fait qu'on peut invertir l'intégrale avec 
la somme définissante $k_{\Phi, P, \chi}(x,y)$. L'argument que c'est possible est standard 
et repose sur la dite astuce de Selberg qui permet de se ramener à une fonction 
$\Phi$ de type $h \ast h^{*}$ où $h^{*} (x) = \overline{h(x^{-1})}$.
\edem
\elem

\subsection{Quelques majorations}\label{par:resGen}

Soit $\calU$ l'algèbre enveloppante 
de $\Lie(G(\rmF_{\infty}))_{\otimes} \C$. 
Elle agit sur $C^{\infty}(G(\rmF_{\infty})) \otimes_{\C} C_{c}^{\infty}(\A_{f})$ à gauche et à droite. 
Pour $f \in C^{\infty}(G(\rmF_{\infty})) \otimes_{\C} C_{c}^{\infty}(\A_{f})$
on note ces actions par $X \ast f$ et $ f \ast Y$ respectivement, où
$X,Y \in \calU$.
On notera aussi $R(X)f$ pour $f \ast X$ et 
si $f$ est à plusieurs variables on note
$R_{i}(X)f$ pour l'action de $X$ à droite 
sur $f$ par rapport à la $i$-ème variable. 

Fixons désormais une $f \in C_{c}^{\infty}(G(\A))$. 
Soient $X,Y \in \calU$ et $\chi \in \calX^{G}$.
On a alors:
\begin{equation}\label{eq:diffOnKernel}
R_{1}(X)R_{2}(Y)K_{f,\chi}(x,y) = 
K_{X \ast f \ast Y,\chi}(x,y).
\end{equation}

\blem[cf. \cite{arthur3}, Corollaire 4.6, \cite{moeWald2} paragraphe I.2.4]\label{lem:arthursCor46}
Il existe des $c, N_{0},N_{1} >0$ qui ne dépendent que de $G$ 
telles que pour
tout $P \in \calF(M_{0})$ on a:
\begin{gather}
\sum_{\chi \in \calX^{G}}|k_{f,P,\chi}(x,y)| \le 
c \|x\|^{N_{0}}\|y\|^{N_{0}}, \quad 
\forall x,y \in G(\A), \\
\int_{A_{G}^{\infty}}
\sum_{\chi \in \calX^{G}}|k_{f,P, \chi}(x,ay)|da \le 
c \|x\|^{N_{1}}\|y\|^{N_{1}}, \quad 
\forall x,y \in G(\A)^{1}.
\end{gather}
\elem

\blem[cf. \cite{arthur5}, lemme 2.3]\label{lem:artLemme23} 
Soient $n \in \N^{*}$ et $P \in \calF(M_{0})$.
La donnée pour tout $1 \le i \le n$ 
des sous-groupes paraboliques $Q_{i} \sps P$, 
des points $z_{i},x_{i} \in G(\A)$ et des nombres 
complexes $c_{i} \in \C$ tels que:
\[
\sum_{i=1}^{n}c_{i}\int_{[N_{P}]}k_{f,Q_{i}}(z_{i},nmx_{i})dn = 0 
\]
pour tout $m \in M_{P}(\rmF)\bsl M_{P}(\A)^{1}$ on a pour tout 
$\chi \in \calX^{G}$ et tout $m \in M_{P}(\rmF)\bsl M_{P}(\A)^{1}$
\[
\sum_{i=1}^{n}c_{i}\int_{[N_{P}]}k_{f,Q_{i},\chi}(z_{i},nmx_{i})dn = 0.
\]
\elem

\bcor[cf. \cite{arthur5}, pages 100-101]\label{cor:artLemme23Cor}
 Il existe un $T_{f} \in \all_{0}$ qui ne dépend que 
du support de $f$ tel que pour tous $z,x \in G(\A)$, $\chi \in \calX^{G}$ 
et $P \in \calF(M_{0})$ on a pour tout $m \in M_{P}(\rmF)\bsl M_{P}(\A)^{1}$:
\begin{gather*}
k_{f,P,\chi}(z,mx) = \htau_{P}(H_{P}(z)-H_{P}(x)-T_{f})k_{f,P,\chi}(z,mx), \\
k_{f,P,\chi}(mz,x) = \htau_{P}(H_{P}(x)-H_{P}(z)-T_{f})k_{f,P,\chi}(mz,x).
\end{gather*}
\ecor

Soient $P_{1}, P_{2} \in \calF(M_{0})$ tels que $P_{1} \sbs P_{2}$. 
On pose $\sigma_{1}^{2}$ la fonction caractéristique de 
$H \in \all_{1}$ tels que
\begin{equation}\label{eq:sigmaDef}
\al(H) > 0 \ \forall \ \al \in \Delta_{1}^{2}, \quad 
\al(H) \le 0 \ \forall \ \al \in \Delta_{1} \smin \Delta_{1}^{2}, \quad 
\varpi(H) >0 \ \forall \ \varpi \in \hDelta_{2}. 
\end{equation}

La fonction $\sigma_{1}^{2}$ a été introduite dans \cite{arthur3}, paragraphe 6. 
Voici ses deux propriétés qu'on utilisera, dont les preuves se trouvent dans loc. cit. 

\blem\label{lem:corArth62}
\begin{enumerate}[i)]
 Soient $P_{1}, P_{2} \in \calF(M_{0})$ tels que $P_{1} \sbs P_{2}$. 
\item On a
\begin{equation*}\label{eq:sigmaHauTau}
\tau_{1}^{2}\htau_{2} = \sum_{Q \sps P_{2}}\sigma_{1}^{Q}.
\end{equation*}
\item Il existe une constante $c >0$ telle que pour tous $H, X \in \all_{1}^{G}$
tels que $\sigma_{1}^{2}(H-X) = 1$ on a
\[
\|H\| \le c(\|H_{1}^{2}\| + \|X\|)
\]
où $H_{1}^{2}$ c'est la projection de $H$ à $\all_{1}^{2}$.
\item  Soit $\epsilon >0$.
Alors pour tout $T \in \all_{1}$ tel que $\al(T) > \epsilon \|T\|$ 
pour tout $\al \in \Delta_{1}$ et tout 
$H_{1} \in \all_{1}^{G}$ tel que $\sigma_{1}^{2}(H-T)=1$ on a
\[
\|H_{1}\| \le c\|H_{1}^{2}\|
\]
où $H_{1}^{2}$ c'est la projection de $H_{1}$ à $\all_{1}^{2}$ et 
la constante $c$ ne dépend que de $\epsilon$.
\end{enumerate}
\elem

On rassemble finalement quelques résultats, et ses conséquences, 
de la section 5 de \cite{arthur3} 
dans le lemme suivant. Dans loc. cit. on suppose que les variables $x,y$
ci-dessous appartiennent à $G(\A)^{1}$, mais les preuves passent sans changement 
pour $x, y \in G(\A)$, en prenant compte de la propriété (\ref{eq:haut8}) 
de la hauteur.

\blem\label{lem:lemmeArth51}
\begin{enumerate}[i)]
\item  Il existe des constantes positives $N_{1}$ et $c$ 
telles que pour tout couple de sous-groupes paraboliques semi-standards 
$P_{1} \subseteq P_{2}$ de 
$G$, tout $x \in G(\A)$ et tout $T \in \all_{P_{1}}$ l'on a:
\[
\sum_{\delta \in P_{1}(\rmF)\bsl G(\rmF)}
F^{1}(\delta x,T)\sigma_{1}^{2}(H_{1}^{2}(x)-T) \le c 
e^{N_{1}\|T\|}\|x\|^{N_{1}}.
\] 
\item 
Soit $P \in \calF(M_{0})$. Il existe
des constantes positives $c$, $N$ telles que 
pour tout $x \in G(\A)$ et tout $X \in \all_{P}$
\[
\sum_{\delta \in \PGf}\htau_{P}(H(\delta x)-X) 
\le c \dsl \|x\|e^{\|X\|}\rb^{N}.
\]
En particulier, la somme est finie.
\item Soient $P \in \calF(M_{0})$, 
$T \in \all_{0}$ et $N \ge 0$. Ils existent 
des constantes positives $c'$, $N'$ telles que 
pour toute fonction 
$\phi$ sur $P(\rmF)\bsl G(\A)$ 
et tout $x,y \in G(\A)$ la somme:
\[
\sum_{\delta \in \PGf}|\phi(\delta x)|
\htau_{P}(H_{P}(\delta x)-H_{P}(y)-T)
\]
est majorée par
\[
c' \|x\|^{N'}\|y\|^{N'}\sup_{x' \in G(\A)}
\dsl |\phi(x')| \|x'\|^{-N} \rb.
\]
\end{enumerate}
\elem 

\bcor\label{cor:art23REALLY} 
Pour tous sous-groupes paraboliques $P \sbs Q$,  
tout $\chi \in \calX^{G}$ et tout $x,y \in G(\A)$ on a
\[
\int_{[N_{P}]}k_{Q,\chi}(nx,y)dn = \sum_{\gamma \in P(\rmF) \bsl Q(\rmF)}
k_{P,\chi}(x, \gamma y)
\] 
la somme étant finie.
\bdem 
Le résultat est clair sans $\chi$. La somme est finie en vertu du corollaire \ref{cor:artLemme23Cor} 
et du lemme \ref{lem:lemmeArth51} \textit{ii)} ci-dessus. 
Le résultat suit maintenant en appliquant le lemme \ref{lem:artLemme23}.
\edem
\ecor

Soient $P, R \in \calF(M_{0})$ tels que $P \sbs R$.
On invoque l'identité due à Arthur \cite{arthur3}, 
proposition 1.1:
\begin{equation}\label{eq:basicidentity}
\sum_{P \sbs Q \sbs R}
(-1)^{d_{Q}^{R}} = 
\begin{cases}
0  \text{ si } P \neq R, \\
1  \text{ sinon}.
\end{cases}
\end{equation}

\subsection{$\Gl_{n} \hrar \Gl_{n+1}$}\label{par:glnGlnplus1Sp}

Soit $W$ un $\rmF$-espace vectoriel de dimension finie $n+1$ 
et soit $V \sbs W$ un sous-espace de dimension $n$, où $n \in \N$.
Notons $\tlG = \Gl(W)$.
Fixons un vecteur 
$e_{0} \in W \smallsetminus V$ et notons $D_{0}$ la droite 
qu'il engendre. 
On a alors $W = V \oplus D_{0}$ ce qui permet 
d'identifier $G = \Gl(V)$ comme un sous-groupe de 
$\tlG$ stabilisant $V$ et fixant $e_{0}$. 
Choisissons $M_{0}$ un sous-groupe de Levi minimal de $G$ et soit 
$M_{\tlzero}$ l'unique sous-groupe de Levi minimal de $\tlG$ contenant $M_{0}$. 
On note $D_{i} \sbs V$, où $i = 1, \ldots, n$, les droites stabilisées par $M_{0}$.

Les résultats des paragraphes précédentes s'appliquent 
aux groupes $G$ et $\tlG$ et leur Levi's minimales $M_{0}$
et $M_{\tlzero}$. Les objets associées à $\tlG$ seront notés toujours avec un tilde. 
Pour le choix du sous-groupe compact maximal, on fixe des vecteurs non-nuls 
$e_{i} \in D_{i}$ pour tout $i =1, \ldots, n$ ce qui avec le choix du vecteur $e_{0}$ 
défini les isomorphismes $\tlG \cong \Gl_{n+1}$ et $G \cong \Gl_{n}$. 
On pose alors $\tlK = \prod_{v}\tlK_{v}$ où, pour une place fini $v$ de $\rmF$ on note 
$\tlK_{v} = \Gl_{n+1}(\calO_{v})$, où $\calO_{v}$ c'est l'anneau des entiers de la complétion de $\rmF$ en $v$, 
pour une place réelle $v$ on pose $\tlK_{v} = O(n+1)$-le groupe orthogonal anisotrope et 
pour une place complexe on met $\tlK_{v} = U(n+1)$- le groupe unitaire anisotrope. 
On pose aussi $K = \tlK \cap G(\A)$. 
Dans ce cas $\tlK$ et $K$ vérifient les conditions 
du paragraphe \ref{par:prelimstraceSp} par rapport à $M_{\tlzero}$ et 
$M_{0}$ respectivement. 
Les inclusions $G \hrar \tlG$  et $M_{0} \hrar M_{\tlzero}$ induisent l'inclusion 
$\Omega^{G} \hrar \Omega^{\tlG}$.
On choisit aussi des représentants du groupe de Weyl 
$\Omega^{\tlG}$ de $\tlG$ comme les éléments permutants les vecteurs $e_{i}$. 
On a alors pour tout $\tls \in \Omega^{\tlG}$ que $w_{\tls} \in \tlG(\rmF) \cap \tlK$ 
et si $s \in \Omega^{G}$ alors $w_{s} \in G(\rmF) \cap K$.

On identifie
$\all_{0}$ et $\all_{0}^{*}$ avec des sous-espaces de $\all_{\tlzero}$ 
et $\all_{\tlzero}^{*}$ respectivement. En particulier la mesure de Haar et la norme 
euclidienne sur $\all_{0}$ sont celles d'un sous-espace de $\all_{\tlzero}$.

Pour tout $\tlP \in \calF(M_{\tlzero})$ on admet la notation:
\begin{equation*}
P := \tlP \cap G \in \calF(M_{0}).
\end{equation*}
Pour un $\rmF$-espace vectoriel $\calV$ on note $\calV^{*} = \Hom_{\rmF}(\calV, \rmF)$.
On a alors $V^{*} \sbs W^{*}$ grâce à la décomposition $W = V \oplus D_{0}$.
Le groupe $G$ (resp. $\tlG$) agit naturellement sur $V^{*}$ (resp. $W^{*}$) donc 
aussi sur $V \times V^{*}$ (resp. $W \times W^{*}$). 
Pour $\tlP \in \calF(M_{\tlzero})$
on note
$\calV_{\tlP} \sbs V \times V^{*}$ 
le plus grand sous-espace de $V \times V^{*}$ stabilisé par $\tlP$ vu comme un sous-espace 
de $W \times W^{*}$.
On note aussi $Z_{\tlP} \sbs V$ le plus petit sous-espace de $V$ 
tel que $M_{\tlP}$ stabilise $Z_{P} \oplus D_{0} \sbs W$.

On note $H_{\tlP}$ le plus grand sous-groupe de $M_{P}$ agissant trivialement sur $Z_{\tlP}$ et 
$G_{\tlP}$ le plus grand sous-groupe de $M_{P}$ agissant trivialement sur $\calV_{\tlP}$. 
On note aussi $\tlG_{\tlP}$ le plus grand sous-groupe de $M_{\tlP}$ agissant trivialement sur $\calV_{\tlP}$. 
On a alors $M_{P} =H_{\tlP} \times G_{\tlP}$ et $M_{\tlP} = H_{\tlP} \times \tlG_{\tlP}$ de façon que 
l'inclusion $M_{P} \hrar M_{\tlP}$ est l'identité sur $H_{\tlP}$ et induit 
l'inclusion $G_{\tlP} \hrar \tlG_{\tlP}$ analogue à $G \hrar \tlG$ mais 
associée à l'inclusion des espaces $Z_{\tlP} \hrar Z_{\tlP} \oplus D_{0}$.
Dans ce contexte on pose 
$A_{\tlP}^{st,\infty} := A_{H_{\tlP}}^{\infty}$, 
$\all_{\tlP}^{st} := \all_{H_{\tlP}}$. Avec nos identifications $\all_{\tlP}^{st} \sbs \all_{\tlzero}$ 
ce qui détermine les mesures de Haar sur $\all_{\tlP}^{st}$ et sur $A_{\tlP}^{st, \infty}$. 
On a $A_{\tlP}^{st,\infty} \cap (H_{\tlP}(\A)^{1} \times G_{\tlP}(\A)) = 1$ et 
$M_{P}(\A) = A_{\tlP}^{st,\infty} (H_{\tlP}(\A)^{1} \times G_{\tlP}(\A))$. 
On fixe donc l'unique mesure de Haar sur $H_{\tlP}(\A)^{1} \times G_{\tlP}(\A)$ 
de façon que la mesure de Haar sur $M_{P}(\A)$ choisie soit produit de cette mesure et celle sur $A_{\tlP}^{\infty}$. 
Soulignons qu'on a les décomposition suivantes:
\begin{equation}\label{eq:MPArelSt}
M_{P}(\A) = A_{\tlP}^{st,\infty} (H_{\tlP}(\A)^{1} \times G_{\tlP}(\A)), \quad 
M_{\tlP}(\A) = A_{\tlP}^{st,\infty} (H_{\tlP}(\A)^{1} \times \tlG_{\tlP}(\A)).
\end{equation}

L'application naturelle $A_{\tlP}^{st, \infty} \rar A_{\tlP}^{\tlG,\infty}$ est un isomorphisme. 
On note $j_{\tlP}$ son Jacobien
et on note $\iota_{\tlP}^{st} : (\all_{\tlP,\C}^{\tlG})^{*} \rar (\all_{\tlP,\C}^{st})^{*} := \Hom_{\R}(\all_{\tlP}^{st}, \C) = \all_{H_{\tlP},\C}^{*}$ 
l'isomorphisme induit.
Si $\tlP = \tlG$ on met $j_{\tlP} = 1$.
Notons qu'on a l'égalité suivante
\[
H_{\tlP}(a) = H_{P}(a) \quad \forall a \in A_{\tlP}^{st,\infty}.
\]
Ainsi, pour toute fonction $\uphi$ sur $\all_{\tlP}$ qui est $\all_{\tlG}$-invariante, 
et tout $\la \in (\all_{\tlP,\C}^{st})^{*}$
on a:
\begin{equation}\label{eq:HQtlQSp}
\int\limits_{A_{\tlP}^{st, \infty}} e^{\la(H_{P}(a))} \uphi(H_{P}(a))da = j_{\tlP}^{-1}
\int\limits_{A_{\tlP}^{\tlG, \infty}}
e^{\iota_{\tlP}^{st}(\la)(H_{\tlP}(a))}
\uphi(H_{\tlP}(a))da = 
j_{\tlP}^{-1}
\int\limits_{\all_{\tlP}^{\tlG}}
e^{\iota_{\tlP}^{st}(\la)(H)}
\uphi(H)dH
\end{equation}

Fixons une $\rmF$-base de $V$. Cette base, couplée avec le vecteur $e_{0} \in D_{0}$
définit un plongement de $\tlG$ dans $\Gl_{n+1}$ ce qui définit, en vertu 
de la discussion dans le paragraphe \ref{par:hauteurs}, une hauteur sur $\tlG(\A)$ 
que l'on fixe. 
On définit la hauteur sur $G(\A)$ comme la restriction de la hauteur 
sur $\tlG(\A)$. On a dans ce cas:

\blem\label{lem:hautComp}
 Il existe des constantes $c_{1}, c_{2},t_{1}, t_{2} > 0$ telles que 
\[
c_{1}\|x\|^{t_{1}} \ge \|\tlx^{1}\| \ge c_{2}\|x\|^{t_{2}}, \quad \forall \ x \in G(\A),
\]
où $\tlx^{1} \in \tlG(\A)^{1}$ c'est la 
projection de $x$ sur $\tlG(\A)^{1}$ selon la décomposition $\tlG(\A) = \tlG(\A)^{1}A_{\tlG}^{\infty}$.
\bdem 
En utilisera la notation du paragraphe \ref{par:hauteurs}. 
Pour une place $v$ de $\rmF$ on note $|\cdot|_{G, v}$ la 
norme sur $G(\rmF_{v})$ définie à partir de la base de $V$ 
que l'on a fixée pour définir la hauteur sur $\tlG(\A)$. 
On note aussi $\|\cdot\|_{G}$ la hauteur sur $G(\A)$ définie comme produit des $|\cdot|_{G,v}$.

Soit $x \in G(\A)$ et soit 
$\tlx^{1} \in \tlG(\A)^{1}$ la 
projection de $x$ sur $\tlG(\A)^{1}$ selon la décomposition $\tlG(\A) = \tlG(\A)^{1}A_{\tlG}^{\infty}$.
L'inégalité $c_{1}\|x\|^{t_{1}} \ge \|\tlx^{1}\|$ découle de la propriété (\ref{eq:haut8}) de la hauteur 
appliquée au groupe $\tlG$.

Soit $x = x^{1}a$ 
la décomposition de $x \in G(\A)$ selon 
$G(\A) = G(\A)^{1}A_{G}^{\infty}$. 
Pour une place $v$ de $\rmF$ on note $x_{v}^{1}$ et $\tlx^{1}_{v}$ la $v$-composante de $x^{1}$ et $\tlx^{1}$
 respectivement. Pour toute place non-archimédienne $v$ on a donc $x_{v}^{1} = \tlx_{v}^{1}$.
 Pour toute place archimédienne $v$ de $\rmF$ on a par contre
\[
|\tlx_{v}^{1}|_{v} = \max ( |a|_{G, v}^{\frac{n}{n+1}},  
|x_{v}^{1}a^{\frac{1}{n+1}}|_{G,v})
\] 
où $a^{\frac{1}{n+1}}$ c'est l'unique élément de $A_{G}^{\infty}$ tel que 
$(a^{\frac{1}{n+1}})^{n+1} = a$.
Puisque $\max(|b|,|c|) \ge \sqrt{|b||c|}$,
en appliquant la propriété (\ref{eq:haut8}) au groupe $G$, on trouve
\[
\prod_{v}
\max ( |a|_{G, v}^{\frac{n}{n+1}},  
|x_{v}^{1}a^{\frac{1}{n+1}}|_{G,v}) \ge 
c\|a\|_{G}^{t_{1}}\|x^{1}\|_{G}^{t_{2}}
\]
pour certains constantes positives $c,t_{1},t_{2} > 0$. 
En prenant $t_{3} = \min(t_{1},t_{2})$ et en utilisant les 
propriétés (\ref{eq:haut1}) et (\ref{eq:haut2}) de la hauteur 
on trouve que cela est plus grand que $c'\|x\|_{G}^{t_{3}}$ 
pour une constante $c' >0$. En utilisant l'équivalence des hauteurs sur $G(\A)$, ceci 
est plus grand que $c''\|x\|^{t_{4}}$ pour certains $c'', t_{4} >0$, 
ce qu'il fallait démontrer.
\edem
\elem

\subsection{Les groupes unitaires}\label{par:unitGp}

Soit $\rmE$ une extension quadratique de $\rmF$ 
et $\sigma$ le générateur du groupe de Galois $\Gal(\rmE/\rmF)$. 
Posons $W_{\rmE} = W \otimes_{\rmF} \rmE$, $V_{\rmE} = V \otimes_{\rmF} \rmE$ 
et $D_{0,\rmE} = D_{0} \otimes_{\rmF} \rmE$. 
Supposons la donnée d'une forme $\sigma$-hermitienne non-dégénérée $\tlPhi$ sur $W_{\rmE}$
de façon que $\nu_{0} := \tlPhi(e_{0},e_{0}) \neq 0$ et que $V_{\rmE}$ soit orthogonal à $D_{0, \rmE}$.
Notons $\tlU = U(W_{\rmE},\tlPhi)$ et $U = U(V_{\rmE},\tlPhi|_{V_{\rmE}})$ les groupes unitaires associés. 
Ces sont des $\rmF$-groupes 
algébriques et l'on voit $U$ comme un sous-groupe de 
$\tlU$ grâce à l'inclusion 
$V_{\rmE} \hrar W_{\rmE} = V_{\rmE} \oplus D_{0,\rmE}$. Comme dans le cas du groupe linéaire, les objets associés à $\tlU$ 
seront notés toujours avec un tilde. 

Choisissons $M_{0}$ un sous-groupe de Levi minimal de $U$ et notons 
$M_{\tlzero}$ le sous-groupe de Levi de $\tlU$ contenant $M_{0}$, 
minimal pour cette propriété. Le groupe $M_{\tlzero}$ n'est pas forcément minimal, mais il 
est uniquement déterminé par $M_{0}$. 
Pour se placer dans le contexte du paragraphe \ref{par:prelimstraceSp} on choisi aussi 
un sous-groupe de Levi minimal $M_{\tlU, min}$ de $\tlU$ contenu dans $M_{\tlzero}$. 
On fixe un compact maximal admissible
$K$ de $U(\A)$ ainsi que $\tlK$ dans $\tlU(\A)$ 
qui contient $K$.

Pour tout $\tlP \in \calF(M_{\tlzero})$ on note $P = U \cap \tlP$. 
L'application $\calF(M_{\tlzero}) \ni \tlP \mapsto P \in \calF(M_{0})$ est alors une bijection. 
On va décrire son inverse.  Tout $P \in \calF(M_{0})$ est défini comme le stabilisateur 
du drapeau de type
\[
0 = V_{0,\rmE}  \sbn V_{1,\rmE} \sbn \cdots \sbn V_{k,\rmE}
\]
où $k \ge 0$ et pour $0 \le i \le k$ les $V_{i,\rmE}$ sont des sous-$\rmE$-espaces 
de $V_{\rmE}$ isotropes, i.e. $\Phi(v,w) = 0$ pour tout $v,w \in V_{\rmE}$. 
On définit donc $\tlP \in \calF(M_{\tlzero})$ comme le stabilisateur du même drapeau 
que $P$ mais vu dans $W_{\rmE}$. 

On va réaliser le groupe $\Omega^{U}$ comme un sous-groupe de $\Omega^{\tlU}$. 
Soit $Z$ le plus petit sous-espace non-isotrope de $V$ stabilisé par $M_{0}$.
Fixons une base $e_{1}, \ldots, e_{d}$ du plus grand sous-espace isotrope de $V$ stabilisé 
par $M_{0}$. Il existe alors des uniques vecteurs isotropes $f_{1}, \ldots, f_{d} \in V$ tels que 
$\Phi(e_{i}, f_{j}) = \delta_{ij}$ pour $i,j \in \{1, \ldots, d\}$ et tels que $M_{0}$ 
stabilise l'espace engendré par $f_{1}, \ldots, f_{d}$. 
On réalise alors le groupe de Weyl de $\Omega^{U}$ comme un sous-groupe de 
$U(\rmF)$ qui agit trivialement sur $Z$ et qui est 
engendré par les permutations de vecteurs $e_{1}, \ldots, e_{d}$ 
et par les involutions $\sigma_{i}$ tels que $\sigma_{i}(e_{i}) = f_{i}$, 
$\sigma_{i}(f_{i}) = e_{i}$, $\sigma_{i}(f_{j}) = f_{j}$ et $\sigma_{i}(e_{j}) = e_{j}$ 
pour $i = 1,\ldots, d$ et $j \in \{1,\ldots, d\} \smin \{i\}$.
Le groupe de Weyl $\Omega^{\tlU}$ admet alors une construction identique 
dans laquelle on peut utiliser les mêmes vecteurs $e_{i}, f_{j}$ et, peut-être, deux vecteurs additionnels.
Il est clair alors qu'on a $\Omega^{U} \sbs \Omega^{\tlU}$.

Pour tout $\tlP \in \calF(M_{\tlzero})$, l'inclusion naturelle $A_{P} \hrar A_{\tlP}$ est un isomorphisme, 
on identifie alors ces groupes, en utilisant toujours $A_{P}$. Cette identification identifie 
$\all_{0}$ avec $\all_{P_{\tlzero}}$ où $P_{\tlzero} \in \calP(M_{\tlzero})$ quelconque. 
On a dans ce cas
\[
H_{\tlP}(x) = H_{P}(x) \quad \forall \ x \in U(\A).
\]
En plus, comme on explique dans le paragraphe 2.4 de \cite{leMoi}, pour tout $P \in \calF(M_{0})$ on a que $\Delta_{P}$ et $\Delta_{\tlP}$ 
ainsi que $\hDelta_{P}$ et $\hDelta_{\tlP}$ sont égaux à $1/2$-près. 
Il s'ensuit que pour tout $\tlP_{1}, \tlP_{2} \in \calF(M_{\tlzero})$, tout 
$x \in U(\A)$ et tout $T \in \all_{0}$ on a
\begin{equation}\label{eq:tautlistau}
\tau_{1}^{2}(H_{1}(x) - T) = \tau_{\tlone}^{\tltwo}(H_{\tlone}(x) - T), \quad 
\htau_{1}^{2}(H_{1}(x) - T) = \htau_{\tlone}^{\tltwo}(H_{\tlone}(x) - T).
\end{equation}

On fixe un $P_{0} \in \calP(M_{0})$ et on note $P_{\tlzero} = \tlP_{0}$ l'unique élément de $\calF(M_{\tlzero})$ 
dont l'intersection avec $U$ égale $P_{0}$. 
On fixe une chambre positive $\all_{0}^{+} := \all_{P_{0}}^{+}$.  
Soit $P \in \calF(M_{0})$ contenant $P_{0}$.
Pour tout $T \in \all_{0}^{+}$
on note $T_{P}$ la projection 
de $T$ à $\all_{P}$. 
Notons que $\all_{P}^{+} = \all_{\tlP}^{+}$. On écrira donc parfois $T_{\tlP}$ au lieu 
de $T_{P}$ et $\all_{\tlzero}^{+}$ au lieu de $\all_{0}^{+}$. 

Fixons $\tlB \in \calP(M_{\tlU, min})$ contenu dans $P_{\tlzero}$. Soit $P \in \calF(M_{0})$ contenant $P_{0}$.
Il est démontré dans \cite{ichYamUn}, lemme 2.1, qu'on peut, et on va, choisir les domaines de Siegel $\Sgl_{\tlB}^{\tlP}$ de $\tlU(\A)$
de façon qu'on ait
\begin{equation}\label{eq:SiegelInSiegel}
\Sgl_{P_{0}}^{P} \sbs \Sgl_{\tlB}^{\tlP}.
\end{equation}
On pose alors $F^{P}(\cdot,\cdot) := F_{P_{0}}^{P}(\cdot,\cdot)$ 
et $F^{\tlP}(\cdot,\cdot) = F_{\tlB}^{\tlP}(\cdot,\cdot)$ 
où $F_{P_{0}}^{P}$ est définie dans le paragraphe \ref{par:SiegelSp}. 
On a alors le résultat suivant.
\blem\label{lem:FPisFtlP}
Il existe un $T_{+} + \in \all_{0}^{+}$ tel que pour tous $T \in T_{+} + \all_{0}^{+}$, $x \in U(\A)$ et $P \in \calF(M_{0})$ contenant $P_{0}$ 
on a
\[
F^{P}(x,T) = F^{\tlP}(x,T).
\]
\bdem
Soit $\tlT_{+} \in \all_{\tlB}^{+}$ tel que 
le lemme 6.4 de \cite{arthur3} appliqué au groupe $\tlU$ et son sous-groupe parabolique minimal $\tlB$, 
soit vrai pour tout $\tlT \in \tlT_{+} \in \all_{\tlB}^{+}$. 
On fixe $T_{+} \in \all_{0}^{+}$ de même façon par rapport au groupe $U$ et son sous-groupe minimal $P_{0}$, 
en demandant en plus que la projection de 
$\tlT_{+} + \all_{\tlB}^{+}$ à $\all_{0}$ contienne $T_{+} + \all_{0}^{+}$. 

Soient alors $T \in T_{+} + \all_{0}^{+}$, $x \in U(\A)$ et $P \in \calF(M_{0})$ contenant $P_{0}$.
En utilisant les propriétés (\ref{eq:tautlistau}) et (\ref{eq:SiegelInSiegel}) on voit que 
$F^{P}(x,T) \le F^{\tlP}(x,T)$. Supposons que $F^{\tlP}(x,T) = 1$. 
En vertu du lemme 6.4 de \cite{arthur3} on a
\[
\sum_{P_{0} \sbs R \sbs P}\sum_{\delta \in R(\rmF) \bsl P(\rmF)}F^{R}(\delta x, T)\tau_{R}^{P}(H_{R}(\delta x) - T) = 1.
\]
Si $F^{P}(x,T) \neq 1$ on a $F^{R}(\delta_{0} x, T)\tau_{R}^{P}(H_{R}(\delta_{0} x) - T) = 1$ pour certains 
$P_{0} \sbs R \sbs P$ et $\delta_{0} \in R(\rmF) \bsl P(\rmF)$. 
Soit alors $\tlT \in \tlT_{+} \in \all_{\tlB}^{+}$ dont la projection à $\all_{0}$ vaut $T$. 
On a pour tout $\varpi \in \hDelta_{\tlB}$ que $\varpi(\tlT) \ge \varpi(T)$. On voit donc que pour tout $\tlQ \sps \tlB$ 
on a $F_{\tlB}^{\tlQ}(\cdot, T) \le F_{\tlB}^{\tlQ}(\cdot, \tlT)$. 
En particulier, d'après ce qu'on a démontré, on trouve que $F_{\tlB}^{\tlR}(\delta_{0} x, \tlT) = F_{\tlB}^{\tlP}(x,\tlT) = 1$. 
En vertu des égalités (\ref{eq:tautlistau}) on a
$\tau_{\tlR}^{\tlP}(H_{\tlR}(\delta_{0} x) - T) = 1$. Mais pour tout $\al \in \Delta_{\tlP}^{\tlR}$ 
on a $\al(\tlT) = \al(T)$ donc on a $\tau_{\tlR}^{\tlP}(H_{\tlR}(\delta_{0} x) - \tlT) = 1$. 
Soit $\tldlt_{0} \in \tlR(\rmF) \bsl \tlP(\rmF)$ la projection de $\delta_{0}$. 
On a donc $F^{\tlR}(\tldlt_{0} x, \tlT)\tau_{\tlR}^{\tlP}(H_{\tlR}(\tldlt_{0} x) - \tlT) = 1$ 
et $F^{\tlP}(x,\tlT) = 1$. 
On obtient donc un absurde en utilisant de nouveau le lemme 6.4 de loc. cit.
\edem 
\elem
On fixe donc $T_{+}$ comme dans le lemme \ref{lem:FPisFtlP} ci-dessus. 

On voit donc que, tant que $x$ appartient à $U(\A)$, on peut utiliser la notation avec le tilde ou sans interchangeablement 
dans le cas du groupe unitaire.
Au début, pour souligner les similarités avec le cas linéaire on va utiliser celle avec un tilde, 
sinon on va opter pour la notation plus économique sans le tilde. 

\subsection{$\Gl_{n} \hrar \Res_{\rmE/\rmF} \Gl_{n}$}\label{par:GinGE}

Pour tout $\rmE$-groupe algébrique $H$ on note $H_{\rmE} = \Res_{\rmE/\rmF}H$. 

On considère les groupes $G_{\rmE} := \Res_{\rmE/\rmF} (\Gl(V_{\rmE}))$ et $\tlG_{\rmE} := \Res_{\rmE/\rmF} (\Gl(W_{\rmE}))$.
L'inclusion des $\rmF$-espaces 
$V \hrar V_{\rmE} = V \otimes_{\rmE} \rmE$ 
induit l'inclusion $g \mapsto g \otimes \Id_{\rmE}$ de 
$G$ dans $G_{\rmE}$ que l'on fixe. On a de même 
$\tlG \hrar \tlG_{\rmE}$ grâce à $W \hrar W_{\rmE} = W \otimes_{\rmF} \rmE$. 
On voit donc tous les groupes $G$, $G_{\rmE}$ et $\tlG$ comme des $\rmF$-sous-groupes 
de $\tlG_{\rmE}$. 
Soit $M_{0, \rmE}$ le sous-$\rmF$-groupe de Levi minimal de $G_{\rmE}$ 
défini comme le stabilisateur des mêmes droites que $M_{0}$ mais tensorisées par $\rmE$.
On définit $M_{\tlzero, \rmE}$ de même façon par rapport à $M_{\tlzero}$. 
Cela nous place dans le cadre du paragraphe \ref{par:prelimstraceSp} ainsi que \ref{par:glnGlnplus1Sp}.
On fixe le compact maximal $\tlK_{\rmE}$ de $\tlG_{\rmE}(\A) = \tlG(\A \otimes_{\rmF} \rmE)$ 
de même façon qu'on l'a fait dans le paragraphe \ref{par:glnGlnplus1Sp} pour le groupe $\tlG(\A)$, 
par rapport à la même base de $W$. On pose $K_{\rmE} = \tlG_{\rmE}(\A) \cap \tlK_{\rmE}$. 
On a donc $K \sbs K_{\rmE}$ et $\tlK \sbs \tlK_{\rmE}$. On identifie 
les groupes de Weyl $\Omega^{G}$ et $\Omega^{G_{\rmE}}$ 
ainsi que $\Omega^{\tlG}$ et $\Omega^{\tlG_{\rmE}}$ ce qui nous 
fixe les représentants des groupes de Weyl $\Omega^{G}$ et $\Omega^{\tlG}$ 
avec les choix qu'on a fait dans le paragraphe \ref{par:glnGlnplus1Sp}.

Notons que l'application $\calF(M_{0, \rmE}) \ni P_{\rmE}  \mapsto P = P_{\rmE} \cap G \in \calF(M_{0})$ 
est une bijection. Pour un $P \in \calF(M_{0})$ on note alors $P_{\rmE} \in \calF(M_{0, \rmE})$ 
le sous-groupe dont intersection avec $G$ vaut $P$. 
Soit $P \in \calF(M_{0})$. L'inclusion 
$M_{P}(\A) \hrar M_{P_{\rmE}}(\A)$ induit un isomorphisme 
$\all_{P} \irar \all_{P_{\rmE}}$ de façon que le diagramme suivant
\[
 \xymatrix{ M_{P}(\A) \ar@{^{(}->}[r] \ar[d]_{H_P} & M_{P_{\rmE}}(\A) \ar[d]^{H_{P_{\rmE}}} \\
               \all_P \ar[r]^{\sim} & \all_{P_{\rmE}} 
}
\]
soit commutatif. On identifie alors $\all_{P_{\rmE}}$ avec $\all_{P}$ et $\all_{P_{\rmE}}^{*}$ avec $\all_{P}^{*}$.
On a donc pour tout $x \in G(\A) \hrar G_{\rmE}(\A)$
\[
H_{P}(x) = H_{P_{\rmE}}(x).
\]
On utilisera donc la notation $H_{P}$ partout. 
En plus, l'inclusion $A_{P} \hrar A_{P_{\rmE}}$ est un isomorphisme. 
On a donc $\Delta_{P} = \Delta_{P_{\rmE}} $ et $\hDelta_{P} = \hDelta_{P_{\rmE}}$ etc.
Il s'ensuit que pour tous $P_{1}, P_{2} \in \calF(M_{0})$, 
$x \in G(\A)$ et tout $T \in \all_{0}$ on a
\begin{equation}\label{eq:tauPE}
\tau_{P_{1}}^{P_{2}}(H_{P}(x) - T) = \tau_{P_{1,\rmE}}^{P_{2,\rmE}}(H_{P,\rmE}(x) - T).
\end{equation}
De même pour les fonctions $\htau$. On voit donc qu'on peut omettre l'indice $\rmE$ 
dans ces situations. Notons pourtant qu'on a 
\[
2\rho_{P} = \rho_{P, \rmE} \in \all_{P}^{*}.
\]

La discussion ci-dessus s'applique mot par mot dans le cas $\tlG \hrar \tlG_{\rmE}$ 
et un $\tlP \in \calF(M_{\tlzero})$.

On choisit pour tout $B \in \calP(M_{0})$ (resp. $\tlB \in \calP(M_{\tlzero})$) 
des domaines de Siegel $\Sgl_{B}^{G} \sbs G(\A)$ 
et $\Sgl_{B_{\rmE}}^{G_{\rmE}} \sbs G_{\rmE}(\A)$ 
(resp. $\Sgl_{\tlB}^{\tlG} \sbs \tlG(\A)$ et $\Sgl_{\tlB_{\rmE}}^{\tlG_{\rmE}} \sbs \tlG_{\rmE}(\A)$ ), 
comme dans le paragraphe \ref{par:SiegelSp}, 
de façon qu'on ait $\Sgl_{B}^{G} \sbs \Sgl_{B_{\rmE}}^{G_{\rmE}}$ 
(resp. $\Sgl_{\tlB}^{\tlG} \sbs \Sgl_{\tlB_{\rmE}}^{\tlG_{\rmE}}$).

On fixe un $\tlBmin \in \calF(M_{\tlzero})$ 
et on note $\all_{\tlzero}^{+} := \all_{\tlBmin}^{+}$.
Soit $T \in \all_{\tlzero}$, pour tout $\tlQ \in \calF(M_{\tlzero})$ 
on note alors $T_{\tlQ} \in \all_{\tlQ}$ la projection de 
$s^{-1}T $ à $\all_{\tlQ}$, où $s \in \Omega^{\tlG}$ 
est tel que $s^{-1}\tlB \sbs \tlQ$. 
Soit $s \in \Omega^{\tlG}$. 
En utilisant le fait que $w_{s} \in \tlG(\rmF) \cap \tlK$ 
on a alors pour tout $\tlP, \tlQ \in \calF(M_{0})$ tels que $\tlP \sbs \tlQ$
et tout $x \in \tlG(\A)$:
\begin{equation*}\label{eq:weyl2Sp}
\htau _{s\tlP}^{s\tlQ}(H_{s\tlP}(x)-T) =
\htau_{\tlP}^{\tlQ}(H_{\tlP}(w_{s}^{-1}x)-s^{-1}T), \quad 
\tau _{s\tlP}^{s\tlQ}(H_{s\tlP}(x)-T) =
\tau_{\tlP}^{\tlQ}(H_{\tlP}(w_{s}^{-1}x)-s^{-1}T). 
\end{equation*}

Pour tout $\tlP \in \calF(M_{0})$ et $s \in \Omega^{\tlG}$ tel que $s\tlBmin \sbs \tlP$, on pose 
alors 
\[
F^{\tlP}(\tlx,T) = F^{\tlP}_{s\tlBmin}(\tlx,T_{s\tlBmin}),  \tlx \in \tlG(\A), \quad 
F^{\tlP_{\rmE}}(\tlx,T) = F^{\tlP_{\rmE}}_{s\tlB_{0,\rmE}}(\tlx,T_{s\tlBmin}),  \tlx \in \tlG_{\rmE}(\A).
\]
La définition ne dépend pas du choix de $s$.

\blem\label{lem:FPisFtlPGl} 
Il existe un $T_{+} \in \all_{\tlzero}^{+}$ tel que pour tout $\tlx \in \tlG(\A)$, tout $T \in T_{+} + \all_{\tlzero}^{+}$  
et tout $\tlP \in \calF(M_{0})$ on a
\[
F^{\tlP}(\tlx, T) = F^{\tlP_{\rmE}}(\tlx, T). 
\]
\bdem 
La preuve est analogue à celle du lemme \ref{lem:FPisFtlP}. 
\edem
\elem
On fixe $T_{+}$ comme dans le lemme \ref{lem:FPisFtlPGl} ci-dessus. 

Vu les relations (\ref{eq:tauPE}) et le lemme \ref{lem:FPisFtlPGl} ci-dessus 
on écrira simplement $\tau_{P}$ au lieu de $\tau_{P_{\rmE}}$, $F^{\tlP}$ 
au lieu de $F^{\tlP_{\rmE}}$ etc. 

Puisqu'on considère les groupes $G$ et $\tlG$ comme des sous-groupes de $G_{\rmE}$ 
et $\tlG_{\rmE}$ respectivement, plongés tous dans $\tlG_{\rmE}$, pour tout 
$\tlP \in \calF(M_{\tlzero})$ et $P \in \calF(M_{0})$ on a 
$\tlP_{\rmE}(\rmF) = \tlP(\rmE)$, $P_{\rmE}(\rmF) = P(\rmE)$, de même pour les parties de Levi 
et sous-groupes unipotentes. On écrira alors $P(\rmE)$ au lieu de $P_{\rmE}(\rmF)$ etc. 
En ce qui concerne les points adéliques, on garde l'écriture $P_{\rmE}(\A)$ etc.

\section{Opérateurs de troncature}\label{sec:opsTronc}

\subsection{Les sous-groupes paraboliques relativement standards}\label{par:relSt}

Soient $H$ et $\tlH$ des $\rmF$-groupes réductif munis d'inclusion $H \hrar \tlH$.
On pense principalement aux 
$H = U \hrar \tlU = \tlH$,  
$H = G_{\rmE}  \hrar \tlG_{\rmE} = \tlH$ et 
$H = \tlG \hrar \tlG_{\rmE} = \tlH$.
Tout ce qu'on a fixé dans le paragraphe \ref{par:prelimstraceSp} s'applique à $H$. 
Soit $M_{H} := M_{0}$-un sous-groupe de Levi minimal de $H$ et on choisit un sous-groupe de Levi 
$M_{\tlH}$ de $\tlH$ qui contient $M_{H}$ et minimal pour cette propriété 
(dans les cas qui nous intéressent il est uniquement déterminé par $M_{H}$).
On fixe $P_{H} \in \calF(M_{H})$ et on note $\calF(P_{H}) = \{P \in \calF(M_{H}) | P \sps P_{H}\}$ 
ainsi que $\relPh = \{\tlP \in \calF(M_{\tlH}) | P_{H} \sbs \tlP\}$. 
Pour tout $\tlP \in \relPh$ on note $P := \tlP \cap H$, alors 
$P \in \calF(P_{H})$.

\textbf{Le cas $G_{\rmE}  \hrar \tlG_{\rmE}$.}
On prend $M_{G_{\rmE}} = M_{0,\rmE}$ 
et $M_{\tlG_{\rmE}} = M_{\tlzero,\rmE}$. On 
fixe un $B \in \calP(M_{0})$ et on pose 
$P_{G_{\rmE}} = B_{\rmE}$. 
La restriction de scalaires induit une bijection entre 
$\relPb$ et $\relPbe$ et il sera plus commode d'utiliser ce premier ensemble dans le pratique. 
On appelle les éléments de $\relPb$ les \textit{sous-groupes paraboliques relativement standards} de $\tlG$ 
et les éléments de $\calF(B)$ les sous-groupes paraboliques standards de $G$.
On note que c'est le seul cas qu'on considère où il n'y a pas de bijection entre 
$\calF(P_{H})$ et $\relPh$, du coup c'est le seul cas en réalité où on utilisera 
l'ensemble $\relPh$.

\textbf{Le cas $U \hrar \tlU$.}
On prend $M_{U} = M_{0}$, $M_{\tlU} = M_{\tlzero}$ et $P_{U} = P_{0}$ que l'on a déjà fixé dans le paragraphe \ref{par:unitGp}. 
On appelle les éléments de $\calF(P_{0})$ les sous-groupes paraboliques standards de $U$. 
Il résulte de la discussion dans le paragraphe  \ref{par:unitGp} que 
$\relPz \ni \tlP \mapsto P = \tlP \cap U \in \calF(P_{0})$ est une bijection. 
Dans le contexte des groupes unitaires on appellera parfois aussi les éléments de $\relPz$ 
\textit{les sous-groupes paraboliques relativement standards} de $\tlU$.
 
\textbf{Le cas $\tlG \hrar \tlG_{\rmE}$.}
On prend $M_{\tlG} = M_{\tlzero}$ et $M_{\tlG_{\rmE}} = M_{\tlzero,\rmE}$. 
Pour tout choix d'un sous-groupe de Borel $\tlB \in \calF(M_{\tlzero})$ 
on a alors $\calF(M_{\tlzero,\rmE}, \tlB) = \calF(\tlB_{\rmE})$ qui est en bijection 
avec $\calF(\tlB)$.
On ne fixe pas le groupe $\tlB$ qu'on lassera varier.

On utilisera alors les groupes $H$ et $\tlH$ 
dans ces trois contextes dans cette section.
Notons qu'on a alors la généralisation du lemme 6.4 de \cite{arthur3}
suivante:
\blem\label{lem:hardNarashimArthur}
Pour tous $\tlQ \in \relPh$, $x \in H(\A)$ et $T \in T_{+} + \alltz^{+}$ 
on a
\[
\sum_{\relPh \ni \tlP \sbs \tlQ}
\sum_{\delta \in P(\rmF) \bsl Q(\rmF)}
F^{\tlP}(\delta x,T_{\tlP})\tau_{\tlP}^{\tlQ}(H_{\tlP}(\delta x) - T_{\tlP}) = 1.
\]
\bdem 
Dans les cas $U \hrar \tlU$ et $\tlG \hrar \tlG_{\rmE}$ c'est une conséquence des lemmes
\ref{lem:FPisFtlP} et \ref{lem:FPisFtlPGl} respectivement. Dans 
le cas $G_{\rmE}  \hrar \tlG_{\rmE}$
il nous faut supposer que pour tout sous-groupe de Borel relativement standard $\tlB_{\rmE}$ de $\tlG$
les compacts $\omega_{\tlB_{\rmE}}$ vérifient $\omega_{B_{\rmE}} \sbs \omega_{\tlB_{\rmE}}$ 
(ce qui est possible car ils contiennent $B$). Le lemme est démontré alors dans \cite{ichYamGl}, lemme 2.3, 
et traduit dans ce langage dans \cite{leMoi2}, proposition 1.3. 
\edem
\elem

Pour $\tlQ \in \calF(M_{\tlH}, P_{H})$ et $\phi \in \calB_{loc}(\tlQ(\rmF) \bsl \tlH(\A))$ 
on notera
\[
\phi_{\tlQ}(x) = \int_{[N_{\tlQ}]}\phi(nx)dn, \quad x \in N_{\tlQ}(\A)M_{\tlQ}(\rmF) \bsl \tlH(\A).
\]

\subsection{Opérateur de troncature de Jacquet-Lapid-Rogawski}

Fixons $\tlB \in \calP(M_{\tlzero})$ un sous-groupe de Borel de $\tlG$. 
Soit $\tlQ \in \calF(\tlB)$ et $\phi$ une fonction localement intégrable sur 
$\tlQ(\rmE) \bsl \tlG_{\rmE}(\A)$. 
On note pour $T \in \all_{\tlzero}$
\[
\La_{JLR}^{\tlQ,T}\phi(x) = \sum_{\tlB \sbs \tlP \sbs \tlQ}(-1)^{d_{\tlP}^{\tlQ}}
\sum_{\delta \in P(\rmF) \bsl Q(\rmF)}\htau_{\tlP}^{\tlQ}(H_{\tlP}(\delta x) - T_{\tlP})
\phi_{\tlP_{\rmE}}(\delta x), \quad x \in  N_{\tlQ}(\A)M_{\tlQ}(\rmF)  \bsl \tlG(\A).
\]
L'opérateur est bien défini car les sommes $\sum_{\delta \in P(\rmF) \bsl Q(\rmF)}$ 
portent sur des ensembles finis en vertu du lemme \ref{lem:corArth62} \textit{ii)}. Il dépend du choix du sous-groupe 
de Borel $\tlB$ qu'on précisera chaque fois.
La notation $\La_{JLR}$ vient du fait que cet opérateur a été introduit dans \cite{jlr}. 
Dans loc. cit. on trouve les démonstrations des assertions présentées dans ce paragraphe. 

\blem\label{lem:formDInvtl}
 Soit $\tlQ \in \calF(\tlB)$. Alors pour tout $x \in \tlG(\A)$ et tout $T \in \all_{\tlzero}$ on a 
\[
\phi_{\tlQ_{\rmE}}(x) = \sum_{\tlB \sbs \tlP \sbs \tlQ}\sum_{\delta \in P(\rmF) \bsl Q(\rmF)}
\tau_{\tlP}^{\tlQ}(H_{\tlP}(\delta x) - T_{\tlP})\La_{JLR}^{\tlP,T}\phi(\delta x).
\]
\elem


On a aussi la propriété de décroissance suivante:
\brop\label{prop:opTronqReltl}
 Soient $\tlQ \in \calF(\tlB)$, $r_{1}, r_{2} \ge 0$ et $\tlK_{0} \sbs M_{\tlQ_{\rmE}}(\A_{f})$ 
un compact ouvert. Il existe alors un sous-ensemble fini $\{X\}$ de l'algèbre enveloppante 
de la complexification de l'algèbre de Lie de $M_{\tlQ_{\rmE}}(\A)^{1} \cap M_{\tlQ_{\rmE}}(\rmF_{\infty})$ 
tel que pour tout espace mesurable $(\Omega, d\omega)$, 
tout $\phi : \Omega  \rar  C^{\infty}(\tlQ(\rmE) \bsl \tlG_{\rmE}(\A) / \tlK_{0})$ mesurable, 
tout $T \in T_{+} + \all_{\tlzero}^{+}$ et
tout $x \in \Sgl^{M_{\tlQ}}_{\tlB \cap M_{\tlQ}} \cap M_{\tlQ}(\A)^{1}$ 
on a que l'expression
\[
\int_{\Omega}|\La_{JLR}^{\tlQ,T}\phi(\omega, x)|d\omega 
\]
où on applique l'opérateur à la variable $x$, est majorée par
 \[
\|x \| ^{-r_{1}} \sum_{X}\sup_{\tly \in M_{\tlQ_{\rmE}}(\A)^{1}}
\dsl 
\|\tly \|^{-r_{2}}\int_{\Omega}|
\int_{[N_{\tlQ_{\rmE}}]}
R(X)\phi(\omega, \tln \tly) d\tln |d \omega
\rb.
 \]
\erop
%

\subsection{Opérateur de troncature relatif sur $G$}\label{par:opTronqm}

Dans ce paragraphe on introduit un opérateur de troncature 
qui tronque les fonctions sur $G_{\rmE}(\A)$ en prenant en compte le centre de $G_{\rmE}$.
Pour $\tlQ \in \relPb$, une fonction localement intégrable $\phi$ sur 
$Q(\rmE) \bsl G_{\rmE}(\A)$
et  $T \in \all_{\tlzero}$ on pose
\[
\Laz^{T, \tlQ}\phi(x) = 
\sum_{B \sbs \tlP \sbs \tlQ}(-1)^{d_{\tlP}^{\tlQ}}
\sum_{\delta \in P(\rmF) \bsl Q(\rmF)}\htau_{\tlP}^{\tlQ}(H_{\tlP}(\delta x) - T_{\tlP})
\phi_{P_{\rmE}}(\delta x), \quad x \in N_{Q}(\A)M_{Q}(\rmF) \bsl G(\A)
\]
où la somme porte sur les sous-groupes paraboliques de $\tlG$ qui contient $B$ 
(i.e. sous-groupes paraboliques relativement standards) qui sont contenus dans $\tlQ$.
Tout comme l'opérateur $\La_{JLR}$, l'opérateur 
$\Laz$ est bien défini et vérifie la formule d'inversion suivante:
\blem\label{lem:formDInv}
 Soit $\tlQ \in \relPb$. Alors pour tout $x \in G(\A)$ et tout $T \in \all_{\tlzero}$ on a 
\[
\phi_{Q_{\rmE}}(x) = \sum_{B \sbs \tlP \sbs \tlQ}\sum_{\delta \in P(\rmF) \bsl Q(\rmF)}
\tau_{\tlP}^{\tlQ}(H_{\tlP}(\delta x) - T_{\tlP})\Laz^{\tlP,T}\phi(\delta x).
\]
\bdem
En utilisant la définition de l'opérateur $\Laz$ on a que la somme ci-dessus égale:
\begin{multline*}
\sum_{B \sbs \tlP \sbs \tlQ}
\sum_{\delta \in P(\rmF) \bsl Q(\rmF)}
\tau_{\tlP}^{\tlQ}(H_{\tlP}(\delta \tlx) - T_{\tlP}) \\
\dsl
\sum_{B \sbs \tlR \sbs \tlP}(-1)^{d_{\tlR}^{\tlP}}
\sum_{\eta \in R(\rmF) \bsl P(\rmF)}\htau_{\tlR}^{\tlP}(H_{\tlR}( \eta \delta \tlx) - T_{\tlR})
\phi_{R_{\rmE}}(\eta \delta x)
\rb \\ =
\sum_{B \sbs \tlR \sbs \tlQ}\sum_{\delta \in R(\rmF) \bsl Q(\rmF)}
\phi_{R_{\rmE}}(\delta x) 
\dsl 
\sum_{\tlR \sbs \tlP \sbs \tlQ} 
(-1)^{d_{\tlR}^{\tlP}}
\htau_{\tlR}^{\tlP}(H_{\tlR}(\delta \tlx) - T_{\tlR})
\tau_{\tlP}^{\tlQ}(H_{\tlP}(\delta \tlx) - T_{\tlP})
\rb.
\end{multline*}
En vertu du lemme combinatoire de Langlands, la somme entre les parenthèses vaut $0$ si $\tlR \neq \tlQ$ 
et $1$ si $\tlR = \tlQ$ d'où le résultat. 
\edem 
\elem

On va démontrer l'analogue de la proposition \ref{prop:opTronqReltl} pour l'opérateur 
$\Laz^{T, \tlQ}$ suivant:
\brop\label{prop:opTronqRel}
 Soient $\tlP \in \relPb$, $r_{1}, r_{2} \ge 0$ et $K_{0} \sbs M_{P_{\rmE}}(\A_{f})$ 
un compact ouvert. Il existe alors un sous-ensemble fini $\{X\}$ de l'algèbre enveloppante 
de la complexification de l'algèbre de Lie de 
$(H_{\tlP_{\rmE}}(\rmF_{\infty}) \cap H_{\tlP_{\rmE}}(\A)^{1}) \times G_{\tlP_{\rmE}}(\rmF_{\infty})$ 
(voir (\ref{eq:MPArelSt}))
tel que pour tout espace mesurable $(\Omega, d\omega)$, 
tout $\phi : \Omega \rar C^{\infty}(P(\rmE) \bsl G_{\rmE}(\A) / K_{0})$ mesurable, 
tout $T \in T_{+} + \all_{\tlzero}^{+}$ et
tout $x \in \Sgl_{B \cap M_{P}}^{M_{P}} \cap (H_{\tlP}(\A)^{1} \times G_{\tlP}(\A))$ 
on a que l'expression
\[
\int_{\Omega}|\Laz^{\tlP,T}\phi(\omega, x)|d\omega 
\]
où on applique l'opérateur $\Laz^{\tlP,T}$ à la variable $x$,
 est majorée par
 \[
\|x \| ^{-r_{1}} \sum_{X}
\sup_{y \in H_{\tlP_{\rmE}}(\A)^{1} \times G_{\tlP_{\rmE}}(\A)}
\dsl 
\|y \|^{-r_{2}}\int_{\Omega}|
\int_{[N_{P_{\rmE}}]}
R(X)\phi(\omega, ny )dn|d \omega
\rb.
 \]
\bdem 

Il suffit de démontrer le cas $\tlP = \tlG$, le cas général étant la composition de ce cas 
et de cas de l'opérateur $\Lajlr$.

On aura besoin quelques résultats préparatoires avant.
Soient $D_{i}$, où $i = 1, \ldots, n$, les droites dans $V$ stabilisés par $M_{0}$ de façon que 
 $V_{i} := \oplus_{j=1}^{i} D_{i}$, où $i = 0, 1,\ldots, n$, soit stabilisé par $B$.
Pour $i = 0, 1, \ldots,n $ soit $e_{i}^{*} \in \all_{\tlzero}^{*}$ 
le caractère par lequel $A_{\tlzero}$ agit sur $D_{i}$. 
Posons aussi $e_{j}^{\vee} \in \all_{\tlzero}$ les 
éléments tels que $e_{i}^{*}(e_{j}^{\vee}) = \delta_{ij}$ où 
$i,j = 0,1, \ldots, n$.
On pose pour $i=1,\ldots, n$
\[
\tlvpi_{i}^{-} = \dfrac{n+1-i}{n+1}(\sum_{j=1}^{i}e_{j}^{*}) - 
\dfrac{i}{n+1}(e_{0}^{*} + \sum_{j=i+1}^{n} e_{i}^{*}), \quad 
\tlvpi_{i}^{+} = \dfrac{n+1-i}{n+1}(e_{0}^{*} + \sum_{j=1}^{i-1}e_{j}^{*}) - 
\dfrac{i}{n+1}(\sum_{j=i}^{n} e_{i}^{*}).
\]
Alors $\tlvpi_{i}^{-}, \tlvpi_{i}^{+} \in (\all_{\tlzero}^{\tlG})^{*}$. 
On pose aussi $\tlvpi_{l}^{-} = \tlvpi_{l}^{+} = 0$ pour 
$l \nin \{1, \ldots, n\}$.

Fixons un sous-groupe parabolique relativement standard $\tlQ$ 
de $\tlG$. Alors $\tlQ$ est le stabilisateur du drapeau de type
\begin{equation*}
0 = V_{i_{0}} \subsetneq \cdots \subsetneq  V_{i_{k-1}}
\subsetneq 
V_{i_{k}} \oplus D_{0} \subsetneq \cdots \subsetneq V_{i_{l}} 
\oplus D_{0} =W
\end{equation*}
où $0 = i_{0} < i_{1} < \cdots < i_{k-1} \le i_{k} < i_{k+1} < \cdots < i_{l} = n$. 
Alors
\begin{equation*}
\hDelta_{\tlQ} = \{\tlvpi_{i_{a}}^{-}, \, \tlvpi_{i_{b}+1}^{+} | \, 
1 \le a \le k-1, \ k \le b \le l-1 \}.
\end{equation*}
Posons 
$\tlvpi_{\tlQ}^{-} := \tlvpi^{-}_{m_{1}}$ et 
$\tlvpi_{\tlQ}^{+} := \tlvpi^{+}_{m}$ où 
$m_{1} = \max(\{ j | \tlvpi_{j}^{-} \in \hDelta_{\tlQ} \} \cup \{0\})$
et
$m_{2} = \min(\{j | \tlvpi_{j}^{+} \in \hDelta_{\tlQ} \} \cup \{0\})$. 
Si $\tlvpi_{\tlQ}^{-} \neq 0$
notons $\tlal_{\tlQ}^{-} \in \Delta_{\tlQ}$ 
la racine simple associée 
à $\tlvpi_{\tlQ}^{-}$
comme dans le paragraphe \ref{par:prelimstraceSp}. 
Si $\tlvpi_{\tlQ}^{-} = 0$ on pose $\tlal_{\tlQ}^{-} = 0$. 
On défini $\\tlal_{\tlQ}^{+} \in \Delta_{\tlQ} \cup \{0\}$ de même façon.  

On pose aussi pour $i=1,\ldots, n-1$
\[
\varpi_{i} = \dfrac{n-i}{n}(\sum_{j=1}^{i}e_{j}^{*}) - 
\dfrac{i}{n}(\sum_{j=i+1}^{n} e_{i}^{*}).
\]
Alors $\hDelta_{B} = \{\varpi_{i} | i = 1,\ldots, n-1\}$. 
Pour $l \nin \{1, \ldots, n-1\}$ on pose $\varpi_{l} = 0$.

On
définit $\hiota : \{\tlvpi_{i}^{-}, \tlvpi_{j}^{+} | i,j \in \{1,\ldots,n \}\} \rar \hDelta_{B} \cup \{0\}$ 
par $\hiota(\tlvpi_{i}^{-}) = \varpi_{i}$ et $\hiota(\tlvpi_{j}^{+}) = \varpi_{j-1}$ 
et on note $\hiota_{\tlQ}$ sa restriction à $\hDelta_{\tlQ}$. 
On a alors $\hDelta_{Q} = \hiota_{\tlQ}(\hDelta_{\tlQ}) \smin \{0\}$. 
On définit alors $\iota_{\tlQ} : \Delta_{\tlQ} \rar \Delta_{Q} \cup \{0\}$ 
de façon suivante. Soient $\tlal \in \Delta_{\tlQ}$ et 
$\tlvpi_{\tlal} \in \hDelta_{\tlQ}$ le poids associé à $\tlal$ comme dans le paragraphe \ref{par:prelimstraceSp}. 
Si $\hiota_{\tlQ}(\tlvpi_{\tlal}) = 0$ on pose $\iota_{\tlQ}(\tlal) = 0$, sinon, on 
a $\hiota_{\tlQ}(\tlvpi_{\tlal}) = \varpi_{\al}$ où $\al \in \Delta_{Q}$ et on pose 
$\iota_{\tlQ}(\tlal) = \al$. On vérifie facilement les assertions suivantes:
\begin{enumerate}
\item[$\tlQ1)$] Si $\al \in \Delta_{Q}$ est tel que $\iota_{\tlQ}^{-1}(\al) = \{\tlal\}$, 
alors pour tout $H \in \all_{\tlQ}^{st}$ on a $\al(H) = \tlal(H)$.
\item[$\tlQ2)$]  Il existe au plus un $\al \in \Delta_{Q}$ tel que 
$\# \iota_{\tlQ}^{-1}(\al) > 1$. Cela arrive si et seulement si 
on a $0< i_{k-1} = i_{k} < n$ auquel cas il existe un unique $\al \in \Delta_{Q}$, noté $\al_{\tlQ}$,
tel que 
$\iota_{\tlQ}^{-1}(\al_{\tlQ}) = \{\tlal_{\tlQ}^{-}, \tlal_{\tlQ}^{+}\}$. 
Dans ce cas $\al_{\tlQ}(H) = (\tlal_{\tlQ}^{-} + \tlal_{\tlQ}^{+})(H)$ 
pour tout $H \in \all_{\tlQ}^{st}$.
\item[$\tlQ3)$]  $0 \in \iota_{\tlQ}(\Delta_{\tlQ})$ si et seulement si $i_{k-1} = i_{k}$ 
et $i_{k} \in \{0,n\}$. Dans ce cas, si $i_{k} = 0$ on a $\tlal_{\tlQ}^{+} \in \Delta_{\tlQ}$, $\tlal_{\tlQ}^{-} = 0$ 
et $\iota_{\tlQ}(\tlal_{\tlQ}^{+}) = 0$; si $i_{k} = n$ on a $\tlal_{\tlQ}^{-} \in \Delta_{\tlQ}$, $\tlal_{\tlQ}^{+} = 0$ 
et $\iota_{\tlQ}(\tlal_{\tlQ}^{-}) = 0$.
\end{enumerate}

On dit que $\tlR \sps \tlQ$ est admissible par rapport à $\tlQ$ s'il vérifie 
$\Delta_{\tlQ}^{\tlR} = \iota_{\tlQ}^{-1}(\Delta_{Q}^{R})$.
On a donc
\begin{itemize}
\item $\tlG$ est admissible par rapport à $\tlQ$ si et seulement si $0 \nin \iota_{\tlQ}(\Delta_{\tlQ})$.
\item Si $i_{k-1} \neq i_{k}$ tout $\tlR \sps \tlQ$ est admissible par rapport à $\tlQ$.
\item $\tlR$ n'est pas admissible par rapport à $\tlQ$ si et seulement si 
$\{\tlal_{\tlQ}^{-}, \tlal_{\tlQ}^{+} \} \cap \Delta_{\tlQ}^{\tlR}$ est de cardinalité $1$. 
\end{itemize}

\blem\label{lem:admissiblePars} Soient $\tlR \sps \tlQ$ et pour tout $R \sps S \sps Q$  
soit $c_{S} \in \C$ un scalaire.
\begin{enumerate}[1)]
\item Supposons que $\tlR$ est admissible par rapport à $\tlQ$. 
Dans ce cas
\begin{enumerate}[i)]
 \item $|\sum_{\tlQ \sbs \tlS \sbs \tlR} (-1)^{d_{\tlS}^{\tlG}}c_{S}| = |\sum_{Q \sbs S \sbs R} (-1)^{d_{S}^{G}}c_{S}|$.
 \item Pour tout $H \in \all_{\tlQ}^{st}$ on a $\tau_{Q}^{R}(H) \ge \tau_{\tlQ}^{\tlR}(H)$.
 \item Pour tout $H \in \all_{\tlQ}^{st}$ on a $\sum_{\tlal \in \Delta_{\tlQ}^{\tlR}} \tlal(H) = \sum_{\al \in \Delta_{Q}^{R}} \al(H)$.
 \end{enumerate}
 \item Si $\tlR$ n'est pas admissible par rapport à $\tlQ$ on a 
 $\sum_{\tlQ \sbs \tlS \sbs \tlR} (-1)^{d_{\tlS}^{\tlG}}c_{S} =0$.
\end{enumerate}
\bdem
Démontrons le point 1). 
Supposons que $\tlR$ est admissible par rapport à $\tlQ$. 
Supposons d'abord que pour tout $\al \in \Delta_{Q}$ on a $\# \iota_{\tlQ}^{-1}(\al) = 1$.
Dans ce cas pour tout $Q \sbs S \sbs R$ il existe un unique $\tlS$ entre $\tlQ$ et $\tlR$ 
tel que $\tlS \cap G = S$. En plus, pour tous $\tlS_{1}, \tlS_{2}$ 
entre $\tlR$ et $\tlQ$ 
on a $d_{\tlS_{1}} - d_{S_{1}} = d_{\tlS_{2}} - d_{S_{2}}$ et le point \textit{i)} suit. 
Les points \textit{ii)} et \textit{iii)} découlent facilement du fait que 
tout $\al \in \Delta_{Q}^{R}$ vérifie la condition $\tlQ1)$ ci-dessus.

Supposons qu'il existe un $\al \in \Delta_{Q}$ tel que $\# \iota_{\tlQ}^{-1}(\al) > 1$. 
Alors $0 < i_{k-1} = i_{k} < n$, un tel $\al $ est unique et égale $\al_{\tlQ}$ et 
$\iota_{\tlQ}^{-1}(\al_{\tlQ}) = \{\tlal_{\tlQ}^{-}, \tlal_{\tlQ}^{+}\} \sbs \Delta_{\tlQ}^{\tlR}$.
Soit $Q \sbs S \sbs R$. Si $\al_{\tlQ} \nin \Delta_{Q}^{S}$, il existe un unique 
$\tlS$ entre $\tlQ$ et $\tlR$ tel que $\tlS \cap G = S$ et en plus 
$d_{\tlS}^{\tlG} = d_{S}^{G}$. 
Si $\al_{\tlQ} \in \Delta_{Q}^{S}$, soit $\tlS$ tel que 
$\Delta_{\tlQ}^{\tlS} = \iota_{\tlQ}^{-1}(\Delta_{Q}^{S})$. 
Alors les sous-groupes $\tlS^{+}$ et $\tlS^{-}$ correspondants 
à $\Delta_{\tlQ}^{\tlS} \smin \{\tlal_{\tlQ}^{-}\}$ et $\Delta_{\tlQ}^{\tlS} \smin \{\tlal_{\tlQ}^{+}\}$ 
vérifient $\tlS^{-} \cap G = \tlS^{+} \cap G = \tlS \cap G = S$ et ces sont tous les 
sous-groupes entre $\tlQ$ et $\tlR$ 
dont l'intersection avec $G$ vaut $S$.
On a donc 
$(-1)^{d_{\tlS}^{\tlG}} + (-1)^{d_{\tlS^{+}}^{\tlG}} + (-1)^{d_{\tlS^{-}}^{\tlG}} = (-1)^{d_{\tlS^{-}}^{\tlG}} = (-1)^{d_{S}^{G}}$ 
et le point \textit{i)} suit. Les points \textit{ii)} et \textit{iii)} sont alors évidents car tout 
$\al \in \Delta_{Q}$ vérifie soit la condition $\tlQ1)$ soit $\tlQ2)$ ci-dessus.

Démontrons le point 2).
Dans ce cas il existe un unique $\tlal_{\tlQ} \in \{\tlal_{\tlQ}^{-}, \tlal_{\tlQ}^{+} \} \cap \Delta_{\tlQ}^{\tlR}$. 
Il vérifie $\iota_{\tlQ}(\tlal_{\tlQ}) \nin \Delta_{Q}^{R}$.
Soit $S$ entre $Q$ et $R$. Notons $\tlS$ le sous-groupe parabolique tel que 
$\Delta_{\tlQ}^{\tlS} = \iota_{\tlQ}^{-1}(\Delta_{Q}^{S})$. Alors 
$\tlal_{\tlQ}  \nin \Delta_{\tlQ}^{\tlS}$ et si l'on note $\tlS_{1}$
le sous-groupe correspondant à 
$\Delta_{\tlQ}^{\tlS} \cup \{\tlal_{\tlQ}\}$ on a $\tlS_{1} \cap G = S$ et 
$\tlS$ et $\tlS_{1}$ sont les seuls sous-groupes paraboliques entre $\tlQ$ et $\tlR$
dont l'intersection avec $G$ vaut $S$. On a donc $d_{\tlS}^{\tlG} -1= d_{\tlS_{1}}^{\tlG}$ et le résultat suit.
\edem
\elem

On est prêt à démontrer la proposition \ref{prop:opTronqRel}. En utilisant le lemme \ref{lem:hardNarashimArthur} on a 
pour $x \in \Sgl_{B}^{G}$
\[
\Laz^{T}\phi(\omega, x) = 
\sum_{B \sbs \tlP_{1} \sbs \tlP_{2}}
\sum_{\delta \in P_{1}(\rmF) \bsl G(\rmF)}
F^{\tlone}(\delta x, T_{\tlone})\sigma_{\tlone}^{\tltwo}(H_{\tlone}(\delta x) - T_{\tlone})
\phi_{\tlP_{1, \rmE}, \tlP_{2, \rmE}}(\omega, \delta x)
\]
où $\phi_{\tlP_{1,\rmE}, \tlP_{2,\rmE}}(\omega, y) = \sum_{\tlP_{1} \sbs \tlP \sbs \tlP_{2}}(-1)^{d_{\tlP}^{\tlG}}\phi_{P_{\rmE}}(\omega, y)$. 
Fixons $\tlP_{1} \sbs \tlP_{2}$. En vertu du lemme \ref{lem:admissiblePars} ci-dessus on 
a $\phi_{\tlP_{1,\rmE}, \tlP_{2,\rmE}} \equiv 0$ sauf si $\tlP_{2}$ est admissible par rapport à $\tlP_{1}$ 
donc on suppose que c'est le cas. En utilisant le lemme \ref{lem:admissiblePars}
on a donc $|\phi_{\tlP_{1,\rmE}, \tlP_{2,\rmE}}| = |\phi_{P_{1,\rmE}, P_{2,\rmE}}|$ où 
$\phi_{P_{1,\rmE}, P_{2,\rmE}}(\omega, y) := \sum_{P_{1} \sbs P \sbs P_{2}}(-1)^{d_{P}^{G}}\phi_{P_{\rmE}}(\omega, y)$.

Soit $\delta \in P_{1}(\rmF) \bsl G(\rmF)$. 
Écrivons $\delta x = n_{2}n_{1}^{2}amk$ où 
$n_{2} \in N_{2}(\A)$,  $n_{1}^{2} \in N_{1}^{2}(\A)$ et $k \in K$ appartiennent aux compacts fixés,
$a \in A_{\tlone}^{st,\infty}$ et 
$m \in H_{\tlone}(\A)^{1} \times G_{\tlone}(\A)$ est tel que $F^{\tlone}(m, T_{\tlone}) = 1$. 
D'après le corollaire 1.5 de \cite{leMoi2} on a que $m$ appartient à un compact fixé dans $M_{1}(\rmF) \bsl M_{1}(\A)$. 
On a alors
\[
\sigma_{\tlone}^{\tltwo}(H_{\tlone}(\delta x) - T_{\tlone}) = 
\sigma_{\tlone}^{\tltwo}(H_{\tlone}(a) + H_{\tlone}(m)- T_{\tlone}) = 
\sigma_{\tlone}^{\tltwo}(H_{\tlone}(a) + T_{\tlone}')
\]
pour un $T_{\tlone}' \in \all_{\tlone}$ dont la norme est majorée par celle de $T$. 
En utilisant la définition de la fonction $\sigma_{\tlone}^{\tltwo}$ donnée dans le paragraphe \ref{par:resGen} et 
le fait que $n_{1}^{2} \in N_{\tlone}^{\tltwo}(\A)$ on obtient:
\[
\phi_{P_{1,\rmE}, P_{2,\rmE}}(\omega,  \delta x) = \phi_{P_{1,\rmE}, P_{2,\rmE}}(\omega,  n_{2}n_{1}^{2}mak) =  
\phi_{P_{1,\rmE}, P_{2,\rmE}}(\omega, a a^{-1}n_{1}^{2}amk) = \phi_{P_{1,\rmE}, P_{2,\rmE}}(\omega, a c) 
\]
où $c$ appartient à un compact fixé de $G(\A)$
 qui dépend de $T$
 
 Il est démontré dans \cite{arthur5}, pages 93-95, qu'il 
existe un nombre fini des $\rmF$-sous-groupes $N_{I}$ de $N_{1,\rmE}$, et pour tout $I$ il existe un 
$\beta_{I} \in (\all_{1}^{2})^{*}$ tel que $\beta_{I} = \sum_{\al \in \Delta_{1}^{2}} a_{I,\al}\al$ où $a_{I,\al} >0$ 
tels que pour tout $n > 0$ il existe un nombre fini d'éléments $X_{i}$ de la complexification de l'algèbre 
enveloppante de l'algèbre de Lie de $G_{\rmE}(\rmF_{\infty})$ tels que $|\phi_{P_{1, \rmE}, P_{2, \rmE}}(\omega, ac)|$ est majorée
par
\begin{equation}\label{eq:sumNi}
 \sum_{I,i} e^{-n\beta_{I}(H_{1}(a))} \int_{[N_{I}]}|R(X_{i})\phi(\omega, u \delta x)|du.
\end{equation}
En vertu du lemme \ref{lem:admissiblePars} point \textit{1) ii)} on a $\tau_{1}^{2}(H_{1}(a) -T') = 1$ 
pour un $T' \in \all_{1}$
 car $\sigma_{\tlone}^{\tltwo}(H_{\tlone}(a) - T_{\tlone}') = 1$ implique $\tau_{\tlone}^{\tltwo}(H_{\tlone}(a)-T_{\tlone}') = 1$.
Puisque $\tau_{1}^{2}(H_{1}(a) - T') = 1$, il existe des constantes $c_{1}', c_{1} > 0$ indépendantes de $n$ telles que 
$-n\beta_{I}(H_{1}(a)) \le c_{1}' + -n c_{1}(\sum_{\al \in \Delta_{1}^{2}}\al(H_{1}(a))$. 
En vertu du lemme \ref{lem:admissiblePars} point \textit{1) iii)} 
ci-dessus, on a $\sum_{\al \in \Delta_{1}^{2}}\al(H_{1}(a)) = \sum_{\tlal \in \Delta_{\tlone}^{\tltwo}}\tlal(H_{\tlone}(a))$. 
On a $ac = \delta x$.
Si l'on note $(\delta x)^{1}$ la projection de $\delta x$ à $\tlG_{\rmE}(\A)^{1}$, on voit que puisque
$\sigma_{\tlone}^{\tltwo}(H_{\tlone}(a) - T_{\tlone}') = 1 $, le lemme \ref{lem:corArth62} \textit{ii)}
implique que $-nc_{1}(\sum_{\tlal \in \Delta_{\tlone}^{\tltwo}}\tlal(H_{\tlone}(a)) \le c_{2}' - nc_{2} \log(\|(\delta x)^{1}\|)$ 
pour des constantes $c_{2}',c_{2} >0$ indépendantes de $n$. En appliquant alors le lemme 
\ref{lem:hautComp} on trouve  
une constante $c_{3} >0$ indépendante de $n$ telle que $-n\beta_{I}(H_{0}(a))\le c_{2}' - nc_{3} \log(\|\delta x\|)$ pour tout $I$.
Puisque $x \in \Sgl_{B}^{G}$, il existe finalement une constante $c_{4} > 0$ telle que $\|\delta x\| \ge c_{4}\|x\|$, cela montre que 
pour tous $t_{1} \in [0,1]$, $I$ et $i$ on a que l'intégrale sur $\Omega$ de l'expression \eqref{eq:sumNi} est majorée par
\begin{multline*}
  \sum_{I,i} \int\limits_{\Omega}\int\limits_{[N_{I}]}|R(X_{i})\phi(u \delta x)|\|u \delta x\|^{-n c_{3}(1-t_{1})}du d\omega \|x\|^{-nc_{4}t_{1}} 
  \le  \\ c_{5} \|x\|^{-nc_{4}t_{1}} 
  \sum_{i} \sup_{y \in G_{\rmE}(\A)}
  \int\limits_{\Omega} |R(X_{i})\phi(\omega, y)|\|y\|^{-n c_{3}(1-t_{1})}d\omega 
\end{multline*}
pour un $c_{5} >0$. 
D'autre côté, en vertu du lemme \ref{lem:lemmeArth51} \textit{i)} 
on a $\sum_{\delta \in P_{1}(\rmF) \bsl G(\rmF)}F^{\tlone}(\delta x, T_{\tlone})\sigma_{\tlone}^{\tltwo}(H_{\tlone}(\delta x) - T_{\tlone}) 
\le c_{6} \|x\|^{N}$ où $N$ ne dépend pas de $T$. En prenant alors $n$ suffisamment grand et $t_{1}$ adapté on conclut la preuve. 
\edem
\erop

\subsection{Opérateur de troncature diagonal}\label{par:opTronqD}

Soient $\tlQ \in \relPb$ et $\phi$ une fonction localement intégrable 
sur  $(Q(\rmF) \times \tlQ(\rmF)) \bsl (G(\A) \times \tlG(\A))$. 
On pose pour $x \in N_{Q}(\A)M_{Q}(\rmF) \bsl G(\A)$, $\tlx \in N_{\tlQ}(\A)M_{\tlQ}(\rmF) \bsl \tlG(\A)$ 
et $T \in \all_{\tlzero}$
\[
\Lad^{\tlQ,T}\phi(x, \tlx) = \sum_{B \sbs \tlP \sbs \tlQ}(-1)^{d_{\tlP}^{\tlQ}}
\sum_{\delta \in P(\rmF) \bsl Q(\rmF)}
\htau_{\tlP}^{\tlQ}(H_{\tlP}(\delta \tlx) - T)
\phi_{P \times \tlP}(\delta x, \delta \tlx).
\]

On a la propriété d'inversion:
\blem\label{lem:formDInvD}
 Soit $\tlQ \in \relPb$. Alors pour tous $x  \in G(\A)$,  $\tlx \in \tlG(\A)$ et 
 $T \in \all_{\tlzero}$ on a 
\[
\phi_{Q \times \tlQ}(x, \tlx) = \sum_{B \sbs \tlP \sbs \tlQ}\sum_{\delta \in P(\rmF) \bsl Q(\rmF)}
\tau_{\tlP}^{\tlQ}(H_{\tlP}(\delta \tlx) - T_{\tlP})\Lad^{\tlP,T}\phi(\delta x, \delta \tlx).
\]
\bdem
La preuve est identique à celle du lemme \ref{lem:formDInv} ci-dessus.
\edem 
\elem

\brop\label{prop:opTronqD}
 Soient $\tlP \in \relPb$, $r_{1}, r_{2} \ge 0$ et $K_{0} \sbs (M_{P} \times M_{\tlP})(\A_{f})$ 
un compact ouvert. Il existe alors un sous-ensemble fini $\{X\}$ de l'algèbre enveloppante 
de la complexification de l'algèbre de Lie de 
\[
(H_{\tlP}(\rmF_{\infty}) \cap H_{\tlP}(\A)^{1}) \times G_{\tlP}(\rmF_{\infty})
\times 
(H_{\tlP}(\rmF_{\infty}) \cap H_{\tlP}(\A)^{1}) \times \tlG_{\tlP}(\rmF_{\infty})
\]
(voir (\ref{eq:MPArelSt})) 
tel que pour tout espace mesurable $(\Omega, d\omega)$, 
tout 
\[
\phi : \Omega \rar C^{\infty}( 
(P(\rmF) \times \tlP(\rmF)) \bsl (G(\A) \times \tlG(\A)) / K_{0})
\]
 mesurable, 
tout $T \in T_{+} + \all_{\tlzero}^{+}$ et
tout 
\[
x \in \Sgl_{B \cap M_{P}}^{M_{P}} \cap (H_{\tlP}(\A)^{1} \times G_{\tlP}(\A))
\]
on a que l'expression
\[
\int_{\Omega}|\Lad^{\tlP,T}\phi(\omega, x, x)|d\omega 
\]
où on applique l'opérateur $\Lad^{\tlP,T}$ à la fonction $\phi(\omega, (\cdot, \cdot))$ 
et on met comme argument $(x,x) \in 
(P(\rmF) \times P(\rmF)) \bsl (G(\A) \times G(\A)) \sbs 
(P(\rmF) \times \tlP(\rmF)) \bsl (G(\A) \times \tlG(\A))$, 
 est majorée par
 \[
\|x \|^{-r_{1}} \sum_{X}
\sup_{\begin{subarray}{c}
y \in H_{\tlP}(\A)^{1} \times G_{\tlP}(\A),\\
\tly \in H_{\tlP}(\A)^{1} \times \tlG_{\tlP}(\A)
\end{subarray}
}
\dsl 
\|y \|^{-r_{2}} \|\tly \|^{-r_{2}}\int_{\Omega}|
\int_{[N_{P}]}
\int_{[N_{\tlP}]}
R(X)\phi(\omega, ny, \tln \tly )dndn |  d \omega
\rb.
 \]
\bdem
On suit de près l'argument de \cite{arthur5}, pages 92-95.

En utilisant le lemme \ref{lem:hardNarashimArthur} pour $G \hrar \tlG$ (il est valable bien sûr dans ce cas aussi), on a 
pour pour $x \in \Sgl_{B \cap M_{P}}^{M_{P}} \cap (H_{\tlP}(\A)^{1} \times G_{\tlP}(\A))$
\[
\Lad^{\tlP,T}\phi(x,x) = 
\sum_{B \sbs \tlP_{1} \sbs \tlP_{2} \sbs \tlP}
\sum_{\delta \in P_{1}(\rmF) \bsl P(\rmF)}
F^{\tlone}(\delta x, T_{\tlone})\sigma_{\tlone}^{\tltwo}(H_{\tlone}(\delta x) - T_{\tlone})
\phi_{\tlP_{1}, \tlP_{2}}(\delta x)
\]
où 
\begin{equation}\label{eq:phip1p2}
\phi_{\tlP_{1}, \tlP_{2}}(y) = \sum_{\tlP_{1} \sbs \tlQ \sbs \tlP_{2}} (-1)^{d_{\tlQ}^{\tlP}}\phi_{Q \times \tlQ}(\omega, y,y)
\end{equation}
Soit $\delta \in P_{1}(\rmF) \bsl P(\rmF)$. 
Selon les décompositions données dans le paragraphe \ref{par:glnGlnplus1Sp}, 
écrivons $\delta x = n_{2}^{P}n_{1}^{2}amk$ où 
$n_{2}^{P} \in N_{2}^{P}(\A)$,  $n_{1}^{2} \in N_{1}^{2}(\A)$ et $k \in K \cap M_{P}(\A)$ appartiennent aux compacts fixés,
$a \in A_{\tlone}^{st, \tlP, \infty} := A_{\tlone}^{st, \infty} \cap (H_{\tlP}(\A)^{1} \times G_{\tlP}(\A))$ et 
$m \in H_{\tlone}(\A)^{1} \times G_{\tlone}(\A)$ est tel que $F^{\tlone}(m, T_{\tlone}) = 1$. 
D'après le corollaire 1.5 de \cite{leMoi2} on a que $m$ appartient à un compact fixé dans $M_{1}(\rmF) \bsl M_{1}(\A)$. 
On a alors
\[
\sigma_{\tlone}^{\tltwo}(H_{\tlone}(\delta x) - T_{\tlone}) = 
\sigma_{\tlone}^{\tltwo}(H_{\tlone}(a) + H_{\tlone}(m)- T_{\tlone}) = 
\sigma_{\tlone}^{\tltwo}(H_{\tlone}(a) + T_{\tlone}')
\]
pour un $T_{\tlone}' \in \all_{\tlone}$ dont la norme est majorée par celle de $T$. 
En utilisant le fait que $n_{1}^{2} \in N_{\tlone}^{\tltwo}(\A)$ on voit donc que
\[
\phi_{\tlP_{1}, \tlP_{2}}(\omega,  \delta x) = \phi_{\tlP_{1}, \tlP_{2}}(\omega,  n_{2}^{P}n_{1}^{2}mak) =  
\phi_{\tlP_{1}, \tlP_{2}}(\omega, a a^{-1}n_{1}^{2}amk) = \phi_{\tlP_{1}, \tlP_{2}}(\omega, a c) 
\]
où $c$ appartient à un compact fixé de $M_{P}(\A)$
 qui dépend de $T$.
 
Au début de la preuve de la proposition \ref{prop:opTronqRel} nous avons introduit une application
$\iota_{\tlQ} : \Delta_{\tlQ} \rar \Delta_{Q} \cup \{0\}$ pour tout 
$\tlQ \in \relPb$. 
Notons alors $\iota_{\tlone}^{\tltwo} : \Delta_{\tlone}^{\tltwo} \rar \Delta_{1}^{2} \cup \{0\}$ l'application 
suivante: pour tout $\tlal \in \Delta_{\tlone}^{\tltwo}$ 
on pose $\iota_{\tlone}^{\tltwo}(\tlal) = \iota_{\tlone}(\tlal)$ 
si $\iota_{\tlone}(\tlal) \in \Delta_{1}^{2}$, sinon on pose $\iota_{\tlone}^{\tltwo}(\tlal) = 0$. 
Pour tout $\tlal \in \Delta_{\tlone}^{\tltwo}$ (resp. $\al \in \Delta_{1}^{2}$) 
soit $N_{\tlal}$ (resp. $N_{\al}$) 
la partie unipotente du sous-groupe de $\tlP_{1}$ (resp. $P_{1}$)
associé à $\Delta_{\tlone}^{\tltwo} \smin \{\tlal\}$ (resp. $\Delta_{1}^{2} \smin \{\al\}$). 

Notons $\all_{\tlone}^{st, \tlP}$ l'image de $ A_{\tlone}^{st, \tlP, \infty}$ par l'application 
$H_{\tlP}$ définie dans le paragraphe \ref{par:prelimstraceSp}. On a alors
\begin{enumerate}
\item[P1)] Pour tout $\tlal \in \Delta_{\tlone}^{\tltwo}$ 
on a $\iota_{\tlone}^{\tltwo}(\tlal) = 0$ si et seulement si $N_{\tlal} \cap N_{1} = N_{2}$.
\item[P2)] Pour tout $\al \in \Delta_{1}^{2}$, l'ensemble $(\iota_{\tlone}^{\tltwo})^{-1}(\{\al\})$ 
n'est pas vide et pour tout $\tlal \in (\iota_{\tlone}^{\tltwo})^{-1}(\{\al\})$ 
on a $N_{\tlal} \cap G = N_{\al}$.
\item[P3)] Il existe au plus un $\al \in \Delta_{1}^{2}$ tel que $\# (\iota_{\tlone}^{\tltwo})^{-1}(\{\al\}) > 1$. 
Si c'est le cas, on a, 
avec la notation de la preuve de la proposition \ref{prop:opTronqRel}, 
$(\iota_{\tlone}^{\tltwo})^{-1}(\{\al\}) = \{\tlal_{\tlone}^{-}, \tlal_{\tlone}^{+}\}$ et
pour tout $H \in \all_{\tlone}^{st, \tlP}$ on a $\al(H) = (\tlal_{\tlone}^{-} + \tlal_{\tlone}^{+})(H)$.
\item[P4)] Si $\al \in \Delta_{1}^{2}$ est tel que $(\iota_{\tlone}^{\tltwo})^{-1}(\{\al\})= \{\tlal\}$ 
on a $\al(H) = \tlal(H)$ pour tout $H \in \all_{\tlone}^{st, \tlP}$. 
\end{enumerate}

Pour des $\tlR, \tlQ \in \relPb$ tels que $\tlR \sbs \tlQ \sbs \tlP$ 
on note $\nl_{\tlR} = \Lie N_{\tlR}$, $\nl_{\tlR}^{\tlQ} = \Lie N_{\tlR}^{\tlQ}$, 
$\nl_{R} = \Lie N_{R}$ etc. 
Pour un $\rmF$-groupe algébrique $H$ on note aussi $H_{\Q} = \Res_{\rmF/\Q}H$.

Soit $\tlal \in \Delta_{\tlone}^{\tltwo}$. Posons $\nl_{\tlal} = \Lie N_{\tlal}$, 
$\nl_{\tlal}^{\tltwo} = \nl_{\tlal} \cap \nl_{\tlone}^{\tltwo}$ et 
soient $n_{\tlal}$ et $m_{\tlal}$ les dimensions sur $\Q$ 
de $\nl_{\tlal, \Q}^{\tltwo}(\Q)$ et de $\nl_{\tlal, \Q}^{\tltwo}(\Q) \cap \nl^{2}_{1,\Q}(\Q)$ respectivement, où 
$\nl_{\tlal, \Q}^{\tltwo} := (\nl_{\tlal}^{\tltwo})_{\Q}$ etc.
Fixons $\{Y_{\tlal,1}, \ldots, Y_{\tlal, n_{\tlal}}\}$
une $\Q$-base de $\nl_{\tlal,\Q}^{\tltwo}(\Q)$ ainsi que $\{X_{\tlal, 1}, \ldots, X_{\tlal, m_{\tlal}}\}$ 
une $\Q$-base de $\nl_{\tlal,\Q}^{\tltwo}(\Q) \cap \nl^{2}_{1,\Q}(\Q)$. On suppose que 
tout $Y_{\tlal,i}$ (resp. $X_{\tlal,j}$) 
est un vecteur propre de la racine $\beta_{\tlal, i}$ \resp{$\gamma_{\tlal, j}$}
pour l'action de $A_{\tlone}$ sur $\nl_{\tlone}^{\tltwo}$ 
\resp{de $A_{1}$ sur $\nl_{1}^{2}$}. 
On suppose aussi que,  pour $k \le l$, la hauteur de  $\beta_{\tlal, k}$ (resp. de $\gamma_{\tlal, k}$) 
est plus grande que celle de $\beta_{\tlal, l}$ (resp. de $\gamma_{\tlal, l}$).

Soit alors $\nl_{\tlal, i}$, pour $0 \le i \le n_{\tlal}$,
la  $\Q$-algèbre de Lie 
engendré par les $\{Y_{\tlal, 1}, \ldots, Y_{\tlal, i}\}$  
et $\nl_{2, \Q} \times \nl_{\tltwo, \Q}$. On a donc
$\nl_{\tlal, n_{\tlal}} = \nl_{2, \Q} \times \nl_{\tlal, \Q}$. 
Notons aussi, pour $0 \le j \le m_{\tlal}$, $\nl_{\tlal, n_{\tlal} + j}$ la somme directe de l'algèbre engendré 
par les $\{X_{\tlal, 1}, \ldots, X_{\tlal, j}\}$ avec 
$\nl_{2, \Q} \times \nl_{\tlal, \Q}$. On a donc 
$\nl_{\tlal, n_{\tlal} + m_{\tlal}} = (\nl_{\tlal, \Q} \cap \nl_{\tlone,\Q}) \times \nl_{\tlal, \Q}$. 
Soit $\exp : \nl_{1} \times \nl_{\tlone} \rar N_{1} \times N_{\tlone}$
un isomorphisme  $A_{1} \times A_{\tlone}$ équivariant.
 On note donc aussi $N_{\tlal, k} = \exp \nl_{\tlal, k}$ où $0 \le k \le n_{\tlal} + m_{\tlal}$. 
C'est un $\Q$-sous-groupe normal de $N_{1,\Q} \times N_{\tlone, \Q}$.

Pour un $\Q$-sous-groupe $N$ de $N_{1, \Q} \times N_{\tlone, \Q}$ on définit 
$\pi(N)$ l'opérateur qui appliqué à une fonction $\upphi$ sur $G(\A) \times \tlG(\A)$ 
donne
\[
G(\A) \times \tlG(\A) \ni y \mapsto \int_{[N]_{\Q}} \upphi(n y)dn
\]
où l'on note $[N]_{\Q} = N(\Q) \bsl N(\A_{\Q})$ où $\A_{\Q}$ c'est l'anneau des adèles de $\Q$. 
On a donc que $\phi_{\tlP_{1}, \tlP_{2}}(\omega, \cdot)$ donné par (\ref{eq:phip1p2}) ci-dessus 
s'obtient, à un signe près, par application de l'opérateur 
\[
\prod_{\tlal \in \Delta_{\tlone}^{\tltwo}}(\pi (N_{2, \Q} \times N_{\tltwo, \Q}) -\pi ((N_{\tlal, \Q} \cap N_{1}) \times N_{\tlal, \Q}))
\]
à $\phi(\omega, \cdot)$. 
Or, on a pour $\tlal \in \Delta_{\tlone}^{\tltwo}$ que $\pi (N_{2, \Q} \times N_{\tltwo, \Q}) - \pi ((N_{\tlal, \Q} \cap N_{1}) \times N_{\tlal, \Q})$ 
égale
\[
\sum_{i = 1}^{n_{\tlal}}
( \pi(N_{\tlal, i-1}) - \pi(N_{\tlal, i})) + 
\sum_{j = 1}^{m_{\tlal}}
( \pi(N_{\tlal, n_{\tlal} + j-1}) - \pi(N_{\tlal, n_{\tlal} + j})).
\]
 
 Soient $\calS_{\tlone}, \calS_{1} \sbs \Delta_{\tlone}^{\tltwo}$ 
 tels que $\calS_{\tlone} \sqcup \calS_{1} = \Delta_{\tlone}^{\tltwo}$ et 
$\iota_{\tlone}^{\tltwo}(\calS_{1}) \sbs \Delta_{1}^{2}$. 
 Pour tout $\tlal \in \calS_{\tlone}$ on choisit un $i_{\tlal}$ entre $1$ et $n_{\tlal}$ 
 et on construit l'ensemble $I_{\tlone} = \{i_{\tlal}\}_{\tlal \in \calS_{\tlone}}$. 
 De même, pour tout $\tlal \in \calS_{1}$ on choisit un $j_{\tlal}$ entre $1$ et $m_{\tlal}$ 
 et on construit $I_{1} = \{j_{\tlal}\}_{\tlal \in \calS_{1}}$.
Pour ces données choisies, on pose
\begin{gather*}
N_{I_{\tlone}} = \prod_{\tlal \in \calS_{\tlone}}N_{\tlal, i_{\tlal}} \cap N_{\tlone}, \quad 
N_{I_{1}} = \prod_{\tlal \in \calS_{1}}N_{\tlal, n_{\tlal} + j_{\tlal}} \cap N_{1}, \\
N_{\brI_{\tlone}} = \prod_{\tlal \in \calS_{\tlone}}N_{\tlal, i_{\tlal} -1} \cap N_{\tlone}, \quad 
N_{\brI_{1}} = \prod_{\tlal \in \calS_{1}}N_{\tlal, n_{\tlal} + j_{\tlal} -1} \cap N_{1},
\end{gather*}
ainsi que $\nl^{I_{\tlone}}$ \resp{$\nl^{I_{1}}$} l'espace engendré 
par $\{Y_{\tlal, i_{\tlal}}\}_{\tlal \in \calS_{\tlone}}$ \resp{$\{X_{\tlal, j_{\tlal}}\}_{\tlal \in \calS_{1}}$}. 
Fixons une base de $\nl^{I_{\tlone}}(\Q)$  \resp{de $\nl^{I_{1}}(\Q)$} parmi les éléments 
de $\{Y_{\tlal, i_{\tlal}}\}_{\tlal \in \calS_{\tlone}}$ \resp{de $\{X_{\tlal, j_{\tlal}}\}_{\tlal \in \calS_{1}}$}.
On note aussi  $\nl^{I_{\tlone}}(\Q)'$ \resp{$\nl^{I_{1}}(\Q)'$} 
le sous-ensemble de $\nl^{I_{\tlone}}(\Q)$ \resp{$\nl^{I_{1}}(\Q)$} 
composé d'éléments dont toutes les coordonnées dans la base choisie sont non-nulles. 

Soit $\psi : \Q \bsl \A_{\Q} \rar \C^{*}$ un caractère continu non-trivial. En utilisant la formule d'inversion de 
Fourier, on voit que $\phi_{\tlP_{1}, \tlP_{2}}(\omega, y)$ égale la somme 
sur tous les choix des $\calS_{\tlone}$, $\calS_{1}$, $I_{\tlone}$, $I_{1}$ comme ci-dessus, de
\[
\sum_{\eta \in \nl^{I_{1}}(\Q)'}
\sum_{\xi \in \nl^{I_{\tlone}}(\Q)'}
\int\limits_{[\nl^{I_{1}}]_{\Q}}
\int\limits_{[\nl^{I_{\tlone}}]_{\Q}}
\int\limits_{[N_{\brI_{1}}]_{\Q}}
\int\limits_{\mathrlap{[N_{\brI_{\tlone}}]_{\Q}}}
\phi(\omega, n \exp(X) y, \tln \exp(\tlX) y)
\psi( \langle X, \eta \rangle )
\psi( \langle \tlX, \xi \rangle )d\tln dn d\tlX dX
\]
où $\bilif$ c'est le produit scalaire défini par les bases respectives. 

On met $y = ac$ comme avant. 
En raisonnant maintenant comme dans les pages 94-95 de \cite{arthur5} 
on trouve que pour tout $n \in \N$ assez grand, il existe 
une constante $C$ qui ne dépend que de $n$, $T$ 
et du compact ouvert $K_{0}$ 
ainsi qu'un nombre fini $\{Y_{j}\}$, 
$\{X_{i}\}$ des éléments des algèbres enveloppantes 
des complexifications des algèbres de Lie de 
$((H_{\tlP}(\rmF_{\infty}) \cap H_{\tlP}(\A)^{1}) \times G_{\tlP}(\rmF_{\infty}))$ et 
$(H_{\tlP}(\rmF_{\infty}) \cap H_{\tlP}(\A)^{1}) \times \tlG_{\tlP}(\rmF_{\infty})$ 
respectivement 
tel que l'expression ci-dessus, quand $y = ac$, 
est majorée par
\[
C
\sum_{i,j}
\int\limits_{[N_{I_{1}}^{P}]_{\Q}}
\int\limits_{[N_{I_{\tlone}}^{\tlP}]_{\Q}}
|
\int\limits_{[N_{P}]}
\int\limits_{\mathrlap{[N_{\tlP}]}}
R_{2}(Y_{j})R_{3}(X_{i})\phi(\omega, n_{P} n_{I_{1}}^{P} ac, n_{\tlP} n_{I_{\tlone}}^{\tlP} ac )dn_{\tlP}dn_{P}|
dn_{I_{\tlone}}^{\tlP} dn_{I_{1}}^{P} 
\]
fois $e$ à la puissance
\[
-n\sum_{\tlal \in \calS_{\tlone}}\beta_{\tlal, i_{\tlal}}(H_{\tlone}(a)) - 
n\sum_{\tlal \in \calS_{1}}\gamma_{\tlal, j_{\tlal}}(H_{1}(a)).
\]
En utilisant la définition de l'ensemble $\calS_{1} \sbs \Delta_{\tlone}^{\tltwo}$ 
et ensuite les propriétés P3) et P4) ci-dessus il résulte du fait 
que $\calS_{\tlone} \sqcup \calS_{1} = \Delta_{\tlone}^{\tltwo}$ 
que la somme ci-dessus égale 
$-n \sum_{\tlal \in \Delta_{\tlone}^{\tltwo}} c_{\tlal} \tlal (H_{\tlone} (a))$ 
pour certains $c_{\tlal} >0 $ qui ne dépendent que de $\calS_{\tlone}$, $\calS_{1}$, $I_{\tlone}$ 
et $I_{1}$. En raisonnant donc maintenant 
de façon analogue comme dans la preuve de la proposition \ref{prop:opTronqRel} après l'équation
(\ref{eq:sumNi}), en s'appuyant sur les lemmes \ref{lem:corArth62} \textit{ii)}, \ref{lem:lemmeArth51} \textit{i)} 
et \ref{lem:hautComp} on trouve que pour tout $r_{1}, r_{2} > 0$ 
il existe des constantes $C', C'' >0 $ telles que pour $\relPb \ni \tlP_{1} \sbs \tlP_{2} \sbs \tlP$ fixés on a
\begin{multline*}
\int_{\Omega}|\sum_{\mathrlap{\delta \in P_{1}(\rmF) \bsl P(\rmF)}}
F^{\tlone}(\delta x, T_{\tlone})\sigma_{\tlone}^{\tltwo}(H_{\tlone}(\delta x) - T_{\tlone}) \ 
\sum_{\mathclap{\tlP_{1} \sbs \tlQ \sbs \tlP_{2} \sbs \tlP}} \
(-1)^{d_{\tlQ}^{\tlP}}\phi_{Q \times \tlQ} (\omega, \delta x, \delta x)| d\omega \le 
C' \sum_{\calS_{\tlone}, \calS_{1}}
\sum_{I_{\tlone}, I_{1}}\sum_{i,j} 
\|x\|^{-r_{1}} \\
\int\limits_{\Omega}
\int\limits_{[N_{I_{1}}^{P}]_{\Q}}
\int\limits_{[N_{I_{\tlone}}^{\tlP}]_{\Q}}
|
\int\limits_{[N_{P}]}
\int\limits_{\mathrlap{[N_{\tlP}]}}
R_{2}(Y_{j})R_{3}(X_{i})\phi(\omega, n n_{I_{1}}^{P}\delta x, \tln n_{I_{\tlone}}^{\tlP} \delta x)d\tln dn|
\|n_{I_{\tlone}}^{\tlP} \delta x\|^{-r_{2}}
\|n_{I_{1}}^{P} \delta x\|^{-r_{2}}
dn_{I_{\tlone}}^{\tlP} dn_{I_{1}}^{P} \\ d\omega \le 
C''
\|x\|^{-r_{1}}
\sum_{i,j}
\sup_{\begin{subarray}{c}
y \in H_{\tlP}(\A)^{1} \times G_{\tlP}(\A),\\
\tly \in H_{\tlP}(\A)^{1} \times \tlG_{\tlP}(\A)
\end{subarray}}
(\int\limits_{\Omega}
|
\int\limits_{[N_{P}]}
\int\limits_{\mathrlap{[N_{\tlP}]}}
R_{2}(X_{j})R_{3}(Y_{i})\phi(\omega, n y, \tln \tly)
d\tln dn
|
\|y\|^{-r_{2}}
\|\tly\|^{-r_{2}}
d\omega).
\end{multline*}
ce qu'il fallait démontrer. 
\edem
\erop

Dans le cas de l'inclusion $U \hrar \tlU$ on introduit aussi l'opérateur de troncature diagonal. 
Soient $Q \in \calF(P_{0})$ et $\phi$ une fonction localement intégrable sur $(Q(\rmF) \times \tlQ(\rmF)) \bsl (U(\A) \times \tlU(\A))$. 
On pose pour $x \in N_{Q}(\A)M_{Q}(\rmF) \bsl U(\A)$, $\tlx \in N_{\tlQ}(\A)M_{\tlQ}(\rmF) \bsl \tlU(\A)$ 
et $T \in \all_{0}$
\[
\Lad^{Q,T}\phi(x, \tlx) = \sum_{P_{0 }\sbs P \sbs Q}(-1)^{d_{P}^{Q}}
\sum_{\delta \in P(\rmF) \bsl Q(\rmF)}
\htau_{P}^{Q}(H_{P}(\delta x) - T)
\phi_{P \times \tlP}(\delta x, \delta \tlx).
\]
Les analogues évidents du lemme \ref{lem:formDInvD} et de la proposition \ref{prop:opTronqD} sont alors 
aussi vrais dans le cas de groupes unitaires. 

\section{Le côté spectral de la formule des traces pour les groupes linéaires}\label{sec:SpecG}
\sectionmark{Le côté spectral pour les groupes linéaires}

Pour l'un des opérateurs $\La$ discutés dans la section \ref{sec:opsTronc} et une fonction 
$\Psi$ de plusieurs variables, par $\La_{n}\Psi$ on entend que l'on applique $\La$ à la $n$-ième variable 
considérant les autres variables fixés. Dans le cas de l'opérateur $\Lad$ introduit dans le paragraphe 
\ref{par:opTronqD}, on écrira $\Ladud \Psi$ pour signifier que l'on l'applique aux premières 2 variables etc.

\subsection{Convergence du noyau spectral tronqué}\label{par:specCvgG}

Soit $\det \in \all_{\tlzero}^{*}$ le déterminant de $\tlG$ pour son action sur $W$. 
On a alors $2\det \in \Hom_{\rmF}(\tlG_{\rmE}, \Gm)$ et pour tout $\tlx \in \tlG_{\rmE}(\A)$ 
et tout $s \in \C$ on définit
\[
|\det \tlx|_{\A}^{s} := e^{s \det H_{\tlG}(\tlx)}.
\]

Soit $\Phi \in C_{c}((G_{\rmE} \times \tlG_{\rmE})(\A))$ et soit 
$k = k_{\Phi}$ son noyau automorphe. 
Pour $\chi \in \calX^{G_{\rmE} \times \tlG_{\rmE}}$, $\tlP \in \relPb$
un sous-groupe parabolique 
relativement standard de $\tlG$, 
$x,y \in N_{P_{\rmE}}(\A)M_{P}(\rmE) \bsl G_{\rmE}(\A)$ et
$\tlx,\tly \in N_{\tlP_{\rmE}}(\A)M_{\tlP}(\rmE) \bsl \tlG_{\rmE}(\A)$
on pose
\[
k_{\tlP,\chi}(x,\tlx, y, \tly) = k_{\Phi, \tlP,\chi}(x,\tlx, y, \tly) := k_{\Phi, P_{\rmE} \times \tlP_{\rmE},\chi}(x,\tlx, y, \tly),
\]
où on utilise la convention du paragraphe \ref{par:prelimstraceSp} 
qu'à $\tlP \in \relPb$ on associe $P := \tlP \cap G \in \calF(B)$.

Pour $T \in \all_{\tlzero}$ et
$(x, \tlx), (y,\tly) \in \GetlGa$ on pose
\begin{equation}\label{eq:kchiGdef}
k_{\chi}^{T}(x,\tlx,y,\tly) = 
k_{\Phi, \chi}^{T}(x,\tlx,y,\tly) = \
\sum_{\mathclap{\tlP \in \relPb}} \ (-1)^{d_{\tlP}^{\tlG}} \
\sum_{\mathclap{
\begin{subarray}{c}
\delta_{1} \in P(\rmE)\bsl G(\rmE)
\\
\delta_{2} \in P(\rmF)\bsl G(\rmF)
\\
\delta_{3} \in \tlP(\rmF)\bsl \tlG(\rmF)
\end{subarray}}} \ 
\htau_{\tlP}(H_{\tlP}(\delta_{2}y)-T_{\tlP})
k_{\tlP,\chi}(\delta_{1} x, \delta_{1} \tlx, \delta_{2} y, \delta_{3} \tly).
\end{equation}

\begin{theo}\label{thm:noyauCVGG} Pour tout $\Phi \in C_{c}^{\infty}(\GetlGa)$, tous $\sigma, \sigma' \in \R$ 
et tout $T \in T_{+} + \all_{\tlzero}^{+}$
on a
\begin{equation}\label{eq:kGcvg}
\sum_{\chi \in \calX^{\GetlG}}
\int_{[G_{\rmE}]}\int_{[G]}\int_{[\tlG]}
|k_{\chi}^{T}(g,g,h,\tlh)| |\det g|_{\A}^{\sigma}
|\det h|_{\A}^{\sigma'}
d\tlh dh dg < \infty. 
\end{equation}

\bdem 
La preuve suivra la route tracée par la preuve du théorème 2.1 de \cite{arthur5}.

Soient $\chi \in \calX^{\GetlG}$, $\tlP \in \relPb$, $g \in G_{\rmE}(\A)$, $h \in G(\A)$ et $\tlh \in \tlG(\A)$ tels que 
$k_{\Phi, \tlP, \chi}(g,g,h,\tlh)$ est non nul. En vertu du lemme \ref{lem:artLemme23}, il existe des
$m \in M_{P_{\rmE}}(\A)^{1}$ et $\tlm \in M_{\tlP}(\A)^{1}$ tels que 
$k_{\Phi, \tlP}(g,\tlm g,mh,\tlh) \neq 0$. Il existe alors $\gamma \in M_{P}(\rmE)$,
$\tlgam \in M_{\tlP}(\rmE)$, $n \in N_{P_{\rmE}}(\A)$ et $\tln \in N_{\tlP_{\rmE}}(\A)$
tels que $g^{-1} n\gamma mh$ et $\tlh^{-1} \tln \tlgam \tlm g$ appartient à 
un compact de $\tlG_{\rmE}(\A) \sps G_{\rmE}(\A)$ fixé qui ne dépend que du support de $\Phi$. 
Puisque $M_{P_{\rmE}}(\A)^{1} \sbs M_{\tlP_{\rmE}}(\A)^{1}$, $N_{P_{\rmE}} \sbs N_{\tlP_{\rmE}}$, $P_{\rmE} \sbs \tlP_{\rmE}$ 
et $K_{\rmE} \sbs \tlK_{\rmE}$, l'argument entre les équations (2.2) et (2.3) de \cite{arthur5}, pages 100-101, 
montre qu'il existe une constante $C \in \R$ 
qui ne dépend que du compact à lequel appartient $g^{-1} n\gamma mh$ telle que
\[
\varpi(H_{\tlP}(g)) \ge C + \varpi(H_{\tlP}(h)), \quad  \forall \ \varpi \in \hDelta_{\tlP}.
\]
Par le même argument appliqué à $\tlh^{-1} \tln \tlgam \tlm g$ on trouve une constante $C'$ telle que
\[
\varpi(H_{\tlP}(\tlh)) \ge C' + \varpi(H_{\tlP}(g)), \quad  \forall \varpi \in \hDelta_{\tlP}.
\]
On voit alors qu'il existe des $T_{\Phi}', T_{\Phi}'' \in \all_{\tlzero}$, 
qui ne dépendent que du support de $\Phi$ tels que si l'on met $T' = T + T_{\Phi}'$ 
et $T'' = T + T_{\Phi}''$ on a que $k_{\chi}^{T}(g,g,h,\tlh)$ égale
\[
\sum_{\tlP}(-1)^{d_{\tlP}^{\tlG}}
\sum_{\mathclap{
\begin{subarray}{c}
\delta_{1} \in P(\rmE)\bsl G(\rmE)
\\
\delta_{2} \in P(\rmF)\bsl G(\rmF)
\\
\delta_{3} \in \tlP(\rmF)\bsl \tlG(\rmF)
\end{subarray}}}
\htau_{\tlP}(H_{\tlP}(\delta_{2}h)-T_{\tlP})
\htau_{\tlP}(H_{\tlP}(\delta_{1} g)-T_{\tlP}')
\htau_{\tlP}(H_{\tlP}(\delta_{3}\tlh)-T_{\tlP}'') \\
k_{\tlP,\chi}(\delta_{1} g, \delta_{1} g, \delta_{2}h, \delta_{3}\tlh).
\]
En particulier, en vertu du lemme \ref{lem:lemmeArth51} \textit{ii)} les sommes dans la définition de 
$k_{\Phi, \chi}^{T}(g,g,h,\tlh)$, où $g \in [G_{\rmE}]$, $h \in [G]$ et $\tlh \in [\tlG]$, sont finies et la fonction est bien définie.

Les fonctions 
$(g,\tlg) \mapsto k_{\tlP,\chi}(g, \tlg, \cdot, \cdot)$ 
et 
$g \mapsto k_{\tlP,\chi}(\cdot, \cdot, g, \cdot)$ égalent leur termes constants le long 
de $P_{\rmE} \times \tlP_{\rmE}$ et $P_{\rmE}$ respectivement. Pour tout $\tlP \in \relPb$, en utilisant le lemme \ref{lem:formDInvD} pour l'opérateur 
$\Ladud^{\tlP, T'}$ (sa version sur $\rmE$, i.e. pour les fonctions définies sur $(G_{\rmE} \times \tlG_{\rmE})(\A)$)
ainsi que le lemme \ref{lem:formDInv} pour l'opérateur $\La_{m,3}^{\tlP,T}$ et ensuite
le lemme \ref{lem:corArth62} \textit{i)} on trouve que
$k_{\chi}^{T}$ égale 
la somme sur les sous-groupes paraboliques relativement standards $\tlP_{1} \sbs \tlP_{4}$ 
et $\tlP_{2} \sbs \tlP_{5}$ de $\tlG$ de
\begin{multline*}
\sum_{
\begin{subarray}{c}
\delta_{1} \in P_{1}(\rmE)\bsl G(\rmE)
\\
\delta_{2} \in P_{2}(\rmF)\bsl G(\rmF)
\end{subarray}}
\sigma_{\tlone}^{\tlqt}(H_{\tlone}(\delta_{1} g)-T_{\tlone}')
\sigma_{\tltwo}^{\tlcq}(H_{\tltwo}(\delta_{2}h)-T_{\tltwo})
\Ladud^{\tlP_{1}, T'} 
\Lazt^{\tlP_{2}, T} \\
\dsl
\sum_{\tlQ_{1}' \sbs \tlP \sbs \tlQ_{2}'}
(-1)^{d_{\tlP}^{\tlG}}
\sum_{
\delta_{3} \in \tlP(\rmF)\bsl \tlG(\rmF)}
\htau_{\tlP}(H_{\tlP}(\delta_{3} \tlh)-T_{\tlP}'')
k_{\tlP,\chi}(\delta_{1} g, \delta_{1} g, \delta_{2}h, \delta_{3}\tlh)
\rb
\end{multline*}
où $\tlQ_{1}'$ c'est le plus petit sous-groupe parabolique de $\tlG$ 
contenant $\tlP_{1} \cup \tlP_{2}$ et $\tlQ_{2}' = \tlP_{4} \cap \tlP_{5}$. 
On considère ces sous-groupes fixés désormais. 
Soit $\tlP$ entre $\tlQ_{1}'$ et $\tlQ_{2}'$.
La fonction $\tlG_{\rmE}(\A) \ni \tly \mapsto k_{\tlP,\chi}(\cdot, \cdot, \cdot, \tly)$ 
égale son terme constant le long de $\tlP_{\rmE}$. On choisit alors un sous-groupe de Borel
$\tlB \in \relPb$ contenu dans $\tlP_{1}$ et on applique le lemme \ref{lem:formDInvtl} 
pour l'opérateur $\Lajlrq^{\tlP, T''}$ et le sous-groupe de Borel $\tlB$ de $\tlG$
et l'on ré-applique le lemme \ref{lem:corArth62} \textit{i)} et l'on trouve 
que l'intégrale \eqref{eq:kGcvg} est majorée par la somme sur les sous-groupes 
paraboliques relativement standards 
$\tlB \sbs \tlP_{1} \sbs \tlP_{4}$, $\tlP_{2} \sbs \tlP_{5}$ 
et $\tlB \sbs \tlP_{3} \sbs \tlP_{6}$
de $\tlG$ de
\begin{multline}\label{eq:AllSigmas}
\sum_{\chi} 
\int\limits_{P_{1}(\rmE)\bsl G_{\rmE}(\A)}
\int\limits_{P_{2}(\rmF)\bsl G(\A)}
\int\limits_{\tlP_{3}(\rmF)\bsl \tlG(\A)}
\sigma_{\tlone}^{\tlqt}(H_{\tlone}(g)-T_{\tlone}')
\sigma_{\tltwo}^{\tlcq}(H_{\tltwo}(h)-T_{\tltwo})
\sigma_{\tltr}^{\tlsix}(H_{\tltr}(\tlh)-T_{\tltr}'')
\\
| 
\Ladud^{\tlP_{1}, T'} 
\Lazt^{\tlP_{2}, T} 
\Lajlrq^{\tlP_{3}, T''}
( \sum_{\tlQ_{1} \sbs \tlP \sbs \tlQ_{2}}
(-1)^{d_{\tlP}^{\tlG}}
k_{\chi, \tlP}(g, g, h, \tlh)) |
|\det g|_{\A}^{\sigma} 
|\det h|_{\A}^{\sigma'}
 d\tlh dh dg
\end{multline}
où $\tlQ_{1}$ c'est le plus petit sous-groupe parabolique de $\tlG$ 
contenant $\tlP_{1} \cup \tlP_{2} \cup \tlP_{3}$ et $\tlQ_{2}= \tlP_{4} \cap \tlP_{5} \cap \tlP_{6}$.

Soient $x,g \in G_{\rmE}(\A)$, $h \in G(\A)$ et $\tly \in \tlG(\A)$.
On a $\Ladud^{\tlP_{1}, T'} \Lazt^{\tlP_{2}, T}
k_{\tlP, \chi}(x,g,h,\tly) = \\
\Ladud^{\tlP_{1}, T'} \Lazt^{\tlP_{2}, T}\int_{[N_{\tlone,\rmE}]}
\int_{[N_{2,\rmE}]}k_{P_{\rmE} \times \tlP_{\rmE}, \chi}(x,n_{\tlone}g,n_{2}h,\tly)dn_{2}dn_{\tlone}$.
Pour tout $\tlP$ entre $\tlQ_{1}$ et $\tlQ_{2}$, en utilisant le corollaire
\ref{cor:art23REALLY} et la décomposition de Bruhat, on a
\[
\int\limits_{[N_{\tlone,\rmE}]} \ \
\int\limits_{\mathclap{[N_{2,\rmE}]}}
k_{P_{\rmE} \times \tlP_{\rmE}, \chi}(x,n_{\tlone}g,n_{2}h,\tly)dn_{2}dn_{\tlone} = 
\sum_{\mathclap{
\begin{subarray}{c}
s \in \Omega^{1} \bsl \Omega^{P} / \Omega^{2} \\
\tls \in \Omega^{\tltr} \bsl \Omega^{\tlP} / \Omega^{\tlone}
\end{subarray}
}} \quad \quad
\sum_{\mathrlap{
\begin{subarray}{c}
\gamma \in (P_{1} \cap sP_{2})(\rmE) \bsl P_{1}(\rmE)\\
\tlgam \in (\tlP_{3} \cap \tls\tlP_{1})(\rmE) \bsl \tlP_{3}(\rmE)
\end{subarray}
}}
k_{P_{2,\rmE} \times \tlP_{1,\rmE}, \chi}(w_{s}^{-1}\gamma x, g, h, w_{\tls}^{-1} \tlgam \tly).
\]
Posons
\[
(\Omega_{\tlQ_{1}}^{\tlP})' = \Omega^{\tlP} \smin \bigcup_{\tlQ_{1} \sbs \tlR \sbn \tlP} \Omega^{\tlR}
\]
et $(\Omega_{\tlQ_{1}, G}^{\tlP})' = (\Omega_{\tlQ_{1}}^{\tlP})' \cap \Omega^{G}$. 
Alors $(\Omega_{\tlQ_{1}}^{\tlP})'$ (resp. $(\Omega_{\tlQ_{1}, G}^{\tlP})'$) est $\Omega^{\tltr}$-stable (resp. $\Omega^{1}$-stable)
à gauche et $\Omega^{\tlone}$-stable (resp. $\Omega^{2}$-stable) à droite.
On a alors les décompositions en parties disjointes suivantes
\[
\Omega^{\tlP} = \coprod_{\tlQ_{1} \sbs \tlR \sbs \tlP}(\Omega_{\tlQ_{1}}^{\tlP})' \quad 
\Omega^{P} = \coprod_{\tlQ_{1} \sbs \tlR \sbs \tlP}(\Omega_{\tlQ_{1}, G}^{\tlP})'.
\]
Pour $\tls \in \Omega^{\tlG}$ soit $\hDelta_{\tlQ_{1}, \tls}= \{\varpi \in \hDelta_{\tlQ_{1}} | \tls\varpi = \varpi\}$. 
Alors $\tls \in (\Omega_{\tlQ_{1}}^{\tlP})'$ si et seulement si $\hDelta_{\tlQ_{1}, \tls}= \hDelta_{\tlP}$. 
Par un argument classique basé sur l'identité (\ref{eq:basicidentity}) 
on obtient pour tous $x,g \in G_{\rmE}(\A)$, $h \in G(\A)$ et $\tly \in \tlG(\A)$.
\begin{equation}\label{eq:alternSumExpl}
\int\limits_{[N_{\tlone,\rmE}]} \ \
\int\limits_{\mathclap{[N_{2,\rmE}]}}
\sum_{\mathrlap{\tlQ_{1} \sbs \tlP \sbs \tlQ_{2}}}
(-1)^{d_{\tlP}^{\tlG}}
k_{\chi, \tlP}(x, n_{\tlone}g, n_{2}h, \tly)dn_{2}dn_{\tlone} \ \mathclap{=}
\sum_{(s, \tls) \in \Upomega_{\tlQ_{1}}'}
\sum_{\mathrlap{
\begin{subarray}{c}
\gamma \in (P_{1} \cap sP_{2})(\rmE) \bsl P_{1}(\rmE)\\
\tlgam \in (\tlP_{3} \cap \tls\tlP_{1})(\rmE) \bsl \tlP_{3}(\rmE)
\end{subarray}
}}
k_{P_{2,\rmE} \times \tlP_{1,\rmE}, \chi}(w_{s}^{-1}\gamma x, g, h, w_{\tls}^{-1} \tlgam \tly)
\end{equation}
où 
\[
\Upomega_{\tlQ_{1}}' = \{(s, \tls) \in (\Omega^{1} \bsl \Omega^{Q_{2}}/\Omega^{2}) \times 
(\Omega^{\tltr} \bsl \Omega^{\tlQ_{2}} / \Omega^{\tlone})| \hDelta_{\tlQ_{2}} = \hDelta_{\tlQ_{1}, s}\cap \hDelta_{\tlQ_{1}, \tls}\}.
\]
Remarquons qu'on regarde donc $s \in \Omega^{G}$ comme un élément de $\Omega^{\tlG}$.

On voit en particulier que l'expression entre la valeur absolue dans (\ref{eq:AllSigmas}) égale
\[
\Ladud^{\tlP_{1}, T'} 
\Lazt^{\tlP_{2}, T} 
\Lajlrq^{\tlP_{3}, T''}
(\sum_{(s, \tls) \in \Upomega_{\tlQ_{1}}'}
\sum_{
\begin{subarray}{c}
\gamma \in (P_{1} \cap sP_{2})(\rmE) \bsl P_{1}(\rmE)\\
\tlgam \in (\tlP_{3} \cap \tls\tlP_{1})(\rmE) \bsl \tlP_{3}(\rmE)
\end{subarray}
}
k_{P_{2,\rmE} \times \tlP_{1,\rmE}, \chi}(w_{s}^{-1}\gamma g, g, h, w_{\tls}^{-1} \tlgam \tlh)
 ).
\]

Écrivons maintenant $g = n_{1}a_{\tlone}^{st}m_{1}k_{1}$, $h = n_{2}a_{\tltwo}^{st}m_{2}k_{2}$ et 
$\tlh = n_{\tltr}a_{\tltr}m_{\tltr}k_{\tltr}$ selon les décompositions suivantes:
\begin{equation*}
\begin{split}
P_{1}(\rmE) \bsl G_{\rmE}(\A) & = [N_{1, \rmE}] \times A_{\tlone}^{st,\infty} 
\times 
(M_{1}(\rmE) \bsl (H_{\tlone, \rmE}(\A)^{1} \times G_{\tlone,\rmE}(\A)))
\times K_{\rmE}, \\
P_{2}(\rmF) \bsl G(\A) & = [N_{2}] \times  A_{\tltwo}^{st,\infty}  \times 
(M_{2}(\rmF) \bsl (H_{\tltwo}(\A)^{1} \times G_{\tltwo}(\A))) \times K, \\
\tlP_{3}(\rmF) \bsl \tlG(\A) & = [N_{\tltr}] \times A_{\tltr}^{\infty}  \times (M_{\tltr}(\rmF) \bsl  M_{\tltr}(\A)^{1}) \times \tlK.
\end{split}
\end{equation*}
Supposons que $\Phi$ est invariante à droite par un compact $K_{0} \sbs (G_{\rmE} \times \tlG_{\rmE})(\A_{f})$. 
Il résulte alors du lemme \ref{lem:artLemme23} que 
pour $(x,\tlx), (y,\tly) \in (G_{\rmE} \times \tlG_{\rmE})(\A)$  fixés
les fonctions 
\begin{gather*}
P_{1}(\rmE) \bsl G_{\rmE}(\A) \times \tlP_{1}(\rmE) \bsl \tlG_{\rmE}(\A) \ni (x_{1},\tlx_{1}) \mapsto \
\sum_{\mathclap{(s, \tls) \in \Upomega_{\tlQ_{1}}'}} \ \
\sum_{\mathrlap{
\begin{subarray}{c}
\gamma \in (P_{1} \cap sP_{2})(\rmE) \bsl P_{1}(\rmE)\\
\tlgam \in (\tlP_{3} \cap \tls\tlP_{1})(\rmE) \bsl \tlP_{3}(\rmE)
\end{subarray}
}}
k_{P_{2,\rmE} \times \tlP_{1,\rmE}, \chi}(w_{s}^{-1}\gamma x_{1} a_{\tlone}^{st}  k_{1}, 
\tlx_{1} a_{\tlone}^{st}  k_{1}, x, w_{\tls}^{-1}\tlgam \tlx),  \\
P_{2}(\rmE) \bsl G_{\rmE}(\A) \ni x_{2} \mapsto 
\sum_{(s, \tls) \in \Upomega_{\tlQ_{1}}'}
\sum_{\mathrlap{
\begin{subarray}{c}
\gamma \in (P_{1} \cap sP_{2})(\rmE) \bsl P_{1}(\rmE)\\
\tlgam \in (\tlP_{3} \cap \tls\tlP_{1})(\rmE) \bsl \tlP_{3}(\rmE)
\end{subarray}
}}
k_{P_{2,\rmE} \times \tlP_{1,\rmE}, \chi}(w_{s}^{-1}\gamma y, \tly, x_{2} a_{\tltwo}^{st}  k_{2}, w_{\tls}^{-1}\tlgam \tlx),  \\ 
\tlP_{3}(\rmE) \bsl \tlG_{\rmE}(\A) \ni \tlx_{3} \mapsto 
\sum_{(s, \tls) \in \Upomega_{\tlQ_{1}}'}
\sum_{\mathrlap{
\begin{subarray}{c}
\gamma \in (P_{1} \cap sP_{2})(\rmE) \bsl P_{1}(\rmE)\\
\tlgam \in (\tlP_{3} \cap \tls\tlP_{1})(\rmE) \bsl \tlP_{3}(\rmE)
\end{subarray}
}}
k_{P_{2,\rmE} \times \tlP_{1,\rmE}, \chi}( w_{s}^{-1}\gamma y, \tly,  x, w_{\tls}^{-1} \tlgam \tlx_{3} a_{\tltr}  k_{\tltr})
\end{gather*}
sont invariantes à droite par  
l'intersection des compacts ouverts
$\bigcap_{k \in K_{\rmE} \hrar (G \times \tlG_{\rmE})(\A)}(kK_{0}k^{-1})$, $\bigcap_{k \in K \hrar G_{\rmE}(\A)}(kK_{0}k^{-1})$ et 
$\bigcap_{\tlk \in \tlK \hrar \tlG_{\rmE}(\A)}(\tlk K_{0}\tlk^{-1})$
avec $(M_{1,\rmE} \times M_{\tlone,\rmE})(\A_{f})$, $M_{2,\rmE}(\A_{f})$ et $M_{\tltr,\rmE}(\A_{f})$ 
respectivement. On applique alors la proposition \ref{prop:opTronqD} pour l'opérateur 
$\Lad^{\tlP_{1},T}$ et la première fonction ci-dessus. 
On applique aussi 
la proposition \ref{prop:opTronqRel} 
pour l'opérateur $\Laz^{\tlP_{2}, T'}$ et la deuxième fonction.
On applique finalement la proposition \ref{prop:opTronqReltl} pour l'opérateur $\Lajlr^{\tlP_{3}, T''}$ et 
la troisième fonction. On trouve 
donc que pour tout $r_{1}, r_{2} \ge 0$ il existe un nombre fini des opérateurs différentiels $X, Y$
tels 
que l'expression (\ref{eq:AllSigmas}) est majorée 
par la somme sur les $X, Y$ de  
\begin{multline}\label{eq:withGreatLemComeGreatExpsSp}
\int\limits_{\calM_{\tlone, \rmE}}
\int\limits_{\calM_{\tltwo}}
\int\limits_{\calM_{\tltr}}
\quad \quad \quad
\sup_{\mathmakebox[0.8cm]{
\begin{subarray}{c}
m_{1}' \in H_{\tlone, \rmE}(\A)^{1} \times G_{\tlone,\rmE}(\A) \\
m_{\tlone}' \in H_{\tlone, \rmE}(\A)^{1} \times \tlG_{\tlone,\rmE}(\A) \\
m_{2}' \in H_{\tltwo, \rmE}(\A)^{1} \times G_{\tltwo, \rmE}(\A) \\
m_{\tltr}' \in M_{\tltr,\rmE}(\A)^{1}
\end{subarray}
}}
(
\|m_{1}'\|
\|m_{\tlone}'\|
\|m_{2}'\|
\|m_{\tltr}'\|)^{-r_{1}} 
\sum_{\chi}
\int\limits_{K_{\rmE}}
\int\limits_{K}
\int\limits_{\tlK}
\int\limits_{A_{\tlone}^{st,\infty}}
\int\limits_{A_{\tltwo}^{st,\infty}}
\int\limits_{A_{\tltr}^{\infty}} 
\sigma_{\tlone}^{\tlqt}(H_{\tlone}(a_{\tlone}^{st}m_{1})-T_{\tlone}')\\
\sigma_{\tltwo}^{\tlcq}(H_{\tltwo}(a_{\tltwo}^{st}m_{2})-T_{\tltwo})
\sigma_{\tltr}^{\tlsix}(H_{\tltr}(a_{\tltr})-T_{\tltr}'')
e^{-2(\rho_{1,\rmE}(H_{1}(a_{\tlone}^{st}m_{1})) + \rho_{2}(H_{2}(a_{\tltwo}^{st}m_{2})) +\rho_{\tltr}(H_{\tltr}(a_{\tltr})))}
\\ | \sum_{(s, \tls) \in \Omega'}
\int\limits_{[N_{1, \rmE}]}
\int\limits_{[N_{\tltr,\rmE}]} 
\sum_{\mathrlap{
\begin{subarray}{l}
\gamma \in (P_{1} \cap sP_{2})(\rmE) \bsl P_{1}(\rmE)\\
\tlgam \in (\tlP_{3} \cap \tls\tlP_{1})(\rmE) \bsl \tlP_{3}(\rmE)
\end{subarray}
}} \
k_{\Phi_{X,Y}, P_{2, \rmE} \times \tlP_{1,\rmE}, \chi}(w_{s}^{-1}\gamma n_{1} a_{\tlone}^{st} m_{1}'  k_{1}, 
a_{\tlone}^{st}m_{\tlone}'k_{1}, 
a_{\tltwo}^{st} m_{2}'  k_{2}, 
w_{\tls}^{-1} \tlgam n_{\tltr} a_{\tltr} m_{\tltr}' k_{\tltr}) \\
dn_{\tltr}dn_{1}
| 
|\det m_{1} a_{\tlone}^{st}|_{\A}^{\sigma}
|\det m_{2} a_{\tltwo}^{st}|_{\A}^{\sigma'}
da_{\tltr}da_{\tltwo}^{st}da_{\tlone}^{st}dk_{\tltr}dk_{2}dk_{1}
)(\|m_{1}\|\|m_{2}\|\|m_{\tltr}\|)^{-r_{2}}
dm_{\tltr}dm_{2}dm_{1}
\end{multline}
où $\calM_{\tlone,\rmE} = \Sgl^{M_{1,\rmE}}_{B_{\rmE} \cap M_{1,\rmE}} \cap (H_{\tlone, \rmE}(\A)^{1} \times G_{\tlone,\rmE}(\A))$,  
$\calM_{\tltwo} = \Sgl^{M_{2}}_{B\cap M_{2}} \cap (H_{\tltwo}(\A)^{1} \times G_{\tltwo}(\A))$ 
$\calM_{\tltr} = \Sgl^{M_{\tltr}}_{ \tlB \cap M_{\tltr}} \cap M_{\tltr}(\A)^{1}$ 
et l'on a remplacé 
la fonction $\Phi$, en passant par la formule (\ref{eq:diffOnKernel}) 
et en utilisant le même raisonnement que sur la page 104 de \cite{arthur5}, 
par une fonction $\Phi_{X,Y} := X \ast \Phi \ast Y$, qui dépend des opérateurs $X$ et $Y$.
 Remarquons que 
le support de $\Phi_{X,Y}$ est contenu dans celui de $\Phi$.

\blem\label{lem:preludeToGreatness}
 Avec la notation ci-dessus, il existe des constantes positives $c, N >0$ telles que
si
\begin{multline*}
\sigma_{\tlone}^{\tlqt}(H_{\tlone}(a_{\tlone}^{st}m_{1})-T_{\tlone}')
\sigma_{\tltwo}^{\tlcq}(H_{\tltwo}(a_{\tltwo}^{st}m_{2})-T_{\tltwo})
\sigma_{\tltr}^{\tlsix}(H_{\tltr}(a_{\tltr})-T_{\tltr}'')
|
\int\limits_{[N_{1, \rmE}]}
\int\limits_{[N_{\tltr,\rmE}]} 
\sum_{(s, \tls) \in \Upomega_{\tlQ_{1}}'} \\
\sum_{\mathrlap{
\begin{subarray}{l}
\gamma \in (P_{1} \cap sP_{2})(\rmE) \bsl P_{1}(\rmE)\\
\tlgam \in (\tlP_{3} \cap \tls\tlP_{1})(\rmE) \bsl \tlP_{3}(\rmE)
\end{subarray}
}} \
k_{\Phi_{X,Y}, P_{2,\rmE} \times \tlP_{1,\rmE}, \chi}(w_{s}^{-1}\gamma n_{1} a_{\tlone}^{st} m_{1}'  k_{1}, 
a_{\tlone}^{st}m_{\tlone}'k_{1}, 
a_{\tltwo}^{st} m_{2}'  k_{2}, 
w_{\tls}^{-1} \tlgam n_{\tltr} a_{\tltr} m_{\tltr}' k_{\tltr}) 
dn_{\tltr}dn_{1}|
\end{multline*}
est non-nulle, alors
\[
\|a_{\tlone}^{st}\|, \|a_{\tltwo}^{st}\|,  \|a_{\tltr}^{\tlG}\| \le c
(\|m_{1}\|  \|m_{2}\| \|m_{1}'\| \|m_{\tlone}'\|\|m_{2}'\|\|m_{\tltr}'\|)^{N}
\]
où $a_{\tltr}^{\tlG}$ c'est la projection de $a_{\tltr}$ à $A_{\tltr}^{\tlG, \infty}$.
\bdem
Si l'expression dans le lemme est non-nulle il existe des 
$(s, \tls) \in \Upomega_{\tlQ_{1}}'$, 
$\gamma \in (P_{1} \cap sP_{2})(\rmE) \bsl P_{1}(\rmE)$, 
$\tlgam \in (\tlP_{3} \cap \tls\tlP_{1})(\rmE) \bsl \tlP_{3}(\rmE)$, $n_{1} \in N_{1, \rmE}(\A)$ et 
$n_{\tltr} \in N_{\tltr, \rmE}(\A)$ tels que
\[
k_{\Phi_{X,Y}, P_{2,\rmE} \times \tlP_{1,\rmE}, \chi}(w_{s}^{-1}\gamma n_{1} a_{\tlone}^{st} m_{1}'  k_{1}, 
a_{\tlone}^{st}m_{\tlone}'k_{1}, 
a_{\tltwo}^{st} m_{2}'  k_{2}, 
w_{\tls}^{-1} \tlgam n_{\tltr} a_{\tltr} m_{\tltr}' k_{\tltr}) \neq 0.
\]
Il résulte alors du lemme \ref{lem:artLemme23} qu'il existe des $m_{\tlone}'' \in M_{\tlone,\rmE}(\A)^{1}$ et 
$m_{2}'' \in M_{2,\rmE}(\A)^{1}$ tels que 
\begin{equation}\label{eq:kSanschinonnul}
k_{\Phi_{X,Y}, P_{2,\rmE} \times \tlP_{1,\rmE}}(w_{s}^{-1}\gamma n_{1} a_{\tlone}^{st} m_{1}'  k_{1}, 
a_{\tlone}^{st}m_{\tlone}'' m_{\tlone}' k_{1}, 
a_{\tltwo}^{st} m_{2}'' m_{2}' k_{2}, 
w_{\tls}^{-1} \tlgam n_{\tltr} a_{\tltr} m_{\tltr}' k_{\tltr}) \neq 0.
\end{equation}
Soit $\tlB \in \calP(M_{\tlzero})$ tel que $\tlB \sbs \tlP_{1}$, que l'on a choisi ci-dessus.
Décomposons $m_{1}'= b_{1}'k$ et $m_{\tltr}' = \tlb_{3}' \tlk$ selon les décompositions d'Iwasawa 
$M_{1,\rmE}(\A) = (B_{\rmE} \cap M_{1,\rmE})(\A) (K_{\rmE} \cap M_{1,\rmE}(\A))$ et
$M_{\tltr,\rmE}(\A) = (\tlB_{\rmE} \cap M_{\tltr,\rmE})(\A) (\tlK_{\rmE} \cap M_{\tltr,\rmE}(\A))$. Décomposons aussi 
$\gamma^{-1}w_{s} = nw_{s'}b$ selon la décomposition de Bruhat, avec $n \in N_{B}(\rmE)$, $s' \in \Omega^{Q_{2}}$ et 
$b \in B(\rmE)$. Notons que $s = s'$ dans $\Omega^{1} \bsl \Omega^{Q_{2}} / \Omega^{2}$ on peut supposer donc que 
$s'= s$. 
Alors, si l'on a (\ref{eq:kSanschinonnul}), on voit qu'il existe 
des $n' \in N_{B,\rmE}(\A)$ et $\xi \in M_{2}(\rmE)N_{2,\rmE}(\A)$ tels que l'expression
\begin{equation}\label{eq:GComp}
(b_{1}')^{-1}(a_{\tlone}^{st})^{-1}n'w_{s} \xi  a_{\tltwo}^{st} m_{2}'' m_{2}' 
\end{equation}
appartient à un compact de $G_{\rmE}(\A)$ qui ne dépend que du support de $\Phi$. 
De même, il existe des $\tln \in N_{\tlB,\rmE}(\A)$ et $\tlxi \in M_{\tlone}(\rmE)N_{\tlone,\rmE}(\A)$ tels que l'expression
\begin{equation}\label{eq:tlGComp}
(\tlb_{3}')^{-1}a_{\tltr}^{-1}\tln w_{\tls} \tlxi a_{\tlone}^{st}m_{\tlone}'' m_{\tlone}'
\end{equation}
appartient dans un compact de $\tlG_{\rmE}(\A)$ qui ne dépend que du support de $\Phi$. 

En utilisant les représentations du plus haut poids $d \varpi$ où $d \in \N$ et $\varpi \in \hDelta_{\tltwo}$ et 
en regardant l'action de l'élément (\ref{eq:GComp}) par des telles représentations, 
comme le fait Arthur dans \cite{arthur5} pages 102-103, on voit qu'il existe une constante $C >0$ telle que pour tout 
$\varpi \in \hDelta_{\tltwo}$ on a
\[
\varpi(H_{\tlB}(a_{\tltwo}^{st})) - s\varpi(H_{\tlB}(a_{\tlone}^{st})) \le C(1 + |\varpi(H_{\tlB}(m_{2}'))| + |s\varpi(H_{\tlB}(b_{1}') )| )
\]
En utilisant alors les propriétés (\ref{eq:haut2}), (\ref{eq:haut6}) et (\ref{eq:haut7}) de la hauteur 
on trouve
\begin{equation}\label{eq:varpiMinsvarpiThm1}
\varpi(H_{\tlB}(a_{\tltwo}^{st})) - s\varpi(H_{\tlB}(a_{\tlone}^{st})) \le C'(1 + \log (\|m_{2}'\| \|m_{1}'\|) ), \quad 
\forall \varpi \in \hDelta_{\tltwo}
\end{equation}
pour un $C' >0$. Par le même raisonnement appliqué à l'expression (\ref{eq:tlGComp}) on trouve
\begin{equation}\label{eq:varpiMinsvarpiThm2}
\varpi(H_{\tlB}(a_{\tlone}^{st})) - \tls\varpi(H_{\tlB}(a_{\tltr})) \le C''(1 + \log (\|m_{\tlone}'\| \|m_{\tltr}'\|) ), \quad 
\forall \varpi \in \hDelta_{\tlone}
\end{equation}
pour un $C'' >0$. 

D'autre côté, en vertu de l'égalité (\ref{eq:alternSumExpl}) ci-dessus, l'expression entre la valeur absolue dans le lemme égale simplement
\begin{multline}\label{eq:altSumReborn}
\int\limits_{[N_{1, \rmE}]}
\int\limits_{[N_{\tltr,\rmE}]} 
\int\limits_{[N_{\tlone, \rmE}]}
\int\limits_{[N_{2,\rmE}]}  \\
\sum_{\tlQ_{1} \sbs \tlP \sbs \tlQ_{2}}
(-1)^{d_{\tlP}}
k_{\Phi_{X,Y},  P_{\rmE} \times \tlP_{\rmE}, \chi}(n_{1} a_{\tlone}^{st} m_{1}'  k_{1}, 
n_{\tlone} a_{\tlone}^{st}m_{\tlone}'k_{1}, 
n_{2} a_{\tltwo}^{st} m_{2}'  k_{2}, 
n_{\tltr} a_{\tltr} m_{\tltr}' k_{\tltr}) 
dn_{2}dn_{\tlone}dn_{\tltr}dn_{1}
\end{multline}
à un signe près. 
En vertu du corollaire \ref{cor:art23REALLY} de nouveau, on a 
pour tout $\tlP$ entre $\tlQ_{1}$ et $\tlQ_{2}$:
\begin{multline*}
\int\limits_{[N_{1, \rmE}]}
\int\limits_{[N_{\tltr,\rmE}]} 
k_{\Phi_{X,Y},  P_{\rmE} \times \tlP_{\rmE}, \chi}(n_{1} a_{\tlone}^{st} m_{1}'  k_{1}, 
n_{\tlone} a_{\tlone}^{st}m_{\tlone}'k_{1}, 
n_{2} a_{\tltwo}^{st} m_{2}'  k_{2}, 
n_{\tltr} a_{\tltr} m_{\tltr}' k_{\tltr})dn_{1}dn_{\tltr} = \\
\sum_{\mathclap{
\begin{subarray}{c}
s' \in \Omega^{2} \bsl \Omega^{P} / \Omega^{1} \\
\tls' \in \Omega^{\tlone} \bsl \Omega^{\tlP} / \Omega^{\tltr}
\end{subarray}
}} \quad \quad
\sum_{\mathrlap{
\begin{subarray}{l}
\gamma \in (P_{2} \cap s'P_{1})(\rmE) \bsl P_{2}(\rmE)\\
\tlgam \in (\tlP_{1} \cap \tls'\tlP_{3})(\rmE) \bsl \tlP_{1}(\rmE)
\end{subarray}
}} \
k_{\Phi_{X,Y}, P_{1,\rmE} \times \tlP_{3,\rmE}, \chi}(a_{\tlone}^{st} m_{1}'  k_{1}, 
w_{\tls'}^{-1}\tlgam n_{\tlone} a_{\tlone}^{st}m_{\tlone}'k_{1}, 
w_{s'}^{-1} \gamma n_{2} a_{\tltwo}^{st} m_{2}'  k_{2}, 
a_{\tltr} m_{\tltr}' k_{\tltr}).
\end{multline*}
Il s'ensuit que l'expression (\ref{eq:altSumReborn}) égale
\[
\int\limits_{[N_{\tlone, \rmE}]}
\int\limits_{[N_{2,\rmE}]} 
\sum_{(s',\tls') \in \Upomega_{\tlQ_{1}}''}
\sum_{\mathrlap{
\begin{subarray}{l}
\gamma \in (P_{2} \cap s'P_{1})(\rmE) \bsl P_{2}(\rmE)\\
\tlgam \in (\tlP_{1} \cap \tls'\tlP_{3})(\rmE) \bsl \tlP_{1}(\rmE)
\end{subarray}
}} \
k_{\Phi_{X,Y}, P_{1,\rmE} \times \tlP_{3,\rmE}, \chi}(a_{\tlone}^{st} m_{1}'  k_{1}, 
w_{\tls'}^{-1}\tlgam n_{\tlone} a_{\tlone}^{st}m_{\tlone}'k_{1}, 
w_{s'}^{-1} \gamma n_{2} a_{\tltwo}^{st} m_{2}'  k_{2}, 
a_{\tltr} m_{\tltr}' k_{\tltr})dn_{2}dn_{\tlone}
\]
où
\[
\Upomega_{\tlQ_{1}}'' = \{(s', \tls') \in (\Omega^{2} \bsl \Omega^{Q_{2}}/\Omega^{1}) \times 
(\Omega^{\tlone} \bsl \Omega^{\tlQ_{1}} / \Omega^{\tltr})| \hDelta_{\tlQ_{2}} = \hDelta_{\tlQ_{1}, s'}\cap \hDelta_{\tlQ_{1}, \tls'}\}.
\]

Donc, si l'expression du lemme est non-nulle,
on obtient, en utilisant le lemme \ref{lem:artLemme23} 
qu'il existe des $(s',\tls') \in \Upomega_{\tlQ_{1}}''$, $\gamma \in (s'P_{1} \cap P_{2})(\rmE) \bsl P_{2}(\rmE)$,
$\tlgam \in (\tls'\tlP_{3} \cap \tlP_{1})(\rmE) \bsl \tlP_{1}(\rmE)$, $n_{2} \in N_{2, \rmE}(\A)$, 
$n_{\tlone} \in N_{\tlone, \rmE}(\A)$ 
$m_{1}'' \in M_{1,\rmE}(\A)^{1}$ et 
$m_{\tltr}'' \in M_{\tltr, \rmE}(\A)^{1}$ tels que
\begin{equation*}
k_{\Phi_{X,Y}, P_{1,\rmE} \times \tlP_{3,\rmE}}(a_{\tlone}^{st} m_{1}'' m_{1}' k_{1}, 
w_{\tls'}^{-1} \tlgam n_{\tlone}  a_{\tlone}^{st} m_{\tlone}'k_{1}, 
w_{s'}^{-1} \gamma n_{2} a_{\tltwo}^{st} m_{2}'  k_{2}, 
a_{\tltr} m_{\tltr}'' k_{\tltr}) \neq 0.
\end{equation*}
En raisonnant comme ci-dessus on trouve une constante $C''' >0$ telle que
\begin{gather}
\varpi(H_{\tlB}(a_{\tlone}^{st})) - s'\varpi(H_{\tlB}(a_{\tltwo}^{st})) \le C'''(1 + \log (\|m_{2}'\| \|m_{1}'\|) ), \quad 
\forall \varpi \in \hDelta_{\tlone}, \label{eq:varpiMinsvarpiThm3}\\
\varpi(H_{\tlB}(a_{\tltr}^{st})) - \tls'\varpi(H_{\tlB}(a_{\tlone})) \le C'''(1 + \log (\|m_{\tlone}'\|) ), \quad 
\forall \varpi \in \hDelta_{\tltr}. \label{eq:varpiMinsvarpiThm4}
\end{gather}

On est en mesure d'appliquer le lemme \ref{lem:fuckYeahLemme} de l'appendice \ref{app:lemme}. On l'applique pour les données suivantes.
Pour le groupe $G$ on prend $\tlG$. 
Pour le sous-groupe parabolique minimal qu'on fixe au début de l'appendice on prend $\tlB$ fixé ci-dessus, 
contenu dans $\tlP_{1} \cap \tlP_{3}$.
Les sous-groupes $\tlP_{1}$, $\tlP_{2}$, $ \tlP_{3}$, $\tlP_{4}$, $ \tlP_{5}$ et $ \tlP_{6}$ 
correspondent donc aux $P_{1}$, $ P_{2}$, $ P_{3}$, $ P_{4}$, $P_{5}$ et $P_{6}$ du lemme. 
Pour les éléments $H_{1}$, $H_{2}$ et $H_{3}$ on prend les projections à $\all_{\tlone}^{\tlG}$, 
$\all_{\tltwo}^{\tlG}$ et $\all_{\tltr}^{\tlG}$ des $H_{\tlB}(\all_{\tlone}^{st})$, 
 $H_{\tlB}(\all_{\tltwo}^{st})$ et  $H_{\tlB}(\all_{\tltr})$ respectivement.
 Pour $X_{1}$, $X_{2}$ et $X_{3}$ on prend $T_{\tlone} - H_{\tlB}(m_{1})$, 
 $T_{\tltwo} - H_{\tlB}(m_{2})$ et $T_{\tltr}$ respectivement. 
 Pour les éléments du groupe de Weyl on prend 
 $s_{1} = s'$, $s_{1}' = s$, $s_{2} = \tls'$ et $s_{2}' = \tls$. 
 Les conditions $(s,\tls) \in \Upomega_{\tlQ_{1}}'$ 
et $(s',\tls') \in \Upomega_{\tlQ_{1}}''$ 
 disent qu'ils vérifient la condition (\ref{eq:hDeltaEses}) du lemme \ref{lem:fuckYeahLemme}.
Les inégalités du lemme \ref{lem:fuckYeahLemme} 
 correspondent aux inégalités qu'on a construit de façon suivant:
 \[
\text{(\ref{eq:varpiMinsvarpi1})} \leftrightarrow  
\text{(\ref{eq:varpiMinsvarpiThm3})}, \quad 
\text{(\ref{eq:varpiMinsvarpi2})} \leftrightarrow  
\text{(\ref{eq:varpiMinsvarpiThm1})}, \quad 
\text{(\ref{eq:varpiMinsvarpi3})} \leftrightarrow  
\text{(\ref{eq:varpiMinsvarpiThm2})}, \quad 
\text{(\ref{eq:varpiMinsvarpi4})} \leftrightarrow  
\text{(\ref{eq:varpiMinsvarpiThm4})}. 
 \] 
Les inégalités déterminent les constantes 
 $M_{1}$, $M_{2}$, $M_{3}$ et $M_{4}$. 
 
 Il nous reste à justifier que la condition (\ref{eq:hDeltaEses2}) soit vérifiée. 
Ceci découle du lemme suivant. 
 \blem\label{lem:relStnotGconj}
  Soient $s \in \Omega^{G}$ et $\tlB, \tlB' \in \calP(M_{\tlzero}) \cap \relPb$. Alors
\begin{equation*}
s(\hDelta_{\tlB} \smin \hDelta_{\tlB'}) \cap \hDelta_{\tlB'} = \varnothing.
\end{equation*}
 \bdem
Raisonnons par absurde. Soient $\tlB$, $\tlB'$ et $s \in \Omega^{G}$ comme ci-dessus 
et soit $\varpi \in \hDelta_{\tlB} \smin \hDelta_{\tlB'}$ tel que $s \varpi \in \hDelta_{\tlB'}$. 
Soit $\tlP \sps \tlB$ associé à $\{ \varpi\} \sbs \hDelta_{\tlB}$ et 
soit $\tlP' \sps \tlB'$ associé à $\{ s\varpi\} \sbs \hDelta_{\tlB'}$. Alors 
$\tlP, \tlP' \in \relPb$ et $s\tlP = \tlP'$. 
Soit $W_{0} \sbs W$ le sous-espace de $W$ tel que $\tlP$ est défini comme le stabilisateur de $W_{0}$. 
Donc, $\tlP'$ est le stabilisateur de $s W_{0}$. Deux cas sont possibles. Soit $W_{0} = V_{0} \sbs V$ 
soit $W_{0} = V_{0} \oplus D_{0} \sbs W$ où $V_{0} \sbs V$ et $D_{0}$ c'est la droite telle  que $W = V \oplus D_{0}$.
Dans le deux cas, le groupe $P = \tlP \cap G$ est le stabilisateur de $V_{0}$ et le groupe $P' = \tlP \cap G$
est le stabilisateur de $sV_{0}$. Les sous-groupes paraboliques $P$ et $P'$ sont standards 
(contiennent $B$) et sont définies comme les stabilisateurs d'un sous-espace de $V$ de dimension $\dim V_{0} = \dim s V_{0}$. 
Ils sont donc égaux. Cela démontre que $sV_{0} = V_{0}$ et, puisque $sV = V$ et $sD_{0} = D_{0}$, par conséquent $sW_{0} = W_{0}$, d'où 
$\tlP = \tlP'$ ce qui est absurde. 
 \edem
 \elem
 
Puisque $s = s_{1}'$ appartient à $\Omega^{G}$, en vertu du lemme \ref{lem:relStnotGconj} ci-desuss, 
il vérifie la propriété 
(\ref{eq:hDeltaEses2}) du lemme \ref{lem:fuckYeahLemme}. Alors, grâce au lemme \ref{lem:fuckYeahLemme}, la propriété (\ref{eq:haut7}) de la hauteur et le lemme \ref{lem:hautComp} on conclut la preuve 
 du lemme \ref{lem:preludeToGreatness}.
 \edem
\elem

Revenons à l'expression (\ref{eq:withGreatLemComeGreatExpsSp}). 
On traite toujours les variables $n_{1}$, $a_{\tlone}^{st}$, $m_{1}$ etc. comme fixées. 
On suppose que les variables $a_{\tlone}^{st}$, $a_{\tltwo}^{st}$ 
et $a_{\tltr}^{\tlG}$ vérifient les conditions du lemme \ref{lem:preludeToGreatness}. 
Posons $g = n_{1} a_{\tlone}^{st} m_{1}'  k_{1}$, $\tlg = a_{\tlone}^{st}m_{\tlone}'k_{1}$, 
$h = a_{\tltwo}^{st} m_{2}'  k_{2}$ et $\tlh = n_{\tltr} a_{\tltr}^{\tlG} m_{\tltr}' k_{\tltr}$,
et notons aussi $\tlg'$ la projection de $\tlg$ 
à $\tlG_{\tlE}(\A)^{1}$. 
Posons pour $x \in G_{\rmE}(\A)$ et $\tlx \tlG_{\rmE}(\A)^{1}$
\[
\Uppsi(x, \tlx) = \int_{A_{\tlG}^{\infty}}k_{\Phi_{X,Y}, P_{2,\rmE} \times \tlP_{1,\rmE}, \chi}(x, \tlg, h, \tlx a_{\tlG})da_{\tlG}.
\]
Par le même argument qu'au début de la preuve on a 
\[
\Uppsi(x, \tlx) = \Uppsi(x, \tlx)\htau_{2}(H_{2}(x) - H_{2}(h) - T_{2, \Phi})
\htau_{\tlone}(H_{\tlone}(\tlx) - H_{\tlone}(\tlg) - T_{\tlone, \Phi})
\]
pour certains $T_{2, \Phi} \in \all_{2}$, $T_{\tlone, \Phi} \in \all_{\tlone}$ qui ne dépendent 
que du support de $\Phi$. 
En faisant un changement de variables on trouve:
\begin{multline*}
\int\limits_{A_{\tlG}^{\infty}} 
\sum_{(s, \tls) \in \Omega'}
|\sum_{\mathrlap{
\begin{subarray}{c}
\gamma \in (P_{1} \cap sP_{2})(\rmE) \bsl P_{1}(\rmE)\\
\tlgam \in (\tlP_{3} \cap \tls\tlP_{1})(\rmE) \bsl \tlP_{3}(\rmE)
\end{subarray}}}
k_{\Phi_{X,Y}, P_{2,\rmE} \times \tlP_{1,\rmE}, \chi}(w_{s}^{-1}\gamma n_{1} a_{\tlone}^{st} m_{1}'  k_{1}, 
a_{\tlone}^{st}m_{\tlone}'k_{1}, 
a_{\tltwo}^{st} m_{2}'  k_{2}, 
w_{\tls}^{-1} \tlgam n_{\tltr} a_{\tltr}^{\tlG}a_{\tlG} m_{\tltr}' k_{\tltr}) | da_{\tlG} \\
\le
\sum_{\mathclap{
\begin{subarray}{c}
\gamma \in P_{2}(\rmE) \bsl G(\rmE)\\
\tlgam \in \tlP_{1}(\rmE) \bsl \tlG(\rmE)
\end{subarray}}}
\Uppsi(\gamma g, \tlgam \tlh)
\htau_{2}(H_{2}(\gamma g) - H_{2}(h) - T_{2, \Phi})
\htau_{\tlone}(H_{\tlone}(\tlgam \tlh) - H_{\tlone}(\tlg) - T_{\tlone, \Phi}).
\end{multline*}
En vertu du lemme  \ref{lem:lemmeArth51} \textit{iii)} et de la propriété 
(\ref{eq:haut6})
pour tout $N_{0}, N_{1} > 0$, il existe 
des constantes positives $c'$, $N'$ telles que l'expression ci-dessus est majorée par 
\[
c' (\|g\| \|h\| \|\tlg\| \|\tlh\|)^{N'} 
\sup_{(x, \tlx) \in G_{\rmE}(\A) \times \tlG_{\rmE}(\A)^{1}}
(\Uppsi( x, \tlx) \|x\|^{-N_{0}} \|\tlx\|^{-N_{1}}).
\]
On prend donc $N_{0}$ et $N_{1}$ comme dans le lemme \ref{lem:arthursCor46} et on trouve que l'expression 
(\ref{eq:withGreatLemComeGreatExpsSp}) est majorée par
\begin{multline*}
\int\limits_{\calM_{\tlone, \rmE}}
\int\limits_{\calM_{\tltwo}}
\int\limits_{\calM_{\tltr}}
\sup_{
\begin{subarray}{c}
m_{1}' \in H_{\tlone, \rmE}(\A)^{1} \times G_{\tlone,\rmE}(\A) \\
m_{\tlone}' \in H_{\tlone, \rmE}(\A)^{1} \times \tlG_{\tlone,\rmE}(\A) \\
m_{2}' \in H_{\tltwo, \rmE}(\A)^{1} \times G_{\tltwo, \rmE}(\A) \\
m_{\tltr}' \in M_{\tltr,\rmE}(\A)^{1}
\end{subarray}
}
(
(\|m_{1}'\|
\|m_{\tlone}'\|
\|m_{2}'\|
\|m_{\tltr}'\|)^{-r_{1}} 
\int\limits_{A_{\tlone}^{st,\infty}}
\int\limits_{A_{\tltwo}^{st,\infty}}
\int\limits_{A_{\tltr}^{\tlG, \infty}}
(\|a_{\tlone}^{st}\|\|a_{\tltwo}^{st}\|\|a_{\tltr}^{\tlG}\|)^{N''} \cdot \\
(\|m_{1}\|\|m_{2}\|\|m_{\tltr}\|
\|m_{1}'\|
\|m_{\tlone}'\|
\|m_{2}'\|
\|m_{\tltr}'\|
)^{N''} 
da_{\tltr}^{\tlG}da_{\tltwo}^{st}da_{\tlone}^{st}
)(\|m_{1}\|\|m_{2}\|\|m_{\tltr}\|)^{-r_{2}}dm_{\tltr}dm_{2}dm_{1}
\end{multline*}
où on a utilise une majoration de type
\[
e^{-2(\rho_{1,\rmE}(H_{1}(a_{\tlone}^{st}m_{1})) + \rho_{2}(H_{2}(a_{\tltwo}^{st}m_{2})) + \rho_{\tltr}(H_{\tltr}(a_{\tltr}))} 
|\det m_{1} a_{\tlone}^{st}|_{\A}^{\sigma}
|\det m_{2} a_{\tltwo}^{st}|_{\A}^{\sigma'}\le 
(\|a_{\tlone}^{st}\|\|a_{\tltwo}^{st}\|\|a_{\tltr}^{\tlG}\|\|m_{1}\|\|m_{2}\|)^{N'''}
\]
et la constante $N''$ ne dépend que de cette constante $N'''$, des constantes fixés $N_{0}$, $N_{1}$ du lemme \ref{lem:arthursCor46}, 
des $t$ et $t'$ de la propriété (\ref{eq:haut6}) de la hauteur et des $t_{1}$, $t_2$ du lemme \ref{lem:hautComp}. 
En prenant alors $r_{1}$ et $r_{2}$ dans le lemme \ref{lem:preludeToGreatness} 
suffisamment grands on obtient la convergence de l'intégrale du théorème.
\edem
\end{theo}

\subsection{Polynômes-exponentielles}\label{par:fonsPolExpSp}
Soit $\calV$ un $\R$-espace vectoriel de dimension finie. 
Par un polynôme-exponentielle sur $\calV$ on entend une fonction 
sur $\calV$ de la forme
\[
f(v) = \sum_{\la \in \calV^{*}}e^{\la(v)}P_{\la}(v), \quad v \in \calV
\]
où
$P_{\la}$ est un polynôme sur $\calV$ 
à coefficients complexes, nul 
pour presque tout $\la \in \calV^{*}$. 
On appelle $\la \in \calV^{*}$ tels que $P_{\la} \neq 0$ les exposants 
de $f$ et le polynôme correspondant à $\la = 0$ le terme 
purement polynomial de $f$.
On a alors le résultat d'unicité suivant: 
si $f$ est comme ci-dessus 
et $g = \sum_{\la \in \calV^{*}}e^{\la}Q_{\la}$ est un 
polynôme-exponentielle sur $\calV$ tel que 
$g(v) = f(v)$ pour tout $v \in \calV$ alors 
pour tout $\la \in \calV^{*}$ on a $P_{\la} = Q_{\la}$.

Pour tout $s \in \C$ 
et tout sous-groupe parabolique relativement standard $\tlQ$ de $\tlG$, avec la notation 
du paragraphe \ref{par:glnGlnplus1Sp},
prenons $s\det - 2\rho_{Q}$ vu comme un élément 
de $(\all_{\tlQ, \C}^{st})^{*}$ et définissons $\upla_{\tlQ,s} \in \all_{\tlQ,\C}^{*}$ 
comme $\iota_{\tlQ}^{st}(s\det - 2\rho_{Q}) + 
2\rho_{\tlQ}$. 
Notons aussi $s_{\tlQ} \in \C$ défini par la propriété:
\begin{equation}\label{eq:sQSp}
e^{- \upla_{\tlQ,s}(H_{\tlQ}(m))} = |\det m|_{\A}^{s_{\tlQ}} \quad \forall \
m \in H_{\tlQ, \rmE}(\A)^{1} \times G_{\tlQ, \rmE}(\A).
\end{equation}
\blem[cf. lemme 3.2 de \cite{leMoi2}]\label{lem:uplasNonNulG}
 Soit $\tlQ$ un sous-groupe parabolique standard 
 de $\tlG$. Pour tout $s \in \C \smin \{-1,1\}$ et tout 
 $\tlvpi^{\vee} \in \hDelta_{\tlQ}^{\vee}$ 
on a $\upla_{\tlQ,s}(\tlvpi^{\vee}) \neq 0$.
\elem

On rappelle la définition de la fonction $\Gamma_{Q}'$ donnée 
par (\ref{eq:GammaQDefSp}) dans le paragraphe \ref{par:prelimstraceSp}, 
ainsi que la définition de la fonction $\hat \theta_{Q}$ 
donnée par (\ref{eq:thetaHatDefSp}) dans le paragraphe \ref{par:haarMes}. 
Le polynôme-exponentielle en question dans le cas des groupes linéaires est donné par le 
lemme suivant:
\blem\label{lem:pQExplicitGSp}
 Soit $\tlQ \in \relPb$.
 \begin{enumerate}[1)]
 \item Pour tout $s \in \C$ 
\[
p_{\tlQ,s}(X) := \int_{A_{\tlQ}^{st,\infty}}e^{(s\det + 2\rho_{\tlQ} -2\rho_{Q})(H_{\tlQ}(a))}\Gamma_{\tlQ}'(H_{\tlQ}(a) ,X)da, \quad 
X \in \all_{\tlQ}
\]
est un polynôme-exponentielle sur $\all_{\tlQ}/\all_{\tlG}$. 
\item Pour tout $m \in H_{\tlQ, \rmE}(\A)^{1} \times G_{\tlQ, \rmE}(\A)$ et tout $T \in \all_{\tlQ}$ on a
\[
\int_{A_{\tlQ}^{st,\infty}}e^{(s\det + 2\rho_{\tlQ} -2\rho_{Q})(H_{\tlQ}(a))}\Gamma_{\tlQ}'(H_{\tlQ}(a) + H_{Q}(m) - T_{\tlQ}, X)da = 
|\det m|_{\A}^{s_{\tlQ}}
e^{\upla_{\tlQ,s}(T_{\tlQ})}
p_{\tlQ,s}(X).
\]
\item Si $s \neq -1,1$, 
pour tout $\tlR \sps \tlQ$ 
il existe un polynôme $P_{\tlQ,\tlR,s}$ de degré au plus $d_{\tlQ}^{\tlG}$ sur $\all_{\tlR}/\all_{\tlG}$ 
tel que
\[
p_{\tlQ,s}(X) := j_{\tlQ}^{-1}\sum_{\tlR \sps \tlQ}e^{\upla_{\tlR,s}(X_{\tlR})}
p_{\tlQ,\tlR,s}(X_{\tlR})
\]
où $p_{\tlQ,\tlG,s}(X_{\tlG}) = (-1)^{d_{\tlQ}^{\tlG}}\hat \theta_{\tlQ}(\upla_{\tlQ,s})^{-1}$ et 
$j_{\tlQ}$ est défini par la propriété \ref{eq:HQtlQSp} dans le paragraphe \ref{par:glnGlnplus1Sp}.
En particulier, si $s \neq -1,1$, la fonction $p_{\tlQ,s}$ est un polynôme-exponentielle
dont le terme purement polynomial est constant 
et égale $(-1)^{d_{\tlQ}^{\tlG}}\hat \theta_{Q}(\upla_{\tlQ,s})^{-1}$.
 \end{enumerate}
\elem

\subsection{Une généralisation du théorème \ref{thm:noyauCVGG}}\label{par:JTforLevisGln}

Soit $\eta : \rmF^{*} \bsl \A^{*} \rar \C$ le caractère quadratique associé 
à l'extension $\rmE/\rmF$ par la théorie de corps de classes. 
On définit le caractère, noté aussi $\eta$, sur $G(\A) \times \tlG(\A)$ 
comme $(h,\tlh) \mapsto \eta(\det (\tlh)^{n} \det (h)^{n+1})$.

Soient $T \in T_{+} + \all_{\tlzero}^{+}$, $s,s' \in \C$ et $\chi \in \calX^{G_{\rmE} \times \tlG_{\rmE}}$.
En vertu du théorème \ref{thm:noyauCVGG}
la distribution $I_{\chi}^{T}(s,s', \cdot)$ qui à $\Phi \in C_{c}^{\infty}(\GetlGa)$ associe
\[
I_{\chi}^{T}(s,s', \Phi) := \int_{[G_{\rmE}]}\int_{[G]}\int_{[\tlG]}
k_{\Phi,\chi}^{T}(g,g,h,\tlh)|\det g|^{s}|\det h|^{s'}\eta(h,\tlh)dxdhd\tlh
\]
est bien définie. 

Soient $W'$ un $\rmF$-espace vectoriel de dimension $m+1$, 
$V' \sbs W'$ un sous-espace de dimension $m$ et $D_{0}' \sbs (W' \smin V') \cup \{0\}$  
une droite, où $m \in \N$.
Soient $H = \prod_{i=1}^{k}\Gl_{n_{i}}$ 
$G' = \Gl(V')$ et $\tlG' = \Gl(W')$ où $k \in \N$, et 
$n_{i} \in \N$ pour $i = 1,\ldots, k$.
On identifie $G'$ avec le sous-groupe de $\tlG'$ qui agit trivialement 
sur $D'_{0}$.
Soient aussi $V'_{\rmE} = V' \otimes_{\rmF} \rmE$ et 
$W'_{\rmE} = W \otimes_{\rmF} \rmE$ des $\rmE$-espaces vectoriels et notons 
$G'_{\rmE} = \Res_{\rmE/\rmF}(\Gl(V'_{\rmE}))$ et $\tlG'_{\rmE} = \Res_{\rmE/\rmF}(\Gl(W'_{\rmE}))$. 
Notons aussi $H_{\rmE} =  \prod_{i=1}^{k}\Res_{\rmE/\rmF}(\Gl_{n_{i}})$. 
Les groupes $H \times G'$, $H \times \tlG'$, $H_{\rmE} \times G'_{\rmE}$ sont alors naturellement 
plongés dans $H_{\rmE} \times \tlG_{\rmE}'$. 
On va généraliser le théorème \ref{thm:noyauCVGG} au cas de l'action naturelle
\[
\Delta (H_{\rmE} \times G'_{\rmE}) \bsl ((H_{\rmE} \times G'_{\rmE}) \times (H_{\rmE} \times \tlG'_{\rmE})) / (H \times G') \times (H \times \tlG').
\]

Soit $B'$ un sous-$\rmF$-groupe de Borel
de $H \times G'$ et fixons aussi $M_{0}$ une partie de Levi de $B'$. 
Soit $M_{\tlzero}$ l'unique sous-groupe de Levi minimal de $H \times \tlG'$ 
tel que $M_{\tlzero} \sps M_{0}$.
On a en particulier l'ensemble $\calF(B')$ de sous-groupes paraboliques 
de $H \times G'$ contenant $B'$ ainsi que $\calF(M_{\tlzero})$ l'ensemble des sous-groupes paraboliques semi-standards de 
$H \times \tlG'$.
Notons $\relPr$ l'ensemble de $\tlP \in \calF(M_{\tlzero})$ 
tels que $\tlP \sps B'$. Pour $\tlP \in \relPr$ on note $P  = (H \times G') \cap \tlP$. 
Soit $M_{0,\rmE}$ (resp. $M_{\tlzero,\rmE}$) l'unique sous-$\rmF$-groupe de Levi minimal 
de $H_{\rmE} \times G'_{\rmE}$ (resp. de $H_{\rmE} \times \tlG'_{\rmE}$) contenant $M_{0}$ 
(resp. $M_{\tlzero}$). Pour tout $P \in \calF(B')$ (resp. $\tlP \in \calF(M_{\tlzero})$) il existe un unique élément de 
$\calF(M_{0, \rmE})$ \resp{$\calF(M_{\tlzero})$}, noté $P_{\rmE}$ \resp{$P_{\rmE}$}, dont l'intersection 
avec $H \times G'$ \resp{$H \times \tlG'$} égale $P$ \resp{$\tlP$}. 
La discussion et conventions du paragraphe \ref{par:GinGE} s'appliquent dans ce cas et va les utiliser sans commentaire. 

Fixons un sous-groupe de Borel $\tlBmin \in \calF(M_{\tlzero})$. 
Soit $\tlP \in \calF(M_{\tlzero})$. Pour tout 
$H \in \all_{\tlBmin}$ on note $H_{\tlP}$ la projection de 
$sH$ à $\all_{\tlP}$ où $s$ est un élément du groupe de Weyl de 
$H \times \tlG'$ tel que $s^{-1}\tlP \sps \tlBmin$.

Pour une fonction 
$f \in C_{c}^{\infty}(\hGpa)$,
un $\tlP \in \relPr$
et une donnée cuspidale $\chi \in \calX^{\hGp }$ on pose
\[
k_{f,\tlP, \chi}(\tlx, \tly) := k_{f, P_{\rmE} \times \tlP_{\rmE}, \chi}(\tlx, \tly),
\]
où $\tlx, \tly \in N_{P_{\rmE} \times \tlP_{\rmE}}(\A) M_{P \times \tlP}(\rmE)\bsl \hGpa$.
Pour un $T \in \mathfrak{a}_{\tlBmin}^{+}$ on pose donc
\[
k^{T}_{f,\chi}(x,\tlx, y, \tly) = 
\sum_{\tlP \in \relPr} (-1)^{d_{\tlP}^{H \times \tlG'}} \quad
\sum_{\mathclap{
\begin{subarray}{c}
\delta_{1} \in P(\rmE)\bsl (H \times G')(\rmE)
\\
\delta_{2} \in P(\rmF)\bsl (H \times G')(\rmF)
\\
\delta_{3} \in \tlP(\rmF)\bsl (H \times \tlG')(\rmF)
\end{subarray}}} \quad
\htau_{\tlP}^{H \times \tlG'}(H_{\tlP}(\delta_{2}y)-T_{\tlP})
k_{f, \tlP,\chi}(\delta_{1} x, \delta_{1} \tlx, \delta_{2} y, \delta_{3} \tly),
\]
où $x \in H_{\rmE}(\A)^{1} \times G'_{\rmE}(\A)$, 
$y \in H(\A)^{1} \times G'(\A)$, $\tly \in H(\A)^{1} \times \tlG'(\A)$ et $\tlx \in H_{\rmE}(\A)^{1} \times \tlG'_{\rmE}(\A)$.
Notons aussi $\det \in \all_{\tlBmin}$ le déterminant de $H \times \tlG'$ et pour $s \in \C$ et 
$x \in H_{\rmE}(\A)^{1} \times G_{\rmE}(\A)$
notons $|\det x|_{\A}^{s}$ l'expression $e^{s\det H_{H \times \tlG}(x)}$. 

\btheo 
Soit $f \in C_{c}^{\infty}( \hGpa)$, 
alors pour tout $T \in  \mathfrak{a}_{\tlBmin}^{+}$ 
suffisamment régulier et tous $\sigma, \sigma' \in \R$ on a
\[
\sum_{\chi \in \calX^{\hGp}} 
\int\limits_{[H_{\rmE}]^{1} \times [G'_{\rmE}]}
\int\limits_{[H]^{1} \times [G']}
\int\limits_{[H]^{1} \times [\tlG']}
|k^{T}_{f,\chi}(x,x, y, \tly)|
|\det x|_{\A}^{\sigma} 
|\det y|_{\A}^{\sigma'}
d\tly dy dx < \infty.
\]
\bdem
La preuve est similaire à celle du théorème \ref{thm:noyauCVGG}.
Les détails sont laissés au lecteur.
\edem 
\etheo

Notons $n'$ le rang du groupe $H \times G'$.
On définit alors le caractère $\eta$ sur $(H \times G')(\A) \times (H \times \tlG')(\A)$ 
par $\eta(y, \tly) = \eta((\det y)^{n'+ 1}(\det \tly)^{n'})$.
Notons alors pour $s, s' \in \C$, $\chi \in \calX^{\hGp}$
et $f \in C_{c}^{\infty}( \hGpa)$
\[
I_{\chi}^{H \times \tlG',T}(s,s', f) =
\int\limits_{[H_{\rmE}]^{1} \times [G'_{\rmE}]}
\int\limits_{[H]^{1} \times [G']}
\int\limits_{[H]^{1} \times [\tlG']}
k^{T}_{f,\chi}(x,x, y, \tly)
|\det x|_{\A}^{s} 
|\det y|_{\A}^{s'}
\eta(y,\tly)
d\tly dy dx.
\]

Revenons dans le contexte de l'inclusion $G \hrar \tlG$.
Soit $\tlQ$ un sous-groupe parabolique relativement standard de 
$\tlG$. 
Comme il est expliqué dans le paragraphe \ref{par:glnGlnplus1Sp}, 
on a les décompositions 
 $M_{\tlQ} \cong H_{\tlQ} \times \tlG_{\tlQ}$ et 
$M_{Q} \cong H_{\tlQ} \times G_{\tlQ}$ où $H_{\tlQ}$, $G_{\tlQ}$ 
et $\tlG_{\tlQ}$ vérifient les conditions de ce paragraphe. 

Soit $\chi \in \calX^{\GetlG}$, et $\{\chi'\} \sbs
\calX^{M_{Q,\rmE} \times M_{\tlQ, \rmE}}$ la préimage de 
$\chi$ par l'application naturelle $\calX^{M_{Q,\rmE} \times M_{\tlQ, \rmE}} \rar \calX^{\GetlG}$.
Pour $T \in \all_{\tlzero}^{+}$ et $s,s' \in \C$, on définit alors la distribution $I_{\chi}^{M_{\tlQ},T}(s, s', \cdot)$ 
sur 
$C_{c}^{\infty}((H_{\tlQ, \rmE}(\A)^{1} \times G_{\tlQ,\rmE}(\A)) \times (H_{\tlQ, \rmE}(\A)^{1} \times \tlG_{\tlQ,\rmE}(\A)))$ par:
\begin{equation*}
I^{M_{\tlQ},T}_{\chi}(s,s', f) = 
\sum_{\chi'}
I^{M_{\tlQ},T}_{\chi'}(s, s' + s'_{\tlQ}+ s_{\tlQ}, f),
\end{equation*}
où $s_{\tlQ}$ est défini par (\ref{eq:sQSp}) et
pour $\chi' \in \calX^{M_{Q,\rmE} \times M_{\tlQ, \rmE}}$, 
$I^{M_{\tlQ},T}_{\chi'}(s,s',  \cdot)$ 
c'est la distribution associée à l'inclusion $M_{Q} \hrar M_{\tlQ}$ 
décrite ci-dessus par rapport au sous-groupe de Levi minimal $M_{0}$ de $M_{Q}$
et aux sous-groupes 
de Borel $B \cap M_{Q}$ de $M_{Q}$ et $\tlBmin \cap M_{\tlQ}$ de $M_{\tlQ}$.

Pour $\Phi \in C_{c}^{\infty}(\GetlGa)$ et $s \in \C$ on 
définit
\begin{multline}\label{eq:fQdefGSp}
\Phi_{\tlQ,s}(x,\tlx) = 
\int\limits_{\mathclap{K_{\rmE} \times K \times \tlK}} \quad \ \
\int\limits_{(N_{Q,\rmE} \times N_{\tlQ,\rmE})(\A)}
\int\limits_{\mathrlap{(A_{\tlQ}^{st,\infty})^{2}}}
e^{\rho_{Q,\rmE}(H_{Q,\rmE}(a_{1}x)) + \rho_{\tlQ,\rmE}(H_{\tlQ}(a_{3}\tlx))} 
\Phi(k_{1}^{-1}a_{1} xn_{Q}k_{2}, k_{1}^{-1}a_{3} \tlx n_{\tlQ}k_{3})\\
|\det a_{1}|_{\A}^{-s}\eta(k_{2},k_{3})da_{1}da_{3}dn_{\tlQ}dn_{Q}dk_{3}dk_{2}dk_{1}, \quad 
x \in H_{\tlQ,\rmE}(\A)^{1} \times G_{\tlQ,\rmE}(\A), \ 
\tlx \in H_{\tlQ,\rmE}(\A)^{1} \times \tlG_{\tlQ,\rmE}(\A)
\end{multline}
alors $\Phi_{\tlQ,s} \in C_{c}^{\infty}((H_{\tlQ,\rmE}(\A)^{1} \times G_{\tlQ,\rmE}(\A)) \times (H_{\tlQ,\rmE}(\A)^{1} \times \tlG_{\tlQ,\rmE}(\A)))$. 

Notons que l'application 
\[
\tlQ \sps \tlP \mapsto M_{\tlQ} \cap \tlP
\]
définit une bijection entre les sous-groupes paraboliques relativement 
standards de $\tlG$
contenus dans $\tlQ$ et les sous-groupes paraboliques 
semi-standards de
$M_{\tlQ}$ contenant $B \cap M_{Q}$. 
Soit donc $\chi \in \calX^{\GetlG}$ et 
soient $\{\chi_{\tlQ}\} \in \calX^{M_{Q_{\rmE}} \times M_{\tlQ, \rmE}}$ qui s'envoient sur $\chi$ par l'application naturelle 
$\calX^{M_{Q_{\rmE}} \times M_{\tlQ, \rmE}} \rar \calX^{\GetlG}$.
En utilisant le lemme \ref{lem:kChiDescente}, 
on s'aperçoit que pour
tout $\tlP \in \relPb$ contenu dans $\tlQ$, 
 tous 
$x,y \in H_{\tlQ,\rmE}(\A)^{1} \times G_{\tlQ,\rmE}(\A)$ 
et tous
$\tlx,\tly \in H_{\tlQ,\rmE}(\A)^{1} \times \tlG_{\tlQ,\rmE}(\A)$ on a 
\begin{multline}\label{eq:kfQisThis}
\sum_{\chi_{\tlQ}}k_{\Phi_{\tlQ,s}, M_{\tlQ} \cap \tlP, \chi_{\tlQ}}(x, \tlx,y, \tly) = 
\int\limits_{K_{\rmE}}
\int\limits_{K}
\int\limits_{\tlK}
\int\limits_{(A_{\tlQ}^{st,\infty})^{2}}
e^{-\rho_{Q,\rmE}(H_{Q,\rmE}(xa_{1}y)) - \rho_{\tlQ,\rmE}(H_{\tlQ}(\tlx a_{3} \tly))} \\
k_{\Phi,\tlP, \chi}(xk_{1}, \tlx k_{1}, a_{1}yk_{2}, a_{3} \tly k_{3})
\eta(k_{2},k_{3})|\det a_{1}|_{\A}^{-s}
da_{1}da_{3}dk_{3}dk_{2}dk_{1}.
\end{multline}

\subsection{Comportement en $T$}\label{par:asymptChaptGSp}

On démontre la proposition suivante.

\brop\label{prop:mainQualitPropGSp}
Soient $\Phi \in C_{c}^{\infty}(\GetlGa)$, $T' \in T_{+} + \all_{\tlzero}^{+}$, 
$s,s' \in \C$, 
$\chi \in \calX^{G_{\rmE} \times \tlG_{\rmE}}$ et $T \in T' + \all_{\tlzero}^{+}$. Alors
\[
I_{\chi}^{T}(s,s',\Phi) = \sum_{\tlQ \in \relPb}
p_{\tlQ,s + s' }(T_{\tlQ} - T'_{\tlQ})e^{\upla_{\tlQ,s + s'}(T'_{\tlQ})}I_{\chi}^{M_{\tlQ}, T'}(s,s',\Phi_{\tlQ,s})
\]
où pour un sous-groupe parabolique $\tlQ$ 
relativement standard de $\tlG$, 
la fonction $p_{\tlQ,s}$ est définie dans le lemme \ref{lem:pQExplicitGSp} 
dans le paragraphe \ref{par:fonsPolExpSp} où 
$\upla_{\tlQ,s } \in (\all_{\tlQ}^{\tlG})^{*}$ est aussi défini,
la distribution $I_{\chi}^{M_{\tlQ}, T'}$ 
est définie dans le paragraphe \ref{par:JTforLevisGln} et 
$\Phi_{\tlQ,s} \in C_{c}^{\infty}((H_{\tlQ,\rmE}(\A)^{1} \times G_{\tlQ,\rmE}(\A)) \times (H_{\tlQ,\rmE}(\A)^{1} \times \tlG_{\tlQ,\rmE}(\A)))$ 
est définie par (\ref{eq:fQdefGSp}) 
dans le même paragraphe.

\bdem
Fixons un \(T' \in T_{+}+\mathfrak{a}_{\tlzero}^{+} \) 
et soit  $T \in T' + \all_{\tlzero}^{+}$.
En utilisant la relation (\ref{eq:GammaRecurr}) 
dans la définition du noyau $k_{\Phi, \chi}^{T}$ (\ref{eq:kchiGdef})
avec $P = \tlP$, 
$H = H_{\tlP}(\delta_{2} h)-T'_{\tlP}$ et
$X = T_{\tlP} - T'_{\tlP}$ 
pour tout $\tlP \in \relPb$ et tout $\delta_{2} \in P(\rmF) \bsl G(\rmF)$ 
on a que $I_{\chi}^{T}(s, s', \Phi)$ égale la somme sur $\tlQ \in \relPb$ de
\begin{multline}\label{eq:firstGlnass}
\int\limits_{Q(\rmE) \bsl G_{\rmE}(\A)}
\int\limits_{Q(\rmF) \bsl G(\A)} \quad \
\int\limits_{\mathclap{\tlQ(\rmF) \bsl \tlG(\A)}} \ 
\Gamma_{\tlQ}'(H_{\tlQ}(h) - T_{\tlQ}', T_{\tlQ} - T_{\tlQ}')
\sum_{\tlP \sbs \tlQ}(-1)^{d_{\tlP}^{\tlQ}} \ 
\sum_{\mathclap{
\begin{subarray}{c}
\delta_{1} \in P(\rmE)\bsl Q(\rmE)
\\
\delta_{2} \in P(\rmF)\bsl Q(\rmF)
\\
\delta_{3} \in \tlP(\rmF)\bsl \tlQ(\rmF)
\end{subarray}}} \ 
\htau_{\tlP}^{\tlQ}(H_{\tlP}(\delta_{2}h)-T_{\tlP})\\
k_{\tlP,\chi}(\delta_{1} g, \delta_{1} g, \delta_{2}h, \delta_{3}\tlh)
|\det g|^{s}_{\A}
|\det h|^{s'}_{\A}
\eta(h,\tlh)d\tlh dh dg.
\end{multline}

Écrivons maintenant $g = n_{1}a_{1}m_{1}k_{1}$, $h = n_{2}a_{2}m_{2}k_{2}$ et 
$\tlh = n_{3}a_{3}m_{3}k_{3}$ selon les décompositions suivantes:
\begin{gather*}
Q(\rmE) \bsl G_{\rmE}(\A) = [N_{Q, \rmE}] \times A_{\tlQ}^{st,\infty} \times 
(M_{Q}(\rmE) \bsl (H_{\tlQ, \rmE}(\A)^{1} \times G_{\tlQ,\rmE})(\A)) \times K_{\rmE}, \\
Q(\rmF) \bsl G(\A) = [N_{Q}] \times A_{\tlQ}^{st,\infty} \times 
(M_{Q}(\rmF) \bsl (H_{\tlQ}(\A)^{1} \times G_{\tlQ})(\A)) \times K, \\
\tlQ(\rmF) \bsl \tlG(\A) = [N_{\tlQ}] \times A_{\tlQ}^{st, \infty} \times (M_{\tlQ}(\rmF) \bsl (H_{\tlQ}(\A)^{1} \times \tlG_{\tlQ})(\A)) \times \tlK.
\end{gather*}

On obtient alors que (\ref{eq:firstGlnass}) égale
\begin{multline*}
\int\limits_{[H_{\tlQ,\rmE}]^{1} \times [G_{\tlQ,\rmE}]}
\int\limits_{[H_{\tlQ}]^{1} \times [G_{\tlQ}]}
\int\limits_{[H_{\tlQ}]^{1} \times [\tlG_{\tlQ}]}
\int\limits_{(A_{\tlQ}^{st,\infty})^{3}}
\int\limits_{K_{\rmE} \times K \times \tlK}
\sum_{\tlP \sbs \tlQ}(-1)^{d_{\tlP}^{\tlQ}} \ \
\sum_{\mathclap{
\begin{subarray}{c}
\delta_{1} \in (M_{Q} \cap P)(\rmE) \bsl M_{Q}(\rmE)
\\
\delta_{2} \in (M_{Q} \cap P)(\rmF)\bsl M_{Q}(\rmF)
\\
\delta_{3} \in (M_{\tlQ} \cap \tlP)(\rmF)\bsl M_{\tlQ}(\rmF)
\end{subarray}}}
\htau_{\tlP}^{\tlQ}(H_{\tlP}(\delta_{2}m_{2})-T_{\tlP})\\
e^{-2(\rho_{Q, \rmE}(H_{Q}(m_{1}a_{1})) + \rho_{Q}(H_{Q}(m_{2}a_{2})) +\rho_{\tlQ}( H_{\tlQ}(m_{3}a_{3}) )) }
\Gamma_{\tlQ}'(H_{\tlQ}(a_{2}m_{2}) - T_{\tlQ}', T_{\tlQ} - T_{\tlQ}')
|\det a_{1}m_{1}|^{s}_{\A}
|\det a_{2}m_{2}|^{s'}_{\A}
\\
\eta(m_{2}k_{2}, m_{3}k_{3})
k_{\tlP,\chi}(\delta_{1} m_{1}a_{1}k_{1}, \delta_{1} m_{1}a_{1}k_{1}, \delta_{2}m_{2}a_{2}k_{2}, 
\delta_{3}m_{3}a_{3}k_{3}) 
dk_{1}dk_{2}dk_{3} da_{1} da_{2} da_{3} dm_{1}dm_{2}dm_{3}.
\end{multline*}
En utilisant le lemme \ref{lem:artLemme23} on voit que
$k_{\tlP,\chi}(\delta_{1} m_{1}a_{1}k_{1}, \delta_{1} m_{1}a_{1}k_{1}, \delta_{2}m_{2}a_{2}k_{2}, \delta_{3}m_{3}a_{3}k_{3})$ 
égale
\[ 
e^{2\rho_{Q,\rmE}(H_{Q}(a_{1})) + 2\rho_{\tlQ,\rmE}(H_{\tlQ}(a_{1}))}
k_{\tlP,\chi}(\delta_{1} m_{1}k_{1}, \delta_{1} m_{1}k_{1}, \delta_{2}m_{2}a_{1}^{-1}a_{2}k_{2}, \delta_{3}m_{3}a_{1}^{-1}a_{3}k_{3}) 
\]
Remarquons que $\rho_{\tlQ,\rmE} = 2\rho_{\tlQ}$ et $\rho_{Q,\rmE} = 2\rho_{Q}$.
En faisant les changements de variable $a_{1}^{-1}a_{3} \mapsto a_{3}$ 
et ensuite $a_{2}^{-1}a_{1} \mapsto a_{1}^{-1}$ on trouve que
\[
\int\limits_{\mathclap{(A_{\tlQ}^{st,\infty})^{2}}}
e^{2\rho_{\tlQ, \rmE}(H_{\tlQ}(a_{1})) -2\rho_{\tlQ}(H_{\tlQ}(a_{3}))}
k_{\tlP,\chi}(\delta_{1} m_{1}k_{1}, \delta_{1} m_{1}k_{1}, \delta_{2}m_{2}a_{1}^{-1}a_{2}k_{2}, \delta_{3}m_{3}a_{1}^{-1}a_{3}k_{3}) 
|\det a_{1}|^{s}_{\A}
 da_{3} da_{1} 
\]
égale
\[
\int\limits_{\mathclap{(A_{\tlQ}^{st,\infty})^{2}}}
e^{\rho_{\tlQ, \rmE}(H_{\tlQ}(a_{2}) - H_{\tlQ}(a_{1})) -2\rho_{\tlQ}(H_{\tlQ}(a_{3}))}
k_{\tlP,\chi}(\delta_{1} m_{1}k_{1}, \delta_{1} m_{1}k_{1}, \delta_{2}m_{2}a_{1}k_{2}, \delta_{3}m_{3}a_{3}k_{3}) 
|\det a_{2}a_{1}^{-1}|^{s}_{\A}
 da_{3} da_{1}. 
\]
En vertu du lemme \ref{lem:pQExplicitGSp} \textit{2)}, l'intégration par rapport à $da_{2}$ devient
\begin{multline*}
\int\limits_{\mathclap{A_{\tlQ}^{st,\infty}}}
e^{((s + s')\det + 2\rho_{\tlQ} - 2\rho_{Q} )(H_{\tlQ}(a_{2})) }
\Gamma_{\tlQ}'(H_{\tlQ}(a_{2}m_{2}) - T_{\tlQ}', T_{\tlQ} - T_{\tlQ}')da_{2} = \\
|\det m_{2}|^{s_{\tlQ} + s'_{\tlQ}}
e^{\upla_{\tlQ, s + s'}(T_{\tlQ}')}
p_{\tlQ, s +s'}(T_{\tlQ} - T_{\tlQ}').
\end{multline*}
Le reste d'intégration c'est alors
$e^{\upla_{\tlQ, s + s'}(T_{\tlQ}')}
p_{\tlQ, s +s'}(T_{\tlQ} - T_{\tlQ}')$ fois
\begin{multline*}
\int\limits_{[H_{\tlQ,\rmE}]^{1} \times [G_{\tlQ,\rmE}]}
\int\limits_{[H_{\tlQ}]^{1} \times [G_{\tlQ}]}
\int\limits_{[H_{\tlQ}]^{1} \times [\tlG_{\tlQ}]}
\int\limits_{(A_{\tlQ}^{st,\infty})^{2}}
\int\limits_{K_{\rmE} \times K \times \tlK}
\sum_{\tlP \sbs \tlQ}(-1)^{d_{\tlP}^{\tlQ}} \quad
\sum_{\mathclap{
\begin{subarray}{c}
\delta_{1} \in (M_{Q} \cap P)(\rmE) \bsl M_{Q}(\rmE)
\\
\delta_{2} \in (M_{Q} \cap P)(\rmF)\bsl M_{Q}(\rmF)
\\
\delta_{3} \in (M_{\tlQ} \cap \tlP)(\rmF)\bsl M_{\tlQ}(\rmF)
\end{subarray}}}
\htau_{\tlP}^{\tlQ}(H_{\tlP}(\delta_{2}m_{2})-T_{\tlP})\\
e^{-\rho_{Q, \rmE}(H_{Q}(m_{1}a_{1}m_{2}) -\rho_{\tlQ, \rmE}(H_{\tlQ}(m_{1}a_{3}m_{3}) ) }
|m_{1}|^{s}_{\A}
|m_{2}|^{s' + s'_{\tlQ} + s_{\tlQ}}_{\A}
\eta(m_{2}, m_{3})\eta(k_{2}, k_{3})|\det a_{1}|_{\A}^{-s} \\
k_{\tlP,\chi}(\delta_{1} m_{1}k_{1}, \delta_{1} m_{1}k_{1}, \delta_{2}m_{2}a_{1}k_{2}, 
\delta_{3}m_{3}a_{3}k_{3}) 
dk_{1}dk_{2}dk_{3} da_{1} da_{3} dm_{1}dm_{2}dm_{3}
\end{multline*}
où on utilise le fait que $\rho_{Q, \rmE}(H_{Q}(m_{1}) = \rho_{\tlQ, \rmE}(H_{Q}(m_{1})$ 
si $m_{1} \in H_{\tlQ}(\A)^{1} \times G_{\tlQ}(\A)$. 
En utilisant la discussion du paragraphe \ref{par:JTforLevisGln}, l'équation 
(\ref{eq:kfQisThis}) en l'occurrence, on trouve que l'expression ci-dessus égale
$I_{\chi}^{M_{\tlQ}, T}(s, s', \Phi_{\tlQ,s})$ où $\Phi_{\tlQ,s}$ est définie 
par (\ref{eq:fQdefGSp}), ce qui termine la preuve de la proposition. 
\edem
\erop

En utilisant la proposition \ref{prop:mainQualitPropGSp} démontrée 
ci-dessus et le lemme \ref{lem:pQExplicitGSp} qui 
décrit les fonctions $p_{\tlQ}$ explicitement on obtient 
le comportement en $T$ de la distribution $I_{\chi}^{T}$.

\begin{theo}\label{thm:mainQualitThmG}
Soient $\Phi \in C_{c}^{\infty}(\GetlGa)$ et $\chi \in \calX^{G_{\rmE} \times \tlG_{\rmE}}$.
La fonction
$T \mapsto I^{T}_{\chi}(s, s',\Phi)$ 
 où
$\chi \in \calX^{G_{\rmE} \times \tlG_{\rmE}}$, $s,s' \in \C$ et
$T$ 
parcourt $T_{+} + \mathfrak{a}_{\tlzero}^{+}$ 
est un polynôme-exponentielle. 
De plus, si $s + s' \neq -1,1$ sa
partie purement polynomiale est constante
et donnée par  
\begin{gather*}
I_{\chi}(s, s',\Phi) := \sum_{\tlQ \in \relPb}(-1)^{d_{\tlQ}^{\tlG}}j_{\tlQ}^{-1}\hat \theta_{\tlQ}(\upla_{\tlQ,s+s'})^{-1}
e^{\upla_{\tlQ,s+s'}(T'_{\tlQ})}I_{\chi}^{M_{\tlQ},T'}(s,s',\Phi_{\tlQ,s})
\end{gather*}
pour tout $T' \in T_{+} + \all_{\tlzero}^{+}$. 
En particulier, la distribution $I_{\chi}$
ne dépend pas de $T'$.
\end{theo}

\brem\label{rem:JMQisPolExpGSp}
 Soit $\tlQ$ un sous-groupe parabolique relativement standard 
de $\tlG$. 
Par le même raisonnement que dans la proposition \ref{prop:mainQualitPropGSp}
 on obtient que pour tous $s,s' \in \C$ la
distribution $I_{\chi}^{M_{\tlQ},T}(s,s', \cdot)$  
définie dans le paragraphe \ref{par:JTforLevisGln}, 
est un polynôme-exponentielle en $T$ qui ne dépend pas
de $T_{\tlQ} \in \all_{\tlQ}$. Cependant, 
si $\tlQ \neq \tlG$ le terme purement polynomial n'est pas constant. 
\erem

\subsection{Équivariance}\label{par:compConjugGSp}

Soient $\Phi \in C_{c}^{\infty}(\GetlGa)$, $x \in G_{\rmE}(\A)$ et $(y,\tly) \in G(\A) \times \tlG(\A)$. 
Notons $\Phi^{x,(y,\tly)} \in C_{c}^{\infty}(\GetlGa)$ la fonction 
définie par $\Phi^{x,(y,\tly)}(g,\tlg) = \Phi(xgy^{-1}, x\tlg\tly^{-1})$.

Soit $\chi \in \calX^{\GetlG}$.
On voit que $I_{\chi}^{T}(s,s',\Phi^{x,(y,\tly)})$ pour 
$T \in T_{+} + \all_{\tlzero}^{+}$ et $s,s' \in \C$ égale
\begin{multline*}
\int_{[G_{\rmE}]}\int_{[G]}\int_{[\tlG]}
\sum_{\tlP \in \relPb}(-1)^{d_{\tlP}^{\tlG}}
\sum_{
\begin{subarray}{c}
\delta_{1} \in P(\rmE)\bsl G(\rmE)
\\
\delta_{2} \in P(\rmF)\bsl G(\rmF)
\\
\delta_{3} \in \tlP(\rmF)\bsl \tlG(\rmF)
\end{subarray}}
\htau_{\tlP}(H_{\tlP}(\delta_{2}hy)-T_{\tlP})
k_{\tlP,\chi}(\delta_{1} g, \delta_{1} g, \delta_{2} h, \delta_{3} \tlh)\\
|\det gx|_{\A}^{s}
|\det hy|_{\A}^{s'}
\eta(hy, \tlh \tly)d\tlh dh dg.
\end{multline*}

Pour $h \in G(\A)$ et $P \in \calF(B)$ soit $k_{P}(h)$ un élément 
de $K$ tel que $hk_{P}(h)^{-1} \in P(\A)$. 
Alors, en utilisant l'égalité (\ref{eq:GammaRecurr})
on a pour tout $\tlP \in \relPb$:
\begin{equation*}
\htau_{\tlP}(H_{\tlP}(\delta_{2} h y)-T_{\tlP}) =
\sum_{\tlQ \supseteq \tlP}(-1)^{d_{\tlQ}^{\tlG}}
\htau_{\tlP}^{\tlQ}(H_{\tlP}(\delta_{2} h)-T_{\tlP})
\Gamma_{\tlQ}'(H_{\tlQ}(\delta_{2} h)-T_{\tlQ},
-H_{\tlQ}(k_{Q}(\delta_{2} h)y))
\end{equation*}
d'où on obtient que $I_{\chi}^{T}(s,s',\Phi^{x,(y,\tly)})$ égale la 
somme sur $\tlQ \in \relPb$ de 
\begin{multline}\label{eq:equivG1}
\int\limits_{Q(\rmE) \bsl G_{\rmE}(\A)}
\int\limits_{Q(\rmF) \bsl G(\A)}
\int\limits_{\tlQ(\rmF) \bsl \tlG(\A)}
\Gamma_{\tlQ}'(H_{\tlQ}(h)-T_{\tlQ}, -H_{\tlQ}(k_{Q}(h)y))
\sum_{\tlP \sbs \tlQ}(-1)^{d_{\tlP}^{\tlQ}} \quad
\sum_{\mathclap{
\begin{subarray}{c}
\delta_{1} \in (P \cap M_{Q})(\rmE)\bsl M_{Q}(\rmE)
\\
\delta_{2} \in (P \cap M_{Q})(\rmF)\bsl M_{Q}(\rmF)
\\
\delta_{3} \in (\tlP \cap M_{\tlQ})(\rmF)\bsl M_{\tlQ}(\rmF)
\end{subarray}}} \\
\htau_{\tlP}^{\tlQ}(H_{\tlP}(\delta_{2}h)-T_{\tlP})
k_{\tlP,\chi}(\delta_{1} g, \delta_{1} g, \delta_{2}h, \delta_{3}\tlh)
|\det gx|_{\A}^{s}
|\det hy|_{\A}^{s'}
\eta(hy, \tlh \tly)d\tlh dh dg.
\end{multline}

Écrivons maintenant $g = n_{1}a_{1}m_{1}k_{1}$, $h = n_{2}a_{2}m_{2}k_{2}$ et 
$\tlh = n_{3}a_{3}m_{3}k_{3}$ selon les décompositions suivantes:
\begin{gather*}
Q(\rmE) \bsl G_{\rmE}(\A) = [N_{Q, \rmE}] \times A_{\tlQ}^{st,\infty} \times 
(M_{Q}(\rmE) \bsl (H_{\tlQ, \rmE}(\A)^{1} \times G_{\tlQ,\rmE}(\A))) \times K_{\rmE}, \\
Q(\rmF) \bsl G(\A) = [N_{Q}] \times A_{\tlQ}^{st,\infty} \times 
(M_{Q}(\rmF) \bsl (H_{\tlQ}(\A)^{1} \times G_{\tlQ}(\A))) \times K, \\
\tlQ(\rmF) \bsl \tlG(\A) = [N_{\tlQ}] \times A_{\tlQ}^{st, \infty} \times (M_{\tlQ}(\rmF) \bsl (H_{\tlQ}(\A)^{1} \times \tlG_{\tlQ}(\A))) \times \tlK.
\end{gather*}
On a alors: 
\begin{equation*}
\Gamma_{\tlQ}'(H_{\tlQ}(n_{2}a_{2}m_{2}k_{2})-T_{\tlQ},-H_{\tlQ}(k_{Q}(n_{2}a_{2}m_{2}k_{2})y))
 = 
\Gamma_{\tlQ}'(H_{\tlQ}(a_{2}) + H_{\tlQ}(m_{2})-T_{\tlQ},-H_{\tlQ}(k_{2}y)).
\end{equation*}

On fait les mêmes changements de variable que dans la
preuve
de la proposition \ref{prop:mainQualitPropGSp}. 
On obtient, grâce à la relation (\ref{eq:HQtlQSp}) et le lemme \ref{lem:pQExplicitGSp} \textit{2)}, 
que l'intégration par rapport à $da_{2}$ devient 
\begin{multline*}
\int_{A_{\tlQ}^{st,\infty}}
e^{((s + s')\det + 2\rho_{\tlQ} - 2\rho_{Q} )(H_{\tlQ}(a_{2})) }
\Gamma_{\tlQ}'(H_{\tlQ}(a_{2}) + H_{\tlQ}(m_{2})-T_{\tlQ},-H_{\tlQ}(k_{2}y))da_{2} = \\
|\det m_{2}|_{\A}^{s_{\tlQ} + s'_{\tlQ}} 
e^{\upla_{\tlQ, s+ s'}(T_{\tlQ})}u_{\tlQ,s +s',y}(k_{2})
\end{multline*}
où 
\begin{equation*}
u_{\tlQ,s,y}(k) = \int_{\all_{\tlQ}^{\tlG}}
e^{\upla_{\tlQ,s}(H)}\Gamma_{\tlQ}'(H,-H_{\tlQ}(ky))dH, \quad k \in K, \ s \in \C
\end{equation*}
est une fonction continue de la variable $k \in K$ en vertu du lemme \ref{lem:GammaIsComp}.

Le reste de l'intégration dans (\ref{eq:equivG1}) c'est
$e^{\upla_{\tlQ, s+ s'}(T_{\tlQ})}|\det x|_{\A}^{s}|\det y|_{\A}^{s'} \eta(y,\tly)$ fois
\begin{multline}\label{eq:equivG2}
\int\limits_{[H_{\tlQ,\rmE}]^{1} \times [G_{\tlQ,\rmE}]}
\int\limits_{[H_{\tlQ}]^{1} \times [G_{\tlQ}]}
\int\limits_{[H_{\tlQ}]^{1} \times [\tlG_{\tlQ}]}
\int\limits_{(A_{\tlQ}^{st,\infty})^{2}}
\int\limits_{K_{\rmE} \times K \times \tlK}
\sum_{\tlP \sbs \tlQ}(-1)^{d_{\tlP}^{\tlQ}} \quad
\sum_{\mathclap{
\begin{subarray}{c}
\delta_{1} \in (P \cap M_{Q})(\rmE)\bsl M_{Q}(\rmE)
\\
\delta_{2} \in (P \cap M_{Q})(\rmF)\bsl M_{Q}(\rmF)
\\
\delta_{3} \in (\tlP \cap M_{\tlQ})(\rmF)\bsl M_{\tlQ}(\rmF)
\end{subarray}}}
\htau_{\tlP}^{\tlQ}(H_{\tlP}(\delta_{2}m_{2})-T_{\tlP})\\
e^{-\rho_{Q_{\rmE}}(H_{Q}(m_{1}a_{1}m_{2}) -\rho_{\tlQ}( H_{\tlQ}(m_{1}a_{3}m_{3}) ) }
|m_{1}|^{s}_{\A}
|m_{2}|^{s' + s'_{\tlQ} + s_{\tlQ}}_{\A}
\eta(m_{2}, m_{3})\eta(k_{2}, k_{3})|\det a_{1}|_{\A}^{-s} \\
u_{\tlQ,s + s',y}(k_{2})k_{\tlP,\chi}(\delta_{1} m_{1}k_{1}, \delta_{1} m_{1}k_{1}, \delta_{2}m_{2}a_{1}k_{2}, 
\delta_{3}m_{3}a_{3}k_{3}) 
dk_{1}dk_{2}dk_{3} da_{1} da_{3} dm_{1}dm_{2}dm_{3}.
\end{multline}
Si l'on pose alors
\begin{multline*}
\Phi_{\tlQ,s,s',y}(x,\tlx) \ \ \mathclap{=} \,
\int\limits_{\mathclap{K_{\rmE} \times K \times \tlK}} \quad \ 
\int\limits_{(N_{Q} \times N_{\tlQ})_{\rmE}(\A)}
\int\limits_{\mathrlap{(A_{\tlQ}^{st,\infty})^{2}}}
e^{\rho_{Q,\rmE}(H_{Q,\rmE}(a_{1}x)) + \rho_{\tlQ,\rmE}(H_{\tlQ}(a_{3}\tlx))} 
\Phi(k_{1}^{-1}a_{1} xn_{Q}k_{2}, k_{1}^{-1}a_{3} \tlx n_{\tlQ}k_{3})\\
|\det a_{1}|_{\A}^{-s}
u_{\tlQ,s +s',y}(k_{2})
\eta(k_{2},k_{3})da_{1}da_{3}dn_{\tlQ}dn_{Q}dk_{3}dk_{2}dk_{1},
\end{multline*}
où $x \in H_{\tlQ,\rmE}(\A)^{1} \times G_{\tlQ,\rmE}(\A)$ et $\tlx \in H_{\tlQ,\rmE}(\A)^{1} \times G_{\tlQ,\rmE}(\A)$,
on a bien $\Phi_{\tlQ,s,s',y} \in C_{c}^{\infty}((H_{\tlQ,\rmE}(\A)^{1} \times G_{\tlQ,\rmE}(\A)) \times (H_{\tlQ,\rmE}(\A)^{1} \times \tlG_{\tlQ,\rmE}(\A)))$
et on voit, 
en se basant sur la discussion du paragraphe \ref{par:JTforLevisGln} 
et en utilisant le lemme \ref{lem:kChiDescente} du paragraphe \ref{par:repRegul},
 que (\ref{eq:equivG2}) c'est juste $I_{\chi}^{M_{\tlQ},T}(s, s', \Phi_{\tlQ,s,s',y})$.
On obtient donc le théorème suivant.

\begin{theo}\label{thm:invarianceTheoG} 
Soient $x \in G_{\rmE}(\A)$, $(y,\tly) \in G(\A) \times \tlG(\A)$, $\chi \in \calX^{G_{\rmE} \times \tlG_{\rmE}}$ et
$s,s' \in \C$. La distribution $I_{\chi}^{T}(s,s', \cdot)$ vérifie pour tout 
$\Phi \in C_{c}^{\infty}(\GetlGa)$ la propriété suivante:
\begin{multline*}
I_{\chi}^{T}(s,s',\Phi^{x,(y, \tly)}) - |\det x|_{\A}^{s}|\det y|_{\A}^{s'} \eta(y,\tly) I_{\chi}^{T}(s,s',\Phi) = \\ 
|\det x|_{\A}^{s}|\det y|_{\A}^{s'} \eta(y,\tly) 
\sum_{\tlQ \in \relPb \smin \{\tlG\}} 
e^{\upla_{\tlQ,s + s'}(T_{\tlQ})}
I_{\chi}^{M_{\tlQ},T}(s,s',\Phi_{\tlQ,s, s', y})
\end{multline*}
où les distributions $I_{\chi}^{M_{\tlQ},T}(s, s', \cdot)$ 
sont définies dans la section \ref{par:JTforLevisGln}. 
En particulier, pour $s +s' \neq -1,1$ on a
\[
I_{\chi}(s,s',\Phi^{x,(y, \tly)})= |\det x|_{\A}^{s}|\det y|_{\A}^{s'} \eta(y,\tly) I_{\chi}^{T}(s,s',\Phi).
\]

\bdem 
La formule pour la différence 
\[
I_{\chi}^{T}(s,s',\Phi^{x,(y, \tly)}) - |\det x|_{\A}^{s}|\det y|_{\A}^{s'} \eta(y,\tly) I_{\chi}^{T}(s,s',\Phi) 
\]
 est claire 
après les calculs qu'on a faits. Si $s + s' \neq -1,1$, cette formule-ci démontre aussi l'équivariance 
de la distribution $I_{\chi}^{T}(s,s',\cdot)$, car
si $\tlQ \subsetneq \tlG$, d'après la remarque \ref{rem:JMQisPolExpGSp}, le terme 
$I_{\chi}^{M_{\tlQ},T}(s,s',\Phi_{\tlQ,s, s', y})$ est un polynôme-exponentielle en $T$ qui ne dépend pas 
de $T_{\tlQ} \in \all_{\tlQ}$. En outre, $\upla_{\tlQ,s + s'}$ est non-trivial
sur $\all_{\tlQ}^{\tlG}$ en vertu du lemme \ref{lem:uplasNonNulG}. 
Il en découle que l'expression $e^{\upla_{\tlQ,s + s'}(T_{\tlQ})}
I_{\chi}^{M_{\tlQ},T}(s,s',\Phi_{\tlQ,s, s', y})$ n'a pas de terme 
constant dans ce cas et par conséquent les termes constants de $I_{\chi}^{T}(s,s',\Phi^{x,(y, \tly)})$ 
et de $ |\det x|_{\A}^{s}|\det y|_{\A}^{s'} \eta(y,\tly) I_{\chi}^{T}(s,s',\Phi)$ coïncident. 
\edem
\end{theo}

\subsection{Données cuspidales non-induites}\label{par:exemplDistG}

Dans ce paragraphe on étudie la distribution $I_{\chi}$ 
pour des classes $\chi$ distinguées. 

Soit $H$ un $\rmF$-groupe réductif quelconque. On dit que $\chi \in \calX^{H}$ est non-induite
si $\chi = \{(H,\sigma)\}$ où $\sigma$ est une représentation cuspidale de $H(\A)^{1}$.

\brop\label{prop:DistSimplG}
  Soient $\Phi \in C_{c}^{\infty}(\GetlGa)$, $\chi^{G_{\rmE}} \in \calX^{G_{\rmE}}$ et $\chi^{\tlG_{\rmE}} \in \calX^{\tlG_{\rmE}}$. 
  Notons $\chi := \chi^{G_{\rmE}} \times \chi^{\tlG_{\rmE}} \in \calX^{\GetlG}$.
\begin{enumerate}[1)]
\item Supposons que $\chi^{\tlG_{\rmE}}$ est non-induite. 
Alors, pour tous $s,s' \in \C$ et $T \in T_{+} + \all_{\tlzero}^{+}$ on a
\begin{multline*}
I_{\chi}^{T}(s,s', \Phi) = 
\int\limits_{\mathclap{[G_{\rmE}]}} \ \int\limits_{[G]}\int\limits_{\mathclap{[\tlG]}}
k_{\Phi, \chi}(g,g, h,\tlh)
|\det g|_{\A}^{s}|\det h|_{\A}^{s'}\eta(h,\tlh)d\tlh dh dg = 
\sum_{\mathclap{\tlQ \in \relPb}}
\\
\int\limits_{[G_{\rmE}]}\int\limits_{Q(\rmF) \bsl G(\A)}\int\limits_{\mathclap{[\tlG]}}
\tau_{\tlQ}(H_{\tlQ}(h) - T) 
\Lazt^{\tlQ, T} 
k_{\Phi, \chi}(g,g, h,\tlh)
|\det g|_{\A}^{s}|\det h|_{\A}^{s'}\eta(h,\tlh)d\tlh dh dg.
\end{multline*}
\item Supposons que $\chi^{G_{\rmE}}$ est non-induite et que pour tout $(\tlM_{\rmE}, \sigma) \in \chi^{\tlG_{\rmE}}$ 
le sous-groupe de Levi $\tlM_{\rmE}$ ne stabilise aucune droite. 
Alors, pour tous $s,s' \in \C$ et $T \in T_{+} + \all_{\tlzero}^{+}$ on a
\begin{multline*}
I_{\chi}^{T}(s,s', \Phi) = 
\int_{[G_{\rmE}]}\int_{[G]}\int_{[\tlG]}
k_{\Phi, \chi}(g,g, h,\tlh)
|\det g|_{\A}^{s}|\det h|_{\A}^{s'}\eta(h,\tlh)d\tlh dh dg = \\
\int_{[G_{\rmE}]}\int_{[G]}\int_{[\tlG]}
\La_{d,12}^{T} 
k_{\Phi, \chi}(g,g, h,\tlh)
|\det g|_{\A}^{s}|\det h|_{\A}^{s'}\eta(h,\tlh)d\tlh dh dg.
\end{multline*}
\end{enumerate}
En plus, si $s + s' \neq -1,1$, dans les deux cas ci-dessus on a 
$I_{\chi}(s,s', \Phi) = I_{\chi}^{T}(s,s', \Phi)$ 
pour tout $T \in T_{+} + \all_{\tlzero}^{+}$.
\bdem 
Démontrons le point \textit{1)}. Soit $\tlP \in \relPb$. On a alors $k_{\Phi, \tlP, \chi} \equiv 0$ 
si $\tlP \neq \tlG$ donc $k_{\Phi, \tlP, \chi}  = k_{\Phi, \tlP, \chi}^{T}$ et la première égalité 
découle du théorème \ref{thm:noyauCVGG}.
Pour la deuxième égalité, on applique la formule d'inversion donnée par le lemme \ref{lem:formDInv} pour l'opérateur 
$\La_{m,3}^{T}$
ce qui donne formellement l'égalité voulue. La convergence découle de la preuve du 
théorème \ref{thm:noyauCVGG}, où l'expression (\ref{eq:AllSigmas}) est un cas particulier des intégrales considérées.

Démontrons le point \textit{2)}. Soit $\tlP \in \relPb$ et $P = G \cap \tlP$. 
Si $P \neq G$, on a $k_{\Phi, \tlP, \chi} \equiv 0$ car $\chi^{G_{\rmE}}$ est non-induite. 
Supposons donc que $P = G$. Si $\tlP \neq \tlG$ on a que $M_{\tlP}$ stabilise la droite $D_{0}$. La condition imposée 
sur $\chi^{\tlG_{\rmE}}$ implique alors $k_{\Phi, \tlP, \chi} \equiv 0$. On obtient alors 
la première égalité. Pour la deuxième, on utilise la formule d'inversion 
donnée dans le lemme \ref{lem:formDInvD} 
et l'on obtient que 
$I_{\chi}^{T}(s,s', \Phi)$ égale la somme sur $\tlP \in \relPb$ de
\[
\int_{P(\rmE) \bsl G_{\rmE}(\A)}\int_{[G]}\int_{[\tlG]}
\tau_{\tlP}(H_{\tlP}(g) - T)
\La_{d,12}^{\tlP, T} 
k_{\Phi, \chi}(g,g, h,\tlh)
|\det g|_{\A}^{s}|\det h|_{\A}^{s'}\eta(h,\tlh)d\tlh dh dg.
\]
La convergence des intégrales découle de nouveau de la convergence de l'expression (\ref{eq:AllSigmas}) 
du théorème \ref{thm:noyauCVGG}. Soit donc $\tlP \in \relPb$ et $P = \tlP \cap G$. 
Si $P \neq G$ la fonction $[G_{\rmE}] \ni g \mapsto \La_{d,12}^{\tlP, T} 
k_{\Phi, \chi}(g,\cdot, \cdot,\cdot)$ est nulle par hypothèse sur $\chi^{G_{\rmE}}$. 
Si $P = G$ et $\tlP \neq \tlG$ le groupe 
$M_{\tlP}$ stabilise la droite $D_{0}$ et 
$[\tlG_{\rmE}] \ni \tlg \mapsto \La_{d,12}^{\tlP, T} 
k_{\Phi, \chi}(\cdot,\tlg, \cdot,\cdot) \equiv 0$ d'où le résultat.

Si $s + s' \neq -1,1$ on voit que $I_{\chi}^{T}(s,s', \Phi)$ ne dépend pas de 
$T$ dans les deux cas qu'on a considéré, donc c'est bien $I_{\chi}(s,s', \Phi)$.
\edem
\erop
\section{Le côté spectral de la formule des traces pour les groupes unitaires}\label{sec:SpecU}
\sectionmark{Le côté spectral pour les groupes unitaires}

\renewcommand{\Laiyd}{\La_{A,1}}
\subsection{Convergence du noyau spectral tronqué}\label{par:specCvgU}
On démontre le théorème \ref{thm:noyauCVGU} dans ce paragraphe. 
La preuve est similaire à celle du théorème \ref{thm:noyauCVGG} et sera donc un peu moins détaillée.

Soient $F \in C_{c}^{\infty}((U \times \tlU)(\A))$, $k = k_{F}$ 
son noyau automorphe. 
Pour un sous-groupe parabolique 
standard $P$ de $U$,
$\chi \in \calX^{U \times \tlU}$, $x,y \in N_{P}(\A)M_{P}(\rmF) \bsl U(\A)$ et
$\tlx, \tly \in N_{\tlP}(\A)M_{\tlP}(\rmF) \bsl \tlU(\A)$ posons
\[
k_{P, \chi}(x ,\tlx, y ,\tly) =
k_{F, P, \chi}(x ,\tlx, y ,\tly) = 
k_{F,  P \times \tlP, \chi}(x ,\tlx, y ,\tly)
\]
où on utilise le fait que pour tout sous-groupe parabolique $P \sps P_{0}$
il existe un unique $\tlP \in \calF(M_{\tlzero})$ tel que $\tlP \cap U = P$, comme on l'a déjà remarqué dans le paragraphe 
\ref{par:relSt}. 

On pose alors pour $x,y \in U(\rmF) \bsl U(\A)$, $\tlx, \tly \in \tlU(\rmF) \bsl \tlU(\A)$ 
et $T \in \all_{0}$
\begin{equation}\label{eq:kchiUdef}
k_{\chi}^{T}(x ,\tlx, y ,\tly) =
k_{F, \chi}^{T}(x ,\tlx, y ,\tly) 
= \sum_{P \sps P_{0}}(-1)^{d_{P}} \ \ 
\sum_{\mathclap{\delta_{1},\delta_{2} \in \PUf}} \ \ 
k_{F, P, \chi}(\delta_{1}x, \delta_{1}\tlx, \delta_{2}y, 
\delta_{2}\tly)
\htau_{P}(H_{P}(\delta_{2}y)-T).
\end{equation}

\begin{theo}\label{thm:noyauCVGU} On a pour tout $F \in C_{c}^{\infty}((U \times \tlU)(\A))$
et tout $T \in T_{+} + \all_{0}^{+}$
\begin{equation}\label{eq:kUcvg}
\sum_{\chi \in \calX^{U \times \tlU}} \int_{[U \times U]}
|k_{F, \chi}^{T}(x,x,y,y)|dxdy < \infty.
\end{equation}

\bdem 
En raisonnant comme au début de la preuve du théorème \ref{thm:noyauCVGG} 
on trouve un $T_{F} \in \all_{0}$ qui ne dépend que du support de $F$ tel que 
si l'on pose $T' = T + T_{F}$ on a 
que $k_{F, \chi}^{T}(x,x,y,y)$ égale:
\[
\sum_{P \sps P_{0}}(-1)^{d_{P}^{U}}
\sum_{\delta_{1},\delta_{2} \in \PUf} 
k_{P,\chi}(\delta_{1}x,\delta_{1}x,\delta_{2}y, \delta_{2}y)
\htau_{P}(H_{P}(\delta_{1}x)-T')
\htau_{P}(H_{P}(\delta_{2}y)-T).
\]
En particulier, en vertu du lemme \ref{lem:lemmeArth51} \textit{ii)} les sommes dans la définition de 
$k_{F, \chi}^{T}$ sont finies et la fonction est bien définie.

Les fonctions  
$P(\rmF) \bsl U(\A) \times \tlP(\rmF) \bsl \tlU(\A) \ni (x, \tlx) \mapsto k_{F, P,\chi}(x, \tlx, \cdot, \cdot)$ et 
$(y,\tly) \mapsto k_{F, P,\chi}(\cdot, \cdot, y, \tly)$
égalent ses termes constants 
le long de $P \times \tlP$.
En utilisant alors le lemme 
\ref{lem:formDInvD} (sa version pour les groupes unitaires)
et ensuite le lemme \ref{lem:corArth62} \textit{i)} on trouve que 
l'intégrale (\ref{eq:kUcvg}) est majorée par la somme 
sur $P_{1} \sbs P_{3}$ et $P_{2} \sbs P_{4}$ de 
\begin{equation*}
\int\limits_{P_{1}(\rmF) \bsl U(\A)} \quad \ \
\int\limits_{\mathclap{P_{2}(\rmF) \bsl U(\A)}} \
\sigma_{1}^{3}(H_{1}(x)-T')
\sigma_{2}^{4}(H_{2}(y)-T)
|
\Ladud^{T',P_{1}}
\Latqd^{T,P_{2}}
(\ \sum_{\mathclap{Q_{1} \sbs P \sbs Q_{2}}} \
(-1)^{d_{P}}
k_{P,\chi}(x, x, y, y))
| dydx
\end{equation*}
où $Q_{1}$ c'est le plus petit sous-groupe parabolique 
standard de $U$ contenant $P_{1} \cup P_{2}$ 
et $Q_{2} = P_{3} \cap P_{4}$. On les fixe désormais. L'intégrale ci-dessus égale aussi
\begin{multline}\label{eq:AllSigmasU}
\int\limits_{P_{1}(\rmF) \bsl U(\A)}
\int\limits_{P_{2}(\rmF) \bsl U(\A)}
\sigma_{1}^{3}(H_{1}(x)-T')
\sigma_{2}^{4}(H_{2}(y)-T) \\
|
\Ladud^{T',P_{1}}
\Latqd^{T,P_{2}}
( \ \sum_{\mathclap{Q_{1} \sbs P \sbs Q_{2}}} 
(-1)^{d_{P}}
\int\limits_{[N_{1}]}
\int\limits_{[N_{\tltwo}]}
k_{P,\chi}(n_{1}x, x, y, n_{\tltwo}y)dn_{\tltwo}dn_{1})
| dydx.
\end{multline}
Pour tout $P$ entre $Q_{1}$ et $Q_{2}$, en utilisant le corollaire
\ref{cor:art23REALLY} et la décomposition de Bruhat, on a
\[
\int\limits_{[N_{1}]}
\int\limits_{[N_{\tltwo}]}
k_{P,\chi}(n_{1}x, x, y, n_{\tltwo}y)dn_{\tltwo}dn_{1} = 
\sum_{
\begin{subarray}{c}
s \in \Omega^{2} \bsl \Omega^{P} / \Omega^{1} \\
\tls \in \Omega^{\tlone} \bsl \Omega^{\tlP} / \Omega^{\tltwo}
\end{subarray}
}
\sum_{
\begin{subarray}{c}
\gamma \in (P_{2} \cap sP_{1})(\rmF) \bsl P_{2}(\rmF)\\
\tlgam \in (\tlP_{1} \cap \tls\tlP_{2})(\rmF) \bsl \tlP_{1}(\rmF)
\end{subarray}
}
k_{P_{1} \times \tlP_{2}, \chi}( x, w_{\tls}^{-1} \tlgam x, w_{s}^{-1}\gamma y,  y).
\]
Posons
\[
(\Omega_{\tlQ_{1}}^{\tlP})' = \Omega^{\tlP} \smin \bigcup_{Q_{1} \sbs R \sbn P} \Omega^{\tlR}
\]
et $(\Omega_{\tlQ_{1}, U}^{\tlP})' = (\Omega_{\tlQ_{1}}^{\tlP})' \cap \Omega^{U}$. 
Notons que $(\Omega_{\tlQ_{1}, U}^{\tlP})' = \Omega^{P} \smin \bigcup_{Q_{1} \sbs R \sbn P} \Omega^{R}$. 
Alors $(\Omega_{\tlQ_{1}}^{\tlP})'$ (resp. $(\Omega_{\tlQ_{1}, U}^{\tlP})'$) est $\Omega^{\tlone}$-stable (resp. $\Omega^{2}$-stable)
à gauche et $\Omega^{\tltwo}$-stable (resp. $\Omega^{1}$-stable) à droite.
On a alors les décompositions en parties disjointes suivantes
\[
\Omega^{\tlP} = \coprod_{Q_{1} \sbs R \sbs P}(\Omega_{\tlQ_{1}}^{\tlP})' \quad 
\Omega^{P} = \coprod_{Q_{1} \sbs R \sbs P}(\Omega_{\tlQ_{1}, U}^{\tlP})'.
\]
Pour $\tls \in \Omega^{\tlQ_{2}}$ soit $\hDelta_{\tls}^{\tlQ_{1}}= \{\varpi \in \hDelta_{\tlQ_{1}} | \tls\varpi = \varpi\}$. 
Alors $\tls \in (\Omega_{\tlQ_{1}}^{\tlP})'$ si et seulement si $\hDelta_{\tls}^{\tlQ_{1}} = \hDelta_{\tlP}$. 
Par un argument classique basé sur l'identité (\ref{eq:basicidentity})
 on a que l'expression entre la valeur absolue dans (\ref{eq:AllSigmasU}) égale
\[
\Ladud^{T',P_{1}}
\Latqd^{T,P_{2}}
(\sum_{(s, \tls) \in \Upomega_{\tlQ_{1}}'}
\sum_{
\begin{subarray}{c}
\gamma \in (P_{2} \cap sP_{1})(\rmF) \bsl P_{2}(\rmF)\\
\tlgam \in (\tlP_{1} \cap \tls\tlP_{2})(\rmF) \bsl \tlP_{1}(\rmF)
\end{subarray}
}
k_{P_{1} \times \tlP_{2}, \chi}( x, w_{\tls}^{-1} \tlgam x, w_{s}^{-1}\gamma y,  y))
\]
où 
\[
\Upomega_{\tlQ_{1}}' = \{(s, \tls) \in (\Omega^{2} \bsl \Omega^{Q_{2}}/\Omega^{1}) \times 
(\Omega^{\tlone} \bsl \Omega^{\tlQ_{2}} / \Omega^{\tltwo})| \hDelta_{\tlQ_{2}} = \hDelta_{s}^{\tlQ_{1}} \cap \hDelta_{\tls}^{\tlQ_{1}}\}.
\]

Écrivons maintenant $x = n_{1}a_{1}m_{1}k_{1}$, $y = n_{2}a_{2}m_{2}k_{2}$ et 
selon les décompositions suivantes:
\begin{equation*}
P_{1}(\rmF) \bsl U(\A) = [N_{1}]   \times A_{1}^{\infty}  \times [M_{1}]^{1} \times K, \quad 
P_{2}(\rmF) \bsl U(\A) = [N_{2}]  \times A_{2}^{\infty}  \times [M_{2}]^{1} \times K.
\end{equation*}
Supposons que $F$ est invariante à droite par un compact $K_{0} \sbs (U \times \tlU)(\A_{f})$. 
Comme il est expliqué dans \cite{arthur5} entre les pages 93 et 94, 
pour $(x,\tlx), (y,\tly) \in (U \times \tlU)(\A)$ et $(s, \tls) \in \Upomega_{\tlQ_{1}}'$ fixés
les fonctions 
\begin{gather*}
P_{1}(\rmF) \bsl U(\A) \ni x_{1} \mapsto 
\sum_{
\begin{subarray}{c}
\gamma \in (P_{2} \cap sP_{1})(\rmF) \bsl P_{2}(\rmF)\\
\tlgam \in (\tlP_{1} \cap \tls\tlP_{2})(\rmF) \bsl \tlP_{1}(\rmF)
\end{subarray}
}
k_{P_{1} \times \tlP_{2}, \chi}(x_{1} a_{1}k_{1} , w_{\tls}^{-1} \tlgam \tlx, w_{s}^{-1}\gamma y,  \tly), \\
\tlP_{2}(\rmF) \bsl \tlU(\A) \ni \tlx_{2} \mapsto 
\sum_{
\begin{subarray}{c}
\gamma \in (P_{2} \cap sP_{1})(\rmF) \bsl P_{2}(\rmF)\\
\tlgam \in (\tlP_{1} \cap \tls\tlP_{2})(\rmF) \bsl \tlP_{1}(\rmF)
\end{subarray}
}
k_{P_{1} \times \tlP_{2}, \chi}( x, w_{\tls}^{-1} \tlgam \tlx, w_{s}^{-1}\gamma y, \tlx_{2}  a_{2} k_{2}) 
\end{gather*}
sont invariantes à droite par l'intersection
du compact ouvert $\bigcap_{k \in K}(kK_{0}k^{-1})$  
avec $M_{1}(\A_{f})$ et $M_{\tltwo}(\A_{f})$ respectivement.
 On applique alors la version pour les groupes unitaires de la proposition
\ref{prop:opTronqD} pour les opérateurs $\Ladud^{T', P_{1}}$ et $\Latqd^{T, P_{2}}$. On trouve 
donc que pour tout $r_{1}, r_{2} \ge 0$ il existe un nombre fini des opérateurs différentiels $X, Y$
tels 
que l'expression \ref{eq:AllSigmasU} est majorée 
par la somme sur les $X, Y$ de  
\begin{multline*}
\int\limits_{\calM_{1}}
\int\limits_{\calM_{2}}\qquad \quad
\sup_{\mathmakebox[0.8cm]{
\begin{subarray}{c}
(m_{1}', m_{\tlone}') \in (M_{1} \times M_{\tlone})(\A)^{1} \\\
(m_{2}', m_{\tltwo}') \in (M_{2} \times M_{\tltwo})(\A)^{1}
\end{subarray}
}}
(\|m_{1}'\| \|m_{\tlone}'\|
\|m_{2}'\|\|m_{\tltwo}'\|)^{-r_{1}} 
\sum_{\chi}
\int\limits_{K \times K}
\int\limits_{A_{1}^{\infty}}
\int\limits_{A_{2}^{\infty}}
\int\limits_{[N_{1}]}
\int\limits_{[N_{2}]}
\\
\sigma_{1}^{3}(H_{1}(a_{1})-T')
\sigma_{2}^{4}(H_{\tltwo}(a_{2})-T)
e^{-2\rho_{1}(H_{1}(a_{1})) -2 \rho_{2}(H_{2}(a_{2}))  }
\sum_{(s, \tls) \in \Omega'}
\\
|
\sum_{
\begin{subarray}{c}
\gamma \in (P_{2} \cap sP_{1})(\rmF) \bsl P_{2}(\rmF)\\
\tlgam \in (\tlP_{1} \cap \tls\tlP_{2})(\rmF) \bsl \tlP_{1}(\rmF)
\end{subarray}
}
k_{F_{X,Y}, P_{1} \times \tlP_{2}, \chi}( a_{1}m_{1}'k_{1}, w_{\tls}^{-1} \tlgam n_{1} a_{1}m_{\tlone}'k_{1}, 
w_{s}^{-1}\gamma n_{2} a_{2}m_{2}' k_{2}, a_{2} m_{\tltwo}' k_{2}))
| \\
dn_{2}dn_{1}da_{2}da_{1}dk_{2}dk_{1}
)(\|m_{1}\|\|m_{2}\|)^{-r_{2}}
dm_{2}dm_{1}
\end{multline*}
où $\calM_{1} = \Sgl^{M_{1} \cap P_{0}}_{P_{0}} \cap M_{1}(\A)^{1}$, 
$\calM_{2} = \Sgl^{M_{2} \cap P_{0}}_{P_{0}} \cap M_{2}(\A)^{1}$
et l'on a remplacé 
la fonction $F$ par $F_{X,Y} := X \ast F \ast Y$ grâce à la formule (\ref{eq:diffOnKernel}) 
et en utilisant le même raisonnement que sur la page 104 de \cite{arthur5}. Remarquons que 
le support de $F_{X,Y}$ est contenu dans celui de $F$.

\blem\label{lem:preludeToGreatnessU}
 Avec la notation ci-dessus, il existe des constantes positives $c, N >0$ telles que
si pour un $(s, \tls) \in \Upomega_{\tlQ_{1}}'$ on a
\begin{multline*}
\sigma_{1}^{3}(H_{1}(a_{1})-T')
\sigma_{2}^{4}(H_{\tltwo}(a_{2})-T) \\
|
\sum_{
\begin{subarray}{c}
\gamma \in (P_{2} \cap sP_{1})(\rmF) \bsl P_{2}(\rmF)\\
\tlgam \in (\tlP_{1} \cap \tls\tlP_{2})(\rmF) \bsl \tlP_{1}(\rmF)
\end{subarray}
}
k_{F_{X,Y}, P_{1} \times \tlP_{2}, \chi}( a_{1}m_{1}'k_{1}, w_{\tls}^{-1} \tlgam n_{1} a_{1}m_{\tlone}'k_{1}, 
w_{s}^{-1}\gamma n_{2} a_{2}m_{2}' k_{2}, a_{2} m_{\tltwo}' k_{2}))
| \neq 0
\end{multline*}
alors
\[
\|a_{1}\|, \|a_{2}\| \le c
(\|m_{\tlone}'\|
\|m_{2}'\|)^{N}.
\]
\bdem
Le même raisonnement qu'au début de la preuve du lemme \ref{lem:preludeToGreatness}
montre qu'il existe une constante $C' >0$ telle que 
\begin{gather}
\varpi(H_{1}(a_{1})) - s\varpi(H_{2}(a_{2})) \le C'(1 + \log (\|m_{2}'\|) ), \quad 
\forall \varpi \in \hDelta_{1} \label{eq:varpiMinsvarpiThmU0}, \\
\varpi(H_{2}(a_{2})) - \tls\varpi(H_{1}(a_{1})) \le C'(1 + \log (\|m_{\tlone}'\|) ), \quad 
\forall \varpi \in \hDelta_{\tltwo}\label{eq:varpiMinsvarpiThmU1}.
\end{gather}
Puisque les éléments de $\hDelta_{1}$ égales ceux de $\hDelta_{\tlone}$ à $1/2$-près 
on remplace les inégalités (\ref{eq:varpiMinsvarpiThmU0}) par
\begin{equation}
\varpi(H_{1}(a_{1})) - s\varpi(H_{2}(a_{2})) \le C''(1 + \log (\|m_{2}'\|) ), \quad 
\forall \varpi \in \hDelta_{\tlone} \label{eq:varpiMinsvarpiThmU2}.
\end{equation}
pour une constante $C'' > 0$.

On est en mesure d'appliquer le corollaire \ref{cor:fuckYeah} de l'appendice \ref{app:lemme}. On l'applique pour les données suivantes. 
Pour le groupe $G$ on prend $\tlU$. Dans l'appendice \ref{app:lemme} on fixe un sous-groupe parabolique minimal 
de $\tlU$, on prend donc un $\tlB \in \calP(M_{\tlU, min})$ contenant $\tlP_{0}$. 
Pour les groupes $P_{1}$, $P_{2}$, $ P_{3}$, $P_{4}$ du corollaire
on prend $\tlP_{1}$, $\tlP_{2}$, $\tlP_{3}$ et $\tlP_{4}$ dans cet ordre.  
Pour les éléments $H_{1}$ et $H_{2}$ on prend 
$H_{1}(a_{1})$ et $H_{2}(a_{2})$ respectivement.
 Pour $X_{1}$ et $X_{2}$ on prend $T_{1}'$ et
 $T_{2}$ respectivement. 
 Pour les éléments du groupe de Weyl on prend 
 $s = s$ et  $s' = \tls$.
 La condition $(s,\tls) \in \Upomega_{\tlQ_{1}}'$ dit qu'ils vérifient la condition (\ref{eq:hDeltaEsesU}) du corollaire \ref{cor:fuckYeah}.
Les inégalités du corollaire \ref{cor:fuckYeah} 
 correspondent aux inégalités qu'on a construit de façon suivant:
 \[
\text{(\ref{eq:varpiMinsvarpi1U})} \leftrightarrow  
\text{(\ref{eq:varpiMinsvarpiThmU2})}, \quad 
\text{(\ref{eq:varpiMinsvarpi2U})} \leftrightarrow  
\text{(\ref{eq:varpiMinsvarpiThmU1})}.
 \] 
Les inégalités déterminent les constantes 
 $M_{1}$, $M_{2}$, $M_{3}$ et $M_{4}$. 
 Grâce au corollaire \ref{cor:fuckYeah} et la propriété (\ref{eq:haut7}) de la hauteur on conclut la preuve 
 du lemme \ref{lem:preludeToGreatness}.
\edem
\elem
La suite de la preuve est maintenant complètement analogue à 
la partie de la preuve du théorème \ref{thm:noyauCVGG} qui suit la preuve du lemme \ref{lem:preludeToGreatness}. 
\edem
\etheo

\subsection{Polynômes-exponentielles}\label{par:fonsPolExpU}
On utilisera le langage du paragraphe \ref{par:fonsPolExpSp}.

Pour tout 
sous-groupe parabolique standard $Q \in \calF(P_{0})$ de $U$ on pose 
$\upla_{Q} := 2\rho_{\tlQ} - 2\rho_{Q} \in \all_{\tlQ}^{*} = \all_{Q}^{*}$
(on rappelle qu'un sous-groupe parabolique standard $Q$ de $U$ définit le 
sous-groupe parabolique relativement standard $\tlQ$ uniquement).
On a alors le résultat suivant démontré dans \cite{leMoi}, lemme 4.2.
\blem\label{lem:positiveLinearU}
 Soit $Q$ un sous-groupe parabolique standard 
 de $U$. Alors, pour tout  $\varpi^{\vee} \in  \hDelta_{Q}^{\vee}$ 
on a $\upla_{Q}(\varpi^{\vee}) > 0$.
\elem

Le polynôme-exponentielle en question est donné par le lemme suivant.
\blem[cf. \cite{leMoi}, lemme 4.3]\label{lem:pQExplicitU}
 Soit $Q$ un sous-groupe parabolique standard de $U$.
Alors, pour tout $R \sps Q$ 
il existe un polynôme $P_{Q,R}$ de degré au plus $d_{Q}$ sur $\all_{R}$ 
tel que la fonction
\[
p_{Q}(X) := \int_{A_{Q}^{\infty}}e^{\upla_{Q}(H_{Q}(a))}\Gamma_{Q}'(H_{Q}(a),X)da, \quad 
X \in \all_{Q}
\]
égale
\[
\sum_{R \sps Q}e^{\upla_{R}(X_{R})}
p_{Q,R}(X_{R})
\]
où $p_{Q,U}(X_{U}) = (-1)^{d_{Q}}\hat \theta_{Q}(\upla_{Q})^{-1}$.
En particulier, la fonction $p_{Q}$ est un polynôme-exponentielle
dont le terme purement polynomial est constant 
et égale $(-1)^{d_{Q}}\hat \theta_{Q}(\upla_{Q})^{-1}$.
\elem

\subsection{Une généralisation du théorème \ref{thm:noyauCVGU}}\label{par:JTforLevisU}

Dans le cas unitaire le théorème \ref{thm:noyauCVGU} dit que 
pour $T \in T_{+} + \all_{0}^{+}$ et $\chi \in \calX^{U \times \tlU}$, 
la distribution 
\[
C_{c}^{\infty}((U \times \tlU)(\A)) \ni F \mapsto J_{\chi}^{T}(F) := 
\int_{[U]}\int_{[U]}
k_{F,\chi}^{T}(x,x,y,y)dxdy
\]
est bien définie. 

Soit $H = \prod_{i=1}^{k}\Res_{\rmE/\rmF}\Gl_{n_{i}}$ où $k \in \N$ et 
$n_{i} \in \N^{*}$ pour $i = 1, \ldots ,k$. Soient aussi $U' = U(V', \Phi')$ et 
$\tlU = U(W', \tlPhi')$ un couple de groupes unitaires munis de l'inclusion 
$U' \hrar \tlU'$ comme dans le paragraphe 
\ref{par:unitGp}. 
On va généraliser le théorème \ref{thm:noyauCVGU} au cas de l'inclusion 
$H \times U' \hrar H \times \tlU'$.

Fixons $P_{0}$ un sous-$\rmF$-groupe parabolique minimal 
de $H \times U'$ et fixons aussi $M_{0}$ une partie de Levi de $P_{0}$ définie sur $\rmF$. 
Notons $M_{\tlzero}$ le sous-$\rmF$-groupe de Levi de $H \times \tlU'$ 
contenant $M_{0}$, minimal pour cette propriété. Alors, $M_{\tlzero}$ est uniquement déterminé 
par $M_{0}$. 
Soit $P \in \calF(P_{0})$ un sous-groupe parabolique standard de $G \times U'$. 
Il résulte des discussions dans les paragraphes \ref{par:unitGp} 
et \ref{par:GinGE} qu'il existe un unique sous-groupe parabolique 
$\tlP \in \calF(M_{\tlzero})$ tel que $\tlP \cap (H \times U') = P$. 
On le note justement $\tlP$.

Pour une fonction 
$f \in C_{c}^{\infty}((H \times U')(\A)^{1} \times (H \times \tlU')(\A)^{1})$,  
un sous-groupe parabolique standard $P$ de $H \times U'$ 
et une donnée cuspidale $\chi \in \calX^{(H \times U') \times (H \times \tlU')}$ on pose
\[
k_{f,P, \chi}(\tlx, \tly) := 
k_{f,P \times \tlP, \chi}(\tlx, \tly), \quad 
\tlx, \tly \in N_{P \times \tlP}(\A)M_{P \times \tlP}(\rmF) \bsl (H \times U')(\A)^{1} \times (H \times \tlU')(\A)^{1}.
\]
Pour un $T \in \mathfrak{a}_{P_{0}}^{+}$ on pose
\[
k^{T}_{f,\chi}(x,\tlx, y, \tly) = 
\sum_{P \sps P_{0}} (-1)^{d_{P}^{H \times U'}} \ 
\sum_{\mathclap{\delta_{1}, \delta_{2} \in P(\rmF) \bsl (H \times U')(\rmF)}} \ 
\htau_{P}^{G \times U'}(H_{P}(\delta_{2} y)-T)
k_{f,P,\chi}(\delta_{1} x, \delta_{1} \tlx, \delta_{2} y, \delta_{2} \tly)
\]
où $x,y \in [H \times U']^{1}$ et $\tlx, \tly \in [H \times \tlU']^{1}$.

\btheo Soit $f \in C_{c}^{\infty}((H \times U')(\A)^{1} \times (H \times \tlU')(\A)^{1})$, 
alors pour tout $T \in  \mathfrak{a}_{P_{0}}^{+}$ assez régulier on a
\[
\sum_{\chi \in \calX^{(H \times U') \times (H \times \tlU')}}
\int_{[H \times U']^{1}} \int_{[H \times U']^{1}}
|k_{f,\chi}^{T}(x,x,y,y)|dxdy < \infty.
\]
\bdem
La preuve est similaire à celle du théorème \ref{thm:noyauCVGU}.
Les détails sont laissés au lecteur. 
\edem 
\etheo

Notons alors pour $\chi \in \calX^{(H \times U') \times (H \times \tlU')}$ 
et $f \in C_{c}^{\infty}((H \times U')(\A)^{1} \times (H \times \tlU')(\A)^{1})$
\[
J_{\chi}^{H \times U',T}(f) = \int_{[H \times U']^{1}} \int_{[H \times U']^{1}}
k_{f,\chi}^{T}(x,x,y,y)dxdy.
\]

Revenons au groupe $U$.
Soit $Q$ un sous-groupe parabolique standard de $U$. 
L'inclusion $M_{Q} \hrar M_{\tlQ}$ vérifie alors les conditions de ce paragraphe.

Soit $\chi \in \calX^{U \times \tlU}$, et $\{\chi'\} \sbs
\calX^{M_{Q} \times M_{\tlQ}}$ la préimage de 
$\chi$ par l'application naturelle $\calX^{M_{Q} \times M_{\tlQ}} \rar \calX^{U \times \tlU}$.
Pour $T \in \all_{0}^{+}$ assez régulier on définit 
alors la distribution $J_{\chi}^{M_{Q},T}(\cdot)$ sur 
$C_{c}^{\infty}(M_{Q}(\A)^{1} \times M_{\tlQ}(\A)^{1})$ par:
\begin{equation*}
J^{M_{Q},T}_{\chi}(f) = 
\sum_{\chi'}
J^{M_{Q},T}_{\chi'}(f),
\end{equation*}
où 
pour $\chi' \in \calX^{M_{Q} \times M_{\tlQ}}$, 
$J^{M_{\tlQ},T}_{\chi'}(\cdot)$ 
c'est la distribution associée à l'inclusion $M_{Q} \hrar M_{\tlQ}$ 
décrite ci-dessus par rapport au sous-groupe de Levi minimal $M_{0}$ de $M_{Q}$
et au sous-groupe parabolique minimal $P_{0} \cap M_{Q}$ de $M_{Q}$.

Pour $F \in C_{c}^{\infty}((U \times \tlU)(\A))$ on
définit pour $x \in M_{Q}(\A)^{1}$ et $\tlx \in M_{\tlQ}(\A)^{1}$
\begin{equation}\label{eq:fQdefU}
F_{Q}(x,\tlx) = 
\int\limits_{K \times K}
\int\limits_{(N_{Q} \times N_{\tlQ})(\A)}
\int\limits_{A_{Q}^{\infty}}
e^{2\rho_{Q}(H_{Q}(a_{1}))} 
F(k_{1}^{-1}a_{1} xn_{Q}k_{2}, k_{1}^{-1}a_{1} \tlx n_{\tlQ}k_{2})
da_{1}dn_{\tlQ}dn_{Q}dk_{2}dk_{1}.
\end{equation}
Alors $F_{Q} \in C_{c}^{\infty}(M_{Q}(\A)^{1} \times M_{\tlQ}(\A)^{1})$. 

Notons que l'application 
\[
Q \sps P \mapsto M_{Q} \cap P
\]
définit une bijection entre les sous-groupes paraboliques
standards de $U$
contenus dans $Q$ et les sous-groupes paraboliques 
de
$M_{Q}$ contenant $P_{0} \cap M_{Q}$. 
Soit donc $\chi \in \calX^{U \times \tlU}$ et 
soient $\{\chi_{Q}\} \in \calX^{M_{Q} \times M_{\tlQ}}$ qui s'envoient sur $\chi$ par l'application naturelle 
$\calX^{M_{Q} \times M_{\tlQ}} \rar \calX^{U \times \tlU}$.
En utilisant le lemme \ref{lem:kChiDescente}, 
et en remarquant que si l'on écrit $2\rho_{Q} = (\rho_{Q} + \rho_{\tlQ}) - (\rho_{\tlQ} - \rho_{Q})$ 
on a $ - (\rho_{Q} + \rho_{\tlQ}) - (\rho_{\tlQ} - \rho_{Q})  =-2\rho_{\tlQ}$, 
on s'aperçoit que pour tout $P \in \calF(P_{0})$ contenu dans $Q$, 
tous
$x,y \in M_{Q}(\A)^{1}$ 
et tous
$\tlx,\tly \in M_{\tlQ}(\A)^{1}$ on a 
\begin{equation}\label{eq:kfQisThisU}
\sum_{\chi_{Q}}k_{F_{Q}, M_{Q} \cap P, \chi_{Q}}(x, \tlx ,y, \tly) = 
\int\limits_{K \times K}
\int\limits_{A_{Q}^{\infty}}
e^{-2\rho_{\tlQ}(H_{Q}(a_{1}))} 
k_{F, P, \chi}(x k_{1}, \tlx k_{1}, a_{1}y k_{2}, a_{1} \tly k_{2})
da_{1}dk_{2}dk_{1}.
\end{equation}

\subsection{Comportement en $T$}\label{par:asymptChaptU}

On démontre la proposition suivante.

\brop\label{prop:mainQualitPropU}
Soient $F \in C_{c}^{\infty}((U \times \tlU)(\A)$, $T' \in T_{+} + \all_{0}^{+}$, 
$\chi \in \calX^{U \times \tlU}$ et $T \in T' + \all_{0}^{+}$. Alors
\[
J_{\chi}^{T}(F) = \sum_{Q \sps P_{0}}
p_{Q}(T_{Q} - T'_{Q})e^{\upla_{Q}(T'_{Q})}J_{\chi}^{M_{Q}, T'}(F_{Q})
\]
où pour un sous-groupe parabolique $Q$ 
standard de $U$, 
la fonction $p_{Q}$ est définie dans le lemme \ref{lem:pQExplicitU} 
dans le paragraphe \ref{par:fonsPolExpU} où 
$\upla_{Q} \in \all_{Q}^{*}$ est aussi défini,
la distribution $J_{\chi}^{M_{Q}, T'}$ 
est définie dans le paragraphe \ref{par:JTforLevisU} et 
$F_{Q} \in C_{c}^{\infty}((M_{Q} \times M_{\tlQ})(\A)^{1})$ 
est définie par (\ref{eq:fQdefU}) 
dans le même paragraphe.

\bdem
Fixons un \(T' \in T_{+}+\mathfrak{a}_{0}^{+} \) 
et soit  $T \in T' + \all_{0}^{+}$.
En utilisant la relation (\ref{eq:GammaRecurr}) 
dans la définition du noyau $k_{F, \chi}^{T}$ (\ref{eq:kchiUdef})
avec $P = P$, 
$H = H_{P}(\delta_{2} y)-T'_{P}$ et
$X = T_{P} - T'_{P}$ 
pour tout $P \in \calF(P_{0})$ et tout $\delta_{2} \in P(\rmF) \bsl U(\rmF)$ 
on a que $J_{\chi}^{T}(F)$ égale la somme sur $Q \in \calF(P_{0})$ de
\begin{multline}\label{eq:firstUass}
\int\limits_{Q(\rmF) \bsl U(\A)}
\int\limits_{Q(\rmF) \bsl U(\A)}
\Gamma_{Q}(H_{Q}(y) - T_{Q}', T_{Q} - T_{Q}') \\
\sum_{P \sbs Q}(-1)^{d_{P}^{Q}}
\sum_{\delta_{1}, \delta_{2} \in (P \cap M_{Q})(\rmF)\bsl M_{Q}(\rmF)}
\htau_{P}^{Q}(H_{P}(\delta_{2} y)-T_{\tlP})
k_{P,\chi}(\delta_{1} x, \delta_{1} x, \delta_{2}y, \delta_{2}y)
dxdy.
\end{multline}

Écrivons maintenant $x = n_{1}a_{1}m_{1}k_{1}$ et $y = n_{2}a_{2}m_{2}k_{2}$
selon la décomposition 
\begin{equation*}
Q(\rmF) \bsl U(\A) = [N_{Q}] \times A_{Q}^{\infty} \times [M_{Q}]^{1} \times K.
\end{equation*}

On obtient alors que (\ref{eq:firstUass}) égale
\begin{multline*}
\int\limits_{[M_{Q}]^{1} \times [M_{Q}]^{1}}
\int\limits_{(A_{Q}^{\infty})^{2}}
\int\limits_{\mathrlap{K \times K}}
e^{-2\rho_{Q}(H_{Q}(a_{1}a_{2}))}
\Gamma_{Q}'(H_{Q}(a_{2}) - T_{Q}', T_{Q} - T_{Q}')
\sum_{P \sbs Q}(-1)^{d_{P}^{Q}} \quad
\sum_{\mathclap{\delta_{1}, \delta_{2} \in (P \cap M_{Q})(\rmF)\bsl M_{Q}(\rmF)}}\\
\htau_{\tlP}^{\tlQ}(H_{\tlP}(\delta_{2}m_{2})-T_{\tlP})
k_{P,\chi}(\delta_{1} m_{1}a_{1}k_{1}, \delta_{1} m_{1}a_{1}k_{1}, \delta_{2}m_{2}a_{2}k_{2}, \delta_{2}m_{2}a_{2}k_{2}) 
dk_{1}dk_{2} da_{1} da_{2}dm_{1}dm_{2}.
\end{multline*}
En utilisant le lemme \ref{lem:artLemme23} on voit que
$k_{P,\chi}(\delta_{1} m_{1}a_{1}k_{1}, \delta_{1} m_{1}a_{1}k_{1}, \delta_{2}m_{2}a_{2}k_{2}, \delta_{2}m_{2}a_{2}k_{2})$ 
égale
\[ 
e^{2(\rho_{Q} + \rho_{\tlQ})(H_{Q}(a_{1}))}
k_{P,\chi}(\delta_{1} m_{1}k_{1}, \delta_{1} m_{1}k_{1}, \delta_{2}m_{2}a_{1}^{-1}a_{2}k_{2}, \delta_{2}m_{2}a_{1}^{-1}a_{2}k_{2}).
\]
En faisant le changement de variable 
$a_{2}^{-1}a_{1} \mapsto a_{1}$ et en utilisant le lemme \ref{lem:pQExplicitU}, on trouve
que l'intégration par rapport à $da_{2}$ devient
\[
\int\limits_{A_{Q}^{\infty}}
e^{2(\rho_{\tlQ} - \rho_{Q})(H_{Q}(a_{2}))}
\Gamma_{Q}'(H_{Q}(a_{2}) - T_{Q}', T_{Q} - T_{Q}')da_{2}= 
e^{\upla_{Q}(T_{Q}')}
p_{Q}(T_{Q} - T_{Q}').
\]
Le reste d'intégration c'est alors $e^{\upla_{Q}(T_{Q}')}p_{Q}(T_{Q} - T_{Q}')$ fois
\begin{multline*}
\int\limits_{[M_{Q}]^{1} \times [M_{Q}]^{1}}
\int\limits_{(A_{Q}^{\infty})}
\int\limits_{K \times K}
\sum_{P \sbs Q}(-1)^{d_{P}^{Q}}
\sum_{\delta_{1}, \delta_{2} \in (P \cap M_{Q})(\rmF)\bsl M_{Q}(\rmF)}
\htau_{\tlP}^{\tlQ}(H_{\tlP}(\delta_{2}m_{2})-T_{\tlP})\\
e^{-2\rho_{\tlQ}(H_{Q}(a_{1}))}
k_{P,\chi}(\delta_{1} m_{1}k_{1}, \delta_{1} m_{1}k_{1}, \delta_{2}m_{2}a_{1}k_{2}, \delta_{2}m_{2}a_{1}k_{2}) 
dk_{1}dk_{2} da_{1} dm_{1}dm_{2}.
\end{multline*}
En utilisant la discussion du paragraphe \ref{par:JTforLevisU}, l'équation 
(\ref{eq:kfQisThisU}) en l'occurrence, on trouve que l'expression ci-dessus égale
$J_{\chi}^{M_{\tlQ}, T}(F_{Q})$ où $F_{Q}$ est définie 
par (\ref{eq:fQdefU}), ce qui termine la preuve de la proposition. 
\edem
\erop

En utilisant la proposition \ref{prop:mainQualitPropU} démontrée 
ci-dessus et le lemme \ref{lem:pQExplicitU} qui 
décrit les fonctions $p_{Q}$ explicitement on obtient 
le comportement en $T$ de la distribution $J_{\chi}^{T}$.

\begin{theo}\label{thm:mainQualitThmU}
Soient $F \in C_{c}^{\infty}((U \times \tlU)(\A))$ et $\chi \in \calX^{U \times \tlU}$.
La fonction
$T \mapsto J^{T}_{\chi}(F)$  
 où
$T$ 
parcourt $T_{+} + \mathfrak{a}_{0}^{+}$ 
est un polynôme-exponentielle
dont la partie purement polynomiale est constante
et donnée par  
\begin{equation*}
J_{\chi}(F) := \sum_{Q \sps P_{0}}(-1)^{d_{Q}}\hat \theta_{Q}(\upla_{Q})^{-1}
e^{\upla_{Q}(T_{Q}')}J_{\chi}^{M_{Q},T'}(F_{Q})
\end{equation*}
pour tout $T' \in T_{+} + \all_{0}^{+}$. 
En particulier, la distribution $J_{\chi}$ 
ne dépend pas de $T'$.
\end{theo}

\brem\label{rem:JMQisPolExpU}
 Soit $Q$ un sous-groupe parabolique standard 
de $U$. 
Par le même raisonnement que dans la proposition \ref{prop:mainQualitPropU}
 on obtient que la 
distribution $J_{\chi}^{M_{Q},T}$ 
définie dans le paragraphe \ref{par:JTforLevisU}, 
est un polynôme-exponentielle en $T$ qui ne dépend pas
de $T_{Q} \in \all_{Q}$. Cependant, 
si $Q \neq U$ le terme purement polynomial n'est pas constant. 
\erem

\subsection{Invariance}\label{par:compConjugU}

Soient $F \in C_{c}^{\infty}((U \times \tlU)(\A))$ et $g,h  \in U(\A)$. 
Notons $F^{g,h} \in C_{c}^{\infty}((U \times \tlU)(\A))$ la fonction 
définie par $F^{g,h}(x,\tlx) = F(gxh^{-1}, g\tlx h^{-1})$. 

Soit $\chi \in \calX^{U \times \tlU}$.
On voit que $J_{\chi}^{T}(F^{g,h})$ pour 
$T \in T_{+} + \all_{\tlzero}^{+}$ égale
\begin{equation*}
\int_{[U]}\int_{[U]}
\sum_{P \sps P_{0}}(-1)^{d_{P}}
\sum_{\delta_{1},\delta_{2} \in P(\rmF)\bsl U(\rmF)}
\htau_{P}(H_{P}(\delta_{2}yh)-T_{P})
k_{P,\chi}(\delta_{1} x, \delta_{1} x, \delta_{2} y, \delta_{2} y)dxdy.
\end{equation*}

Pour $y \in U(\A)$ et $P \in \calF(P_{0})$ soit $k_{P}(y)$ un élément 
de $K$ tel que $yk_{P}(y)^{-1} \in P(\A)$. 
Alors, en utilisant l'égalité (\ref{eq:GammaRecurr})
on a:
\begin{equation*}
\htau_{P}(H_{P}(\delta_{2} y h)-T_{P}) =
\sum_{Q \supseteq P}(-1)^{d_{Q}}
\htau_{P}^{Q}(H_{P}(\delta_{2} y)-T_{P})
\Gamma_{Q}'(H_{Q}(\delta_{2} y)-T_{Q},
-H_{Q}(k_{Q}(\delta_{2} y)h))
\end{equation*}
d'où on obtient que $J_{\chi}^{T}(F^{g,h})$ égale la 
somme sur $Q \in \calF(P_{0})$ de 
\begin{multline}\label{eq:equivU1}
\int\limits_{Q(\rmF) \bsl U(\A)}
\int\limits_{Q(\rmF) \bsl U(\A)}
\Gamma_{Q}'(H_{Q}(y)-T_{Q}, -H_{Q}(k_{Q}(y)h))
\sum_{P \sbs Q}(-1)^{d_{P}^{Q}}
\sum_{ \delta_{1}, \delta_{2} \in (P \cap M_{Q})(\rmF)\bsl M_{Q}(\rmF)} \\
\htau_{P}^{Q}(H_{P}(\delta_{2}y)-T_{P})
k_{P,\chi}(\delta_{1} x, \delta_{1} x, \delta_{2}y, \delta_{2}y)dxdy.
\end{multline}

Écrivons maintenant $x = n_{1}a_{1}m_{1}k_{1}$ et $y = n_{2}a_{2}m_{2}k_{2}$ 
selon la décomposition 
\begin{equation*}
Q(\rmF) \bsl U(\A) = [N_{Q}] \times A_{Q}^{\infty} \times [M_{Q}]^{1} \times K.
\end{equation*}
On a alors:
\begin{equation*}
\Gamma_{Q}'(H_{Q}(n_{2}a_{2}m_{2}k_{2})-T_{Q}, -H_{Q}(k_{Q}(n_{2}a_{2}m_{2}k_{2})h))
 = 
\Gamma_{Q}'(H_{Q}(a_{2})-T_{Q}, -H_{Q}(k_{2}h)).
\end{equation*}

On fait le même changement de variable que dans la
preuve
de la proposition \ref{prop:mainQualitPropU}. 
On obtient, en vertu du lemme \ref{lem:pQExplicitU}, 
que l'intégration par rapport à $da_{2}$ devient:
\[
\int_{A_{Q}^{\infty}}
e^{\upla_{Q}(H_{Q}(a_{2})) }
\Gamma_{Q}'(H_{Q}(a_{2})-T_{Q}, -H_{Q}(k_{2}h))da_{2}= 
e^{\upla_{Q}(T_{Q})} 
u_{Q,h}(k_{2})
\]
où 
\begin{equation*}
u_{Q,h}(k) = \int_{\all_{Q}}
e^{\upla_{Q}(H)}\Gamma_{Q}'(H,-H_{Q}(kh))dH, \quad k \in K
\end{equation*}
est une fonction continue de la variable $k \in K$ en vertu du lemme \ref{lem:GammaIsComp}.

Le reste de l'intégration dans (\ref{eq:equivG1}) c'est
$e^{\upla_{Q}(T_{Q})}$ fois
\begin{multline}\label{eq:equivU2}
\int\limits_{[M_{Q}]^{1}}
\int\limits_{[M_{Q}]^{1}}
\int\limits_{A_{Q}^{\infty}}
\int\limits_{K \times \tlK}
\sum_{P \sbs Q}(-1)^{d_{P}^{Q}}
\sum_{ \delta_{1}, \delta_{2} \in (P \cap M_{Q})(\rmF)\bsl M_{Q}(\rmF)} 
\htau_{P}^{Q}(H_{P}(\delta_{2}y)-T_{P})\\
e^{-2\rho_{\tlQ}(H_{Q}(a_{1}))}
k_{P,\chi}(\delta_{1} m_{1}k_{1}, \delta_{1} m_{1}k_{1}, \delta_{2}m_{2}a_{1}k_{2}, 
\delta_{2}m_{2}a_{1}k_{2})u_{Q,h}(k_{2})
dk_{1}dk_{2}da_{1} dm_{1}dm_{2}.
\end{multline}
Si l'on pose alors
\begin{multline*}\label{eq:fQdefGy}
F_{Q,h}(x,\tlx) = 
\int\limits_{K \times \tlK}
\int\limits_{(N_{Q} \times N_{\tlQ})(\A)}
\int\limits_{A_{\tlQ}^{\infty}}
e^{2\rho_{Q}(H_{Q}(a_{1}))} 
F(k_{1}^{-1}a_{1} x n_{Q}k_{2}, k_{1}^{-1}a_{1} \tlx n_{\tlQ} k_{2})\\
u_{\tlQ,h}(k_{2})
da_{1}dn_{\tlQ}dn_{Q}dk_{2}dk_{1}, \quad 
x \in M_{Q}(\A)^{1}, \ 
\tlx \in M_{\tlQ}(\A)^{1}
\end{multline*}
on a bien $F_{Q,h} \in C_{c}^{\infty}(M_{Q}(\A)^{1} \times M_{\tlQ}(\A)^{1})$
et on voit, 
en se basant sur la discussion du paragraphe \ref{par:JTforLevisU} 
et en utilisant le lemme \ref{lem:kChiDescente} du paragraphe \ref{par:repRegul},
que (\ref{eq:equivU2}) c'est juste $J_{\chi}^{M_{Q},T}(F_{Q,h})$.
On trouve alors le théorème suivant.

\begin{theo}\label{thm:invarianceTheoU} Soient $g,h \in U(\A)$ et $\chi \in \calX^{U\times \tlU}$. 
La distribution $J_{\chi}^{T}(\cdot)$ vérifie pour tout 
$F \in C_{c}^{\infty}((U \times \tlU)(\A))$ la propriété suivante:
\begin{equation*}
J_{\chi}^{T}(F^{g,h}) - J_{\chi}^{T}(F)  = 
\sum_{Q \in \calF(P_{0}) \smin \{U\}} 
e^{\upla_{Q}(T_{Q})}
J_{\chi}^{M_{Q},T}(F_{Q,h})
\end{equation*}
où les distributions $J_{\chi}^{M_{Q},T}$ 
sont définies dans la section \ref{par:JTforLevisU}. 
En particulier, on a
\[
J_{\chi}(F^{g,h}) = J_{\chi}(F).
\]

\bdem 
L'argument est complètement analogue a celui de la preuve du théorème \ref{thm:invarianceTheoG} ci-dessus
et repose sur le théorème \ref{thm:mainQualitThmU}, sa remarque \ref{rem:JMQisPolExpU} 
et la propriété des exposants $\upla_{Q}$ donnée par le lemme \ref{lem:positiveLinearU}.
\edem
\end{theo}

\subsection{Données cuspidales non-induites}\label{par:exemplDistU}

Dans ce paragraphe on étudie la distribution $J_{\chi}$
pour des classes $\chi$ distinguées. On utilisera le langage du paragraphe \ref{par:exemplDistG}.

\brop\label{prop:DistSimplU} Soient $F \in C_{c}^{\infty}((U \times \tlU)(\A))$,
$\chi^{U} \in \calX^{U}$ et $\chi^{\tlU} \in \calX^{\tlU}$. Notons $\chi = \chi^{U} \times \chi^{\tlU} \in \calX^{U \times \tlU}$.
Supposons que $\chi^{U}$ ou $\chi^{\tlU}$ est non-induite. Alors, pour tout 
$T \in T_{+} + \all_{0}^{+}$ on a
\[
J_{\chi}(F) = J_{\chi}^{T}(F) = \int\limits_{\mathclap{[U]}}\int\limits_{\mathrlap{[U]}}
k_{F, \chi}(x,x, y,y)dxdy = 
\int\limits_{\mathclap{[U]}}\int\limits_{\mathrlap{[U]}} \La_{d,12}^{T}\La_{d, 34}^{T}k_{F, \chi}(x,x, y,y)dxdy.
\]
\bdem 
Soit $P \in \calF(P_{0})$. Si $P \neq U$ on a 
$k_{F, P, \chi} \equiv 0$ 
d'où $k_{F, \chi} = k_{F, \chi}^{T}$ 
et les deux premières égalités découlent du théorème \ref{thm:noyauCVGU}.

Démontrons la dernière égalité. 
En utilisant la formule du lemme \ref{lem:formDInvD}
on trouve que $J_{\chi}(F)$ égale la somme sur $P, Q \in \calF(P_{0})$ de
\[
\int\limits_{P(\rmF) \bsl U(\A)}\int\limits_{Q(\rmF) \bsl U(\A)}
\tau_{P}(H_{P}(x) - T)\tau_{Q}(H_{Q}(y) - T)\La_{d, 12}^{P, T}\La_{d, 34}^{Q, T}k_{F, \chi}(x,x, y,y)dxdy.
\]
La convergence absolue des intégrales découle de la convergence de l'expression (\ref{eq:AllSigmasU}) 
dans la preuve du théorème \ref{thm:noyauCVGU}. Or, 
si $P \neq U$ ou $Q \neq U$ on a 
$\La_{d, 12}^{P, T}\La_{d, 34}^{Q, T}k_{F, \chi}\equiv 0$ 
par hypothèse sur $\chi$, d'où le résultat.
\edem
\erop

\subsection{Indépendance des choix}

Soient $s, s' \in \C$ tels que $s + s' \neq 1$.
Dans ce paragraphe on constate que les distributions $I_{\chi}(s,s', \cdot)$  
et $J_{\chi}$
ne dépendent d'aucun choix, sauf le choix d'une mesure de Haar 
sur $G(\A)$, $G_{\rmE}(\A)$, $\tlG(\A)$, $\tlG_{\rmE}(\A)$, $U(\A)$, 
$\tlU(\A)$ et les choix des mesures sur les points adéliques des parties unipotentes $N_{P}(\A)$ des $\rmF$-sous-groupes 
paraboliques de tous les groupes en question, le choix étant que 
le volume $N_{P}(\rmF) \bsl N_{P}(\A)$ soit $1$.

En effet, l'indépendance vient des théorèmes \ref{thm:mainQualitThmG} 
et \ref{thm:mainQualitThmU} et se démontre de même façon que 
dans le paragraphe 4.5 de \cite{leMoi} ou bien le paragraphe 3.5 de \cite{leMoi2} 
où démontre la même propriété pour les distributions 
géométriques. Par exemple, pour démontrer l'invariance de la distribution 
$I_{\chi}(s,s', \cdot)$ du choix d'un sous-groupe de Levi minimal $M_{0}$ de $G$, on prend 
$\gamma \in G(\rmF)$ et $M_{0}' = \gamma M_{0} \gamma^{-1}$ et on 
note $I_{M_{0}', \chi}$ la distribution associée
à $M_{0}'$ et les autres données conjugués par $\gamma$. 
À l'aide du théorème \ref{thm:mainQualitThmG} on trouve donc que 
pour tout $\Phi \in C_{c}^{\infty}(\GetlGa)$ 
on a $I_{M_{0}', \chi}(s,s', \Phi) = I_{\chi}(s,s', \Phi^{\gamma, (\gamma, \gamma)})$
et le résultat suit du théorème \ref{thm:invarianceTheoG} car 
$|\det \gamma|_{\A}^{s}|\det \gamma|_{\A}^{s'} \eta(\gamma, \gamma) = 1$. 

\section{Formule des traces relative de Jacquet-Rallis pour les groupes linéaires}\label{sec:RTFG}

\subsection{Les invariants}\label{par:invpar}

Soit $\tlgl_{\rmE} = \Res_{\rmE/\rmF}(\End_{\rmE}(W_{\rmE}))$. 
On note dans ce paragraphe $\A_{\rmE}$  - la $\rmE$-droite affine. 
Notons 
\begin{equation}\label{eq:invMap}
Q : \tlgl_{\rmE} \rar \Res_{\rmE/\rmF}(\A_{\rmE})^{2n+1}
\end{equation}
 l'application suivante: 
pour tout $X \in \tlgl_{\rmE}$ on a $Q(X) = (Q(X)_{i})_{i=1,\ldots, 2n+1}$ 
où $Q(X)_{i} = \Tr \bigwedge^{i} X$ pour $i = 1, \ldots n+1$ et 
$Q(X)_{i} = e_{0}^{*}(X^{i - (n+1)}e_{0})$ pour $i = n+2, \ldots, 2n+1$ 
où $e_{0}^{*} \in \Res_{\rmE/\rmF}(W_{\rmE}^{*})$ est définie par $e_{0}^{*}(V_{\rmE}) = 0$ 
et $e_{0}^{*}(e_{0}) = 1$. Alors, $Q$ est $G_{\rmE}$-invariante pour l'action de 
$G_{\rmE}$ sur $\tlgl_{\rmE}$ par adjonction. 
Soit $\calO$ l'ensemble des sous-$\rmF$-variétés de $\tlgl_{\rmE}$ définies comme des pré-images 
par $Q$
des points $\rmF$-rationnels dans $\Res_{\rmE/\rmF}(\A_{\rmE})^{2n+1}$. En particulier, 
à tout $\ol \in \calO$ on peut associer, et on le fait, 
un polynôme $P_{\ol} \in \rmE[T]$ qui est le polynôme caractéristique commun à tous les éléments de $\ol$.

Pour tout $\xi \in \rmE$ soit $\xi I \in \tlgl_{\rmE}(\rmF)$ l'homothétie de rapport le scalaire $\xi$. 
On note alors
\[
D_{\xi} = \Res_{\rmE/\rmF}(\{ X \in \End_{\rmE}(W_{\rmE})| \det (X - \xi I) = 0\}).
\]
Pour une sous-$\rmF$-variété $\calV$ de $\tlgl_{\rmE}$, par $\calV \smin D_{\xi}$ on entend toujours l'ouvert complément de $\calV \cap D_{\xi}$ 
dans $\calV$ au sens de schémas algébriques. On a donc 
pour toute $\rmF$-algèbre $R$:
\[
(\calV \smin D_{\xi})(R) = \{X \in \calV(R) \sbs \End_{\rmE}(W_{\rmE})(\rmE \otimes_{\rmF} R) | \det (X - \xi I) \in (\rmE \otimes_{\rmF} R)^{*}\}.
\]
Soit $\ol \in \calO$.
Remarquons que, puisque tous les éléments de $\ol$ ont le même polynôme caractéristique, pour tout $\xi \in \rmE$,  
on a soit $\ol \sbs D_{\xi}$ soit $\ol \sbs \tlgl_{\rmE} \smin D_{\xi}$.

Pour tout $\xi \in \rmE$ on définit l'application, dite de Cayley, 
$\kappa_{\xi} : \tlgl_{\rmE} \smin D_{1} \rar  \tlgl_{\rmE} \smin D_{\xi}$ de façon suivante:
\begin{equation}\label{eq:cayleyDef}
\kappa_{\xi}(X) = -\xi(1+X)(1-X)^{-1}.
\end{equation}
L'application $\kappa_{\xi}$ a été introduite dans le contexte de la formule des traces relative de Jacquet-Rallis 
dans \cite{zhang2}, section 3. Voici quelques-une des ses propriétés.
\blem\label{lem:cayleyDef} Soit $\xi \in \rmE$.
\begin{enumerate}[1)]
\item L'application $\kappa_{\xi} : \tlgl_{\rmE} \smin D_{1} \rar  \tlgl_{\rmE} \smin D_{\xi}$ définie par 
(\ref{eq:cayleyDef}) ci-dessus est bien définie et est un isomorphisme des $\rmF$-variétés algébriques 
son inverse étant $\kappa_{\xi}^{-1}(X) = -(\xi + X)(\xi - X)^{-1}$.
\item Pour tout $\ol \in \calO$ disjoint avec $D_{1}$ 
il existe un unique $\ol' \in \calO$ disjoint avec $D_{\xi}$ 
tel que la restriction de $\kappa_{\xi}$ à $\ol$   
induit un isomorphisme 
entre $\ol$ et $\ol'$.
\end{enumerate}
\bdem
La première assertion découle du fait que $\kappa_{\xi}$ est algébrique 
et du fait qu'en effet son inverse est donnée par la formule $X \mapsto -(\xi + X)(\xi - X)^{-1}$, ce qu'on vérifie facilement.
La deuxième assertion c'est le lemme 3.5 de \cite{zhang2}. 
\edem
\elem
  
Notons $\Gamma_{\rmE/\rmF}$ le groupe de Galois 
de l'extension $\rmE/\rmF$ et soit $\sigma$ son générateur. 
On considère à la fois $\tlgl_{\rmE}$ comme une variété contenant $\tlG_{\rmE}$ et 
comme l'algèbre de Lie de $\tlG_{\rmE}$. Les deux points de vue sont complètement 
compatibles. 
Le groupe $\Gal(\rmE/\rmF)$ agit sur la variété $\tlgl_{\rmE}$ grâce à son action naturelle sur $W_{\rmE} = W \otimes_{\rmF} \rmE$. 
Pour $X \in \tlgl_{\rmE}$ on note donc $\brX\in \tlG_{\rmE}$ défini 
pour tout $v \in W_{\rmE}$ comme $\brX v = \sigma(X \sigma(v))$. 
On a donc $\tlG = \tlG_{\rmE}^{\Gamma_{\rmE/\rmF}}$.
Soit $S_{W}$ la sous-$\rmF$-variété de $\tlG_{\rmE}$ définie comme 
\[
S_{W} = \{g \in \tlG_{\rmE}| g \brg = 1\}.
\]
Alors $S_{W}$ est une $\rmF$-variété algébrique lisse isomorphe, en vertu du théorème de Hilbert 90, 
à l'espace homogène $\tlG_{\rmE}/ \tlG$ via l'isomorphisme:
\[
\tlG_{\rmE}/ \tlG \ni g \mapsto g \brg^{-1} \in S_{W}.
\]
Notons aussi $\sn_{W}$ l'espace tangent à l'identité de $S_{W}$ défini comme
\[
\sn_{W} := \{X \in \tlgl_{\rmE} | X + \brX = 0\}.
\]

Le groupe $G$ agit sur $S_{W}$ préservant l'identité par la restriction de son action par conjugaison sur $\tlG_{\rmE}$, il agit donc aussi sur $\sn_{W}$. 
Il est clair que pour tout $\xi \in \rmE$, ces action préservent les sous-variétés 
$\tlG_{\rmE} \smin D_{\xi}$ et $\sn_{W} \smin D_{1}$.
Notons $\rmE^{1} = \{ \xi \in \rmE| \, \xi \bar \xi = 1\}$. On a alors
\blem[cf. \cite{zhang2}, lemme 3.4]\label{lem:kappaIsoG} Pour tout $\xi \in \rmE^{1}$ la restriction de 
l'application $\kappa_{\xi}$ à $\sn_{W} \smin D_{1}$ 
est à valeurs dans $S_{W}$ est induit un isomorphisme $G$-invariant
entre $\sn_{W} \smin D_{1}$ et $S_{W} \smin D_{\xi}$.
\elem

Fixons un $\tau \in \rmE$ tel que $\sigma(\tau) = - \tau$. Notons aussi $\tau$ la multiplication 
par la matrice $\tau I$ dans $\tlgl_{\rmE}$. 
Alors, $\tau$ induit une application bijective entre $\calO$ et lui même.
Soit $\tlgl$ l'algèbre de Lie de $\tlG$. 
Le groupe $G$ agit sur $\tlG$ par conjugaison, et la multiplication par $\tau$ 
induit un isomorphisme $G$-équivariant des variétés $\sn_{W}$ et $\tlgl$. 
Il sera plus commode de travailler avec $\tlgl$ qu'avec $\sn_{W}$.

\subsection{Compatibilité de l'application de Cayley avec la décomposition de Levi}\label{par:cayleyLevi}

Soit $\tlP \in \calF(M_{\tlzero})$. On note $\ml_{\tlP}$ l'algèbre de Lie 
de $M_{\tlP}$ et $\nl_{\tlP}$ l'algèbre de Lie de $N_{\tlP}$. On fixe 
une mesure de Haar sur $\nl_{\tlP}(\A)$ de façon que 
le volume du quotient $\nl_{\tlP}(\rmF) \bsl \nl_{\tlP}(\A)$ soit $1$.

On a alors
\blem\label{lem:caleyLevi}
 Soient $\xi \in \rmE^{1}$, $\tlP \in \calF(M_{\tlzero})$, 
$\gamma \in ((S_{W} \cap M_{\tlP_{\rmE}}) \smin D_{\xi})(\rmF)$ 
et notons $\zeta = \tau \kappa_{\xi}^{-1}(\gamma)$.
Choisissons un $\delta \in M_{\tlP}(\rmE)$ tel que $\gamma = \delta \bar \delta^{-1}$. 
On a donc
\begin{enumerate}[1)]
\item $\zeta \in (\ml_{\tlP} \smin D_{\tau})(\rmF)$.
\item Pour tout $n \in N_{\tlP_{\rmE}}(\A)$ 
on a $\delta n \brn^{-1} \brdel^{-1} \in (S_{W} \smin D_{\xi})(\A)$.
\item Pour tout $n \in N_{\tlP_{\rmE}}(\A)$ il existe un unique $N \in \nl_{\tlP}(\A)$ tel que 
\[
\tau \kappa_{\xi}^{-1}(\delta n \brn^{-1} \brdel^{-1}) = \zeta + N.
\]
\item L'application $n \mapsto N$ décrite dans le point ci-dessus induit un homéomorphisme
\[
N_{\tlP_{\rmE}}(\A) /N_{\tlP}(\A) \rar \nl_{\tlP}(\A)
\]
qui préserve les mesures de Haar. En particulier, pour tout $\upphi \in C_{c}^{\infty}(S_{W}(\A))$ 
on a
\[
\int_{N_{\tlP_{\rmE}}(\A)/ N_{\tlP}(\A)}\upphi (\delta n \brn^{-1} \brdel^{-1})dn = 
\int_{\nl_{\tlP}(\A)} \upphi(\kappa_{\xi}(\tau^{-1}( \zeta + N)))dN
\]
et l'intégrale à gauche ne dépend pas du choix de $\delta$. 
\item Soit $\ol \in \calO$ tel que $\gamma \in (S_{W} \cap M_{\tlP_{\rmE}} \cap \ol)(\rmF)$. 
Alors, pour tout $n \in N_{\tlP_{\rmE}}(\A)$ on a 
$\delta n \brn^{-1} \brdel^{-1} \in (S_{W}  \cap \ol)(\A)$.
\end{enumerate}

\bdem 
Le premier point est clair. Pour le deuxième, on a bien sûr 
$\delta n \brn^{-1} \brdel^{-1} \in S_{W}(\A)$ et en plus 
$\delta n \brn^{-1} \brdel^{-1} = \gamma \brdel n \brn^{-1} \brdel^{-1}$. 
Puisque $\gamma, \brdel \in M_{\tlP}(\rmE)$ et $n \brn^{-1} \in N_{\tlP_{\rmE}}(\A)$ 
on a $ \brdel n \brn^{-1} \brdel^{-1} \in N_{\tlP_{\rmE}}(\A)$ et donc 
les polynômes caractéristiques de $\delta n \brn^{-1} \brdel^{-1}$ 
et $\gamma$ coïncident, d'où le résultat. 

Démontrons le point \textit{3)}. Soit $n_{0} = \brdel n \brn^{-1} \brdel^{-1} \in N_{\tlP_{\rmE}}(\A)$ 
et pour $\eta \in M_{\tlP}(\rmE)$ notons $n_{0}^{\eta} = \eta^{-1} n_{0} \eta$. 
On a alors que $\kappa_{\xi}^{-1}(\delta n \brn^{-1} \brdel^{-1} )$ égale
\begin{multline*}
-(\xi + \gamma n_{0})(\xi - \gamma n_{0})^{-1} = 
-(\xi + \gamma + \gamma (n_{0}-I))(\xi - \gamma - \gamma (n_{0} - I))^{-1} = 
-(\xi+ \gamma)(\xi - \gamma) \cdot \\
(I + (\xi + \gamma)^{-1}\gamma (n_{0}^{\xi - \gamma} - I))
(I - (\xi - \gamma)^{-1}\gamma(n_{0}-I))^{-1} = 
\tau^{-1}\zeta (I + \al_{1}(n_{0}^{\al_{2}} - I))(I + \sum_{k=1}^{\infty} (\al_{3}(n_{0}-I))^{k})
\end{multline*}
pour certains $\al_{1}, \al_{2}, \al_{3} \in M_{\tlP}(\rmE)$.
La somme qui apparaît ci-dessus est bien sûr finie
car $\al_{3}(n_{0}-I) \in \nl_{\tlP_{\rmE}}(\A)$.
On voit donc que 
$N := \zeta (I + \al_{1}(n_{0}^{\al_{2}} - I))(I + \sum_{k=1}^{\infty} (\al_{3}(n_{0}-I))^{k}) - \zeta \in \nl_{\tlP_{\rmE}}(\A)$. 
Puisque $\tau^{-1} \zeta \in (\sn_{W} \smin D_{1})(\rmF)$, il résulte du lemme \ref{lem:kappaIsoG} 
que $\tau^{-1} N \in (\sn_{W} \cap \nl_{\tlP_{\rmE}}) (\A)$. 
Le même lemme montre que l'application 
$n \mapsto N$ est un isomorphisme entre $N_{\tlP_{\rmE}}(\A) /N_{\tlP}(\A)$ 
et $\nl_{\tlP}(\A)$. Le fait que la mesure soit préservée découle du fait 
que cette application est polynomiale à coefficients rationnels, ce qui démontre 
aussi le point \textit{4)}. 

Démontrons le point \textit{5)}. Soit $\ol' \in \calO$ tel que 
$\tau \kappa_{\xi}^{-1}((S_{W} \smin D_{\xi}) \cap \ol) \sbs \ol'$. En vertu du lemme 
\ref{lem:kappaIsoG}, $\ol'$ existe et est unique. En vertu du point \textit{3)} 
il suffit de montrer que $\zeta + N \in \ol'(\A)$ pour tout $N \in \nl_{\tlP}(\A)$. 
Cela n'est pas difficile et est démontré dans le cas des groupes unitaires 
dans \cite{leMoi}, proposition 2.5. 
\edem
\elem

\subsection{Fonctions localisées}\label{par:locfonsG}

Pour toute place finie $v$ de $\rmF$ (resp. $w$ de $\rmE$) soit 
$\rmF_{v}$ le complété de $\rmF$ à $v$ 
(resp. $\rmE_{w}$ le complété de $\rmE$ à $w$) 
et soit $\calO_{\rmF_{v}}$ (resp. $\calO_{\rmE_{w}}$) 
l'anneau des entiers de $\rmF_{v}$ (resp. $\rmE_{w}$). 
Soit $\calS_{\rmF}$ l'ensemble des toutes places de $\rmF$ et notons 
$\calS_{\infty}$ et $\calS_{f}$ les sous-ensembles de $\calS_{\rmF}$ composés 
des places infinies et finies respectivement.

Pour une $\rmF$-variété $X$ lisse on note $C_{c}^{\infty}(X(\A)) := C_{c}^{\infty}(X(\rmF_{\infty})) \otimes_{\C} C_{c}^{\infty}(X(\A_{f}))$ 
où $C_{c}^{\infty}(X(\rmF_{\infty}))$ c'est l'espace de fonctions 
lisses à support compact sur la variété différentiable $X(\rmF_{\infty})$ 
et $C_{c}^{\infty}(X(\A_{f}))$ c'est l'espace de fonctions localement constantes à support compact 
sur $X(\A_{f})$.
Pour tout $v \in \calS_{f}$ (resp. $v \in \calS_{\infty}$) on note aussi dans ce contexte 
$C^{\infty}(X(\rmF_{v}))$ l'espace de fonctions localement constantes (resp. lisses) sur $X(\rmF_{v})$ 
et $C_{c}^{\infty}(X(\rmF_{v}))$ le sous-espace de $C^{\infty}(X(\rmF_{v}))$ composé de fonctions 
à support compact. 

Fixons un isomorphisme $\tlG_{\rmE} \cong \Res_{\rmE/\rmF}(\Gl_{n+1})$. On a alors que
les variétés $\tlG_{\rmE}$, $S_{W}$, $\sn_{W}$, $\tlgl_{\rmE}$ et $\tlgl$ 
sont définies sur l'anneau des entiers de $\rmF$ 
et on peut en particulier parler de ses $\calO_{\rmF_{v}}$-points pour tout $v \in \calS_{f}$.
Soit $X$ l'une des ces variétés. 
Pour tout $v \in \calS_{f}$ on note $\ind_{X(\calO_{\rmF_{v}})} \in C_{c}^{\infty}(X(\rmF_{v}))$ 
la fonction caractéristique du compact $X(\calO_{\rmF_{v}})$.
On dit que $f \in C_{c}^{\infty}(X(\A))$ est décomposable si $f$ s'écrit 
$f_{\infty} \otimes_{v \in \calS_{f}}f_{v}$ où 
$f_{\infty} \in C_{c}^{\infty}(X(\rmF_{\infty}))$ et pour tout $v \in \calS_{f}$ on a 
$f_{v} \in C_{c}^{\infty}(X(\rmF_{v}))$ 
et $f_{v} = \ind_{X(\calO_{\rmF_{v}})}$ pour presque toute place $v \in \calS_{f}$. On a alors que 
toute fonction $f \in C_{c}^{\infty}(X(\A))$ est une somme finie des fonctions décomposables.

On note finalement $supp(f)$ le support d'une fonction $f$.

\blem\label{lem:finiOls} Pour tout $f \in C_{c}^{\infty}(S_{W}(\A))$ il existe un nombre fini des $\ol \in \calO$ 
tels que $supp(f) \cap \ol(\A) \neq \varnothing$.
\bdem
Par définition de l'application $Q$, définie par (\ref{eq:invMap}), l'ensemble 
\[
\bigcup_{\ol \in \calO} Q(\ol(\A))
\]
est contenu dans le réseau discret $\rmE^{2n+1}$ de $(\A \otimes_{\rmE} \rmE)^{2n+1}$. 
En plus, pour tout $\ol \in \calO$, l'ensemble $Q(\ol(\A))$ est composé d'un seul point et si $\ol' \neq \ol$ on a 
$Q(\ol'(\A)) \neq Q(\ol(\A))$. 
Puisque l'image $Q(supp(f))$ du support de la fonction $f$ par $Q$ est compact 
dans $(\A \otimes_{\rmE} \rmE)^{2n+1}$, on obtient  
qu'il existe un nombre fini des $\ol \in \calO$ tels que 
$Q(\ol(\A)) \cap Q(supp(f)) \neq \varnothing$ d'où le résultat. 
\edem
\elem 


\bdefi\label{def:olxiAdapt}
Soient $\ol \in \calO$, $\xi \in \rmE^{1}$ et $f \in \calC_{c}^{\infty}(S_{W}(\A))$.
On dit que $f$ est $\ol$-$\xi$-adaptée si elle vérifie les conditions suivantes:
\begin{enumerate}[1)]
\item La fonction $f$ est décomposable. On écrit alors $f = f_{\infty} \otimes_{v \in \calS_{f}} f_{v}$.
\item Pour tout $\ol' \in \calO$ si $\ol'(\A) \cap supp(f) \neq \varnothing$ alors $\ol'= \ol$.
\item $\ol \cap D_{\xi} = \varnothing$.
\item Il existe un ensemble fini $\calS \sbs \calS_{\rmF}$ tel que
\begin{enumerate}[a)]
\item $\calS \sps \calS_{\infty}$.
\item Pour tout $v \in \calS_{\rmF} \smin \calS$ on a $f_{v} = \ind_{S_{W}(\calO_{\rmF_{v}})}$.
\item 
Soit $\ol' = \tau \kappa_{\xi}^{-1}(\ol) \in \calO$. Alors, pour tout 
$v \in \calS_{\rmF} \smin \calS$ et toute place $w$ de $\rmE$ au dessus de $v$ 
on a que $P_{\ol'}(\tau)$, $P_{\ol'}(-\tau)$, $P_{\ol}(\xi)$, $\tau$, et $\xi$ appartiennent à $\calO_{\rmE_{w}}^{*}$.
\item On a $supp(f_{\infty}) \sbs (S_{W} \smin D_{\xi})(\rmF_{\infty})$ et pour tout 
$v \in \calS \smin \calS_{\infty}$ on a $supp(f_{v}) \sbs (S_{W} \smin D_{\xi})(\rmF_{v})$.
\end{enumerate}
\end{enumerate}
\edefi

\blem\label{lem:redToAdapt}
 Soit $f \in C_{c}^{\infty}(S_{W}(\A))$ décomposable et $\ol_{1}, \ldots, \ol_{m} \in \calO$ tous les éléments de $\calO$ 
qui vérifient $\ol_{i}(\A) \cap supp(f) \neq \varnothing$
donnés par le lemme \ref{lem:finiOls}. 
 Pour tout $i \in \{1,\ldots, m\}$ 
il existe un $\xi_{i} \in \rmE^{1}$ et une fonction $f_{i} \in C_{c}^{\infty}(S_{W}(\A))$ 
qui est $\ol_{i}$-$\xi_{i}$-adaptée telle que 
pour tout $s \in (\calS_{W} \cap \ol_{i})(\A)$ 
on a $f(s) = f_{i}(s)$.
\bdem 
Soit $i \in \{1, \ldots, m\}$. 
Prenons un $\xi_{i} \in \rmE^{1}$ tel que $\ol_{i} \cap D_{\xi_{i}} = \varnothing$. 
Écrivons $f = \otimes_{v} f_{v}$ et soit $\ol_{i}' = \tau \kappa_{\xi_{i}}^{-1}(\ol_{i})$.
Soit $\calS$ un sous-ensemble fini de $\calS_{\rmF}$ qui vérifie les conditions suivantes:
\begin{itemize}
\item $\calS \sps \calS_{\infty}$.
\item Pour tout $v \in \calS_{\rmF} \smin \calS$ et toute place $w$ de $\rmE$ au dessus de $v$ 
on a $P_{\ol_{i}'}(\tau)$, $P_{\ol_{i}'}(-\tau)$, $P_{\ol_{i}}(\xi_{i})$, $\tau$, $\xi_{i} \in \calO_{\rmE_{w}}^{*}$.
\item Pour tout $v \in \calS_{\rmF} \smin \calS$ on a $f_{v} = \ind_{S_{W}(\calO_{\rmF_{v}})}$.
\end{itemize}
Pour tout $v \in \calS$ il est clair qu'il existe un $\phi_{v} \in C^{\infty}(S_{W}(\rmF_{v}))$ 
tel que $supp(\phi_{v}) \sbs (S_{W} \smin (D_{\xi_{i}} \cup \bigcup_{j \in \{1, \ldots, m\} \smin \{i\}} \ol_{j}))(\rmF_{v})$ 
et la restriction de $\phi_{v}$ à $(S_{W} \cap \ol_{i})(\rmF_{v})$ égale $1$. Prenons une telle fonction $\phi_{v}$ 
pour tout $v \in \calS$ et posons $f_{i} := f \cdot \prod_{v \in \calS} \phi_{v}$. 
Alors, $f_{i} \in C_{c}^{\infty}(S_{W}(\A))$ et $f_{i}$ est $\ol_{i}$-$\xi_{i}$-adaptée. 
En plus, pour tout $s \in (S_{W} \cap \ol_{i})(\A)$ on a $f_{i}(s) = f(s)$ ce qu'il fallait démontrer. 
\edem
\elem

\subsection{Descente vers l'algèbre de Lie}\label{par:Liedescente}

\blem\label{lem:Liedescente} Soient $\ol \in \calO$, $\xi \in \rmE^{1}$ et 
$f \in C_{c}^{\infty}(S_{W}(\A))$ une fonction $\ol$-$\xi$-adaptée. 
Il existe une fonction décomposable $\Phi \in C_{c}^{\infty}(\tlgl(\A))$ telle que
\[
f(s) = \Phi(\tau \kappa_{\xi}^{-1}(s)), \quad \forall s \in (S_{W} \cap \ol)(\A). 
\]
\bdem 
On va construire $\Phi \in C_{c}^{\infty}(\tlgl(\A))$ de type 
$\Phi_{\infty} \otimes_{v \in \calS_{f}}\Phi_{v}$. Soit $\calS$ l'ensemble de places de $\rmF$ 
associé à $f$ comme dans le point \textit{2)} de la définition \ref{def:olxiAdapt} 
et écrivons $f = f_{\infty} \otimes_{v \in \calS_{f}} f_{v}$.
Pour tout $v \in \calS_{\rmF} \smin \calS$ on pose alors $\Phi_{v} = \ind_{\tlgl(\calO_{\rmF_{v}})}$.
Soit $v \in \calS \smin \calS_{\infty}$. On définit $\Phi_{v}$ de façon suivante. Si $X \in (\tlgl \smin D_{\tau})(\rmF_{v})$ 
on pose $\Phi_{v}(X) = f_{v}(\kappa_{\xi}(\tau^{-1}X))$ et 
si $X \in D_{\tau}(\rmF_{v})$ on pose $\Phi_{v}(X) = 0$. 
En utilisant le fait que $supp(f_{v}) \sbs (S_{W} \smin D_{\xi})(\rmF_{v})$ et le lemme \ref{lem:kappaIsoG} 
on voit que $\Phi_{v} \in C_{c}^{\infty}(\tlgl(\rmF_{v}))$.
De même on définit $\Phi_{\infty}(X)$ par 
$f_{\infty}(\kappa_{\xi}(\tau^{-1}X))$  si $X \in (\tlgl \smin D_{\tau})(\rmF_{\infty})$ et par $0$ sinon. 
On trouve alors que $\Phi_{\infty} \in C_{c}^{\infty}(\tlgl(\rmF_{\infty}))$.
 On a donc bien 
que $\Phi := \Phi_{\infty} \otimes_{v \in \calS_{f}} \Phi_{v} \in C_{c}^{\infty}(\tlgl(\A))$. 

Soit alors $s = \prod_{v}s_{v} \in (S_{W} \cap \ol)(\A)$. 
On a $f_{\infty}(\prod_{v \in \calS_{\infty}}s_{v}) = \Phi_{\infty}(\tau \kappa_{\xi}^{-1}(\prod_{v \in \calS_{\infty}}s_{v}))$ 
et $f_{v}(s_{v}) = \Phi_{v}(\tau \kappa_{\xi}^{-1}(s_{v}))$ 
pour $v \in \calS \smin \calS_{\infty}$ par construction. 
Supposons que $v \in \calS_{\rmF} \smin \calS$. 
Soit $\ol' = \tau \kappa_{\xi}^{-1}(\ol)$.
En vertu du point \textit{4) c)} de la définition \ref{def:olxiAdapt}, 
l'application $\tau \kappa_{\xi}^{-1}$ 
induit une bijection entre 
$S_{W}(\calO_{\rmF_{v}}) \cap \ol(\rmF_{v})$ et
$\tlgl(\calO_{\rmF_{v}}) \cap \ol'(\rmF_{v})$.  
On a donc $s_{v} \in S_{W}(\calO_{\rmF_{v}})$ 
si et seulement si $\tau \kappa_{\xi}^{-1} (s_{v}) \in \tlgl(\calO_{\rmF_{v}})$
d'où $f_{v}(s_{v}) = \ind_{S_{W}(\calO_{\rmF_{v}})}(s_{v}) = \Phi(\tau \kappa_{\xi}^{-1}(s_{v})) = \ind_{\tlgl(\calO_{\rmF_{v}})}(\tau \kappa_{\xi}^{-1}(s_{v}))$.
\edem
\elem

\subsection{Le côté géométrique de la formule des traces pour les groupes linéaires}\label{par:geomGroupesG}

Soient $f \in C_{c}^{\infty}(S_{W}(\A))$, $\tlP \in \relPb$ et $\ol \in \calO$.
Posons
\[
k_{f,\tlP, \ol}(x) = \sum_{\gamma \in (M_{\tlP_{\rmE}} \cap S_{W} \cap \ol)(\rmF)}
\int_{N_{\tlP_{\rmE}}(\A) / N_{\tlP}(\A)}f(x^{-1}\delta n \brn^{-1} \brdel^{-1} x)dn, \quad x \in M_{\tlP}(\rmF)N_{\tlP}(\A) \bsl \tlG(\A)
\]
où $\delta \in M_{\tlP}(\rmE)$ est tel que $\delta \brdel^{-1} = \gamma$. La définition ne dépend pas du choix de $\delta$ 
en vertu du lemme \ref{lem:caleyLevi} \textit{4)}. On définit aussi pour $\Phi \in C_{c}^{\infty}(\tlgl(\A))$ 
la fonction
\[
k_{\Phi, \tlP, \ol}(x) = \sum_{\xi \in (\ml_{\tlP} \cap \ol)(\rmF)}
\int_{\nl_{\tlP}(\A)}\Phi(x^{-1}(\xi + N)x)dN, \quad x \in M_{\tlP}(\rmF)N_{\tlP}(\A) \bsl \tlG(\A).
\]
Dans le paragraphe 2.1 de \cite{leMoi2} nous avons défini une relation 
d'équivalence sur $\tlgl(\rmF)$, notons $\calO_{\tlgl(\rmF)}$ l'ensemble de classes d'équivalence 
pour cette relation. Il existe alors une bijection entre $\calO_{\tlgl(\rmF)}$ et 
l'ensemble de classes $\ol \in \calO$ tels que $\ol \cap \tlgl \neq \varnothing$, 
la bijection étant: à $\ol' \in \calO_{\tlgl(\rmF)}$ on associe l'unique $\ol \in \calO$ tel que 
$(\ol \cap \tlgl)(\rmF) = \ol'$. 
Avec la notation du paragraphe 2.2 de loc. cit, si $\ol' \in \calO_{\tlgl(\rmF)}$ et $\ol \in \calO$ se correspondent de cette façon, 
pour tout $\tlP \in \relP$ et tout $\Phi \in C_{c}^{\infty}(\tlgl(\A)$ on a:
\[
k_{\Phi, \tlP, \ol'} = k_{\Phi, \tlP, \ol}.
\]

On a alors le corollaire suivant du lemme \ref{lem:Liedescente}.

\bcor\label{cor:kpGiskptlG}
Soient $\ol \in \calO$, $\xi \in \rmE^{1}$ et 
$f \in C_{c}^{\infty}(S_{W}(\A))$ une fonction $\ol$-$\xi$-adaptée. Notons $\ol' = \tau k_{\xi}^{-1}(\ol)$.
Il existe alors une fonction décomposable $\Phi \in C_{c}^{\infty}(\tlgl(\A))$ telle que
$supp(\Phi) \sbs \tlgl(\A) \smin D_{\tau}(\A)$ et telle que pour tout $\tlP \in \relPb$ et tout $x \in G(\A)$ on a 
\[
k_{f,\tlP,\ol}(x) = k_{\Phi, \tlP, \ol'}(x).
\]
\bdem 
On prend $\Phi$ donné par le lemme \ref{lem:Liedescente}. Il résulte alors directement de ce lemme-là
et du point \textit{4)} 
du lemme \ref{lem:caleyLevi} que $k_{f,\tlP,\ol}(1) = k_{\Phi, \tlP, \ol'}(1)$. 
En utilisant le fait que $\tau k_{\xi}^{-1}$ est $G$-équivariant, on obtient 
le résultat voulu. 
\edem
\ecor

Pour $f \in C_{c}^{\infty}(S_{W}(\A))$, $T \in \all_{\tlzero}$ et $\ol \in \calO$ posons
\[
k_{f, \ol}^{T}(x) = \sum_{\tlP \in \relPb}(-1)^{d_{\tlP}^{\tlG}}
\sum_{\delta \in P(\rmF) \bsl G(\rmF)}\htau_{\tlP}(H_{\tlP}(\delta x) - T_{\tlP})
k_{f,\tlP, \ol}(\delta x), \quad x \in G(\rmF) \bsl G(\A).
\]

\btheo\label{thm:noyauGeomCVGG}
 Soit $f \in C_{c}^{\infty}(S_{W}(\A))$.
\begin{enumerate}[1)]
\item Pour tout  $\sigma \in \R$ et tout $T \in T_{+} + \all_{\tlzero}^{+}$ on a
\[
\sum_{\ol \in \calO}\int_{[G]}|k_{f, \ol}^{T}(x)| |\det x|_{\A}^{\sigma}dx < \infty.
\]
\item Pour $s \in \C$, $T \in T_{+} + \all_{\tlzero}^{+}$ et $\ol \in \calO$ on pose
\[
I_{\ol}^{T}(s,f) = \int_{[G]}k_{f, \ol}^{T}(x) |\det x|_{\A}^{s} \eta(\det x)dx.
\]
Alors, la fonction $T \mapsto I_{\ol}^{T}(s,f)$ est un polynôme-exponentielle, dont le terme purement polynômial
est constant si $s \neq -1,1$.
\item Si $s \in \C \smin \{-1,1\}$ on note $I_{\ol}(s, \cdot)$ le terme constant de la distribution 
$I_{\ol}^{T}(s, \cdot)$. Alors, la distribution $I_{\ol}(s, \cdot)$ 
ne dépend que du choix de la mesure de Haar 
sur $G(\A)$ et des choix des mesures sur les points adéliques des parties unipotentes $N_{\tlP}$ et $N_{\tlP_{\rmE}}$ des $\rmF$-sous-groupes 
paraboliques de $\tlG$ et $\tlG_{\rmE}$, le choix étant que 
les volumes $N_{\tlP}(\rmF) \bsl N_{\tlP}(\A)$ et 
$N_{\tlP_{\rmE}}(\rmF) \bsl N_{\tlP_{\rmE}}(\A)$ soient $1$.
\item Soit $y \in G(\A)$. On note $f^{y} \in C_{c}^{\infty}(S_{W}(\A))$ 
la fonction donnée par $f^{y}(s) = f(\Ad(y)s)$. On a alors
\[
I_{\ol}(s, f^{y}) = \eta (\det y) |\det y|_{\A}^{s}I_{\ol}(s,f), \quad \forall \, s \in \C \smin \{-1,1\}, \ 
\forall \, \ol \in \calO.
\]
\end{enumerate}
\bdem 
Toute fonction $f \in C_{c}^{\infty}(S_{W}(\A))$ est une somme finie des fonctions décomposables, on peut alors supposer 
que $f$ est décomposable. Soient $\ol_{1}, \ldots, \ol_{m} \in \calO$ tous les éléments de $\calO$ 
qui vérifient $\ol_{i}(\A) \cap supp(f) \neq \varnothing$
donnés par le lemme \ref{lem:finiOls}.   
Pour tout $i \in \{1,\ldots, m\}$ soient $\xi_{i} \in \rmE^{1}$ et
$f_{i} \in C_{c}^{\infty}(S_{W}(\A))$ la fonction $\ol_{i}$-$\xi_{i}$-adaptée associés à $f$ 
par le lemme \ref{lem:redToAdapt}.
Pour tout $i \in \{1,\ldots, m\}$ et tout $T \in \all_{\tlzero}$ 
on a $k_{f,\ol_{i}}^{T} = k_{f_{i},\ol_{i}}^{T}$ 
et 
\[
\sum_{\ol \in \calO}|k_{f,\ol}^{T}| = \sum_{i=1}^{m}
\sum_{\ol \in \calO}|k_{f_{i},\ol}^{T}| = 
\sum_{i=1}^{m}|k_{f_{i},\ol_{i}}^{T}|.
\]
En plus, pour tout $y \in G(\A)$ et tout $i \in \{1,\ldots, m\}$ 
on a $k_{f^{y},\ol_{i}}^{T} = k_{f_{i}^{y},\ol_{i}}^{T}$ ce qui démontre qu'on peut supposer 
que $f$ est $\ol$-$\xi$-adaptée pour certains $\ol \in \calO$ et $\xi \in \rmE^{1}$. 
Soit donc $\Phi \in C_{c}^{\infty}(\gl(\A))$ 
une fonction associée à $f$ comme dans le lemme \ref{lem:Liedescente}. 
En vertu de ce lemme-là et du corollaire \ref{cor:kpGiskptlG}, on a pour tout $y \in G(\A)$ 
et tout $T \in \all_{\tlzero}$
que $k_{f^{y},\ol}^{T} = k_{\Phi^{y}, \ol'}^{T}$ où 
$\ol' = \tau \kappa_{\xi}(\ol)$ et $\Phi^{y}(X) = \Phi(\Ad(y)X)$. 

Le théorème découle alors des résultats suivants démontrés dans \cite{leMoi2}.
Le point \textit{1)} découle du théorème 2.6, 
le point \textit{2)} résulte de la proposition 3.7 et du théorème 3.8, 
le point \textit{3)} c'est la conséquence du paragraphe 3.5
et finalement le point \textit{4)} suit du théorème 3.11.
\edem
\etheo

\subsection{La formule des traces relative}\label{par:RTFG}

Soit $\Phi \in C_{c}^{\infty}(\GetlGa)$
Notons alors pour $y \in S_{W}(\A)$:
\[
f_{\Phi}(y) := \int_{G_{\rmE}(\A)} \int_{\tlG(\A)} \Phi(x, x y \tlh) d\tlh dx.
\]

On définit alors dans ce cas pour $s \in \C$, $\ol \in \calO$ et $T \in T_{+} + \all_{\tlzero}^{+}$
\[
I_{\ol}^{T}(s,\Phi) := I_{\ol}^{T}(s, f_{\Phi})
\]
et si $s \neq -1,1$:
\[
I_{\ol}(s,\Phi) := I_{\ol}(s, f_{\Phi}).
\]

On utilise la notation de la section \ref{sec:SpecG}.

\btheo\label{thm:RTFG}[La formule des traces relative de Jacquet-Rallis, cas des groupes linéaires]
Pour tout $\Phi \in C_{c}^{\infty}(\GetlGa)$,
tout $T \in T_{+} + \all_{\tlzero}^{+}$ et tous $s,s' \in \C$ on a
\begin{equation}\label{eq:rtfGT}
\sum_{\chi \in \calX^{G_{\rmE} \times \tlG_{\rmE}}}I_{\chi}^{T}(s,s', \Phi) = 
\sum_{\ol \in \calO}I_{\ol}^{T}(s + s',\Phi).
\end{equation}
En plus, si $s + s' \neq -1,1$ on a 
\begin{equation}
\sum_{\chi \in \calX^{G_{\rmE} \times \tlG_{\rmE}}}I_{\chi}(s,s', \Phi) = 
\sum_{\ol \in \calO}I_{\ol}(s + s',\Phi).
\end{equation}
\bdem 
Pour $\tlP \in \relPb$ soient 
$k_{\Phi, \tlP} := \sum_{\chi \in \calX^{G_{\rmE} \times \tlG_{\rmE}}}k_{\Phi, \tlP, \chi}$ 
et $k_{f_{\Phi}, \tlP} = \sum_{\ol \in \calO}k_{f_{\Phi}, \tlP, \ol}$. 
On vérifie alors facilement que pour tous $h \in G(\A)$ et $s,s' \in \C$ on a
\begin{equation}\label{eq:rtfHelp}
\int\limits_{[G_{\rmE}]} \int\limits_{[\tlG]}
\sum_{\mathrlap{
\begin{subarray}{c}
\delta_{1} \in P(\rmE)\bsl G(\rmE)
\\
\delta_{3} \in \tlP(\rmF)\bsl \tlG(\rmF)
\end{subarray}}}
k_{\Phi, \tlP}(\delta_{1} x, \delta_{1} x,h,\delta_{3} \tlh)
|\det x|_{\A}^{s}|\det h|_{\A}^{s'}\eta(h,\tlh)
d\tlh dx =  
k_{f_{\Phi}, \tlP}(h)|\det h|_{\A}^{s+s'}\eta(\det h).
\end{equation}
Si l'on pose alors pour $T \in \all_{\tlzero}$, 
$k_{\Phi}^{T} := \sum_{\chi \in \calX^{G_{\rmE} \times \tlG_{\rmE}}} k_{\Phi,\chi}^{T}$ 
et $k_{f_{\Phi}}^{T} := \sum_{\ol \in \calO}k_{f_{\Phi},\ol}^{T}$, on voit 
que l'égalité (\ref{eq:rtfGT}) se traduit en
\[
\int_{[G_{\rmE}]}\int_{[G]}\int_{[\tlG]}
k_{\Phi}^{T}(x,x,h,\tlh) |\det x|_{\A}^{s}
|\det h|_{\A}^{s'}\eta(h,\tlh)d \tlh dh dx = 
\int_{[G]}k_{f_{\Phi}}^{T}|\det h|_{\A}^{s + s'}\eta( \det h)dh,
\]
ce qui est vrai en vertu de l'égalité \eqref{eq:rtfHelp} et de la convergence absolue 
des intégrales en question, donnée par les théorèmes \ref{thm:noyauCVGG} 
et \ref{thm:noyauGeomCVGG} \textit{1)}.

La deuxième partie découle maintenant de la première et des théorèmes 
\ref{thm:mainQualitThmG} et \ref{thm:noyauGeomCVGG} \textit{2)}, car égalité 
des polynômes-exponentielles implique égalité de leur termes constantes. 
\edem
\etheo
\section{Formule des traces relative de Jacquet-Rallis pour les groupes unitaires}\label{sec:RTFU}

Les résultats de cette section sont analogues aux résultats de la section \ref{sec:RTFG} et pour la plupart des résultats on 
envoie le lecteur à la section \ref{sec:RTFG}. 
Soit $\tlul = \Lie(\tlU)$. On regarde $U$, $\tlU$ et $\tlul$ 
comme des sous-variétés de $\tlgl_{\rmE}$.

\subsection{Préparations pour le côté géométrique}\label{par:invparU}

Soit $P \in \calF(M_{0})$. On note $\ml_{\tlP}$ l'algèbre de Lie 
de $M_{\tlP}$ et $\nl_{\tlP}$ l'algèbre de Lie de $N_{\tlP}$. On fixe 
une mesure de Haar sur $\nl_{\tlP}(\A)$ de façon que 
le volume du quotient $\nl_{\tlP}(\rmF) \bsl \nl_{\tlP}(\A)$ soit $1$.

On a alors:
\blem\label{lem:caleyLeviU}
 Soient $\xi \in \rmE^{1}$, $P \in \calF(M_{0})$, 
$\gamma \in ((\tlU \cap M_{\tlP}) \smin D_{\xi})(\rmF)$ 
et notons $\zeta = \kappa_{\xi}^{-1}(\gamma)$.
On a donc:
\begin{enumerate}[1)]
\item $\zeta \in (\ml_{\tlP} \smin D_{1})(\rmF)$.
\item Pour tout $n \in N_{\tlP}(\A)$ 
on a $\gamma n \in (\tlU \smin D_{\xi})(\A)$.
\item Pour tout $n \in N_{\tlP}(\A)$ il existe un unique $N \in \nl_{\tlP}(\A)$ tel que 
\[
\kappa_{\xi}^{-1}(\gamma n) = \zeta + N.
\]
\item L'application $n \mapsto N$ décrite dans le point ci-dessus induit un homéomorphisme
\[
N_{\tlP}(\A) \rar \nl_{\tlP}(\A)
\]
qui préserve les mesures de Haar. En particulier, pour tout $\upphi \in C_{c}^{\infty}(\tlU(\A))$ 
on a
\[
\int_{N_{\tlP}(\A)}\upphi (\gamma n)dn = 
\int_{\nl_{\tlP}(\A)} \upphi(\kappa_{\xi}( \zeta + N))dN.
\]
\item Soit $\ol \in \calO$ tel que $\gamma \in (\tlU \cap M_{\tlP} \cap \ol)(\rmF)$. 
Alors, pour tout $n \in N_{\tlP}(\A)$ on a 
$\gamma n \in (\tlU  \cap \ol)(\A)$.
\end{enumerate}

\bdem 
Pour le premier point, pour $X \in \tlgl_{\rmE}(\A)$ soit $X^{*}$ l'endomorphisme 
adjoint par rapport à la forme hermitienne $\tlPhi$ définissante $\tlU$. 
Il suffit de monter que $\zeta \in \tlul(\rmF)$, autrement dit que $\zeta + \zeta^{*} = 0$.
On a alors:
\[
\zeta + \zeta^{*} = 
-(\xi + \gamma)(\xi - \gamma)^{-1} - (\xi^{*} + \gamma^{*})(\xi^{*} - \gamma^{*})^{-1}.
\]
En utilisant le fait que $\xi^{*} = \xi^{-1}$ et $\gamma^{*}  = \gamma^{-1}$ on obtient le résultat voulu. 
Le deuxième point est clair. La preuve des points \textit{3)}, \textit{4)} et \textit{5)} est analogue 
à celle du lemme \ref{lem:caleyLevi}.
\edem
\elem

Dans le paragraphe \ref{par:locfonsG} nous avons fixé l'isomorphisme $\tlG_{\rmE} \cong \Res_{\rmE/\rmF} \Gl_{n+1}$. 
On regarde alors $\tlU$ comme une sous-$\rmF$-variété de $\Res_{\rmE/\rmF} \Gl_{n+1}$. 
La variété $\tlU$ est alors défini sur une localisation de l'anneau des entiers de $\rmF$ et 
notons 
$\calS_{\tlU}$ le sous-ensemble fini des places finies de $\rmF$ tel que pour tout $v \in \calS_{f} \smin \calS_{\tlU}$ 
l'anneau $\calO_{\rmF_{v}}$ contient cette localisation. Pour tout $v \in \calS_{f} \smin \calS_{\tlU}$  
on peut alors parler des $\calO_{\rmF_{v}}$-points de $\tlU$ et de $\tlul$.

\blem[cf. lemme \ref{lem:finiOls}]\label{lem:finiOlsU} Pour tout $f \in C_{c}^{\infty}(\tlU(\A))$ il existe un nombre fini des $\ol \in \calO$ 
tels que $supp(f) \cap \ol(\A) \neq \varnothing$.
\elem 

\bdefi\label{def:olxiAdaptU}
Soient $\ol \in \calO$, $\xi \in \rmE^{1}$ et $f \in \calC_{c}^{\infty}(\tlU(\A))$.
On dit que $f$ est $\ol$-$\xi$-adaptée si elle vérifie les conditions suivantes:
\begin{enumerate}[1)]
\item La fonction $f$ est décomposable. On écrit alors $f = f_{\infty} \otimes_{v \in \calS_{f}} f_{v}$.
\item Pour tout $\ol' \in \calO$ si $\ol'(\A) \cap supp(f) \neq \varnothing$ alors $\ol'= \ol$.
\item $\ol \cap D_{\xi} = \varnothing$.
\item Il existe un ensemble finie des places $\calS$ tel que
\begin{itemize}
\item $\calS \sps \calS_{\infty} \cup \calS_{\tlU}$.
\item Soit $\ol' \in \calO$ tel que $\kappa_{\xi}^{-1}(\ol) = \ol'$ donné par le lemme \ref{lem:cayleyDef}.
Alors, pour tout $v \in \calS_{\rmF} \smin \calS$ et toute place $w$ de $\rmE$ au dessus de $v$ 
on a que $P_{\ol'}(1)$, $P_{\ol'}(-1)$, $P_{\ol}(\xi)$ et $\xi$ appartiennent à $\calO_{\rmE_{w}}^{*}$.
\item Pour tout $v \in \calS_{\rmF} \smin \calS$ on a $f_{v} = \ind_{\tlU(\calO_{\rmF_{v}})}$.
\item On a $supp(f_{\infty}) \sbs (\tlU \smin D_{\xi})(\rmF_{\infty})$ et pour tout 
$v \in \calS \smin \calS_{\infty}$ on a $supp(f_{v}) \sbs (\tlU \smin D_{\xi})(\rmF_{v})$.
\end{itemize}
\end{enumerate}
\edefi

\blem[cf. lemme \ref{lem:redToAdapt}] 
Soit $f \in C_{c}^{\infty}(\tlU(\A))$ décomposable et $\ol_{1}, \ldots, \ol_{m} \in \calO$ tous les éléments de $\calO$ 
qui vérifient $\ol_{i}(\A) \cap supp(f) \neq \varnothing$
donnés par le lemme \ref{lem:finiOlsU}. 
 Pour tout $i \in \{1,\ldots, m\}$ 
il existe un $\xi_{i} \in \rmE^{1}$ et une fonction $f_{i} \in C_{c}^{\infty}(\tlU(\A))$ 
qui est $\ol_{i}$-$\xi_{i}$-adaptée telle que 
pour tout $x \in (\tlU \cap \ol_{i})(\A)$ 
on a $f(x) = f_{i}(x)$.
\elem

\blem[cf. lemme \ref{lem:Liedescente}]\label{lem:LiedescenteU} 
Soient $\ol \in \calO$, $\xi \in \rmE^{1}$ et 
$f \in C_{c}^{\infty}((\A))$ une fonction $\ol$-$\xi$-adaptée. 
Il existe une fonction décomposable $\Phi \in C_{c}^{\infty}(\tlul(\A))$ telle que
\[
f(x) = \Phi(\kappa_{\xi}^{-1}(x)), \quad \forall x \in (\tlU \cap \ol)(\A). 
\]
\elem

\subsection{Le côté géométrique de la formule des traces pour les groupes unitaires}\label{par:geomGroupesU}

Soient $f \in C_{c}^{\infty}(\tlU(\A))$, $P \in \calF(P_{0})$ et $\ol \in \calO$.
Posons
\[
k_{f, P, \ol}(x) = \sum_{\gamma \in (M_{\tlP} \cap \ol)(\rmF)}
\int_{ N_{\tlP}(\A)}f(x^{-1}\gamma n x)dn, \quad x \in M_{\tlP}(\rmF)N_{\tlP}(\A) \bsl \tlU(\A).
\]
On définit aussi pour $\Phi \in C_{c}^{\infty}(\tlul(\A))$ 
la fonction
\[
k_{\Phi, P, \ol}(x) = \sum_{\xi \in (\ml_{\tlP} \cap \ol)(\rmF)}
\int_{\nl_{\tlP}(\A)}\Phi(x^{-1}(\xi + N)x)dN, \quad x \in M_{\tlP}(\rmF)N_{\tlP}(\A) \bsl \tlU(\A).
\]
Dans le paragraphe 2.3 de \cite{leMoi} nous avons défini une relation d'équivalence 
sur $\tlul(\rmF)$. Notons $\calO_{\tlul}'$ l'ensemble de classes d'équivalence pour cette relation.  
Notons aussi $\calO_{\tlul}$ l'ensemble de $\ol \in \calO$ tels que $\ol \cap \tlul \neq \varnothing$. 
On a alors une application surjective naturelle $p : \calO_{\tlul}' \rar \calO_{\tlul}$ 
qui à $\ol' \in \calO_{\tlul}'$ associe $\ol \in \calO_{\tlul}$ tel que $\ol' \sbs \ol(\rmF)$. 
On a alors, avec la notation du début de la section 3 de \cite{leMoi}, 
que pour tout $P \in \calF(P_{0})$, tout $\Phi \in C_{c}^{\infty}(\tlul(\A))$ 
et tout $\ol \in \calO_{\tlul}$:
\[
k_{\Phi, P, \ol} = \sum_{\ol' \in p^{-1}(\ol)}k_{\Phi, P, \ol'}.
\]

Pour $f \in C_{c}^{\infty}(\tlU(\A))$, $T \in \all_{0}$ et $\ol \in \calO$ posons
\[
k_{f, \ol}^{T}(x) = \sum_{P \sbs P_{0}}(-1)^{d_{P}}
\sum_{\delta \in P(\rmF) \bsl U(\rmF)}\htau_{P}(H_{P}(\delta x) - T_{P})
k_{f, P, \ol}(\delta x), \quad x \in U(\rmF) \bsl U(\A).
\]

\btheo\label{thm:noyauGeomCVGU}
 Soit $f \in C_{c}^{\infty}(\tlU(\A))$.
\begin{enumerate}[1)]
\item Pour tout $T \in T_{+} + \all_{0}^{+}$ on a
\[
\sum_{\ol \in \calO}\int_{[U]}|k_{f, \ol}^{T}(x)|dx < \infty.
\]
\item Pour $T \in T_{+} + \all_{0}^{+}$ et $\ol \in \calO$ on pose
\[
J_{\ol}^{T}(f) = \int_{[U]}k_{f, \ol}^{T}(x)dx.
\]
Alors, la fonction $T \mapsto J_{\ol}^{T}(f)$ est un polynôme-exponentielle, dont le terme purement polynômial
est constant.
\item On note $J_{\ol}(\cdot)$ le terme constant de la distribution 
$J_{\ol}^{T}(\cdot)$. Alors, la distribution $J_{\ol}(\cdot)$ 
ne dépend que du choix de la mesure de Haar 
sur $U(\A)$ et des choix des mesures sur les points adéliques des parties unipotentes $N_{\tlP}$ des $\rmF$-sous-groupes 
paraboliques de $\tlU$, le choix étant que 
le volume $N_{\tlP}(\rmF) \bsl N_{\tlP}(\A)$  soit $1$.
\item Soit $y \in U(\A)$. On note $f^{y} \in C_{c}^{\infty}(\tlU(\A))$ 
la fonction donnée par $f^{y}(x) = f(\Ad(y)x)$. On a alors
\[
J_{\ol}(f^{y}) = J_{\ol}(f), \quad
\forall \, \ol \in \calO.
\]
\end{enumerate}
\bdem 
Tout comme le théorème \ref{thm:noyauGeomCVGG}, les 
lemmes préparatoires \ref{lem:caleyLeviU}, 
\ref{lem:finiOlsU} et
\ref{lem:LiedescenteU} ramènent la preuve aux cas des algèbres de Lie 
qui est traité en détail dans \cite{leMoi}.
\edem
\etheo

\subsection{La formule des traces relative}\label{par:RTFU}

Soit $\Phi \in C_{c}^{\infty}((U \times \tlU)(\A))$.
Notons alors pour $x \in \tlU(\A)$
\[
f_{\Phi}(x) := \int_{U(\A)}\Phi(y, x y) dy.
\]
On voit que $f_{\Phi} \in C_{c}^{\infty}(\tlU(\A))$. On définit alors dans ce cas pour $\ol \in \calO$ et $T \in T_{+} + \all_{0}^{+}$
\[
J_{\ol}^{T}(\Phi) := J_{\ol}^{T}(f_{\Phi}), \quad 
\]
et 
\[
J_{\ol}(\Phi) := J_{\ol}(f_{\Phi}).
\]

On utilise la notation de la section \ref{sec:SpecU}. 

\btheo[La formule des traces relative de Jacquet-Rallis, cas des groupes unitaires]\label{thm:RTFU}
Pour tout $\Phi \in C_{c}^{\infty}((U \times \tlU)(\A))$ et
tout $T \in T_{+} + \all_{0}^{+}$ on a
\begin{equation*}
\sum_{\chi \in \calX^{U \times \tlU}}J_{\chi}^{T}(\Phi) = 
\sum_{\ol \in \calO}J_{\ol}^{T}(\Phi)
\end{equation*}
et
\begin{equation*}
\sum_{\chi \in \calX^{U \times \tlU}}J_{\chi}(\Phi) = 
\sum_{\ol \in \calO}J_{\ol}(\Phi).
\end{equation*}
\bdem
La preuve du théorème est analogue à la preuve de son homologue, le théorème 
\ref{thm:RTFG} du paragraphe \ref{par:RTFG}.
\edem
\etheo
\appendix
\newcommand{\enref}[1]{\textit{\ref{#1})}}
\section{Un lemme}\label{app:lemme}

Soit $G$ un $\rmF$-groupe réductif comme dans le paragraphe \ref{par:prelimstraceSp}. Fixons 
$P_{0} \in \calF(M_{0})$ et notons $\calF(P_{0})$ l'ensemble des 
sous-$\rmF$-groupes paraboliques de $G$ contenant $P_{0}$. 

\blem\label{lem:projsArePositif} Soit $P \in \calF(P_{0})$ et soit $\varpi \in \hDelta_{0} \smin \hDelta_{P}$. Alors
\begin{enumerate}[a)]
\item\label{itm:Appa} Pour tout $\al^{\vee} \in \Delta_{P}^{\vee}$ on a $\varpi(\al^{\vee}) \ge 0$.
\item\label{itm:Appb} Pour tout $\varpi^{\vee} \in (\hDelta_{0}^{P})^{\vee}$ on a $\varpi(\varpi^{\vee}) \ge 0$.
\item\label{itm:Appc} Soit $\varpi^{\vee} \in \hDelta_{0}^{\vee}$ un élément correspondant à $\varpi$ et 
soit $\brvpi^{\vee} \in (\hDelta_{0}^{P})^{\vee}$ sa projection sur $\all_{0}^{P}$. 
Alors $\varpi(\brvpi^{\vee}) > 0$.
\end{enumerate}
\bdem
Le lemme 1.7.1 de \cite{labWal} montre que $\htau_{0} \ge \tau_{0}^{P} \htau_{P}$. Il est facile 
de voir que cela entraîne les points \enref{itm:Appa} et \enref{itm:Appb}.  
En ce qui concerne le point \enref{itm:Appc}, on a la forme de Killing $\bilif$ 
sur $\all_{0}^{G}$ qui est bilinéaire, symétrique, non-dégénérée et euclidienne. 
Elle identifie alors $(\all_{0}^{G})^{*}$ avec $\all_{0}^{G}$. 
Par cette identification $\varpi$ correspond à un multiple positif de $\varpi^{\vee}$. 
En plus, les espaces $\all_{0}^{P}$ et $\all_{P}^{G}$ sont orthogonaux pour la forme $\bilif$ 
ce qui démontre le dernier point. 
\edem
\elem

Soient $s \in \Omega^{G}$ et $\varpi \in \hDelta_{0}$. On note $\varpi^{s} = \varpi - s\varpi$ 
et $a_{\al,\varpi}^{s}$ les réels tels que $\varpi^{s} = \sum_{\al \in \Delta_{0}} a_{\al, \varpi}^{s} \al$. 
On note quelques propriétés standards utiles dans la suite
\begin{enumerate}
\item[\textbf{P1)}] Pour tout $\al \in \Delta_{0}$ on a $a_{\al, \varpi}^{s} \ge 0$ et 
si $\gamma \in \Delta_{0}$ est la racine correspondant à $\varpi$ on a 
$a_{\gamma, \varpi}^{s} > 0$ si et seulement si $\varpi^{s} \neq 0$.
\item[\textbf{P2)}]  Soient $P, Q \in \calF(P_{0})$ tels que $P \sbs Q$. Si $s \in \Omega^{Q}$ on a 
$a_{\al, \varpi}^{s} = 0$ pour tout $\al \in \Delta_{0} \smin \Delta_{0}^{Q}$. Dans ce cas, pour 
tout $H \in \all_{P}$ on a $\varpi^{s} (H) = \sum_{\al \in \Delta_{P}^{Q}} a_{\al} \al(H)$ où $a_{\al} \ge 0$ 
pour tout $\al \in \Delta_{P}^{Q}$.
\end{enumerate}

\blem\label{lem:notCrucBut}
 Soient $P, Q, P', Q' \in \calF(P_{0})$, $M, M' \ge 0$, $H, X \in \all_{P}^{G}$, $H', X' \in \all_{P'}^{G}$  
et $s, s' \in \Omega^{Q \cap Q'}$
tels que $P \cup P' \sbs Q \cap Q'$, $\tau_{P}^{Q}(H - X) = 1$, $\tau_{P'}^{Q'}(H' - X') = 1$ 
et
\begin{gather}
\varpi(H) - s'\varpi(H') \le M, \ \forall \ \varpi \in \hDelta_{P},  \label{eq:varpiMinsvarpi}\\
\varpi(H') - s\varpi(H) \le M', \ \forall \ \varpi \in \hDelta_{P'}. \label{eq:varpiMinsPrimvarpi}
\end{gather}
Alors, il existe une constante $C >0$ indépendante des $H, H', X, X', M, M'$ telle que 
\begin{gather}
\al(H) - \al(H') \le C(M + M' + \|X\| + \|X'\|), \ \forall \ \al \in \Delta_{0}^{Q} \smin \Delta_{0}^{P}, \label{eq:alMoinsal}\\
\al(H') - \al(H) \le C(M + M' + \|X\| + \|X'\|), \ \forall \ \al \in \Delta_{0}^{Q'} \smin \Delta_{0}^{P'}. \notag
\end{gather}
\bdem
Puisque l'énoncé est symétrique il suffit de montrer les inégalités (\ref{eq:alMoinsal}).
Soit $\gamma \in \Delta_{0}^{Q} \smin \Delta_{0}^{P}$  
et soit $\upvarpi \in \Delta_{0}^{Q} \smin \Delta_{0}^{P}$ le poids lui correspondant. 
Écrivons $\upvarpi$ dans la base $(\hDelta_{P'} \smin \{\upvarpi\}) \cup \Delta_{0}^{P'} \cup \{\gamma\}$ de 
$(\all_{0}^{G})^{*}$:
\[
\upvarpi = \sum_{\varpi \in \hDelta_{P'} \smin \{\upvarpi\}} c_{\varpi} \varpi+ 
\sum_{\al \in \Delta_{0}^{P'} \cup \{\gamma\}} c_{\al} \al.
\]
En vertu du lemme \ref{lem:projsArePositif} on a $c_{\varpi}, c_{\al} \ge 0$ 
et $c_{\gamma} >0$. Notons $\upvarpi_{0} = \sum_{\al \in \Delta_{0}^{P'} \cup \{\gamma\}} c_{\al} \al$. 

Pour tout $\varpi \in \hDelta_{P'} \smin \{\upvarpi\}$ on multiplie l'équation \eqref{eq:varpiMinsPrimvarpi} 
par $c_{\varpi}$ et ensuite on ajoute toutes ces équations 
à l'équation \eqref{eq:varpiMinsvarpi} avec $\varpi = \upvarpi$. On retrouve alors:
\begin{equation}\label{eq:majorIneq}
\upvarpi_{0}(H) + \sum_{\varpi \in \hDelta_{P'} \smin \{\upvarpi\}} c_{\varpi} \varpi^{s}(H) + 
\upvarpi^{s'}(H') -  \upvarpi_{0}(H') \le M' + CM
\end{equation}
où $C = \sum_{\varpi \in \hDelta_{P'} \smin \{\upvarpi\}} c_{\varpi}$. 

Soit $\varpi \in \hDelta_{P'} \smin \{\upvarpi\}$.
Puisque $Q \sps Q \cap Q'$ on a $\varpi^{s}(H) = \sum_{\al \in \Delta_{P}^{Q}} a_{\al, \varpi}^{s}\al(H)$ 
par la propriété \textbf{P2)} ci-dessus. 
On a $\al(H) \ge \al(X)$ pour tout $\al \in \Delta_{P}^{Q}$ par hypothèse $\tau_{P}^{Q}(H - X) = 1$. 
Grâce à la propriété \textbf{P1)} on 
a $a_{\al, \varpi}^{s} \ge 0$ pour tout $\al \in \Delta_{P}^{Q}$. On voit donc qu'il existe une constante $c'_{\varpi} >0$ telle 
que $-c_{\varpi}\varpi^{s}(H) \le c_{\varpi}'\|X\|$. Par le même raisonnement, il existe une constante $c_{0} >0$ 
telle que $-\upvarpi^{s'}(H') \le c_{0} \|X'\|$.
On a démontré alors que l'inégalité \eqref{eq:majorIneq} ci-dessus, entraîne
\begin{equation}\label{eq:varpiMinvarpiOnly}
\upvarpi_{0}(H) - \upvarpi_{0}(H') \le C_{0}(M' + M + \|X\| + \|X'\|)
\end{equation}
pour une constante $C_{0} >0$. 
 
La condition $P \cup P' \sbs Q \cap Q'$ est équivalente 
à $\Delta_{0}^{P} \cup \Delta_{0}^{P'} \sbs \Delta_{0}^{Q} \cap \Delta_{0}^{Q'}$ 
ce qui démontre que $\Delta_{0}^{P'} \sbs \Delta_{0}^{Q}$. 
En utilisant l'hypothèse $\tau_{P}^{Q}(H - X) = 1$ de nouveau 
et le fait que $\al(H) = 0$ pour tout $\al \in \Delta_{0}^{P}$, on trouve 
\[
\upvarpi_{0}(H) = \qquad 
\sum_{\mathclap{\al \in (\Delta_{0}^{P'} \cup \{\gamma\}) \cap (\Delta_{0}^{Q} \smin \Delta_{0}^{P})}}
\qquad
c_{\al} \al (H) \ge 
c_{\gamma} \gamma (H) \quad + \quad
\sum_{\mathclap{\al \in (\Delta_{0}^{P'} \cap (\Delta_{0}^{Q} \smin \Delta_{0}^{P})) \smin \{\gamma\}}} 
\qquad
c_{\al} \al(X)  
\ge c_{\gamma} \gamma (H)  - C_{0}'\|X\|
\]
pour une constante $C_{0}' >0$. 
En utilisant ceci, l'inégalité \eqref{eq:varpiMinvarpiOnly} et les faits que 
$\upvarpi_{0}(H') = c_{\gamma} \gamma(H')$ et $c_{\gamma} >0$ on conclut la preuve.
\edem
\elem

\bcor\label{cor:notCrucBut}
 Soient $P, Q, Q' \in \calF(P_{0})$, $P' \in \calF(M_{0})$, $M, M' \ge 0$, $H, X \in \all_{P}^{G}$, $H', X' \in \all_{P'}^{G}$  
et $s, s' \in \Omega^{Q \cap Q'}$
tels que $P \cup P' \sbs Q \cap Q'$, $\tau_{P}^{Q}(H - X) = 1$, $\tau_{P'}^{Q'}(H' - X') = 1$ 
et
\begin{gather}
\varpi(H) - s'\varpi(H') \le M, \ \forall \ \varpi \in \hDelta_{P},  \label{eq:varpiMinsvarpic}\\
\varpi(H') - s\varpi(H) \le M', \ \forall \ \varpi \in \hDelta_{P'}. \label{eq:varpiMinsPrimvarpic}
\end{gather}
Alors, il existe une constante $C >0$ indépendante des $H, H', X, X', M, M'$ telle que 
si l'on note $P_{1}$ le plus petit sous-groupe parabolique contenant $P$ et $P'$ 
et $s_{0} \in \Omega^{P_{1}}$ tel que $s_{0}^{-1}P' \in \calF(P_{0})$ on a
\begin{gather*}
\al(H) - \al(s_{0}^{-1}H') \le C(M + M' + \|X\| + \|X'\|), \ \forall \ \al \in \Delta_{0}^{Q} \smin \Delta_{0}^{P},\\
\al(H') - \al(s_{0}H) \le C(M + M' + \|X\| + \|X'\|), \ \forall \ \al \in \Delta_{s_{0}P_{0}}^{Q'} \smin \Delta_{s_{0}P_{0}}^{P'}.
\end{gather*}
\bdem 
Soit $s_{0}$ comme dans l'énoncé. Notons $P'' = s_{0}^{-1}P'$. On a alors, que les inégalités 
(\ref{eq:varpiMinsvarpic}) et (\ref{eq:varpiMinsPrimvarpic}) s'écrivent comme
\begin{gather*}
\varpi(H) - s_{0}^{-1}s'\varpi(s_{0}^{-1}H') \le M, \ \forall \ \varpi \in \hDelta_{P},  \\
\varpi(s_{0}^{-1}H') - s s_{0}\varpi(H) \le M', \ \forall \ \varpi \in \hDelta_{P''}.
\end{gather*}
Remarquons que $s_{0}Q' = Q'$.
On a donc bien $P \cup P'' \sbs Q \cap Q'$ et $s_{0}^{-1}s', s s_{0} \in \Omega^{Q \cap Q'}$. 
En plus on a $\tau_{P''}^{Q'}(s_{0}^{-1}H' - s_{0}^{-1}X') = 1$. En appliquant alors le 
lemme \ref{lem:notCrucBut} on trouve le résultat voulu.
\edem
\ecor

Pour $s \in \Omega$ on pose $\hDelta_{s} = \{\varpi \in \hDelta_{0} | \varpi^{s} = 0\}$.

\blem\label{lem:fuckYeahLemme} Soient $P_{1}, P_{3}, P_{4}, P_{5}, P_{6} \in \calF(P_{0})$, 
$P_{2} \in \calF(M_{0})$,
$M_{1}, M_{2}, M_{3}, M_{4} \ge 0$
$H_{1}, X_{1} \in \all_{1}^{G}$, $H_{2}, X_{2} \in \all_{2}^{G}$, $H_{3}, X_{3} \in \all_{3}^{G}$, $s_{1}, s_{2}, s_{1}', s_{2}' \in \Omega^{P_{4} \cap P_{5} \cap P_{6}}$. 
Soit $Q_{1}$ le plus petit sous-groupe parabolique contenant $P_{1} \cup P_{2} \cup P_{3}$. 
Supposons qu'on a $ Q_{1} \sbs P_{4} \cap P_{5} \cap P_{6} =: Q_{2}$ et que
\begin{gather}
\hDelta_{Q_{2}} = \hDelta_{Q_{1}} \cap \hDelta_{s_{1}} \cap \hDelta_{s_{2}} = 
\hDelta_{Q_{1}} \cap \hDelta_{s_{1}'} \cap \hDelta_{s_{2}'},  \label{eq:hDeltaEses}\\
s_{1}'(\hDelta_{2} \smin \hDelta_{0}) \cap \hDelta_{0} = \varnothing  \label{eq:hDeltaEses2}\\
\sigma_{1}^{4}(H_{1} - X_{1}) = \sigma_{2}^{5}(H_{2} - X_{2}) = \sigma_{3}^{6}(H_{3} - X_{3}) = 1, \notag \\
\varpi(H_{1}) - s_{1}\varpi(H_{2}) \le M_{1}, \ \forall \ \varpi \in \hDelta_{1},  \label{eq:varpiMinsvarpi1}\\
\varpi(H_{2}) - s_{1}'\varpi(H_{1}) \le M_{2}, \ \forall \ \varpi \in \hDelta_{2},  \label{eq:varpiMinsvarpi2}\\
\varpi(H_{1}) - s_{2}'\varpi(H_{3}) \le M_{3}, \ \forall \ \varpi \in \hDelta_{1},  \label{eq:varpiMinsvarpi3}\\
\varpi(H_{3}) - s_{2}\varpi(H_{1}) \le M_{4}, \ \forall \ \varpi \in \hDelta_{3}.  \label{eq:varpiMinsvarpi4}
\end{gather}
Alors, il existe une constante $C >0$ indépendante des éléments
$H_{1}$, $H_{2}$, $H_{3}$, $X_{1}$, $X_{2}$, $X_{3}$, $M_{1}$, $M_{2}$, $M_{3}$, $M_{4}$
telle que
\begin{equation}\label{eq:whatWeWant}
\|H_{1}\|, \|H_{2}\|, \|H_{3}\| \le C(M_{1} + M_{2} + M_{3} + M_{4} + \|X_{1}\| + \|X_{2}\| + \|X_{3}\|).
\end{equation}
\bdem
Dans la preuve on ne prendra pas garde de constantes qui apparaissent. Il est clair que pour obtenir 
une constante $C$ uniforme il suffit de prendre le maximum de toutes les constantes en question. 

Remarquons que si l'on démontre 
\begin{equation}\label{eq:whatWeReallyWant}
\al(H_{1}) \le C(M_{1} + M_{2} + M_{3} + M_{4} + \|X_{1}\| + \|X_{2}\| + \|X_{3}\|)
\end{equation}
pour tout $\al \in \Delta_{1}^{4}$
on a aussi l'inégalité (\ref{eq:whatWeWant}) pour $H_{1}$.
En effet, par définition de $\sigma_{1}^{4}$ (voir (\ref{eq:sigmaDef}) dans le paragraphe \ref{par:resGen}) on a 
$\al(H_{1}) \ge \al(X_{1})$ pour tout $\al \in \Delta_{1}^{4}$. Puisque $\Delta_{1}^{4}$ 
est une base de $\all_{1}^{4}$ on a que la projection de $H$ à $\all_{1}^{4}$ 
vérifie (\ref{eq:whatWeWant}) et on obtient le résultat pour $H_{1}$ en invoquant le lemme \ref{lem:corArth62} \textit{ii)}.

Supposons maintenant qu'on a démontré l'inégalité dans la proposition pour $H_{1}$ seulement. 
On veut en déduire la même inégalité pour $H_{2}$. 
En remarquant qu'on a l'inégalité des fonctions caractéristiques $\tau_{P}^{Q} \ge \sigma_{P}^{Q}$, 
on applique
le corollaire \ref{cor:notCrucBut} avec $P = P_{1}$, $Q = P_{4}$, $P' = P_{2}$ 
et $Q' = P_{5}$ avec les inégalités (\ref{eq:varpiMinsvarpi1}) et (\ref{eq:varpiMinsvarpi2}) 
ci-dessus et on obtient une constante $C' >0 $ telle que
\[
\al(H_{2}) \le C'(M_{1} + M_{2} + \|X_{1}\| + \|X_{2}\| + \|H_{1}\|)
\]
pour tout 
$\al \in \Delta_{2}^{5}$. 
D'après ce qu'on a dit ci-dessus cela implique l'inégalité (\ref{eq:whatWeWant}) 
pour $H_{2}$ si $H_{1}$ la vérifie. 
De même façon, en passant directement par le lemme \ref{lem:notCrucBut}, 
on montre que si $H_{1}$ vérifie 
(\ref{eq:whatWeWant}) alors $H_{3}$ la vérifie aussi.

Fixons un $\gamma \in \Delta_{0}^{4} \smin \Delta_{0}^{1}$. Il suffit donc de montrer que $\gamma(H_{1})$ 
vérifie l'inégalité (\ref{eq:whatWeReallyWant}). 
Soit $s_{0} \in \Omega^{Q_{1}}$ tel que $s_{0}^{-1}P_{2} \in \calF(P_{0})$.
Posons
\begin{gather*}
\calS_{1} \ \mathclap{=} \ (\hDelta_{1} \ \mathclap{\smin} \ \hDelta_{4}) \ \mathclap{\cap} \ \hDelta_{6}, \ 
\calS_{2} \ \mathclap{=} \  (\hDelta_{1} \ \mathclap{\smin} \ \hDelta_{4}) \ \mathclap{\cap} \  (\hDelta_{0} \ \mathclap{\smin} \ \hDelta_{3}), \ 
\calS_{3} \ \mathclap{=} \ (\hDelta_{1} \ \mathclap{\smin} \ \hDelta_{4}) \ \mathclap{\cap} \   \hDelta_{5}, \
\calS_{4} \ \mathclap{=} \  (\hDelta_{1} \ \mathclap{\smin} \ \hDelta_{4}) \ \mathclap{\cap} \  (\hDelta_{0} \ \mathclap{\smin} \ \hDelta_{s_{0}^{-1}P_{2}}), \\
\calS_{5} \ \mathclap{=} \ (\hDelta_{1} \ \mathclap{\smin} \ \hDelta_{4}) \ \mathclap{\cap} \  (\hDelta_{3} \ \mathclap{\smin} \ \hDelta_{6}) \ \mathclap{\cap} \ 
 (\hDelta_{Q_{1}} \ \mathclap{\smin} \ \hDelta_{5}), \ 
 \calS_{6} \ \mathclap{=} \  (\hDelta_{1} \ \mathclap{\smin} \ \hDelta_{4}) \ \mathclap{\cap} \  (\hDelta_{3} \ \mathclap{\smin} \ \hDelta_{6}) \ \mathclap{\cap} \ 
 (\hDelta_{s_{0}^{-1}P_{2}} \ \mathclap{\smin} \ \hDelta_{Q_{1}}).
\end{gather*}
On a donc $(\hDelta_{1} \smin \hDelta_{4}) = \bigcup_{i=1}^{6}\calS_{i}$.
En utilisant le corollaire \ref{cor:notCrucBut} avec $P = P_{1}$, $Q = P_{4}$, $P' = P_{2}$ 
et $Q' = P_{5}$ avec les inégalités (\ref{eq:varpiMinsvarpi1}) et (\ref{eq:varpiMinsvarpi2}) 
ainsi que le lemme \ref{lem:notCrucBut} avec 
$P = P_{1}$, $Q = P_{4}$, $P' = P_{3}$ 
et $Q' = P_{6}$ avec les inégalités (\ref{eq:varpiMinsvarpi3}) et (\ref{eq:varpiMinsvarpi4}) 
on trouve une constante $C' >0$ telle que
\begin{gather}
\gamma(H_{1}) \le \gamma(s_{0}^{-1}H_{2}) + C'(M_{1} + M_{2} + \|X_{1}\| + \|X_{2}\|), \label{eq:almostThereP2} \\
\gamma(H_{1}) \le \gamma(H_{3}) + C'(M_{3} + M_{4} + \|X_{1}\| + \|X_{3}\|). \label{eq:almostThereP3}
\end{gather}
Soit $\upvarpi \in \hDelta_{1} \smin \hDelta_{4}$ le poids correspondant à $\gamma \in \Delta_{0}^{4} \smin \Delta_{0}^{1}$ qu'on a choisi. 
Alors
\begin{itemize}
\item Si $\upvarpi \in \calS_{1}$, on a $\gamma \in \Delta_{0} \smin \Delta_{0}^{6}$ 
et donc $\gamma(H_{3}) \le \gamma(X_{3}) \le C''\|X_{3}\|$ par définition de $\sigma_{3}^{6}$ 
et le résultat suit de l'inégalité (\ref{eq:almostThereP3}). 
\item Si $\upvarpi \in \calS_{2}$ on a $\gamma \in \Delta_{0}^{3}$ et donc $\gamma(H_{3}) = 0$ et le résultat suit aussi de (\ref{eq:almostThereP3}). 
\item Si $\upvarpi \in \calS_{3}$ on a $\gamma(s_{0}^{-1}H_{2}) = s_{0}\gamma(H_{2})$ 
et $s_{0}\gamma \in \Delta_{s_{0}P_{0}} \smin \Delta_{s_{0}P_{0}}^{s_{0}P_{5}}$. 
Mais $s_{0}P_{5} = P_{5}$ car $s_{0} \in \Omega^{Q_{1}} \sbs \Omega^{P_{5}}$. 
Le résultat suit donc de l'inégalité (\ref{eq:almostThereP2}) et de la définition de 
la fonction $\sigma_{2}^{5}$. 
\item  Si $\upvarpi \in \calS_{4}$ on a $s_{0} \gamma \in \Delta_{s_{0}P_{0}}^{P_{2}}$ 
donc $\gamma(s_{0}^{-1}H_{2}) = 0$ et le résultat suit de (\ref{eq:almostThereP2}) de nouveau.
\item Supposons que $\upvarpi \in \calS_{5}$. 
On a alors $\gamma \in (\Delta_{0}^{4} \cap \Delta_{0}^{5} \cap \Delta_{0}^{6}) \smin (\Delta_{0}^{1} \cup \Delta_{0}^{Q_{1}} \cup \Delta_{0}^{3})$. 
En plus, $\upvarpi \nin \hDelta_{4} \cup \hDelta_{5} \cup \hDelta_{6} = \hDelta_{Q_{2}}$ 
et $\upvarpi  \in \hDelta_{Q_{1}}$.
En vertu de la condition \eqref{eq:hDeltaEses} on a donc que soit $\upvarpi  \nin \hDelta_{s_{1}'}$ 
soit $\upvarpi  \nin \hDelta_{s_{2}'}$.  
\begin{itemize}
\item Supposons $\upvarpi  \nin \hDelta_{s_{2}'}$.
 Ajoutons les inégalités 
(\ref{eq:varpiMinsvarpi3}) et (\ref{eq:varpiMinsvarpi4}) avec $\varpi = \upvarpi$, ce qu'on peut faire car 
$\upvarpi \in (\hDelta_{1} \smin  \hDelta_{4}) \cap (\hDelta_{3} \smin \hDelta_{6}) \sbs \hDelta_{1} \cap \hDelta_{3}$.
On obtient donc:
\[
\upvarpi^{s_{2}}(H_{1}) + \upvarpi^{s_{2}'}(H_{3}) \le M_{3} + M_{4}.
\]
D'après les propriétés \textbf{P1)} et \textbf{P2)} ci-dessus, on a $\upvarpi^{s_{2}}(H_{1}) = 
\sum_{\al \in \Delta_{0}^{4} \smin \Delta_{0}^{1}}a_{\al} \al(H_{1})$ 
et $\upvarpi^{s_{2}'}(H_{3}) = \sum_{\al \in \Delta_{0}^{6} \smin \Delta_{0}^{3}}a'_{\al} \al(H_{3})$
où $a_{\al}, a'_{\al} \ge 0$ et $a'_{\gamma} > 0$. En utilisant la définition 
des fonctions $\sigma_{1}^{4}$ et $\sigma_{3}^{6}$, pour tout $\al \in \Delta_{0}^{4} \smin \Delta_{0}^{1}$ 
on a $\al(H_{1}) \ge \al(X_{1})$ et pour tout $\al \in \Delta_{0}^{6} \smin \Delta_{0}^{3}$ on a $\al(H_{3}) \ge \al(X_{3})$.
On trouve donc une constante 
$C' >0$ telle que $a'_{\gamma}\gamma(H_{3})-\upvarpi^{s_{2}'}(H_{3}) \le C' \|X_{3}\|$ 
et $-\upvarpi^{s_{2}}(H_{1}) \le C' \|X_{1}\|$. 
D'où 
\[
\gamma(H_{3}) \le \dfrac{1}{a'_{\gamma}}(C'\|X_{1}\| + C'\|X_{3}\| + M_{3} + M_{4}).
\]
En utilisant alors l'inégalité (\ref{eq:almostThereP3}) on obtient le résultat voulu. 
\item Supposons $\upvarpi  \nin \hDelta_{s_{1}'}$. 
On a alors $\upvarpi \in \calS_{5} \sbs \hDelta_{Q_{1}}$ 
et $s_{0} \in \Omega^{Q_{1}}$ donc $s_{0}\upvarpi = \upvarpi$. 
Donc $\upvarpi \in \hDelta_{1} \cap \hDelta_{2}$. 
On ajoute alors les inégalités 
(\ref{eq:varpiMinsvarpi1}) et (\ref{eq:varpiMinsvarpi2}) avec $\varpi = \upvarpi$
ce qui donne
\[
\upvarpi^{s_{1}'}(H_{1}) + \upvarpi^{s_{1}}(H_{2}) \le M_{1} + M_{2}.
\]
Puisque $s_{1} \in \Omega^{Q_{1}} \sbs  \Omega^{P_{5}}$, 
on a $ \upvarpi^{s_{1}}(H_{2}) = \sum_{\al \in \Delta_{2}^{5}}a_{\al}\al(H_{2})$ 
où $a_{\al} \ge 0$ en vertu de la propriété \textbf{P2)} ci-dessus. En utilisant donc 
la définition de la fonction $\sigma_{2}^{5}$ on trouve une constante $C' >0$ telle que 
$-\upvarpi^{s_{1}}(H_{2}) \le C' \|X_{2}\|$. 
Donc, puisque $\upvarpi^{s_{1}'} \neq 0$, on a 
$\upvarpi^{s_{1}'}(H_{1}) = \sum_{\al \in \Delta_{0}^{4} \smin \Delta_{0}^{1}}a_{\al} \al(H_{1})$ 
où $a_{\al} \ge 0$ pour tout $\al \in \Delta_{0}^{4} \smin \Delta_{0}^{1}$ 
et $a_{\gamma} >0$. Par définition de la fonction $\sigma_{1}^{4}$ on obtient alors:
\[
\gamma(H_{1}) \le c''(C''\|X_{1}\| + C'\|X_{2}\| + M_{1} + M_{2}).
\]
pour certaines constantes $c'', C'' >0$.
\end{itemize}
\item Supposons que $\upvarpi \in \calS_{6}$. 
On a $s_{0}\upvarpi \in \hDelta_{2}$. 
 Ajoutons donc l'inégalité
(\ref{eq:varpiMinsvarpi1}) avec $\varpi = \upvarpi$ et 
(\ref{eq:varpiMinsvarpi2}) avec $\varpi = s_{0}\upvarpi$. 
On trouve:
\[
\upvarpi^{s_{1}'s_{0}}(H_{1}) + (s_{0}\upvarpi)^{s_{1}s_{0}^{-1}}(H_{2}) \le M_{1} + M_{2}.
\]
On a $\hDelta_{Q_{1}} = \hDelta_{1} \cap \hDelta_{3} \cap \hDelta_{s_{0}^{-1}P_{2}} \cap \hDelta_{s_{0}}$. 
Puisque $\upvarpi \nin \hDelta_{Q_{1}}$ et $\upvarpi \in \hDelta_{1} \cap \hDelta_{3} \cap \hDelta_{s_{0}^{-1}P_{2}}$, 
on a bien $s_{0}\upvarpi \neq \upvarpi$. 
Puisque différents éléments de $\hDelta_{0}$ ne sont pas conjugués sous $\Omega^{G}$ on obtient $s_{0}\upvarpi \in \hDelta_{2} \smin \hDelta_{0}$. 
En vertu de la condition (\ref{eq:hDeltaEses2}) on alors $\upvarpi^{s_{1}'s_{0}} \neq 0$ 
et on conclut en raisonnant de même façon que dans le deuxième sous-point 
du point précédent. 
\end{itemize}
\edem
\elem

\bcor\label{cor:fuckYeah} Soient $P_{1}, P_{2}, P_{3}, P_{4} \in \calF(P_{0})$, 
$M_{1}, M_{2}  \ge 0$
$H_{1}, X_{1} \in \all_{1}^{G}$, $H_{2}, X_{2} \in \all_{2}^{G}$, 
$s, s' \in \Omega^{P_{3} \cap P_{4}}$. 
Soit $Q_{1}$ le plus petit sous-groupe parabolique contenant $P_{1} \cup P_{2}$. 
Supposons qu'on a $ Q_{1} \sbs P_{3} \cap P_{4} =: Q_{2}$ et que
\begin{gather}
\hDelta_{Q_{2}} = \hDelta_{Q_{1}} \cap \hDelta_{s} \cap \hDelta_{s'}  \label{eq:hDeltaEsesU}\\
\sigma_{1}^{3}(H_{1} - X_{1}) = \sigma_{2}^{4}(H_{2} - X_{2}) = 1, \notag \\
\varpi(H_{1}) - s\varpi(H_{2}) \le M_{1}, \ \forall \ \varpi \in \hDelta_{1},  \label{eq:varpiMinsvarpi1U}\\
\varpi(H_{2}) - s'\varpi(H_{1}) \le M_{2}, \ \forall \ \varpi \in \hDelta_{2},  \label{eq:varpiMinsvarpi2U},
\end{gather}
Alors, il existe une constante $C >0$ indépendante des éléments
$H_{1}$, $H_{2}$, $X_{1}$, $X_{2}$, $M_{1}$, $M_{2}$ telle que
\begin{equation}\label{eq:whatWeWantU}
\|H_{1}\|, \|H_{2}\| \le C(M_{1} + M_{2} \|X_{1}\| + \|X_{2}\|\|).
\end{equation}
\bdem
Il s'agit de se ramener au lemme \ref{lem:fuckYeahLemme} ci-dessus. 
On l'utilise alors avec les données suivantes. 
Pour les groupes $P_{1}$, $P_{2}$, $ P_{3}$, $P_{4}$, $P_{5}$ et $P_{6}$ du lemme 
on prend $P_{1}$, $P_{2}$, $P_{2}$, $P_{3}$, $P_{4}$ et $P_{4}$ dans cet ordre.  
Pour les éléments $H_{1}$, $H_{2}$ et $H_{3}$ on prend 
$H_{1}$, $H_{2}$ et $H_{2}$ respectivement.
 Pour $X_{1}$, $X_{2}$ et $X_{3}$ on prend $X_{1}$, 
 $X_{2}$ et $X_{2}$ respectivement. 
 Pour les éléments du groupe de Weyl on prend 
 $s_{1} = s$, $s_{1}' = s'$, $s_{2} = s'$ et $s_{2}' = s$. 
 Il est clair que la condition (\ref{eq:hDeltaEsesU}) ci-dessus 
 donne la condition (\ref{eq:hDeltaEses}) du lemme. 
Les inégalités du lemme \ref{lem:fuckYeahLemme} 
 correspondent aux inégalités suivantes:
 \[
\text{(\ref{eq:varpiMinsvarpi1})} \leftrightarrow  
\text{(\ref{eq:varpiMinsvarpi1U})}, \quad 
\text{(\ref{eq:varpiMinsvarpi2})} \leftrightarrow  
\text{(\ref{eq:varpiMinsvarpi2U})}, \quad 
\text{(\ref{eq:varpiMinsvarpi3})} \leftrightarrow  
\text{(\ref{eq:varpiMinsvarpi1U})}, \quad 
\text{(\ref{eq:varpiMinsvarpi4})} \leftrightarrow  
\text{(\ref{eq:varpiMinsvarpi2U})}. 
 \] 
Les inégalités déterminent les constantes 
 $M_{1}$, $M_{2}$, $M_{3}$ et $M_{4}$. Finalement, la condition (\ref{eq:hDeltaEses2}) est trivialement vérifiée 
 car $P_{2} \in \calF(P_{0})$.
\edem
\ecor

\bibliographystyle{alpha-fr}
 \bibliography{bibliography}
\end{document}